\documentclass[letterpaper,12pt]{article}
\usepackage[english]{babel}
\usepackage[utf8x]{inputenc}
\usepackage[T1]{fontenc}
\usepackage{url}
\usepackage{setspace}
\usepackage[title, titletoc]{appendix}
\usepackage{yfonts}

\usepackage[letterpaper,top=3cm,bottom=2cm,left=3cm,right=3cm,marginparwidth=1.75cm]{geometry}

\usepackage{amsmath}
\usepackage{amssymb} 
\usepackage{csquotes}
\usepackage[colorinlistoftodos]{todonotes}
\usepackage[colorlinks=true, allcolors=blue]{hyperref}
\usepackage{setspace}
\usepackage{graphicx}
\usepackage{listings}
\usepackage{epigraph}
\usepackage{wasysym}

\numberwithin{equation}{section}

\title{Analytic Continuation of $\zeta(s)$ Violates the Law of Non-Contradiction (LNC)}

\author{Ayal Sharon \thanks{Patent Examiner, U.S. Patent and Trademark Office (USPTO). ayal.sharon@uspto.gov. This research received no funding, and was conducted in the author's off-duty time. The opinions expressed herein are solely the author's, and do not reflect the views of the USPTO, the U.S. Dept. of Commerce, or the U.S. Government.} 
} 

\begin{document}
\maketitle
\begin{abstract}

The Dirichlet series of $\zeta(s)$ was long ago proven to be divergent throughout half-plane $\text{Re}(s)\le1$. If also Riemann's proposition is true, that there exists an "expression" of $\zeta(s)$ that is convergent at all $s$ (except at $s=1$), then $\zeta(s)$ is both divergent and convergent throughout half-plane $\text{Re}(s)\le1$ (except at $s=1$). 

This result violates all three of Aristotle's "Laws of Thought": the Law of Identity (LOI), the Law of the Excluded Middle (LEM), and the Law of Non-Contradition (LNC). In classical and intuitionistic logics, the violation of LNC also triggers the "Principle of Explosion" / \textit{Ex Contradictione Quodlibet} (ECQ).

In addition, the Hankel contour used in Riemann's analytic continuation of $\zeta(s)$ violates Cauchy's integral theorem, providing another proof of the invalidity of Riemann's $\zeta(s)$. Riemann's $\zeta(s)$ is one of the $L$-functions, which are all invalid due to analytic continuation. This result renders unsound all theorems (e.g. Modularity, Fermat's last) and conjectures (e.g. BSD, Tate, Hodge, Yang-Mills) that assume that an $L$-function (e.g. Riemann's $\zeta(s)$) is valid.

We also show that the Riemann Hypothesis (RH) is not "non-trivially true" in classical logic, intuitionistic logic, or three-valued logics (3VLs) that assign a third truth-value to paradoxes (Bochvar's 3VL, Priest's $LP$).

\end{abstract}

\pagebreak
\tableofcontents
\pagebreak

\epigraph{Tutte le verità sono facili da capire una volta che sono state rivelate. Il difficile è scoprirle.}{Galileo Galilei}

\onehalfspacing

\section{Introduction}

\subsection{The Crux of the Problem}

Riemann \cite{riemann1859number} states the following on the first page of his famous paper, "On the Number of Prime Numbers less than a Given Quantity": 
\footnote{In the original German: "Ueber die Anzahl der Primzahlen unter einer gegebenen Grösse".}
\begin{quotation}
For this investigation my point of departure is provided by the observation of Euler that the product
\begin{equation}
\prod \frac{1}{1-1/p^s} = \sum \frac{1}{n^s}   
\end{equation}
if one substitutes for $p$ all prime numbers, and for $n$ all whole numbers.
The function of the complex variable $s$ which is represented by these two expressions, wherever they converge, I denote by $\zeta(s)$. Both expressions converge only when the real part of $s$ is greater than 1; at the same time an expression for the [zeta] function can easily be found which always remains valid.
\end{quotation}
Riemann's proposition is that there exists an "expression" for $\zeta(s)$ which is convergent for all values of $s$, \textit{in addition to} the Dirichlet series of $\zeta(s)$, which is proven to be \textit{divergent} throughout half-plane $\text{Re}(s)\le1$ (as admitted by Riemann in the cited text). Riemann's proposition violates all three of Aristotle's "Laws of Thought": the Law of Identity (LOI), the Law of the Excluded Middle (LEM), and the Law of Non-Contradition (LNC). 
\footnote{The one-to-two mapping from domain to range also violates the set theory definition of "a function".}
In classical and intuitionistic logics, this violation of LNC triggers the "Principle of Explosion" (\textit{Ex Contradictione Quodlibet}, or "ECQ").

Unfortunately, Riemann confused the mathematical concept of "convergence" with the logical concept of "validity." The Dirichlet series of $\zeta(s)$ is proven to be divergent (i.e. \textit{not convergent}), throughout half-plane $\text{Re}(s)\le1$, and this proof is logically \textit{valid}. Riemann's alternative "expression" for $\zeta(s)$ claims to be \textit{convergent} throughout that same half-plane.
\footnote{Except for a solitary pole at $s=1$.}
If true, $\zeta(s)$ would be both \textit{convergent} and \textit{not convergent} throughout half-plane $\text{Re}(s)\le1$.
\footnote{Again, except for a solitary pole at $s=1$.}
So Riemann's "analytic continuation" of $\zeta(s)$ violates the LOI, LEM, and LNC; and is \textit{not valid} in logics that have any of these as axioms. \footnote{And as stated above, the one-to-two mapping of domain-to-range also violates the set theory definition of a "function".}

The remainder of this paper provides commentary on this fact, including a discussion of the errors in the derivation of Riemann's $\zeta(s)$, and a discussion of the resulting implications for the Riemann Hypothesis and other conjectures that assume that the "analytic continuation" of $\zeta(s)$ is valid. 

Riemann's $\zeta(s)$ is one of the $L$-functions, and all $L$-functions invalid, for the same reason: Riemann's version of analytic continuation results in a contradiction. The invalidity of $L$-functions renders unsound all theorems (e.g. Modularity theorem, Fermat's last theorem) and conjectures (e.g. BSD, Tate, Hodge, Yang-Mills) that falsely assume that an $L$-function (e.g. Riemann's $\zeta(s)$) or "zeta function regularization" is valid.

\subsection{Attacking the Most Specific Problem with the Oldest Ideas}

In "Problems of the Millennium: the Riemann Hypothesis", Bombieri \cite{Bombieri} states:
\begin{quotation}
Not a single example of validity or failure of a Riemann hypothesis for an $L$-function is known up to this date. The Riemann hypothesis for $\zeta(s)$ does not seem to be any easier than for Dirichlet $L$-functions (except possibly for non-trivial Real zeros), leading to the view that its solution may require attacking much more general problems, by means of entirely new ideas.
\footnote{See Bombieri \cite{Bombieri}, p.5.}
\end{quotation}

The present paper takes the opposite approach: attacking the most specific problem by means of the oldest possible ideas. More specifically, the present paper attacks the Riemann hypothesis (RH) by proving the invalidity of the most famous of the $L$-functions: Riemann's $\zeta(s)$. This is done by means of logics of the early 20th century (i.e. "classical logic",  intuitionistic logic, and three-valued logics (3VLs)), using concepts inherited from Aristotelian and medieval logic. 
These inherited concepts are certainly not new. For example, the three "Laws of Thought" 
\footnote{The three "Laws of Thought" are: the Law of Non-Contradiction (LNC), Law of the Excluded Middle (LEM), and Law of Identity (LOI).}
and syllogism are discussed in Aristotle's \textit{Organon} \cite{Aristotle}. 
\footnote{Aristotle's life is dated as 384–322 BCE. See Wikipedia \cite{Aristotle3}, citing Boeckh \cite{Boeckh}, vol.VI, p.195, Jacoby \cite{Jacoby}, FGrHist 244 F 38, and Düring \cite{Duering}, p.253.}
\footnote{The LNC pre-dates Aristotle. See, e.g. Plato's Socratic dialogue \textit{Euthyphro} at Plato \cite{Plato}, p. 264-265, \S 8: "Then the same things are hated by the gods and loved by the gods, and are both hateful and dear to them? ... And upon this view the same things, Euthyphro, will be pious and also impious?"}
\footnote{See also Cohen \cite{Cohen2}, p.328: "Actually, the laws go back well before Aristotle, who was essentially summarizing the views of the pre-Socratic philosophers, most notably Parmenides. It was he, in the 5th century BCE, who had formulated the [LNC] as 'Never will this prevail, that what is not, is.'"} 
\footnote{See also Boole's \cite{Boole} discussion of the LOI in Chapter II, pp.34-36, Para.12-13; the LNC in Chapter III, p.49, Prop. IV; and the LEM in Chapter III, p.48, Prop. II, and also in pp.8 and 99-100.}
The Square of Opposition "shows up already in the second century CE", and "Boethius incorporated it into his writing". 
\footnote{See Parsons \cite{Parsons}, \S2.1: "The diagram accompanying and illustrating the doctrine shows up already in the second century CE; Boethius incorporated it into his writing, and it passed down through the dark ages to the high medieval period, and from thence to today. Diagrams of this sort were popular among late classical and medieval authors, who used them for a variety of purposes."}
The "Principle of Explosion" (\textit{Ex Contradictione Quodlibet}, or "ECQ") is a theorem that dates back to the 12th century. 
\footnote{The first proof of this principle is attributed to 12th century French philosopher William of Soissons. See Wikipedia \cite{Explosion}, citing Priest \cite{Priest9}, p.25, which in turn cites Priest \cite{Priest10}, vol.6, ch.4.}

The focus in this paper is on the Law of Non-Contradiction (LNC) and the "Principle of Explosion" (\textit{Ex Contradictione Quodlibet}, or "ECQ"). LNC is an axiom, and ECQ is a theorem, in the "classical logic" of Whitehead and Russell's \textit{Principia Mathematica} \cite{Whitehead2}, and also in Intuitionistic logics (e.g. Brouwer's and Heyting's). However, LNC and ECQ both fail in the three-valued logics (3VLs) discussed in this paper (e.g. Bochvar's 3VL, and Priest's "Logic of Paradox", which itself is a version of Kleene's 3VL), due to these logics having more than two truth-values. 

Łukasiewicz's Three-Valued Logic (3VL) (c. 1920),
\footnote{See e.g. Łukasiewicz \cite{Lukasiewicz4}.}
Brouwer's intuitionistic logic (c. 1921),
\footnote{See Brouwer \cite{Brouwer4}.}
and Whitehead and Russell's "classical logic" (c. 1925), 
\footnote{See Whitehead et al. \cite{Whitehead2}.}
are all approximately hundred years old. Heyting's intuitionistic logic (c. 1930); 
\footnote{See Heyting \cite{Heyting2} (in German) and \cite{Heyting3} (in English). }
and Bochvar's 3VL (c. 1938),
\footnote{See Bochvar \cite{Bochvar}.}
are slightly more recent.
Kleene's 3VL (c. 1952), 
\footnote{See Kleene \cite{Kleene2}, first published in 1952.}
is over 65 years old. Priest's "Logic of Paradox" ($LP$) (c. 1979)
\footnote{See Priest \cite{Priest5} and Priest \cite{Priest7}.}
is by far the youngest logic of this group, at "only" 40 years old. 

But in spite of the ages of these logics, the author has not found \textit{any} discussion of the RH in the context of \textit{any} of these logics. It appears that logicians do not apply logic to specific examples, and that mathematicians lost interest in foundational questions after the development of Zermelo-Fraenkel set theory.
\footnote{See Wikipedia \cite{ZFC}:" Today, Zermelo–Fraenkel set theory, with the historically controversial axiom of choice (AC) included, is the standard form of axiomatic set theory and as such is the most common foundation of mathematics."}

\subsection{RH is a Problem in Logic, Regardless of Logic's Relationship to Math}

On a related note, Oort et al. \cite{Oort} states:
\begin{quotation}
Historically as well as mathematically, the real conundrum is: where do the Riemann Hypothesis and its avatars belong in the vast and changing landscape of mathematics? The day we will see a proof of the Riemann Hypothesis, this will root and place the statement for the first time.
\footnote{See Oort et al. \cite{Oort}, p.596.}
\end{quotation}

The present paper asserts that Turing \cite{Turing} erred when he classified the RH as a "number-theoretic" problem.
\footnote{See Turing \cite{Turing}, p.165: "It is easy to show that a number of unsolved problems, such as the problem of the truth of Fermat's last theorem, are number-theoretic. There are, however, also problems of analysis which are number-theoretic. The Riemann hypothesis gives us an example of this."} 
Instead, RH is a problem in logic (a.k.a. "foundations of mathematics", or "meta-mathematics"). 

In the logics discussed in this paper, the Riemann hypothesis is false, or a paradox, or "trivially true" (depending on the logic, and whether the analytic continuation of $\zeta(s)$ is true or false). In none of the evaluated scenarios is RH "non-trivially true". So even in the unlikely event that the RH is "non-trivially true" in some other logic, the RH cannot be a "logical truth", because it is not "necessarily true".
\footnote{See Wikipedia \cite{LogicalTruth}: "Logical truths (including tautologies) are truths which are considered to be necessarily true. This is to say that they are considered to be such that they could not be untrue and no situation could arise which would cause us to reject a logical truth. It must be true in every sense of intuition, practices, and bodies of beliefs. However, it is not universally agreed that there are any statements which are \textit{necessarily} true."}
\footnote{See also Gómez-Torrente  \cite{Gomez-Torrente}:"It is typical to hold that, in some sense or senses of 'could', a logical truth could not be false or, alternatively, that in some sense or senses of 'must', a logical truth must be true. But there is little if any agreement about how the relevant modality should be understood."}

Moreover, classifying RH as a problem in logic does not definitively place it within Oort et al.'s \cite{Oort} "landscape of mathematics". The relationship between logic and mathematics is a matter of dispute between the classical school and the intuitionistic school.
\footnote{See Haack \cite{Haack}, pp.216-217: "[T]he Intuitionists think of logic as secondary to mathematics, as a collection of principles which are discovered, \textit{a posteriori}, to govern mathematical reasoning. This obviously challenges the 'classical' conception of logic as the study of principles applicable to all reasoning regardless of subject-matter, as the most fundamental and general of theories, to which even mathematics is secondary."}
\footnote{See Vafeiadou et al. \cite{Vafeiadou}, p.2, citing Brouwer \cite{Brouwer}, p.61: "In direct opposition to Russell and Whitehead's logicism, Brouwer asserted in 1907 that mathematics cannot be considered a part of logic. 'Strictly speaking the construction of intuitive mathematics in itself is an \textit{action} and not a \textit{science}; it only becomes a science, i.e. a totality of causal sequences, repeatable in time, in a mathematics of the second order [metamathematics], which consists of the \textit{mathematical consideration of mathematics} or \textit{of the language of mathematics} ... But there, as in the case of theoretical logic, we are concerned with an \textit{application of mathematics}, that is, with an \textit{experimental science}."}
This debate is discussed in greater detail in Section \ref{math_logic} of this paper.  

\subsection{Other Important Results Discussed in This Paper}

Also discussed in this paper is the invalidity in logics with LNC of the derivation of Riemann's $\zeta(s)$. Riemann used Cauchy's integral theorem to find the limit of the Hankel contour as the Hankel contour approaches the branch cut of $f(s)=\log(-s)$ on half-axis $s\ge0$. But $\log(-s)$ has no defined value (and so is non-holomorphic) on the branch cut. The Hankel contour is either an open path, or closed at $s = +\infty$ and encloses non-holomorphic points. Either way, prerequisites of Cauchy's integral theorem are violated. 

Also, the analytic continuation of $\zeta(s)$ violates the LNC and triggers ECQ, so all conjectures in classical or intuitionistic logic that falsely assume that AC of $\zeta(s)$ is true are rendered unsound. All generalizations of Riemann's $\zeta(s)$, such as Dirichlet $L$-functions, are unsound, because they falsely assume that Riemann's analytic continuation of $\zeta(s)$ is valid. So the BSD conjecture is unsound, because it assumes that Dirichlet $L$-functions are valid. 
Also, it is thereby proven that $\zeta(1)\ne0$, because $\zeta(s)$ is exclusively defined by the Dirichlet series. At $s=1$, $\zeta(s)$ is the famous "harmonic series"  which is \textit{proven} to be divergent. (Coincidentally, also Riemann's $\zeta(s)$ is divergent at $s=1$). So $\zeta(1)\ne0$.

This result also resolves other conjectures that are directly or indirectly equivalent to the BSD Conjecture (e.g. Tate's and Hodge's, as discussed in works by B. Totaro \cite{Totaro} and \cite{Totaro2}, and J.S. Milne \cite{Milne4}). In addition, it also invalidates the Modularity Theorem, and thus also Fermat's last theorem.

Several physics theories (including Yang-Mills theory) falsely assume that the analytic continuation of $\zeta(s)$ (a.k.a. "Riemann zeta function regularization") is true. This false assumption renders these theories unsound. We also apply Venn's "Modern" Square of Opposition to prove that $P\ne NP$.

Also included is discussion of Aristotelian logic, Classical logic, and Non-Classical logics, and a detailed discussion of the invalid derivation of Riemann's $\zeta(s)$. 

\section{Convergence and Divergence are the Two Values of a Bivalent System}

When discussing the RH, and the LNC, it is important to note that the terms "convergent" and "divergent" are defined as the two values of a bivalent logic. Therefore, a series cannot be both simultaneously: 

\begin{quotation}
[A] series 
\begin{equation} \label{eq:3.5}
\sum_{n=0}^{\infty} a_n = a_0+a_1+a_2+\ldots
\end{equation}
is said to be \textit{convergent}, to the sum $s$, if the 'partial sum' 
\begin{equation} \label{eq:3.6}
s_n = a_0+a_1+\ldots+a_n
\end{equation}
tends to a finite limit $s$ when $n \to \infty$; and a series which is not convergent is said to be \textit{divergent}. 
\footnote{See Hardy \cite{Hardy}, p.1.}
\end{quotation}

Therefore, by definition, an infinite series is \textit{either} convergent \textit{or} divergent. The series either converges to a value, or it does not. For example, an oscillating series such as $1-1+1-1+\ldots$ does not converge to any value, and therefore by definition is divergent.
\footnote{See Hardy \cite{Hardy}, p.1.}  

Moreover, if proposition ($P$) states that a given series is convergent, then the negation of that proposition, ($\neg P$), states that the given series is divergent. The converse is also true: if the proposition ($P$) states that the series is divergent, then the negation of that proposition, ($\neg P$), states that the series is convergent. 

According to the LNC, ($P$) and ($\neg P$) cannot both be true simultaneously. So the infinite series $\zeta(s)$ cannot be both convergent and divergent at \textit{any} value of $s$.

\section{The Dirichlet Series $\zeta(s)$ is Proven Divergent Throughout Half-Plane Re(s)<=1}

Riemann \cite{riemann1859number} concedes the divergence of Dirichlet series $\zeta(s)$ in the half-plane $\text{Re}(s)\le1$ as a \textit{given fact}.
\footnote{See Riemann \cite{riemann1859number}, p.1 (emphasis added): "The function of the complex variable s which is represented by these two expressions [the Euler product and the Dirichlet series], wherever they converge, I denote by $\zeta(s)$. Both expressions converge only when the Real part of $s$ is greater than 1".} 
Euler was the first to prove the first to prove that the Dirichlet series $\zeta(s)$ equals the Euler product.  
\footnote{See Wikipedia \cite{Euler} (citing  Derbyshire \cite{Derbyshire}, ch.7): "Leonhard Euler proved the Euler product formula for the Riemann zeta function in his thesis \textit{Variae observationes circa series infinitas} (Various Observations about Infinite Series), published by St Petersburg Academy in 1737."}
\footnote{See the Euler Archive \cite{Euler2} for the original publication (in Latin) of Euler's \textit{Variae observationes circa series infinitas}, and also English and German translations.}
Euler also proved that both the Euler product of $\zeta(s)$ and the Dirichlet series of $\zeta(s)$ are divergent along the half-line $s\le1$, for $s \in \mathbb{R}$.
\footnote{See also Calinger \cite{Calinger}, ch.4, p.136: "His 'Variae observationes' also introduces his famous product decomposition formula also introduces his famous product decomposition formula $p$ for the set of primes,
\begin{equation}
\prod_{p \in P} (1-p)^{-1} = \sum_{n=1}^{\infty} n^{-s}
\end{equation}
Multiplying the right side of the equation yields $\sum_{n=1}^{\infty} n^{-s} = \zeta(s)$. When $s=1$, $\zeta(1)$ is the harmonic series, which diverges to $\infty$. By applying the divergence of the harmonic series to the occurrence of primes, Euler proved indirectly their infinitude, a fact known since antiquity. The corresponding product must have infinitely many factors."}
This is easily confirmed by use of the "Integral test for convergence",
\footnote{See, e.g., Guichard et al.'s \cite{Guichard1}, discussion of the \textit{Integral test for convergence}, at Theorems 13.3.3 and 13.3.4 and their proofs.}
which is simplified in the "p-series" test for convergence.
\footnote{See, e.g. Department of Mathematics Website, Oregon State University \cite{p-series}, and Birdsong \cite{Birdsong}.}

Moreover, the Dirichlet series $\zeta(s)$ is \textbf{\textit{proven}} to be divergent at all values of $s$ in half-plane $\text{Re}(s)\le1$. \footnote{See Hardy et al. \cite{Hardy2}, p.4, Theorem 3: "The series may be convergent for all values of $s$, or for none, or for some only. In the last case there is a number $\sigma_s$ such that the series is convergent for $\sigma>\sigma_s$, and divergent or oscillatory for $\sigma<\sigma_s$. In other words \textit{the region of convergence is a half-plane}." (Citing Jensen \cite{Jensen} for the proof).} 
There exist proofs that the Dirichlet series $\zeta(s)$ is divergent in a portion of that half-plane, and there also exist proofs that it is divergent throughout the entire half-plane: 

\subsection{Dirichlet Series $\zeta(s)$ is Divergent at Real Half-Axis $\{\text{Re}(s)\le1, \text{Im}(s)=0\}$}

Dirichlet series $\zeta(s)$ is divergent for all values of $s$ on the Real half-axis $\{\text{Re}(s)\le1, \text{Im}(s)=0\}$. See the "Integral Test for convergence" (a.k.a. the Maclaurin–Cauchy test for convergence). This is commonly taught in introductory calculus textbooks. For example:
\footnote{See, e.g., Guichard  et al. \cite{Guichard1}, the \textit{Integral test for convergence}, discussed at Theorems 13.3.3 and 13.3.4 and their proofs.}
\footnote{See also the \textit{P-series test for convergence} at the Oregon State Univ. Dept. of Mathematics website \cite{p-series}, and at Birdsong \cite{Birdsong}. The \textit{P-series test for convergence} is the same as Guichard's \cite{Guichard1} Thm. 13.3.4.}
\begin{quotation}
\textit{Theorem 13.3.3}: Suppose that $f(x)>0$ and is decreasing on the infinite interval $[k,\infty)$ (for some $k\ge1$) and that $a_n=f(n)$. Then the series $\sum_{n= 1}^{\infty} a^n$ converges if and only if the improper integral $\int_{1}^{\infty} f(x)\,dx$ converges.

[A] $p$-series is any series of the form $\sum 1/n^p$. If $p\le 0$, [and $p \in \mathbb{R}$, then] $\lim_{n\to \infty} 1/n^p \ne 0$, so the series diverges. For positive values of $p$ we can determine precisely which series converge.

\textit{Theorem 13.3.4}: A $p$-series with $p>0$ converges if and only if $p>1$.

\textit{Proof.} We use the integral test; we have already done $p=1$, so assume that $p\ne1$.
\[\int_1^{\infty} \frac{1}{x^p}\,dx=\lim_{D\to\infty} \left.\frac{x^{1−p}}{1−p}\right|_{1}^D=\lim_{D\to\infty} \frac{D^{1-p}}{1-p} - \frac{1}{1-p}\]
If $p>1$ then $1−p<0$ and $\lim_{D \to \infty} D^{1−p}=0$, so the integral converges. If $0<p<1$, then $1−p>0$ and $\lim_{D\to \infty} D^{1−p} =\infty$, so the integral diverges.
\end{quotation}

\subsection{Dirichlet Series $\zeta(s)$ is Divergent at "Line of Convergence" $\{\text{Re}(s)=1\}$}

Dirichlet series $\zeta(s)$ is divergent for all values of $s$ on the misleadingly-named "line of convergence" $\text{Re}(s)=1$, which is parallel to the Imaginary axis $\text{Re}(s)=0$, and which is the border between the Dirichlet $\zeta(s)$ half-plane of convergence and its half-plane of divergence. At the point $s=1$, where $\zeta(s)$ is the famous "harmonic series", the function $\zeta(s)$ is divergent. At all other values of $s$ on the "line of convergence", $s = 1 + ti$, the function $\zeta(s)$ is "oscillating", which by definition is divergent. 
\footnote{See Hardy et al. \cite{Hardy2}, p.5, Example (iii), citing Bromwich \cite{Bromwich}: "The series $\sum n^{-s}$ has $\sigma=1$ as its line of convergence. It is not convergent at any point of the line of convergence, diverging to $+\infty$ for $s=1$, and oscillating finitely at all other points of the line." (Hardy \cite{Hardy}, p.1. defines "oscillating" as "divergent".  But note that Dirichlet series $\zeta(s)$ is "Cesàro summable" at all points on the "line of convergence" (except at $s=1$, where $\zeta(s)$ is the harmonic series, which is divergent).}  

\subsection{Dirichlet Series $\zeta(s)$ is Divergent Throughout Half-Plane $\{\text{Re}(s)\le 1\}$}

Dirichlet series $\zeta(s)$ is divergent for all values of $s$ in the half-plane $\{\text{Re}(s)\le1, \text{Im}(s)\}$. Hildebrand \cite{Hildebrand} states: 
\footnote{See Hildebrand \cite{Hildebrand}, pp.117-118, Theorem 4.6 (Convergence of Dirichlet series). See also Clark \cite{Clark}, p.11, Theorem 11, regarding the half-plane of divergence.}
\footnote{See Hardy et al. \cite{Hardy2}, p.3, fn. $\ddagger$, citing Jensen \cite{Jensen} and Cahen \cite{Cahen} for proofs of Hardy et al. \cite{Hardy2}, Theorems 1 and 2, respectively. Theorem 1 is: "\textit{If the [Dirichlet] series is convergent for $s=\sigma+ti$, then it is convergent for any value of $s$ whose Real part is greater than $\sigma$.}" (In other words, for $\text{Re}(s)>\sigma$.)}
\footnote{See also Conrad \cite{Conrad}, p.1, Example 1: "If $a_{n} = 1$ for all $n$ then $f(s) = \zeta(s)$, which converges for $\sigma > 1$. It does not converge at $s = 1$". See also Conrad \cite{Conrad}, pp.2-3, Theorems 8 and 9, and: "The contribution of Jensen \cite{Jensen} to Theorem 9 was a proof that convergence at $s_0$ implies convergence on the half-plane to the right of $s_0$." }
\begin{quotation}
For every Dirichlet series there exists a number $\sigma_c \in \mathbb{R} \cup \{\pm \infty\}$, called the abscissa of convergence, such that the series converges in the half-plane $\sigma > \sigma_c$ (the “half-plane of convergence”), and diverges in the half-plane $\sigma < \sigma_c$.
\end{quotation}
Hildebrand \cite{ Hildebrand} also states: 
\footnote{See Hildebrand \cite{Hildebrand}, p.126, Remark regarding Theorem 4.11.}
\footnote{See also Overholt \cite{Overholt}, pp. 65: "Hence every Dirichlet series has an \textit{abscissa of convergence} $\sigma_c$ such that it converges to the right of the line $\sigma=\sigma_c$ and diverges to the left of this line, which is called the \textit{line of convergence} for the series." (Overholt's \cite{Overholt} analytic continuation in pp.157-158, 162, including Proposition 5.1, violates the LNC.)} 
\footnote{See e.g. the discussion regarding the "abscissa of convergence" at Wikipedia \cite{Dirichlet_series} citing Hardy et al. \cite{Hardy2}: "In general the abscissa of convergence of a Dirichlet series is the intercept on the Real axis of the vertical line in the complex plane such that there is convergence to the right of it, and divergence to the left" and "Hence, for every $s$ such that $\sum _{n=1}^{\infty }a_{n}n^{-s}$ diverges, we have $ \sigma \geq \text{Re}(s)$, and this finishes the proof."}
\footnote{The author's proof (not peer reviewed) substitutes Euler's formula into the Dirichlet series, and then performs integration by parts (based on the assumption that $\int$ for $n \in \mathbb{R}$ is an acceptable approximation of $\Sigma$ for $n \in \mathbb{Z}$). See Sharon \cite{Sharon}, Version 4, Appendices A-F, pp.19-34.}
\footnote{Note: The partial sums of $f(t) = \sum_{n=1}^{\infty} \sin (t \cdot \ln(n))$ are not bounded, nor are they for $f(t) = \sum_{n=1}^{\infty} \cos (t \cdot \ln(n))$. Therefore, Sharon \cite{Sharon}, Version 4, Appendices G and H are wrong. However, this means that Dirichlet series test for convergence and Abel's lemma contradict the claim of convergence for Riemann's $\zeta(s)$ throughout the "critical strip" (except for the Imaginary axis $\text{Re}(s)=0$, where $f(\sigma, n)= n^{\sigma}$ is not monotonically decreasing as $n\to \infty$ for constant $\sigma=0$). (This is more evidence that Riemann's $\zeta(s)$ violates the LNC.)}
\begin{quotation}
For example, the Dirichlet series representation (4.11) of the zeta function diverges at every point in the half-plane $\sigma < \sigma_c = 1$ (and even at every point on the line $\sigma = 1$, as one can show by Euler’s summation).
\end{quotation}

Moreover, another method to prove that Dirichlet series $\zeta(s)$ is divergent throughout the half-plane $\text{Re}(s)\le1$ is to rewrite the Dirichlet series into its trigonometric form, 
\footnote{Because $s=\sigma + it$, the Dirichlet series $\zeta(s)$ can be written as:
\begin{equation}
\zeta(s) = \sum n^{-s} = \sum n^{-\sigma-it} = \sum n^{-\sigma}n^{-it}   
\end{equation}
and
\begin{equation}
\begin{split}
n^{-it} & = \exp(-it \cdot \ln(n)) \\
& = \cos(-it \cdot \ln(n)) + i\cdot \sin(-it \cdot \ln(n))
\end{split}
\end{equation}
and therefore $\zeta(s)$ can be written as:
\begin{equation}
\text{Re}\Big[\zeta(s)\Big] = \sum \Big[ n^{-\sigma} \cdot \cos(-it \cdot \ln(n)) \Big]    
\end{equation}
\begin{equation}
\text{Im}\Big[\zeta(s)\Big] = i \cdot \sum \Big[ n^{-\sigma} \cdot \sin(-it \cdot \ln(n)) \Big]       
\end{equation}}
and then performing integration by parts. 
\footnote{See, e.g., Sharon \cite{Sharon}, Appendix E, pp.26-32. (Note: not peer reviewed).}

\section{Riemann's $\zeta(s)$ Claims to be Convergent Throughout Half-Plane Re(s)<=1 (Except at s=1)}

Directly contradicting these proofs of the divergence of Dirichlet series $\zeta(s)$ throughout half-plane $\text{Re}(s)\le1$ is Riemann's \cite{riemann1859number} famous claim to have derived an alternative "expression" of $\zeta(s)$ that is \textit{convergent} at \textit{all} values of $s$, 
\footnote{See Riemann \cite{riemann1859number}, p.1: "The function of the complex variable $s$ which is represented by these two expressions [the Euler product and the Dirichlet series], wherever they converge, I denote by $\zeta(s)$. Both expressions converge only when the Real part of $s$ is greater than $1$; at the same time an expression for the function can easily be found which always remains valid." This is Riemann's $\zeta(s)$. }
except at the point $s=1$ (where, coincidentally, the Dirichlet series $\zeta(s)$ is the divergent  "harmonic series"). 
\footnote{See also Edwards \cite{Edwards}, pp.10-11: "Thus, formula
\begin{equation}
\zeta(s) = \frac{\Pi(-s)}{2\pi i} \int_{+\infty}^{+\infty} \frac{(-x)^s}{e^x - 1} \cdot \frac{dx}{x}
\end{equation}
defines a function $\zeta(s)$ which is analytic at all points of the complex s-plane except for a simple pole at $s=1$. This function coincides with $\sum n^{-s}$ for real values of $s>1$ and in fact, by analytic continuation, throughout the half-plane $\text{Re}((s)>1$. The function $\zeta(s)$ is known as the Riemann zeta function."}
\footnote{In Riemann's definition of $\zeta(s)$, the Euler product is superfluous.}
Regarding this analytic continuation of $\zeta(s)$ to half-plane $\text{Re}(s) \le 1$, Hildebrand \cite{ Hildebrand} states: 
\begin{quotation}
Strictly speaking, we should use a different symbol, say $\overline{\zeta}(s)$, for the analytic continuation ... However, to avoid awkward notations, it has become standard practice to denote the analytic continuation of a Dirichlet series by the same symbol as the series itself, and we will usually follow this practice.
\footnote{See Hildebrand \cite{ Hildebrand}, p.126, Remark to Theorem 4.11.}
\end{quotation}

However, if the analytic continuation $\overline{\zeta}(s)$ were indeed true at \textit{any} value of $s$ in the half-plane $\text{Re}(s)\le1$, then it would contradict the \textit{proven} divergence of Dirichlet series $\zeta(s)$ at that value of $s$. The "standard practice" of using the same symbol $\zeta(s)$ for the two contradictory "expressions" of the Zeta function (in the half-plane $\text{Re}(s)\le1$) is merely an explicit confirmation of the violation of of the LOI and the LNC. In logics that have both LNC and ECQ, the assumption that $\overline{\zeta}(s)$ is true violates the LNC, and triggers ECQ.

Included in Riemann's claim, that his "expression" of $\zeta(s)$ is convergent throughout half-plane $\text{Re}(s)\le1$, is a claim that Riemann's "expression" of $\zeta(s)$ is \textit{convergent} throughout the Real half-axis $\{\text{Re}(s)<1, \text{Im}(s)=0\}$.
\footnote{Riemann's functional equation of $\zeta(s)$ even claims that $\zeta(s)$ has "trivial zeros" on this Real half-axis.}
This directly contradicts the results of the "Integral test for convergence" (a.k.a. the Maclaurin-Cauchy test for convergence) for all values of $s$ on this Real half-axis.

Also included in Riemann's claim, that his "expression" of $\zeta(s)$ is convergent in half-plane $\text{Re}(s)\le1$, is a claim that his "expression" of $\zeta(s)$ is \textit{convergent} throughout the Dirichlet $\zeta(s)$ "line of convergence" at $\text{Re}(s)=1$ (except at the point $s=1$). This directly contradicts Hardy et al.'s \cite{Hardy2} theorem (citing Bromwich \cite{Bromwich}) that:
\begin{quotation}
The series $\sum n^{-s}$ has $\text{Re}(s)=1$ as its line of convergence. It is not convergent at any point of the line of convergence, diverging to $+\infty$ for $s=1$, and oscillating finitely at all other points of the line.
\footnote{See Hardy et al. \cite{Hardy2}, p.5, Example (iii). According to Hardy, the Dirichlet series is not convergent at any point on the "line of convergence", instead diverging to +$\infty$ at $s=1$, and oscillating finitely at all other points of the line [citing Bromwich \cite{Bromwich}].} 
\footnote{See also Sharon \cite{Sharon}, Version 4, Equation E.42 in Appendix E, p.31, which confirms Hardy's comment that Dirichlet series $\zeta(s)$ oscillates at all points on the line of convergence $\text{Re}(s)=1$, except at $(\text{Re}(s)=1, \text{Im}(s)=0)$ where it diverges to infinity. (However, the discussion in Sharon \cite{Sharon}, Version 4, p.32 overlooks this fact.)}
\end{quotation}
As stated in the preceding section, at all values of $s$ on the "line of convergence" (except at $s=1$), the Dirichlet series $\zeta(s)$ is a finitely oscillating series, that does not converge to any value, and therefore by definition is divergent.
\footnote{See Hardy \cite{Hardy}, p.1.}  
\footnote{Riemann concedes that his "expression" of $\zeta(s)$ is divergent at $s=1$, where Dirichlet series $\zeta(s)$ is the "harmonic series".} 

\section{A Third Version of $\zeta(s)$}

Ash \cite{Ash} discloses a third version of $\zeta(s)$, that contradicts both Dirichlet’s version of $\zeta(s)$ and Riemann’s version of $\zeta(s)$:

\begin{quotation}
We start with a trick. Multiply the sum for $\zeta(s)$ by $1/2^{s}$, and we get:
\begin{equation}
\frac{1}{2^{s}}\cdot \zeta(s) = \frac{1}{2^{s}} + \frac{1}{4^{s}} + \frac{1}{6^{s}} + \frac{1}{8^{s}} + \cdots.
\end{equation}
Let's line this up - the trick is to do it twice - underneath the sum for $\zeta(s)$, and subtract:

\begin{table}[h] \label{tab:trick_zeta}
\begin{tabular}{rlllllll}
$\zeta(s)$ = & 1 & + $1/2^{s}$ & + $1/3^{s}$ & + $1/4^{s}$ & + $1/5^{s}$ & + $1/6^{s} + \cdots$ \\
$1/2^{s} \cdot \zeta(s)$ = & & + $1/2^{s}$ & & + $1/4^{s}$ & & + $1/6^{s} + \cdots$ \\
$1/2^{s}\cdot \zeta(s)$ = & & + $1/2^{s}$ & & + $1/4^{s}$ & & + $1/6^{s} + \cdots$ \\  
\hline
$(1 - 1/2^{s} - 1/2^{s})\cdot \zeta(s)$ = & 1 & - $1/2^{s}$ & + $1/3^{s}$ & - $1/4^{s}$ & + $1/5^{s}$ & - $1/6^{s} + \cdots$   \\  
\end{tabular}
\end{table}
The result is:
\begin{equation} \label{cond_RH}
\Big(1-\frac{1}{2^{s-1}}\Big)\cdot \zeta(s) = 1 - \frac{1}{2^{s}} + \frac{1}{3^{s}} - \frac{1}{4^{s}} + \frac{1}{5^{s}} - \frac{1}{6^{s}} + \frac{1}{7^{s}} - \frac{1}{8^{s}} + \frac{1}{9^{s}} - \frac{1}{10^{s}} + \cdots
\end{equation}
In equation Eq. \ref{cond_RH}, the right-hand side is a Dirichlet series in which the coefficients are $1, -1, 1, -1, 1, -1, \cdots$. Notice that $|a_1 + a_2 + \cdots + a_n| < 2$ for any value of $n$. Theorem 11.7 now tells us that the right-hand side of equation Eq. \ref{cond_RH} can be summed provided that $\sigma<0$. The formula
\begin{equation} \label{cond2_RH}
\zeta(s) = \Big(1 - \frac{1}{2^{s-1}}\Big)^{-1} \cdot  \Big( 1 - \frac{1}{2^{s}} + \frac{1}{3^{s}} - \frac{1}{4^{s}} + \frac{1}{5^{s}} - \frac{1}{6^{s}} + \frac{1}{7^{s}} - \frac{1}{8^{s}} + \frac{1}{9^{s}} - \frac{1}{10^{s}} + \cdots \Big)
\end{equation}
therefore can be evaluated provided that $\sigma>0$, with the sole exception of the value at $s=1$
\footnote{See Ash \cite{Ash}, pp.170-171.}
\end{quotation}
Moreover, Ash's \cite{Ash} cited "Theorem 11.7" states the following:
\footnote{See Ash \cite{Ash}, p.169.}
\begin{quotation}
THEOREM 11.7 Suppose that there is some constant $K$ so that $|a_1 + \cdots + a_n|<K$ for all $n$. Then the Dirichlet series $\sum a_n n^{-s}$ converges if $\sigma>0$.
\end{quotation}

This third version of $\zeta(s)$ is convergent in half-plane $\text{Re}(s)>0$, and divergent in half-plane $\text{Re}(s)\le0$. Therefore, this third version contradicts the Dirichlet series version of $\zeta(s)$ throughout the “critical strip” $0<Re(s)\le1$, by being convergent there, and also contradicts Riemann’s version of $\zeta(s)$ throughout half-plane $\text{Re}(s)\le0$, by being divergent there. 

Clearly, in a logic with LNC, only one of these three versions of $\zeta(s)$ can be true. It is impossible for all three versions of $\zeta(s)$ to be true, or even for two of the three to be true. The "Riemann series theorem" provides an explanation as to why this third version of $\zeta(s)$ is invalid as an indicator of convergence/divergence. According to the Riemann series theorem:  
\begin{quotation}
By a suitable rearrangement of terms, a conditionally convergent series may be made to converge to any desired value, or to diverge.
\footnote{See Weisstein \cite{Weisstein}, citing: Bromwich et al. \cite{Bromwich2}, p.74; Gardner \cite{Gardner2}, p.171; and Havil \cite{Havil}, p.102. }
\end{quotation}
So any conditionally convergent series can be rearranged to be convergent to any finite number, and also can be rearranged to be divergent. This is not the case for an absolutely convergent series, which is one of convergent or divergent, regardless of how the terms are rearranged.

The third version of $\zeta(s)$ is created by transforming the Dirichlet series $\zeta(s)$, which is an \textit{unconditionally} convergent series throughout half-plane $\text{Re}(s)>1$, to a series which is \textit{conditionally} convergent throughout half-plane $\text{Re}(s)>0$.
Therefore, it is due to the \textit{conditional} convergence of the third version of $\zeta(s)$ that the third version can be manipulated to have a different zone of convergence than the original Dirichlet series version. According to the Riemann series theorem, the region of convergence of the Dirichlet series version remains constant when the series is rearranged, but the region of convergence of the third version changes when rearranged (thus violating the LNC).

Ash \cite{Ash} concludes the description of the third version of $\zeta(s)$ by stating that it can be analytically continued even further:
\begin{quotation}
However, we need to do still more, and find a way to evaluate $\zeta(s)$ for \textit{all} values of $s$ except $s=1$. That requires a discussion of \textit{functional equations}. 
\footnote{See Ash \cite{Ash}, p.171.}
\end{quotation}
This proposition violates the LNC, because contrary to the LNC, it holds that the third version of $\zeta(s)$ can be convergent for every value of $s$ in half-plane $\text{Re}(s)\le 0$ (except $s=1$), even though the Dirichlet series of $\zeta(s)$ is proven to be divergent throughout said half-plane.

\section{One-to-Two "Functions" Violate the Definition of a Function}

As discussed above, the Dirichlet series "expression" of $\zeta(s)$ is \textit{proven} to be divergent throughout this half-plane $\text{Re}(s)\le 1$. Riemann claims to have derived an additional "expression" of $\zeta(s)$, that is convergent throughout half-plane $\text{Re}(s)\le1$ (except at $s=1$). However, this violates the definition of a "function" in set theory:
\begin{quotation}
A function is a relation that uniquely associates members of one set with members of another set. More formally, a function from $A$ to $B$ is an object $f$ such that every $a \in A$ is uniquely associated with an object $f(a) \in B$. A function is therefore a many-to-one (or sometimes one-to-one) relation.
\footnote{See Stover et al. \cite{Stover}.} 
\end{quotation}
Therefore, if both "expressions" of $\zeta(s)$ were true in half-plane $\text{Re}(s)\le 1$, they would violate this definition of a mathematical function, due to the \textit{one-to-two} mapping from $s$ to the two different values of $\zeta(s)$, throughout the half-plane $\text{Re}(s)\le1$ (except at $s=1$). Each $s$ value would map to two $\zeta(s)$ values: convergent and divergent.

\section{The Law of Non-Contradiction (LNC)}

The central thread of this paper is that Riemann-style "analytic continuation" of the Zeta function violates the Law of Non-Contradiction (LNC). LNC holds that a proposition $p$ cannot be simultaneously both true and false. In logic notation: $\vdash \neg (p \land \neg p)$. The LNC is one of Aristotle's three "Laws of Thought". It is included in classical logic and intuitionistic logic, either as an axiom or as a theorem. In contrast, 3VL rejects the LNC, either implicitly or explicitly. \footnote{See Haack \cite{Haack2}, p.5: "In Łukasiewicz's 3-valued logic (motivated by the idea, already suggested by Aristotle in \textit{De interpretatione} \S 9, in \textit{Organon}), that future contingent sentences are neither true nor false but 'indeterminate') both the Law of the Excluded Middle ('LEM;' 'p or not p') and the Law of Non-Contradiction ('LNC;' 'not both p and not-p') fail."}
\footnote{But see Decker \cite{Decker}, p.69, \S 3.3 Precursors of Paraconsistent Logic: "the introduction of more than two truth -values opens up the possibility that some formulas which are classically interpreted as contradictions no longer evaluate to false."
}
\footnote{By definition, there is no third category other than convergent and divergent. For example, Cesàro summable sequences are classified as "divergent" (e.g. range-bound sine and cosine functions).}

In regards to the RH, in half-plane $\text{Re}(s)\le1$ there are two conflicting definitions of $\zeta(s)$: the Dirichlet series $\zeta(s)$, and Riemann's $\zeta(s)$. \textbf{The LNC holds that a convergent "expression" of $\zeta(s)$ and a divergent "expression" of $\zeta(s)$ cannot both be true, at any value of $s$.}

In Riemann's famous paper, he refers to the Dirichlet series $\zeta(s)$ as an "expression" of $\zeta(s)$ that "converge[s] only when the Real part of $s$ is greater than $1$", and then claims that there is another "expression" that "always remains valid".
\footnote{See Riemann \cite{riemann1859number}, p.1: "The function of the complex variable $s$ which is represented by these two expressions [the Euler product and the Dirichlet series], wherever they converge, I denote by $\zeta(s)$. Both expressions converge only when the Real part of $s$ is greater than $1$; at the same time an expression for the function can easily be found which always remains valid."}
This latter proposition is where Riemann violates the LNC in his paper. 

\section{\textit{Ex Contradictione Quodlibet} (ECQ)}

Classical and intuitionistic logics also contain the "Principle of Explosion" \textit{Ex Contradictione Quodlibet} (ECQ). 
\footnote{In MVL, paradoxes can be assigned a 3rd truth-value, thereby avoiding both LNC and ECQ.} 
ECQ states that a violation of LNC materially implies that every other statement is "trivially true".  
\footnote{See, e.g. Kleene \cite{Kleene2}, p.101. according to which ECQ ($A, \neg A \vdash B$) is valid in both classical and intuitionistic logics.}
In classical logic notation: $\vdash (p \land \neg p)\rightarrow q$. 

\textbf{In regards to the RH, ECQ holds that if both the Riemann $\zeta(s)$ and the Dirichlet series $\zeta(s)$ are both true in half-plane $\text{Re}(s)\le1$, or are falsely assumed to both be true, then the contradiction "trivially" implies the truth of any other statement.}

ECQ is a direct result of the definition of "material implication" in classical logic, which holds that a false statement implies any statement. (In logic notation: $\vdash (p \rightarrow q) = (\neg p \lor q)$). 
\footnote{See Wikipedia \cite{Material_implication}, citing Hurley \cite{Hurley}, pp.364–5; Copi et al. \cite{Copi}, p.371; and Math StackExchange \cite{StackExchange2}.} 
\footnote{See also Priest \cite{Priest8}, p.45: "Recall that a conditional is a sentence of the form 'if $a$ then $c$', which we are writing as $a \rightarrow c$", "If you know that $a \rightarrow c$, it would seem that you can infer that $\neg (a \& \neg c)$ (it is not the case that $a$ and not $c$)", and "Conversely, if you know that $\neg (a \& \neg c)$, it would seem that you can infer $a \rightarrow c$ from this."}
\footnote{But see Priest \cite{Priest8}, p.46, alternative name and notation: "$\neg (a \& \neg c)$ is often written as $a \supset c$, and called the \textit{material conditional}."}
The proof of "material implication" in classical logic relies upon the Law of the Excluded Middle (LEM), 
\footnote{See Wikipedia \cite{Material_implication}, citing  Math StackExchange \cite{StackExchange2}: "Suppose we are given that $P\to Q$. Then, since we have $ \neg P\lor P$ by the law of excluded middle, it follows that $ \neg P\lor Q$. Suppose, conversely, we are given $\neg P\lor Q$. Then if P is true that rules out the first disjunct, so we have Q. In short, $P\to Q$."}
and therefore is not valid in logics that do not have the LEM as an axiom, such as intuitionistic logic and 3VL. Heyting's intuitionism has an alternative version of ECQ, 
\footnote{See Decker \cite{Decker}, pp.67-68, \S 3.2.3 Intuitionism: "As opposed to logicists and formalists, Brouwer ... rejected the use of LEM and the law of double negation (LDN, formally: $p \Leftrightarrow \neg\neg p$, one half of which is axiom 10 above). It was Kolmogorov \cite{Kolmogorov2} (English version: \cite{Kolmogorov}) and Heyting \cite{Heyting2}
\cite{Heyting3} who proposed axiomatic systems for making intuitionism accessible to formal treatment (at the expense of some of Brouwer's basic philosophical beliefs)."}
\footnote{Also, at Decker \cite{Decker}, pp.67-68: "Replacing $\neg\neg p \to p$ (axiom 10) by $\neg p \to (p \to q)$ (ECQ) in fig.1 yields Heyting's system. In the resulting schema, neither LDN nor LEM are derivable any more. The non-validity of LDN and LEM also invalidates the law of material implication (LMI) and effectively coerces proofs to be constructive. ... However, LoC [Law of (Non-)Contradiction] continues to hold."}
but minimal logic (a version of intuitionistic logic) does not have ECQ as a theorem. 
\footnote{See Decker \cite{Decker}, p.68, \S 3.2.3 Intuitionism: "The axiom $\neg p \rightarrow (p \rightarrow q)$ (i.e., ECQ) in Heyting's system is abandoned in Johansson's \textit{minimal logic} \cite{Johansson}. The latter consists of MP and axioms 1-9 (fig. 1) and thus essentially is the same as Kolmogorov's system. In minimal logic, only each negated sentence can be deduced in the presence of contradiction, but not necessarily each sentence whatsoever. In particular, $\neg p \rightarrow (p \rightarrow \neg q)$ can be deduced from axioms 1 - 9, but $\neg p  \rightarrow (p \rightarrow q)$ (ECQ) cannot."}

Three-Valued Logic (3VL) has ECQ,
\footnote{See Decker \cite{Decker}, p.69, \S 3.3 Precursors of Paraconsistent Logic: "[T]he derivability of a contradiction still entails trivialization." (Citing Urquhart \cite{Urquhart}.)}
but "the introduction of more than two truth -values opens up the possibility that some formulas which are classically interpreted as contradictions no longer evaluate to false."
\footnote{See Decker \cite{Decker}, p.69, \S 3.3 Precursors of Paraconsistent Logic.}
Both Bochvar's 3VL
\footnote{See Urquhart \cite{Urquhart2}, pp.252-253, \S 1.6: "The work of the Russian logician Bochvar \cite{Bochvar} represents a new philosophical motivation for many-valued logic; its use as a means of avoiding the logical paradoxes. His system introduces the intermediate value \textit{I} in addition to the classical values \textit{T} and \textit{F}. His idea is to avoid logical paradoxes such as Russell's and Grelling's by declaring the crucial sentences involving them to be meaningless (having the value \textit{I})."}
and Priest's "Logic of Paradox" (\textit{LP}) (which uses the truth tables of Kleene’s "strong" 3-valued logic)
\footnote{See Priest \cite{Priest5} and \cite{Priest7}. See also Hazen et al. \cite{Hazen}, p.2: "Truth values of compound formulas are derived from those of their subformulas by the familiar “truth tables” of Kleene’s (strong) 3-valued logic [\cite{Kleene2}, §64], but whereas for Kleene (thinking of the “middle value” as truth-valuelessness) only the top value (True) is designated, for Priest the top two values are both designated."}
expressly assign the third truth-vale to paradoxes, so paradoxes do not evaluate to "false".
\footnote{See Urquhart \cite{Urquhart2}, pp.252: "[Bochvar's] idea is to avoid logical paradoxes such as Russell's and Grelling's by declaring the crucial sentences involving them to be meaningless (having the value \textit{I})."}
\footnote{See also Panti \cite{Panti1998}, p.48, \S 2.5.1 Bochvar's and Kleene's systems: "In addition to 0 and 1 for \textit{false} and \textit{true}, they have a third value 2. While for Łukasiewicz the third value stands for \textit{possible}, or \textit{not yet detennined}, from Bochvar's point of view it stands for \textit{paradoxical}, or \textit{meaningless}. Any compound proposition that includes a meaningless part is meaningless itself, and hence the [following] truth tables ...", and "Bochvar's systems was proposed in [Bochvar \cite{Bochvar}] as a way for avoiding the logical paradoxes, notably Russell's paradox. We refer to [Rescher, \cite{Rescher}, \S 2.4] and [Urquhart, \cite{Urquhart}, \S 1.6] for a deeper analysis and further references."}
\footnote{See also Smith \cite{Smith3}, pp.17-18: "Consideration of the paradoxes — set-theoretic (e.g. Russell’s) and/or semantic (e.g. the Liar, where it seems impossible to assign either truth value $1$ or $0$ to ‘This sentence is false’) — was a motivation for Bochvar, Moh Shaw-Kwei and others (see Rescher \cite{Rescher}, pp. 13, 29, and 207 for additional references) ... Kleene [\cite{Kleene2}, 335] also considers a different interpretation of his three values: “t, f, u must be susceptible of another meaning besides (i) ‘true’, ‘false’, ‘undefined’, namely (ii) ‘true’, ‘false’, ‘unknown (or value immaterial)’."}
\footnote{See also Visser \cite{Visser}, p.181: "This paper interweaves various themes. Two main themes are four-value logic and the Liar Paradox."}

\section{The LNC and Analytic Continuation of $\zeta(s)$ Cannot Both be True}

\subsection{If Both are True, in Logics with LNC and ECQ, this Triggers ECQ}

If we assume that analytic continuation of $\zeta(s)$ is true, then there exist two true, but contradictory, definitions of $\zeta(s)$ throughout the half-plane $\text{Re}(s)\le1$. One definition says that $\zeta(s)$ is divergent throughout that half-plane. The other definition says that $\zeta(s)$ is convergent there. This situation violates the Law of Non-Contradiction (LNC) at all values of $s$ in the half-plane $\text{Re}(s)\le1$ (except at $s=1$). 
Even if we limit ourselves to the so-called "line of convergence" at $\text{Re}(s)=1$, or to half of the Real-axis $\{\text{Re}(s)<1, \text{Im}(s)=0\}$, analytic continuation of $\zeta(s)$ still violates the Law of Non-Contradiction (LNC).

In any other logic that contains both the LNC and the "principle of explosion" (ECQ), such as classical logic, or intuitionistic logic, any "proof" with a contradiction is "trivially true" due to ECQ.
\footnote{See, e.g. Whitehead and Russell \cite{Whitehead2}, Th. *2.21 on p.99 is their version of ECQ: "$\vdash: \sim p . \supset . p \supset q $", which is described as: "I.e. a false proposition implies any proposition."}
So if we assume such a logic to be the foundation for our mathematics (and thereby assume LNC and ECQ to be true), and if we assume that Riemann's "expression" of $\zeta(s)$ is true (and thereby assume that $\zeta(s)$ is both divergent and convergent in half-plane $\text{Re}(s)\le1$), then any "proof" in such mathematics that relies upon $\zeta(s)$ in half-plane $\text{Re}(s)\le1$ is "trivially true" due to ECQ.

\subsection{If AC of Zeta is False, Then in Logics with LNC, Falsely Assuming it is True Renders a Proof Unsound}

So analytic continuation of $\zeta(s)$ is \textit{false} in half-plane $\text{Re}(s)\le1$, due to the \textit{proven} divergence of the Dirichlet series $\zeta(s)$ there. All other alleged "proofs" of analytic continuation of $\zeta(s)$ to that half-plane must also be \textit{false}, if the logical foundation of analytic number theory is a logic that includes LNC and ECQ (e.g. classical logic or intuitionistic logic).
\footnote{See e.g. Titchmarsh et al. \cite{Titchmarsh}, \S 2.1 to \S 2.10, pp.13-27, which lists seven additional alleged "proofs".}
Any "proof" that assumes that Riemann's $\zeta(s)$ is true in half-plane $\text{Re}(s)\le1$ is "trivially true" in classical and intuitionistic logics, due to ECQ.

Moreover, if analytic continuation of $\zeta(s)$ to that half-plane is false, due to LNC and ECQ, then $\zeta(s)$ is \textit{exclusively} defined by its Dirichlet series definition, which is divergent throughout half-plane $\text{Re}(s)\le 1$. So $\zeta(s)$ \textit{has no zeros}, because the Dirichlet series $\zeta(s)$ has no zeros. 
\footnote{As discussed in Hardy et al. \cite{Hardy2}, p.5, Example (iii), the Dirichlet series $\zeta(s)$ is divergent throughout half-plane $\text{Re}(s)\le 1$.}
Euler proved that $\zeta(s)$ has no zeros in the "other" half-plane, $\text{Re}(s)> 1$,
\footnote{In half-plane $\text{Re}(s)>1$, the Dirichlet series $\zeta(s)$ equals the Euler product of the primes. Each factor of the Euler product is a fraction having "1" as the numerator. So the Euler product cannot equal zero, because at least one numerator of "0" is necessary for the product to equal zero. Therefore, the Dirichlet series cannot equal zero either (in this half-plane). So neither the Dirichlet series nor the Euler product have any zeros in that half-plane.} 
and in half-plane $\text{Re}(s)\le 1$, the Dirichlet series $\zeta(s)$ is divergent (and non-zero) throughout.

Therefore, according to classical logic,
\footnote{I.e. The logic of \textit{Principia Mathematica}.} 
in half-plane $\text{Re}(s)\le1$ (except $s=1$), analytic continuation of $\zeta(s)$ violates all three of the "laws of thought" that are inherited from Aristotelian logic: the Law of Identity (LOI), the Law of the Excluded Middle (LEM), and the Law of Non-Contradiction (LNC). 
\footnote{In addition, analytic continuation of $\zeta(s)$ \textit{also} violates the definition of a function, according to set theory, because a function \textit{cannot} have a one-to-two mapping from domain to range.}
The most serious of these "violations" is the violation of the Law of Non-Contradiction (LNC) in that half-plane, because it triggers ECQ.

\section{Riemann's Analytic Continuation of $\zeta(s)$ is Invalid} \label{Riemann-Invalid_1}

Riemann's analytic continuation of Dirichlet series $\zeta(s)$ is violates the LNC. So in logics with the LNC, Riemann's $\zeta(s)$ is false where it contradicts Dirichlet series $\zeta(s)$. 
\footnote{See Edwards \cite{Edwards}, pp.10-11: "Thus, formula
\begin{equation}
\zeta(s) = \frac{\Pi(-s)}{2\pi i} \int_{+\infty}^{+\infty} \frac{(-x)^s}{e^x - 1} \cdot \frac{dx}{x}
\end{equation}
defines a function $\zeta(s)$ which is analytic at all points of the complex s-plane except for a simple pole at $s=1$. This function coincides with $\sum n^{-s}$ for real values of $s>1$ and in fact, by analytic continuation, throughout the half-plane $\text{Re}((s)>1$. The function $\zeta(s)$ is known as the Riemann zeta function."}

In deriving his expression for $\zeta(s)$, Riemann uses the Hankel contour,
\footnote{See Riemann's use of the Hankel contour in Equation \ref{Eq3} of this paper.}
which is taken directly from Hankel's derivation of the Gamma function $\Gamma(s)$. 
\footnote{See Whittaker et al. \cite{Whittaker}, pp.244-245 and 266.} 
Riemann then uses Cauchy's integral theorem to find the limit of the Hankel contour as the Hankel contour approaches the branch cut of $f(s)=\log(-s)$ for $s \in \mathbb{C}$. But by definition, $\log(-s)$ has no value on half-axis $s\ge0$ (and thus is also non-holomorphic on this half-axis).
\footnote{The geometric proof that $\log(-s)$ is non-holomorphic on half-axis $s\ge0$: In the Cartesian plane, the 1st derivative of $f(x)=\log(-x)$, for $x \in \mathbb{R}$ at a value of $x$, is represented by the slope of the line tangent to $f(x)$ at $x$. However, $f(x)$ has no values at $x\ge0$, so it has no 1st derivative values at $x\ge0$.}

The Hankel contour is either open, or closed. In both cases, the Hankel contour violates prerequisites of Cauchy's integral theorem. A first prerequisite is that all points inside a closed contour must be holomorphic. If the Hankel contour is closed (for example, at $s = +\infty$, which is assumed in the derivation of Riemann's $\zeta(s)$),
\footnote{See Whittaker et al. \cite{Whittaker}, p.245: "We shall write $\int_\infty^{(0+)}$ for $\int_C$, meaning thereby that the path of integration starts at 'infinity' on the Real axis, encircles the origin in the positive direction, and returns to the starting point." }
the contour encloses the non-holomorphic points of the branch cut, which violates this first prerequisite. A second prerequisite of the Cauchy integral theorem is that there be two different paths connecting two points. If the Hankel contour is open at $s = +\infty$, it violates this second prerequisite, which requires that the contour be closed. 

So the Hankel contour contradicts prerequisites of Cauchy's integral theorem. Thus, the derivation of Riemann's $\zeta(s)$ violates the LNC. (This is described in greater detail in Chapter \ref{Riemann-Invalid_2} of this paper).

Moreover, given that the Dirichlet series $\zeta(s)$ is \textit{proven} to be divergent throughout half-plane $\text{Re}(s)\le1$, classical and intuitionistic logics (both of which have both the LNC and ECQ) hold that \textbf{\textit{every}} so-called analytic continuation of $\zeta(s)$ into this half-plane \textbf{\textit{must be false}}, because it violates LNC and triggers ECQ. Accordingly, in logics with LNC and ECQ, $\zeta(s)$ is exclusively defined by the Dirichlet series $\zeta(s)$, which has no zeros and no poles.
\footnote{Throughout half-plane $s\le1$, Dirichlet series $\zeta(s)$ is divergent, and throughout half-plane $s>1$ it is equal to the Euler product, whose factors all have non-zero numerators. See e.g. Hildebrand \cite{ Hildebrand}, pp.147, Theorem 5.2(iv): "$\zeta(s)$ has no zeros in the half-plane $\sigma > 1$." }

\section{Weierstrass's Analytic Continuation is Valid, but Riemann's is Not} \label{Weierstrass}

\subsection{Weierstrass's Chain of Disks}

As Edwards \cite{Edwards} states (emphasis added): \begin{quotation}
It is interesting to note that Riemann does \textit{not} speak of the 'analytic continuation' of the function $\sum n^{-s}$ beyond the halfplane $\text{Re}(s) > 1$, but speaks rather of finding a formula for it which '\textit{remains valid for all $s$}.'
\footnote{See also Edwards \cite{Edwards}, p.9}
\end{quotation}
Furthermore, Edwards \cite{Edwards} compares Weierstrass's and Riemann's versions of analytic continuation, as follows:
\begin{quotation}
The view of analytic continuation in terms of chains of disks and power series convergent in each disk descends from Weierstrass and is quite antithetical to Riemann's basic philosophy that analytic functions should be dealt with \textit{globally}, not locally in terms of power series. \footnote{See Edwards \cite{Edwards}, p.9.}
\end{quotation}

For example, Weierstrass's "chains of disks" analytic continuation of $f(s)= 1/(1-s)$ does \textit{not} make the contradictory claim that $f(s)$ is both convergent and divergent at $s=1$. In fact, Weierstrass's analytic continuation method \textit{avoids} all directly contradictory propositions, by forbidding disks from encircling any pole. Each of Weierstrass's disks represents a proposition distinct from all of the other disks, and distinct from the poles. \footnote{See, e.g., Weyl \cite{Weyl}, pp.1-4, and Coleman \cite{Coleman}, pp.1-2.}

This methodology ensures that no two propositions (i.e. disks, poles) are directly contradictory. It also ensures that no disk violates Cauchy's integral theorem (i.e. each disk  exclusively has holomorphic points). Regarding this method, Weyl \cite{Weyl} states:
\begin{quotation}
In its convergence disc ..., such a function represents a regular analytic function in the sense of Cauchy, 
\footnote{See, e.g., Weyl \cite{Weyl}, p.1} [and]

It is not claimed that each of these continuations can be extended to an analytic chain reaching the end ($\lambda=1$); in general that is false. The fact that each analytic chain contains only a finite number of irregular elements makes it possible to \textit{avoid} these irregular elements.
\footnote{See, e.g., Weyl \cite{Weyl}, p.13}
\end{quotation}
In summary, Weierstrass's method avoids direct contradictions. Riemann's does not.

\subsection{Poincaré's "l’Analysis Situs" (1895) Uses Weierstrass's Chain of Disks}

According to Morgan \cite{Morgan}, page 10: 
\begin{quotation}
\textbf{2.2. l’Analysis Situs (1895).}
\footnote{Citing Poincaré \cite{Poincare}. See also the English translation at Poincaré \cite{Poincare2}. See also Wikipedia \cite{Analysis_Situs}.}
This is a long (121 pages), foundational paper. Poincaré begins by defending the study he is about to undertake by saying 
\begin{quotation}
“Geometry in $n$-dimensions has a real goal; no one doubts this today. Objects in hyperspace are susceptible to precise definition  like those in ordinary space, and even if we can’t represent them to ourselves we can conceive of them and study them.”
\end{quotation}
There then follows a discursive introduction to the study of the topology of manifolds. Many of the approaches and techniques that came to dominate 20th century topology are introduced in this paper. It truly is the beginning of Topology as an independent branch of mathematics.
\end{quotation}
Morgan \cite{Morgan}, further discloses at page 10 (emphasis added): 
\begin{quotation}
[Poincaré] also considers manifolds defined by locally closed, one-one immersions from open subsets of Euclidean $n − p$-space. He goes on to consider manifolds covered by overlapping subsets of either type \textit{(though he is considering the real analytic situation where the extensions are given by analytic continuation)}. Having defined manifolds, he considers orientability, orientations, and homology.
\end{quotation}
Poincaré \cite{Poincare2}, pp. 24-25 clarifies that the anaylytic continuation used is an analogue of Weierstrass's chain of disks, not of Riemann's version of anaylytic continuation:
\begin{quotation}
It can happen that the two manifolds have a common part $V''$ also of $m$
dimensions. In that case, in the interior of $V''$, the $y$ will be analytic functions of the $y'$ and conversely. We then say that the two manifolds $V$ and $V'$ are \textit{analytic continuations} of each other. In this way we can form a chain of manifolds
\begin{equation}
V_1, V_2, \ldots V_n    
\end{equation}
such that each is an analytic continuation of its predecessor, and there is a common part between any two consecutive manifolds of the chain. I shall call this a \textit{connected chain}.
\end{quotation}

\section{The "Calculated Zeros" of Riemann's $\zeta(s)$ are of Other Formulas (That Assume AC of $\zeta(s)$ is True)}

\subsection{The Euler-Maclaurin Formula}

The so-called "calculated zeros of Riemann's $\zeta(s)$" are actually zeros of approximations. For example, Odlyzko et al. \cite{Odlyzko} assumes the following before attempting to calculate the "zeros" of Riemann's $\zeta(s)$: 
\begin{quotation}
The Riemann zeta function is defined for $s = \sigma + it$ by
\begin{equation} \label{eq:Dirichlet2}
\zeta(s) = \sum_{n=1}^{\infty}n^{-s}
\end{equation}
for $\sigma > 1$, and by analytic continuation can be extended to an analytic function of $s$ for all $s \ne 1$ [citing Edwards \cite{Edwards}, Ivić \cite{Ivic}, and Titchmarsh \cite{Titchmarsh}]. 
\end{quotation}

However, as discussed above, assuming that the analytic continuation of $\zeta(s)$ is true generates a paradox. In classical and intuitionistic logics, this violation of the LNC triggers ECQ, and thus renders "trivially true" (and \textit{de facto} invalidates) everything that is built on the assumption (that uses the Euler-Maclaurin Formula).

Odlyzko et al. \cite{Odlyzko} then discloses the following in regards to calculating the "zeros" of Riemann's $\zeta(s)$, not by use of Riemann's $\zeta(s)$, but by use of the Euler-Maclaurin formula: 
\footnote{See Odlyzko et al. \cite{Odlyzko}, p.798} 
\begin{quotation} \label{eq:Dirichlet}
The [Equation] (\ref{eq:Dirichlet2}) suggests the idea of using the Euler-Maclaurin summation formula [citing Abramowitz et al.'s \cite{Abramowitz} Equation 23.1.30] to evaluate $\zeta(s)$, and one easily obtains, for any positive integers $m$ and $n$,
\begin{equation} \label{eq:zeta}
\zeta(s) = \sum_{j=1}^{n-1} j^{-s} + \frac{1}{2} n^{-s} + \frac{n^{1-s}}{s-1} + \sum_{k=1}^{m} T_{k,n}(s) + E_{m,n}(s)
\end{equation}
where
\begin{equation}
T_{k,n}(s) = \frac{B_{2k}}{(2k)!} n^{1-s-2k} \prod_{j=0}^{2k-2} (s+j)
\end{equation}
$B_{2} = 1/6, B_{4} = -1/30, \ldots,$ are the Bernoulli numbers, and
\begin{equation} \label{eq:error}
|E_{m,n}(s)| < \Big|\frac{s+2m+1}{\sigma+2m+1} T_{m+1,n}(s) \Big|
\end{equation}
The formula (\ref{eq:zeta}) with the estimate (\ref{eq:error}) can easily be shown to hold for any $\sigma > -(2m + 1)$. By taking $m$ and $n$ large enough (and using sufficient accuracy in basic arithmetic routines), any value of $\zeta(s)$ can be computed to any desired accuracy by this formula. All calculations of zeros of the zeta function that were published before 1930 relied on this method. Its advantages include the ease of estimating the error term. (This is the main reason this formula is still used for very accurate computations of $\zeta(s)$ for $s$ small, cf. [19].) 
\end{quotation}

However, the use of the Euler-Maclaurin summation formula fails at Odlyzko et al.'s \cite{Odlyzko} first sentence: "[The Dirichlet series definition of $\zeta(s)$] suggests the idea of using the Euler-Maclaurin summation formula ... to evaluate $\zeta(s)$". Apostol \cite{Apostol} indicates
\footnote{See Apostol \cite{Apostol}, p.409, "Introduction".}
why the Euler-Maclaurin summation formula cannot be used to calculate "zeros" of the Dirichlet series of $\zeta(s)$ in half-plane $\text{Re}(s)\le1$:
\begin{quotation}
The integral test for convergence of infinite series compares a finite sum $\sum_{k=1}^{n} f(k)$ and an integral $\int_{1}^{n} f(x)\, dx$ where $f$ is positive and strictly decreasing. The difference between a sum and an integral can be represented geometrically, as indicated in Figure 1. In 1736, Euler \cite{Euler3} used a diagram like this to obtain the simplest case of what came to be known as Euler's summation formula, a powerful tool for estimating sums by integrals, and also for evaluating integrals in terms of sums. Later Euler \cite{Euler4} derived a more general version by an analytic method that is very clearly described in [Hairer et al. \cite{Hairer}, pp. 159-161]. Colin Maclaurin \cite{Maclaurin} discovered the formula independently and used it in his \textit{Treatise of Fluxions}, published in 1742, and some authors refer to the result as the Euler-Maclaurin summation formula.
\end{quotation}

As Apostol \cite{Apostol} indicates, "[t]he integral test for convergence of infinite series compares a finite sum $\sum_{k=1}^{n} f(k)$ and an integral $\int_{1}^{n} f(x)\, dx$ where $f$ is positive and strictly decreasing", and "[t]he difference between a sum and an integral can be represented geometrically". As discussed in the present paper, the Dirichlet series of $\zeta(s)$ \textbf{\textit{fails}} the integral test for convergence of infinite series at all values of $s$ in half-plane $\text{Re}(s)\le1$. This is sufficient reason to disqualify the use of the  Euler-Maclaurin summation formula to calculate "zeros" of $\zeta(s)$ in half-plane $\text{Re}(s)\le1$.

\subsection{Riemann-Siegel formula}

Odlyzko et al. \cite{Odlyzko} also discusses the use of Riemann-Siegel formula for calculating the "zeros" of Riemann's $\zeta(s)$: 
\footnote{See Odlyzko et al. \cite{Odlyzko}, p.798.}

\begin{quotation}
A method for computing $\zeta(s)$ that is much more efficient than the Euler-Maclaurin formula (1.2) was discovered around 1932 in Riemann's unpublished papers by C. L. Siegel \cite{Siegel}. This formula [\cite{Siegel}, Equation (32)], now universally referred to as the Riemann-Siegel formula, is presented in \S 2. Roughly speaking, it enables one to compute $\zeta(\sigma+it)$ for $t$ large and $\sigma$ bounded to within $\pm t^{~c}$ for any constant $c$ in about $t^{1/2}$ steps. (Since $\zeta(\overline{s})$ = $\overline{\zeta(s)}$, we will always assume that $t > 0$.) The Riemann-Siegel formula is the fastest method for computing the zeta function to moderate accuracy that is currently known, and has been used for all large scale computations since the 1930s.
\end{quotation}
However, the Riemann–Siegel formula is: 
\begin{quotation}
an asymptotic formula for the error of the approximate functional equation of the Riemann zeta function, an approximation of the zeta function by a sum of two finite Dirichlet series.
\footnote{See Wikipedia \cite{Riemann-Siegel}.}
\end{quotation}
Edwards \cite{Edwards} 
confirms that the functional equation of $\zeta(s)$ is used in the derivation of the Riemann-Siegel formula. 
\footnote{See Edwards \cite{Edwards}, \S7.2 at pp.137-138, citing Edward's \S1.5 at pp.12-15.}

Unfortunately, the functional equation of Riemann's $\zeta(s)$ is \textbf{\textit{not}} valid in logics with LNC, because the analytic continuation of $\zeta(s)$ is not valid in those logics. The Dirichlet series $\zeta(s)$ is proven to be divergent throughout half-plane $\text{Re}(s)\le 1$, so the analytic continuation of $\zeta(s)$ violates the LNC. Also, the sum of two finite series cannot approximate a divergent infinite series. So in logics with LNC, the Riemann–Siegel formula is an approximation of an invalidity.

Moreover, the Riemann-Siegel formula is "an approximation of the [Riemann] zeta function by a sum of two finite Dirichlet series." 
\footnote{See Wikipedia \cite{Riemann-Siegel}.}
But summing two finite series, in order to obtain a finite value, is not a logically valid method of "approximating" a divergent infinite series (Dirichlet series $\zeta(s)$).

\subsection{Other Methods}

Odlyzko et al. \cite{Odlyzko} also discloses other methods for calculating the "zeros" of Riemann's $\zeta(s)$, including a method by Turing \cite{Turing2}, a method using Fast Fourier Transforms, 
\footnote{See Odlyzko et al. \cite{Odlyzko}, p.800, Eq.1.7; and pp.803-804, \S 3 "Application of the fast Fourier transform."}
etc. See also Gourdon et al. \cite{Gourdon} for additional discussion.

However, these other methods share the same problems as the Euler-Maclaurin and Riemann-Siegel formulas. All of these formulas are approximations of Riemann's $\zeta(s)$, which is \textit{invalid} in half-plane $\text{Re}(s)\le1$ in logics with LNC. Therefore, the functional equation of Riemann's $\zeta(s)$ must also be invalid in logics with LNC. All zeros calculated by these "approximations" are \textit{neither} zeros of Riemann's $\zeta(s)$ \textit{nor} zeros of the Dirichlet series $\zeta(s)$. 

\section{If AC of $\zeta(s)$ is False, RH is a Paradox, Due to Lack of Zeros }

\subsection{Material Implication}

Material implication, is the definition of the conditional "if $p$ then $q$" in both classical and intuitionistic logics. It states that the conditional "if $p$ then $q$", $(p\rightarrow q)$ is logically equivalent to $\neg (p\land \neg q)$. 
\footnote{By De Morgan's Laws (which are accepted in classical logic but not in intuitionistic logic) $\neg (p\land \neg q)$ is further equivalent to $(\neg p\lor q)$. }
So in classical and intuitionistic logics, material implication is counter-intuitively always "true" when $p$ is "false".
\footnote{See Tarski \cite{Tarski}, pp.25-26: "The logicians ... adopted the same procedure with respect to the phrase "\textit{if ..., then ...}" as they had done in the caso of the word "\textit{or}". For this purpose, they extended the usage of this phrase, considering an implication as a meaningful sentence even if no connection whatsoever exists between its two members, and they made the truth or falsity of an implication dependent exclusively upon the truth or falsity of the antecedent and consequent. 

To characterize this situation briefly, we say that contemporary logic uses IMPLICATIONS IN MATERIAL MEANING, or simply, MATERIAL IMPLICATIONS; this is opposed to the usage of IMPLICATIONS IN FORMAL MEANING, or simply, FORMAL IMPLICATION, in which case the presence of a certain formal connection between antecedent and consequent is an indispensable condition of the meaningfulness and truth of the implication. The concept of formal implication ... is narrower than that of material implication[.]"}
\footnote{See also Tarski \cite{Tarski}, p.26:
"In order to illustrate the foregoing remarks, let us consider the following four sentences:

\textit{if $2\cdot 2 =4$, then New York is a large city;}

\textit{if $2\cdot 2 =5$, then New York is a large city;}

\textit{if $2\cdot 2 =4$, then New York is a small city;}

\textit{if $2\cdot 2 =5$, then New York is a small city.}

In everyday language, these sentences would hardly be considered as meaningful, and even less true. From the point of view of mathematical logic, on the other hand, they are all meaningful, the third sentence being false, while the remaining three are true."
}
\footnote{See also Grattan-Guinness \cite{Grattan-Guinness}, p.329, describing Hardy's review, in the \textit{Times Literary Supplement}, of Russell's \cite{Russell3} \textit{Principles of Mathematics} (emphasis added): "On the logical aspects, [Hardy] stressed the \textit{unintuitive character} of [material] implication, that 'every false proposition implies every other proposition, true or false'." Did Hardy fail to consider applying Russell's work to the Riemann Hypothesis?}

This counter-intuitive aspect of "material implication" does not exist in the more narrowly defined "formal implication". In "formal implication", if the statement $p \Rightarrow q$ is \textit{true}, then the statement $\neg p \Rightarrow \neg q$ is \textit{false} (a.k.a "Denying the Antecedent").
\footnote{See Davis et al. \cite{Davis2}. p.301: "Denying the Antecedent (Invalid): $p \Rightarrow q, \neg p \therefore \neg q$.

[D]enying the antecedent [is] easily refuted by finding counterinstances, such as:

If whales are fish, then they are aquatic.

Whales are not fish.

$\therefore$ Whales are not aquatic."}

When the material conditional is applied to the RH, it holds that the RH is true, because RH states:
\begin{quotation}
If $\zeta(s)=0$, then all $\zeta(s)=0$ are on the critical line $\text{Re}(s)=0.5$.
\end{quotation}
and because $\zeta(s)$, as defined by Dirichlet series $\zeta(s)$, has no zeros. So RH is true. 

However, according to material implication, a statement we call "anti-RH" (ARH) is true too. It states:
\begin{quotation}
If $\zeta(s)=0$, then all $\zeta(s)=0$ are off the critical line $\text{Re}(s)=0.5$.
\end{quotation}
So if $\zeta(s)$ has no zeros, then this "anti-RH" is true. Yet  it is paradoxical for both RH and this "anti-RH" to be true.

\subsection{The Vacuous Subjects of the Riemann Hypothesis} \label{RH}

The Riemann Hypothesis (RH) states that "all non-trivial zeros of $\zeta(s)$ are on the critical line $\text{Re}(s)=0.5$". 

The Dirichlet series $\zeta(s)$ is proven to be divergent in half-plane $\text{Re}(s)\le1$. So analytic continuation of $\zeta(s)$ to that half-plane violates the LNC (and therefore must false), and thus in classical and intuitionistic logics, triggers ECQ. So $\zeta(s)$ is exclusively defined by the Dirichlet series $\zeta(s)$, \textit{which has no zeros and no poles}. Therefore, \text{none} of the zeros assumed by the RH exist. 
\footnote{Also, Riemann's functional equation of $\zeta(s)$ is invalidated in logics with LNC by the Dirichlet series $\zeta(s)$, which is proven to be divergent throughout half-plane $\text{Re}(s)\le1$. Thus $\zeta(1-s)$ is divergent at $\text{Re}(s)\ge0$. This contradicts Gelbart et al. \cite{Gelbart} at p.60, "our emphasis will be on explaining how we know that $\zeta(s)$ extends meromorphically to the entire complex plane and satisfies the functional equation." }
These non-existent zeros of $\zeta(s)$ constitute \textit{vacuous subjects} of a proposition, just like Russell's \cite{Russell2} famous example of "the present King of France" in the proposition "the present King of France is bald".

So given that the RH is a proposition with vacuous subjects, what is its truth -value? The answer: it depends on the system of logic that is applied.  
\footnote{RH's truth-value also depends upon the formulation of RH. See, e.g. Gelbart et al. \cite{Gelbart}, p.60: "The Riemann Hypothesis: $\zeta(s) \ne 0$ for $\text{Re}(s) > 1/2$." According to Dirichlet series $\zeta(s)$, this version of RH is true. However, Dirichlet series $\zeta(s)$ \textit{has no poles and no zeros}. In contrast, Gelbart et al. \cite{Gelbart}, p.60 falsely assumes that $\zeta(s)$ has both poles and zeros: "Our role here is not so much to focus on the \textit{zeroes} of $\zeta(s)$, but in some sense rather on its \textit{poles}." }

\subsection{"Vacuous Subjects" Generate Paradoxes}

Turing \cite{Turing} argues that the Riemann Hypothesis ("RH") is a "number-theoretic" problem.
\footnote{See Turing \cite{Turing}, p.165: "It is easy to show that a number of unsolved problems, such as the problem of the truth of Fermat's last theorem, are number-theoretic. There are, however, also problems of analysis which are number-theoretic. The Riemann hypothesis gives us an example of this."}
This classification is incorrect, and is the reason why the problem has remained unsolved for so long.
The fact of the matter is that the RH is a logic problem 
\footnote{Turing's error is in falsely assuming that Riemann's analytic continuation of the Dirichlet series $\zeta(s)$ is valid. See Turing \cite{Turing}, p.165: "We denote by $\zeta(s)$ the function defined for $\text{Re}(s) = \sigma > 1$ by the series $\sum n^{-s}$ and over the rest of the complex plane with the exception of the point $s=1$ by analytic continuation."}
In classical logic, the RH is an undecidable paradox. Riemann's $\zeta(s)$ is not valid in the half-plane where it contradicts the Dirichlet $\zeta(s)$. So $\zeta(s)$ is exclusively defined by the Dirichlet $\zeta(s)$, which \textit{has no zeros}. So the zeros of the RH form an empty set.  Therefore, RH is a proposition that suffers from "reference failure".
\footnote{See Haack \cite{Haack2}, pp.14-15: "Another challenge to classical logic derives from the phenomenon of reference failure, i.e., of sentences containing proper names (such as "Mr. Pickwick" or "Odysseus") or definite descriptions (such as "the present king of France" or "the greatest prime number" which have no referent."}
RH has "vacuous subjects", due to $\zeta(s)$ having no zeros. According to classical logic's "material implication", all propositions pertaining to "vacuous subjects" are true, including contradictory propositions.

In regards to the Riemann Hypothesis ("RH"), its traditional phrasing is "\textit{all zeros of $\zeta(s)$ are on the critical line} $\text{Re}(s)=0.5$".
\footnote{See Edwards \cite{Edwards}, \S 1.9, p.19: "Riemann's next statement is even more baffling. He states that the number of roots [$\rho$ of $\xi(\rho)=0$] \textit{on the line} $\text{Re}(s)=0.2$ is also "about" [$T/2\pi \cdot \log T/2\pi - T/2\pi$] ... He gives no indication of a proof at all, and no one since Riemann has been able to prove (or disprove) this statement ... He says he considers it 'very likely' that the roots all do lie on [the critical line] $\text{Re}(s)=0.5$, but says that he was not able to prove it".  See also Edwards \cite{Edwards}, \S7.8, pp.164-166; and chapter 9, pp.182-202.}
\footnote{Riemann's statement in \cite{riemann1859number}, p.4, as translated by Wilkins, is: "One now finds indeed approximately this number of Real roots [of $\xi(t)=0$] within these limits, and it is very probable that all roots are Real. Certainly one would wish for a stricter proof here ...".} 
The traditional phrasing of the negation of RH ("$\neg$RH") is "\textit{not all zeros of $\zeta(s)$ are on the critical line} $\text{Re}(s)=0.5$". 
However, RH falsely assumes that Riemann's $\zeta(s)$ is valid and \textit{has zeros}. 
\footnote{The Dirichlet series "expression" of $\zeta(s)$ has no zeros, the Euler product "expression" of $\zeta(s)$ has no zeros, and the Riemann "expression" of $\zeta(s)$ is not valid in logics with LNC.} 

If $\zeta(s)$ has no zeros, then according to classical logic, both RH and $\neg$RH are "vacuously true" according to material implication. So in classical logic, RH is a paradox, because it is simultaneously true and false.   
\footnote{See Gardner \cite{Gardner} for many other examples of paradoxes.}
\footnote{See also Scruton \cite{Scruton}, Chapter 27 "Paradox", pp.397-412, and 575.}
This contradiction violates LNC and LEM. and the violation of LNC triggers ECQ.

The results in classical logic are identical to the results in set theory. Because $\zeta(s)$ has no zeros, \textit{both} RH and $\neg$RH are propositions with "vacuous subjects", of the type discussed by Frege, 
\footnote{See Frege's \textit{Über Sinn und Bedeutung} ("On Sense and Denotation") \cite{Frege2}.}
Russell, 
\footnote{See Russell's \textit{On Denoting}, \cite{Russell2}.}
Strawson, 
\footnote{See Strawson's \textit{On Referring}, \cite{Strawson}.}
and others. (The most famous example being "The present King of France is bald"). According to Frege, 
\footnote{See Frege's \textit{Über Sinn und Bedeutung} ("On Sense and Denotation") \cite{Frege2}.}
\textit{neither} RH nor $\neg$RH have any truth-value (i.e. a "truth-value gap"). But according to Russell, \textit{both} RH and  $\neg$RH are \textit{both} true and false (i.e. a "truth-value glut").

Moreover, the same paradoxical results are obtained by rephrasing RH and $\neg$RH to expressly state the assumption that $\zeta(s)$ has zeros. This can be done in two ways: (1) as conditional propositions, or (2) as conjunctions. It turns out that the conditional propositions are negations of the conjunctions (and vice versa).
\footnote{
The negation of the conditional proposition is determined as follows: The sequent for conditional propositions (material implication) is: ($A \supset B) \Leftrightarrow (\lnot A \lor B$). The negation of both sides of this equivalence results in: $\lnot(A \supset B) \Leftrightarrow \lnot(\lnot A \lor B)$, which according to De Morgan's laws and the Law of Double Negation is equivalent to: $\lnot(\lnot A \lor B) \Leftrightarrow (A \land \lnot B)$. So in regards to $RH_1$, its negation is "\textit{$\zeta(s)$ has zeros, and not all zeros of $\zeta(s)$ are on the critical line} $\text{Re}(s)=0.5$", which is $\overline{RH_2}$. When performed on $\overline{RH_1}$, the result is $RH_2$. Negation of the conditional is as follows: $\lnot (A \land \lnot B) \Leftrightarrow  (\lnot A \lor B) \Leftrightarrow (A \supset B)$. So $\lnot RH_1 \Leftrightarrow \overline{RH_2}$ and $\lnot RH_2 \Leftrightarrow \overline{RH_1}$}

If RH is rewritten as a conditional proposition, it becomes $RH_1$: "\textit{if $\zeta(s)$ has zeros, then all zeros of $\zeta(s)$ are on the critical line} $\text{Re}(s)=0.5$". Its sequent is ($A \supset B$). 

Likewise, $\neg$RH becomes $\overline{RH_1}$: "\textit{if $\zeta(s)$ has zeros, then not all zeros of $\zeta(s)$ are on the critical line} $\text{Re}(s)=0.5$". In classical logic, if $\zeta(s)$ actually had zeros, 
one of $RH_1$ and $\overline{RH_1}$ would be true, and the other would be false, because according to material implication, ($A \supset B) \Leftrightarrow (\lnot A \lor B$). 

However, $\zeta(s)$ has no zeros, so material implication holds that the $RH_1$ and $\overline{RH_1}$ are \textit{both true}, because each has an antecedent portion ("\textit{$\zeta(s)$ has zeros}") that is false. 
\footnote{See Carnap \cite{Carnap}, p.8: "The sentence '$(A)\supset(B)$' is an abbreviation for '$[\sim(A)]\lor(B)$'", and "Also, in connection with the conditional '$(A)\supset(B)$' we find it convenient to retain the name '\textit{antecedent}' for the first component '($A$)' and the name '\textit{consequent}' for the second component '($B$)'." So, given that the antecedent of RH is false ($\zeta(s)$ has no zeros), then any consequent is true.}
Therefore, regardless of the truth or falsity of the consequent portion ("\textit{all zeros of $\zeta(s)$ are on the critical line} $\text{Re}(s)=0.5$", or "\textit{not all zeros of $\zeta(s)$ are on the critical line} $\text{Re}(s)=0.5$"), the proposition as a whole is \textit{true}. 

If RH is rewritten as a conjunction ($RH_2$), it becomes: "\textit{$\zeta(s)$ has zeros, and all zeros of $\zeta(s)$ are on the critical line} $\text{Re}(s)=0.5$". Its sequent ia ($A \land B$).
The negation of RH ($\neg$RH) becomes $\overline{RH_2}$: "\textit{$\zeta(s)$ has zeros, and not all zeros of $\zeta(s)$ are on the critical line} $\text{Re}(s)=0.5$". In classical logic, if $\zeta(s)$ actually had zeros, then by conjuction,
one of $RH_2$ and $\overline{RH_2}$ would be true, and the other would be false. 

However, $\zeta(s)$ has no zeros, so conjuction holds that the $RH_2$ and $\overline{RH_2}$ are \textit{both false}, because each has an antecedent portion ("\textit{$\zeta(s)$ has zeros}") that is false. Therefore, regardless of the truth or falsity of the consequent portion ("\textit{all zeros of $\zeta(s)$ are on the critical line} $\text{Re}(s)=0.5$", or "\textit{not all zeros of $\zeta(s)$ are on the critical line} $\text{Re}(s)=0.5$"), according to conjunction, the proposition as a whole is always \textit{false}.

So RH and its negation $\neg$RH are paradoxes, and have either a truth-value glut, or a truth-value gap. These results are impermissible in classical logic, due to the LNC and LEM.

\subsection{The Riemann Hypothesis and Venn's "Modern" Square of Opposition}

Lande \cite{Lande}, p.268: 
\begin{quotation}
It was once thought that what is known as the [Aristotelian] Square of Opposition captured all logically significant sentences, as well as their mutual relations. Although the [Aristotelian] Square of Opposition fails to do justice to the complexity of the sentences that you will soon be encountering, it provides you with a structure that is actually quite helpful for translating increasingly complex sorts of English sentences into logical notation.
\end{quotation}

\textbf{"A" Propositions (Universal Affirmatives): "All S are P".}
"All zeros of $\zeta(s)$ are on the critical line."
"The Riemann Hypothesis (RH)."
$(\forall s) \{(\zeta(s)=0) \rightarrow (\text{Re}(s)=0.5)\}$

\textbf{"E" Propositions (Universal Negations): "No S are P".}
"No zeros of $\zeta(s)$ are on the critical line."
"All zeros of $\zeta(s)$ are off the critical line."
"The Anti-Riemann Hypothesis (ARH)."
$(\forall s) \{(\zeta(s)=0) \rightarrow (\text{Re}(s)\ne0.5)\}$

\textbf{"I" Propositions (Particular/Existential Affirmatives): "Some S are P".}
"Some zeros of $\zeta(s)$ are on the critical line."
"There exists a zero of $\zeta(s)$ on the critical line."
"Negation of the Anti-Riemann Hypothesis ($\neg$ARH)."
$(\exists s) \{(\zeta(s)=0) \land (\text{Re}(s)= 0.5))\}$

\textbf{"O" Propositions (Particular/Existential Negations): "Some S are not P".}
"Some zeros of $\zeta(s)$ are off the critical line."
"There exists a zero of $\zeta(s)$ off the critical line."
"Negation of the Riemann Hypothesis ($\neg$RH)."
$(\exists s) \{(\zeta(s)=0) \land (\text{Re}(s)\ne 0.5)\}$

\begin{figure} [ht]
    \centering
    \includegraphics[scale=1.5]{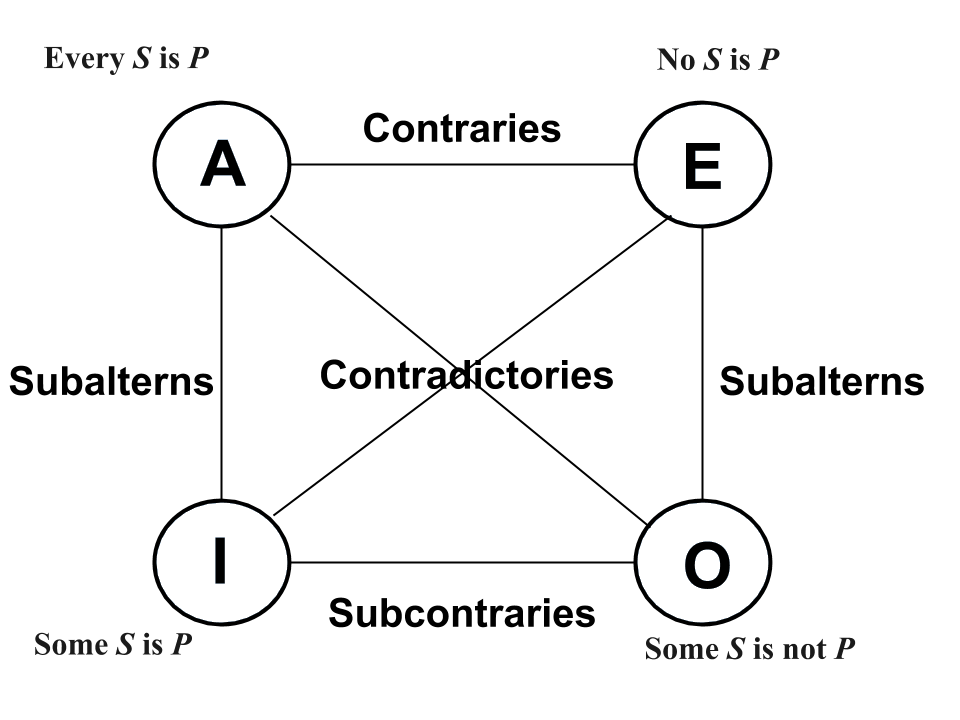}
    \caption{The Traditional Square of Opposition}
    \label{fig:trad_square}
\end{figure}

The Traditional Square of Opposition \footnote{See Parsons \cite{Parsons}.} is shown in Figure \ref{fig:trad_square}. In the Traditional Square of Opposition, RH's contrary is "anti-RH" (ARH). ARH's subaltern is $\neg$RH. RH's subaltern is $\neg$ARH.  RH and $\neg$RH are contradictories, as are ARH and $\neg$ARH. 

The Internet Encyclopedia of Philosophy (IEP) article discloses the following about the Traditional Square of Opposition \cite{Square}:
\begin{quotation}
Given the assumption made within [Aristotelian] categorical logic, that every category contains at least one member, the following relationships, depicted on the [Aristotelian] square, hold:

Firstly, A and O propositions are contradictory, as are E and I propositions. Propositions are contradictory when the truth of one implies the falsity of the other, and conversely. 

Secondly, A and E propositions are contrary. Propositions are contrary when they \textit{cannot} both be true. 

Next, I and O propositions are subcontrary. Propositions are subcontrary when it is impossible for both to be false. 

Lastly, two propositions are said to stand in the relation of subalternation when the truth of the first ("the superaltern") implies the truth of the second ("the subaltern"), but not conversely. 

The presupposition, mentioned above, that all categories have at least one member, has been abandoned by most later logicians.
\end{quotation}

But $\zeta(s)=0$ is an empty set.  According to the truth table of the material implication operator,
both $F \rightarrow T$ and $F \rightarrow F$ are true, SO both $\zeta(s)=0 \rightarrow RH$ and $ \zeta(s)=0  \rightarrow \neg RH$ are true, So in classical logic, both of the contrary A proposition  and E proposition (a.k.a. RH and ARH) are "vacuously true", thereby violating the requirement that they 
\textit{cannot} both be true, and forming an undecidable semantic paradox that triggers ECQ.

According to Davis et al. \cite{Davis2}, p.240: "Aristotelian categoricals and their Venn transforms have the same truth values as long as there is something to which their subjects apply."

\begin{figure} [ht]
    \centering
    \includegraphics[scale=1.5]{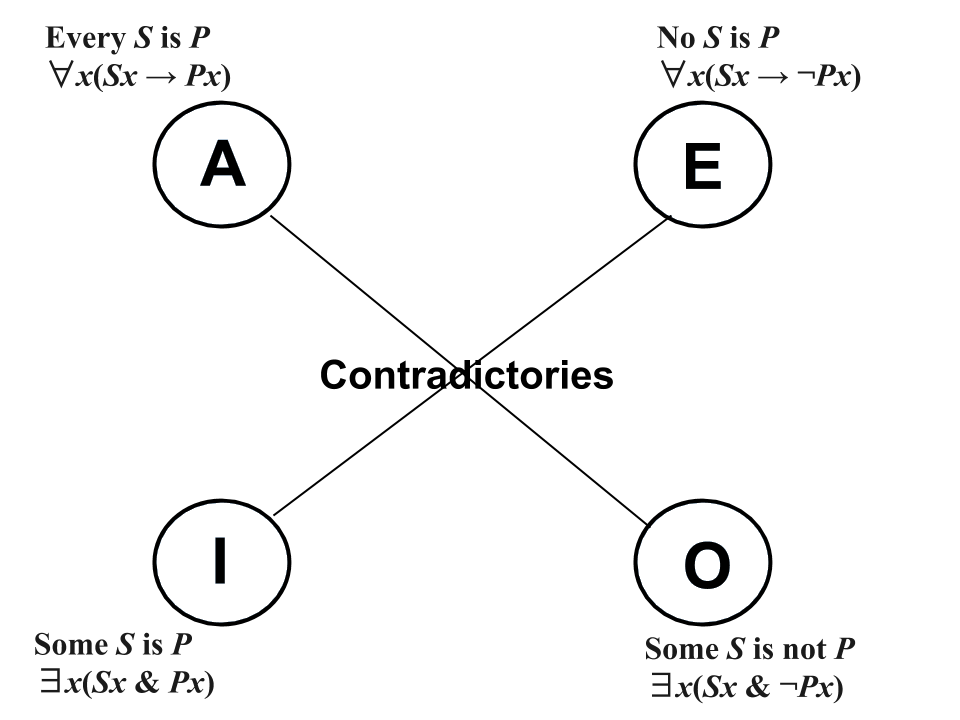}  
    \caption{Venn's "Modern" Square of Opposition}
    \label{fig:venn_square}
\end{figure}

According to the Internet Encyclopedia of Philosophy (IEP) website's entry on the Square of Opposition \cite{Square}, Venn's "Modern" Square of Opposition (see Figure \ref{fig:venn_square}) differs from the Traditional version as follows:
\begin{quotation}
The presupposition [in Aristotelian logic], mentioned above, that all categories have at least one member, has been abandoned by most later logicians. Modern logic deals with uninstantiated terms such as "unicorn" and "ether flow" the same as it does other terms such as "apple" and "orangutan". 

When dealing with "empty categories", the relations of being contrary, being subcontrary and of subalternation no longer hold. Consider, e.g., "all unicorns have horns" and "no unicorns have horns." Within contemporary logic, these are both regarded as true, so strictly speaking, they cannot be contrary, despite the former's status as an A proposition and the latter's status as an E proposition. Similarly, "some unicorns have horns" (I) and "some unicorns do not have horns" (O) are both regarded as false, and so they are not subcontrary.

Obviously then, the truth of "all unicorns have horns" does not imply the truth of "some unicorns have horns," and the subalternation relation fails to hold as well. Without the traditional presuppositions of "existential import", i.e., the supposition that all categories have at least one member, then only the contradictory relation holds.

On what is sometimes called the "modern square of opposition" (as opposed to the traditional square of opposition sketched above) the lines for contraries, subcontraries and subalternation are erased, leaving only the diagonal lines for the contradictory relation.
\end{quotation}

\section{The Truth-Value of RH (a Paradox) Depends on the Logic}

\subsection{In Intuitionistic Logic: RH is False (and Thus Decidable)}

So if Riemann's version of $\zeta(s)$ is false. what about the RH? In intuitionistic logic, the truth-value "true" applies \textit{only} to proven propositions pertaining to objects that have been proven to exist. 
\footnote{See Vafeiadou et al. \cite{Vafeiadou}, p.2, citing Brouwer \cite{Brouwer}, p.79: "Moreover, the '... \textit{existence} of a mathematical system satisfying a set of axioms can never be proved from the consistency of the logical system based on those axioms,' but only by construction."}
\footnote{See also Bridges et al. \cite{Bridges}, \S 2 "The Constructive Interpretation of Logic": "$\exists$ (there exists):	to prove $\exists x P(x)$ we must construct an object $x$ and prove that $P(x)$ holds", and "These \textit{BHK-interpretations} (the name reflects their origin in the work of Brouwer, Heyting, and Kolmogorov) can be made more precise using Kleene’s notion of \textit{realizability}", citing (Dummett \cite{Dummett}, pp.222–234 and Beeson \cite{Beeson}, Chapter VII).}
The truth-value "\textit{false}" applies to objects that have been proven to \textit{not} exist, and the truth-value "\textit{neither true nor false}" (a violation of the LEM) applies to objects whose existence has yet to be proven or disproven. 
\footnote{See Moschovakis \cite{Moschovakis}: "Intuitionistic propositional logic is effectively decidable, in the sense that a finite constructive process applies uniformly to every propositional formula, either producing an intuitionistic proof of the formula or demonstrating that no such proof can exist."}
\footnote{See also Haack \cite{Haack3}, p.92: "[In Intuitionism,] only \textit{constructible} mathematical entities are admitted ... and only constructive proofs of mathematical statements are admitted, so that, for instance, a statement to the effect that there is a number with such-and-such property is provable only if a number with that property is constructible." }

As discussed above, the derivation of Riemann's $\zeta(s)$ is invalid in logics with LNC, and so $\zeta(s)$ is exclusively defined by Dirichlet series $\zeta(s)$, which has no zeros. This means that the RH is directed to non-existent objects. Brouwer's intuitionistic logic refuses to admit propositions regarding mathematical objects, unless the objects have existence proofs. Regarding propositions about objects that have been proven to \textit{not} exist, such propositions are "\textit{false}". 
\footnote{See also Bridges et al. \cite{Bridges}, \S 2 "The Constructive Interpretation of Logic": "$\exists$ (there exists):	to prove $\exists x P(x)$ we must construct an object $x$ and prove that $P(x)$ holds", and "$\forall$  (for each/all):	a proof of $\forall x \in  SP(x)$ is an algorithm that, applied to any object $x$ and to the data proving that $x \in S$, proves that $P(x)$ holds."}
So in Brouwer's intuitionistic logic, RH is "\textit{false}" (and decidable). 
\footnote{By applying the definitions of \S 2 "The Constructive Interpretation of Logic" in Bridges et al. \cite{Bridges}, the absence of zeros turns the RH into a "decision problem" that "can be posed as a yes-no question of the input values". See Wikipedia \cite{DecisionProblem}. The decision question is "Given that there are no zeros of $\zeta(s)$, are they all on the critical line?" The intuitionist answer is "no". }

For intuitionistic logic, the RH's reference to non-existent zeros constitutes a perfect example of its criticisms of classical logic. 

\subsection{In Classical Logic, RH is an Undecidable Paradox}

In classical logic, RH is an undecidable paradox resulting from material implication, because the Dirichlet series $\zeta(s)$ has no zeros, and these non-existent zeros are "vacuous subjects" of the RH. This can be seen most clearly in Venn's "Modern"  Square of Opposition. 

The Riemann Hypothesis (RH) is: "\textit{All} zeros of $\zeta(s)$ are \textit{on} the critical line $\text{Re}(s)=0.5$." According to the original Aristotelian Square of Opposition, RH's contrary (anti-RH or "ARH") is: "\textit{All} zeros of $\zeta(s)$ are \textit{off} the critical line". The RH's subaltern ($\neg$ARH) is: "There exists a zero of $\zeta(s)$ \textit{on} the critical line", and the ARH's subaltern ($\neg$RH) is: "There exists a zero of $\zeta(s)$ \textit{off} the critical line".  

When written in first-order logic notation, Aristotle's Square of Opposition is as follows: RH is $(\forall s) \{(\zeta(s)=0) \rightarrow (\text{Re}(s)=0.5)\}$. RH's contrary is anti-RH (ARH): $(\forall s) \{(\zeta(s)=0) \rightarrow (\text{Re}(s)\ne0.5)\}$. ARH's subaltern ($\neg$RH) is $(\exists s) \{(\zeta(s)=0) \land (\text{Re}(s)\ne 0.5)\}$. RH's subaltern ($\neg$ARH) is $(\exists s) \{(\zeta(s)=0) \land (\text{Re}(s)= 0.5))\}$.

But in the case of "vacuous subjects" (for RH, the absence of $\zeta(s)=0$), the original Aristotelian Square of Opposition \textit{fails}, because it assumes non-vacuous subjects. 
\footnote{See Davis et al. \cite{Davis2}, p.239: "Aristotelian categoricals presuppose that their subjects apply to something."}
In contrast, Venn's "Modern" Square of Opposition is valid \textit{even in the case} of "vacuous subjects". 
\footnote{See Davis et al. \cite{Davis2}, p.240: "The A* and E* propositions here are true: since nothing is a 2006 Edsel, nothing is both a 2006 Edsel and either a four-door or a nonfour-door."}
But Venn's version differs from Aristotle's in that there are no subaltern or contrary relationships - only the contradictory relationships of Aristotle's version. 

When written in first-order logic notation, Venn's "Modern" Square of Opposition is as follows: RH is $(\forall s) \{(\zeta(s)=0) \rightarrow (\text{Re}(s)=0.5)\}$. RH's contrary is anti-RH (ARH): $(\forall s) \{(\zeta(s)=0) \rightarrow (\text{Re}(s)\ne0.5)\}$. $\neg$RH is $(\exists s) \{(\zeta(s)=0) \land (\text{Re}(s)\ne 0.5)\}$. $\neg$ARH is $(\exists s) \{(\zeta(s)=0) \land (\text{Re}(s)= 0.5))\}$. 

RH contradicts $\neg$RH, and ARH contradicts $\neg$ARH. But $\zeta(s)$ has no zeros. So in Venn's "Modern" Square of Opposition, as in classical logic, both $\neg$RH and $\neg$ARH are false. More importantly, their respective contradictory statements (RH and ARH) are both "vacuously true", due to classical logic's material implication (which holds that a false proposition implies any proposition). But according to both classical logic, if RH is true, then ARH \textit{should be} false, and vice versa. Yet both are true, due to their "vacuous subjects" (the non-existent zeros of the zeta function). Therefore, in classical logic, RH is a semantic paradox, is undecidable, violates the LNC, and triggers ECQ. 
\footnote{For another example of an undecidable paradox in classical logic, see the "Liar Paradox" in Gödel's \textit{On Formally Undecidable Propositions of Principia Mathematica and Related Systems} \cite{Godel} (in German) and \cite{Godel2} (in English).}

\subsection{In 3VL, RH has the Third Truth-Value, and is Decidable}

In some three-valued logics (3VL) such as Bochvar's, and Priest's $LP$, paradoxes are the original intended use of the third truth-value. In other 3VLs, such as Łukasiewicz's, paradoxes are \textit{not} the original intended use. But when applied to paradoxes, they assign the third truth-value to paradoxes, and thus the LNC is avoided. 

Moreover, in Łukasiewicz's and Kleene's 3VLs (unlike classical and intuitionist logics), the truth table of material implication shows the 3rd truth-value as \textit{not} resulting in ECQ. 
\footnote{See Urquhart \cite{Urquhart2}, p.260: "The Kleene system does not contain the paradox of material implication $p \vdash q \lor \neg q$; however it contains $p,\neg p \vdash q$, so it is not free of the paradoxes of material implication. The relationship between the two systems can be briefly indicated by noting that while Kleene allows for the possibility "\textit{neither} true nor false", Anderson and Belnap allow for the possibility '\textit{both} true and false'."  See e.g. Belnap \cite{Belnap} for discussion of a 4VL.}
Also, this third truth-value can be assigned the label "indefinite", "undecidable", "unknown", "partially true", or even "paradox".
\footnote{See Stewart \cite{Stewart}, p.242: "We generally assume that an unsolved conjecture, like the Riemann Hypothesis, is either true or false, so either there's a proof or a disproof. ... Classical logic, with its sharp distinction between truth and falsity, with no middle ground, is two-valued. Gödel's discovery suggests that for mathematics, a three-valued logic would be more appropriate: true, false, or undecidable."}

\subsection{Truth-Value and Decidability of RH Depends Upon the Logic Applied}

So, depending on the logic applied to the RH, (or alternatively, the "foundational logic" underlying the RH):
\footnote{See Moschovakis \cite{Moschovakis}: "Philosophically, intuitionism differs from logicism by treating logic as a part of mathematics rather than as the foundation of mathematics;" and "Hilbert’s formalist program, to justify classical mathematics by reducing it to a formal system whose consistency should be established by finitistic (hence constructive) means, was the most powerful contemporary rival to Brouwer’s developing intuitionism. In (\cite{Brouwer3}) Brouwer correctly predicted that any attempt to prove the consistency of complete induction on the natural numbers would lead to a vicious circle."}

\begin{enumerate}
\item \textit{In Intuitionistic logic}: the RH is decidable, and false. 
\item \textit{In Classical logic}: the RH is undecidable: a paradox that violates LNC and triggers ECQ. 
\item \textit{In 3VL (Bochvar's, Priest's)}: the RH is decidable: a paradox that has the 3rd of the three truth-values, does not violate LNC, and thus does not cause ECQ. 
\end{enumerate}
Each of these results is inconsistent with the other results, just as the respective logics are inconsistent with one another.

\subsection{There Exist Many Logics, So What Constitutes Proof?}

These conflicting results for RH, that vary depending on the logic applied, provide support for the criticism against Aristotle's concept of deductive proof. Bertrand Russell attributed this criticism to Timon of Phlius, the Pyrrhonist philosopher: 
\footnote{See Russell \cite{Russell5}, p.234. Russell does not cite any reference for this attribution to Timon. The present author has not found any reference that either supports or contradicts Russell's attribution.}
\begin{quotation}
The only logic admitted by the Greeks was deductive, and all deduction had to start, like Euclid, from general principles regarded as self-evident. Timon denied the possibility of finding such principles. Everything, therefore, will have to be proved by means of something else, and all argument will be either circular or an endless chain hanging from nothing. In either case nothing can be proved.
\end{quotation}

Aristotle believed that his "Three Laws of Thought" - The Law of Identity (LOI), Law of Non-Contraction (LNC), and the Law of the Excluded Middle (LEM) - were "self-evident" general principles. 
\footnote{See Russell \cite{Russell5}, p.234: "The only logic admitted by the Greeks was deductive, and all deduction had to start, like Euclid, from general principles regarded as self evident."}
\footnote{See Cohen \cite{Cohen2}, p.75: "Aristotle's greatest achievement is supposed to have been his 'Laws of Thought,' part of his attempt to put everyday language on a logical footing. His \textit{Prior Analytics} is the first attempt to create a system of formal deductive logic, whereas the  \textit{Posterior Analytics} attempts to use this to systematize scientific knowledge."} 
But given that intuitionistic logic selectively rejects the LEM, and given that the LEM and LNC fail in 3VL,
\footnote{See Haack \cite{Haack2}, p.5: "In Łukasiewicz's 3-valued logic (motivated by the idea, already suggested by Aristotle in \textit{De interpretatione} \S 9, in \textit{Organon}), that future contingent sentences are neither true nor false but 'indeterminate') both the Law of the Excluded Middle ('LEM;' 'p or not p') and the Law of Non-Contradiction ('LNC;' 'not both p and not-p') fail."}
\footnote{Note also that there is disagreement regarding "future contingents" in Aristotle's \cite{Aristotle} \textit{De interpretatione} \S 9, in \textit{Organon}. See also Haack's \cite{Haack3} ch.4 for arguments for and against LEM, due to future contingents.}
Aristotle's general principles are clearly not "self-evident" general principles. 
\footnote{See also Kuznetsov \cite{Kuznetsov}, in the 1974 \textit{Proc. of the ICM},  p.244: "One might also criticize the laws of intuitionistic logic—either from the standpoint of refusing from the so-called 'paradoxes of implication', which lead to different logics of rigorous implication; or from the point of view of accounting for the peculiarities of quantum-mechanical problems (in this case one axiom is doubtful, for the calculus without it see Tolstova \cite{Tolstova}); or in the light of immersion not in $S$4, but in weaker modal logic.}
\footnote{Kuznetsov also argues for what Haack \cite{Haack}, ch.12, \S 1, calls "local pluralism" of logics. See Kuznetsov \cite{Kuznetsov}, p.244: "Moreover, I am keeping to the view that none of fixed logic may be suitable in all the situations, for all cases of life; therefore a general investigation of different large classes of non-classical logics is useful. However, being unable to embrace the nonembraceable, I shall here restrict myself only to the consideration of propositional logics, and from them only the superintuitionistic logics, i.e., classical, intuitionistic, intermediate (between them) and absolutely contradictory."}

Do any "self-evident" general principles exist? LNC, LEM, and LOI were historically accepted as "self-evident" from Aristotle's time, up until Brouwer's intuitionistic logic rejected LEM (for propositions that cannot be either proved or disproved).
\footnote{See Davis \cite{Davis}, p.95.}
Heyting's version rejected the LEM entirely, and Łukasiewicz's 3VL rejected LEM by creating a third truth-value (intended to be used for future contingents, so as to keep LNC). Other 3VLs (e.g. Priest's $LP$) reject LNC. If there are no generally accepted principles, then Timon's argument is correct. 

However, even if there are no "self-evident" general principles, Timon's argument is wrong in the following situation: if a proposition is found to be "exclusively true" (to distinguish it from having "true" as one value in a truth-value glut) in \textit{every} logic, then it \textit{must} be "exclusively true" ("logically true"), 
despite the absence of any "self-evident" general principles. Therefore, the most restrictive logic is the logic that determines "logical truth". A candidate for such a logic would be the most restrictive of the intuitionistic logics.
\footnote{See Wikipedia \cite{LogicalTruth}: "Logical truths (including tautologies) are truths which are considered to be necessarily true. This is to say that they are considered to be such that they could not be untrue and no situation could arise which would cause us to reject a logical truth. It must be true in every sense of intuition, practices, and bodies of beliefs. However, it is not universally agreed that there are any statements which are \textit{necessarily} true."}
\footnote{See also Gómez-Torrente \cite{Gomez-Torrente}: "As we said above, it seems to be universally accepted that, if there are any logical truths at all, a logical truth ought to be such that it could not be false, or equivalently, it ought to be such that it must be true."}

Conversely, if a proposition is found to be "exclusively false" (to distinguish it from having "false" as one value in a truth-value glut) in \textit{every} logic, then it \textit{must} be "exclusively false", despite the absence of any "self-evident" general principles. Therefore, the least restrictive logic is the logic that determines logical falsehoods. 

A candidate for such a logic would be an MVL with an infinite number of truth-values.

\section{Aristotle, the Axiomatic Method, and the LNC}

\subsection{The Different Types of Logic (Including Deductive Logic)}

Logic is the study of the methods and principles used to distinguish between valid and invalid arguments.
\footnote{See Lee \cite{Lee}, p.2} 
A sound argument is valid argument, whose premises are \textit{all} true. 
\footnote{See Lee \cite{Lee}, p.19.}
Unfortunately, logic cannot determine if the premises assumed in an argument are true. Therefore, logic can only identify unsound arguments if they are invalid (regardless of whether or not all premises are true). Logic cannot determine if a \textit{valid} argument is sound or unsound, because it is unable to determine whether premises are true or false. 
\footnote{See Wikipedia \cite{Nyaya}, citing Church \cite{Church}: "Logic is the systematic study of the structure of propositions and of the general conditions of valid inference by a method, which abstracts from the content or matter of the propositions and deals only with their logical form. This distinction between form and matter is made whenever we distinguish between the logical soundness or validity of a piece of reasoning and the truth of the premises from which it proceeds[,] and in this sense is familiar from everyday usage."}
\footnote{Therefore, truth cannot be determined by logic alone. Premises can only be determined by the senses. This undercuts Plato's argument that because the senses are misleading,  truth must be determined by logic alone. See Kline \cite{Kline}, p.48: "Plato stressed the unreliability of sensory perceptions. Empirical knowledge, as Plato put it, yields opinion only."}

Logic is a normative discipline, in that it describes how we \textit{should} argue (i.e. "reason"), not how we \textit{actually} "reason". 
\footnote{See Lee \cite{Lee}, p.19. "[This] is the job of the psychologist."}
There are four main types of arguments: inductive, deductive, abductive, and analogical.
\footnote{See Fontainelle \cite{Fontainelle}, pp.182-183.}

Deductive reasoning "moves from the general to the particular, producing a necessary conclusion whose truth follows from that of premises."  
\footnote{See Fontainelle \cite{Fontainelle}, p.182.}
(The premises are assumed to be true). An example of deductive reasoning is: "Mythical animals do not really exist. Werewolves are mythical animals. Therefore Werewolves do not really exist."
\footnote{See Fontainelle \cite{Fontainelle}, p.182.}

According to Russell \cite{Russell5}: "The only logic admitted by the [ancient] Greeks was deductive, and all deduction had to start, like Euclid, from general principles regarded as self-evident."
\footnote{See Russell \cite{Russell5}, p.234.}
In contrast, Timon of Phlius denied the possibility of finding such self-evident general principles. 
According to Timon:
\footnote{See Russell \cite{Russell5}, p.234.}
 
\begin{quotation}
Everything, therefore, will have to be proved by means of something else, and all argument will be either circular
\footnote{Otherwise known as "begging the question" or \textit{petitio principii} (assuming the principal): the logical fallacy of assuming that the statement under examination is true. In other words, using a premise to support itself.}
\footnote{See e.g. Cameron \cite{Cameron}, citing Cardano \cite{Cardano}, p.246: "Mathematics, however, is, as it were, its own explanation; this, although it may seem hard to accept, is nevertheless true, for the recognition that a fact is so is the cause upon which we base the proof."}
or an endless chain hanging from nothing.
\footnote{Informally referred to as "turtles all the way down".}
In either case nothing can be proved. This argument, as we can see, cut at the root of the Aristotelian philosophy which dominated the Middle Ages.
\end{quotation}
 
The argument is that it is impossible to prove that a proof is sound,
\footnote{These arguments have also been called the "Münchhausen-Trilemma" (\textit{Dogmatismus – unendlicher Regreß – Psychologismus}) attributed to German philosopher Hans Albert. See Wikipedia \cite{Muenchhausen-Trilemma}, citing Westermann \cite{Westermann}, p. 15, in turn citing Albert \cite{Albert}, p. 11.}
As discussed above, a sound argument is both valid \textit{and} its premises are true. \textit{Timon's argument is, paradoxically, a proof by deductive reasoning that there is no proof by deductive reasoning.}
\footnote{Modern day proponents of Pyrrhonist philosophy are called "Fallibilists". Notable proponents of this school of philosophy include Charles Sanders Peirce, Karl Popper, W.V.O. Quine, and Susan Haack.}

This result, in turn. makes deductive reasoning consistent with the philosophy of the ancient Skeptics, and also with other methods of reasoning, by proving that conclusions obtained by deductive reasoning are not necessarily true. 
\footnote{Paradoxically, this consistency with the other main methods of reasoning (regarding the production of uncertain conclusions) addresses the central concern of Aristotelian and classical deductive logic: that of consistency (the LNC).}
Even mathematicians have (inadvertantly) conceded this point, Manin et al. \cite{Manin} stated: "A proof only becomes a proof after the social act of 'accepting it as a proof'."
\footnote{See Cameron \cite{Cameron}, citing Manin et al. \cite{Manin}.} 
Borovik \cite{Borovik} reaffirms this sentiment, stating: "Manin describes the act of acceptance as a social act; however, the importance of its personal, psychological component can hardly be overestimated."
\footnote{See Cameron \cite{Cameron}, citing Borovik \cite{Borovik}, p.35.} 
But as Bertrand Russell \cite{Russell6} perceptively points out: "The fact that an opinion has been widely held is no evidence whatever that it is not utterly absurd; indeed in view of the silliness of the majority of mankind, a widely spread belief is more likely to be foolish than sensible."
\footnote{See Cameron \cite{Cameron}, citing Russell \cite{Russell6}, p.58.} 

In contrast to deductive reasoning, inductive reasoning "begins with the particular and proceeds to the general. Things are observed, [and] then a rule or cause is proposed to account for them."
\footnote{See Fontainelle \cite{Fontainelle}, p.182.}
If the premises are true, then the conclusion is\textit{likely} to be true. The truth of the premises does \textit{not completely determine} the truth of the conclusion. The argument indicates some sort of probability.
\footnote{See Lee \cite{Lee}, p.12.}
"[This] is why, strictly speaking, no scientific theory is regarded as being true."
\footnote{See Fontainelle \cite{Fontainelle}, p.182.}
An example of an inductive argument is: "It has been raining for a month now. So it is likely to rain again tomorrow."
\footnote{See Lee \cite{Lee}, pp.6 and 12.}

The third main type of reasoning, abductive reasoning, "infers the truth of the \textit{best explanation} [out of many,] for a set of facts[,] even if that explanation includes unobserved elements ... Diagnoticians and detectives commonly employ abductive reasoning."
\footnote{See Fontainelle \cite{Fontainelle}, p.183.}
An example of an abductive argument is: "If it rains, the grass becomes wet. The grass is wet. So it is most likely that it rained."
\footnote{See Fontainelle \cite{Fontainelle}, p.183.}
However, "the conclusion is probable but not exclusive; someone might have watered the lawn."
\footnote{See Fontainelle \cite{Fontainelle}, p.183.}

The fourth main type of reasoning, analogical reasoning, "transfers information from a particular source to a particular target ... [and] is always preceded by inductive reasoning".
\footnote{See Fontainelle \cite{Fontainelle}, p.183.}
An example of an analogical argument is: 
 
\begin{quotation}
1. Many objects have been observed to share certain characteristics.

2. We induce a class of these objects by their common characteristics, and name it 'apples'.

3. We observe a target [object] which shares characteristics we have found to be typical of [the class named] 'apples'.

4. [Therefore, we] reason analogically that this [target object] is also an apple.
\footnote{See Fontainelle \cite{Fontainelle}, p.183.}
\end{quotation}
 
However, "the conclusion is only probable; it could be a plastic apple."
\footnote{See Fontainelle \cite{Fontainelle}, p.183.}
Moreover, analogical reasoning is subject to the logical fallacy called "the analogical fallacy", which is the false assumption that "because two or more things are similar in one way, they must be similar in other ways [too]".
\footnote{See Fontainelle \cite{Fontainelle},, pp.183 and 210.}
An example of the analogical fallacy is: 
 
\begin{quotation}
1. The universe is like a watch.

2. A watch can give you an itchy wrist.

3. Therefore the universe can give you an itchy wrist.
\footnote{See Fontainelle \cite{Fontainelle}, p.210.}
\end{quotation}

\subsection{The Axiomatic Method is Deductive Logic}

Courant et al.'s \cite{Courant} definition of "the axiomatic method" (written in 1941) is identical to that of deductive logic:

\begin{quotation}
In general terms the axiomatic point of view can be described as follows: To prove a theorem in a deductive system is to show that the theorem is a necessary logical consequence of some previously proved propositions; these, in turn, must themselves be proved; and so on. The process of mathematical proof would therefore be the impossible task of an infinite regression unless, in going back, one is permitted to stop at some point. Hence there must be number of statements, called \textit{postulates} or \textit{axioms}, which are accepted as true, and for which proof is not required. 
\footnote{See Courant et al. \cite{Courant}, pp.214-215.}
\end{quotation}

Courant's "impossible task of an infinite regression" hints at (but fails to clearly state) Russell's key insight (emphasis added): "all argument will be either circular or an endless chain hanging from nothing. \textit{In either case nothing can be proved.}"

Moreover, Courant adds the following criteria for the axioms:
 
\begin{quotation}
The choice of the propositions selected as axioms is to a large extent arbitrary. But little is gained by the axiomatic method unless the postulates are simple and not too great in number. Moreover, the postulates must be \textit{consistent}, in the sense that no two theorems deductible from them can be mutually contradictory, and \textit{complete}, so that every theorem of the system is deductible from them. 
\footnote{See Courant et al. \cite{Courant}, pp.214-215.}
\end{quotation}
 
Courant's consistency requirement is \textit{itself} an axiom. (We can refer to it either as a \textit{meta-axiom}, or as a \textit{default axiom}), This specific axiom is the Law of Non-Contradiction (LNC), which is discussed in greater detail later in this paper. Courant requires this axiom due to yet another axiom, "explosion" / ECQ, which also is discussed in greater detail later in this paper. Inherent to the LNC is another axiom: that there are only two truth-values (true and false).

Moreover, Courant's requirement that "the postulates must be \textit{consistent}, in the sense that no two theorems deductible from them can be mutually contradictory" can be satisfied only in an intuitionistic logic (because it requires an existence proof for every mathematic object). As discussed in this paper, postulates applied to "vacuous subjects" produce contradictory theorems. 

Furthermore, Courant's completeness requirement ("that every theorem of the system is deductible from [the axioms]") requires clarification. Davis \cite{Davis} defined completeness as follows (emphasis in the original):
 
\begin{quotation}
Hilbert asked for a proof that [Peano arithmetic (PA)] is \textit{complete}, meaning that for any proposition that can be expressed in PA, either it can be proved in PA that the proposition is true or it can be proved in PA that the proposition is false.
\end{quotation}
 
According to Davis's more detailed definition, the completeness requirement is rendered impossible by the consistency requirement. This is because, as shown by Gödel in his second incompleteness theorem, there exist propositions that cannot be proven either true or false. Some propositions have  have both truth-values (e.g. the liar paradox, or a proposition with a vacuous subject).

\subsection{LNC as the First Axiom of Aristotelian and Classical Logic}

Boole states that the LNC is (emphasis added): "... that 'principle of contradiction' which Aristotle has described as \textbf{the fundamental axiom of all philosophy}."
\footnote{See Davis \cite{Davis}, p.33, citing Boole \cite{Boole}, p.49.} 
Boole quotes Aristotle's \textit{Metaphysics} as follows (emphasis added):
\footnote{See Aristotle \cite{Aristotle2}, Book IV, Part 3: "For what a man says, he does not necessarily believe; and if it is impossible that contrary attributes should belong at the same time to the same subject (the usual qualifications must be presupposed in this premise too), and if an opinion which contradicts another is contrary to it, obviously it is impossible for the same man at the same time to believe the same thing to be and not to be; for if a man were mistaken on this point he would have contrary opinions at the same time. It is for this reason that all who are carrying out a demonstration reduce it to this as an ultimate belief; \textbf{for this is naturally the starting-point even for all the other axioms}."}
 
\begin{quotation}
It is impossible that the same quality should both belong and not belong to the same thing ... This is the most certain of all principles ... Wherefore they who demonstrate refer to this as an ultimate opinion. \textbf{For it is by nature the source of all the other axioms.}
\end{quotation}
 
Moreover, as further discussed in Cohen \cite{Cohen} (emphasis added):
\footnote{See Cohen \cite{Cohen}, Part 4: "The Fundamental Principles: Axioms".}
 
\begin{quotation}
... Aristotle goes on in Book $\Gamma$ to argue that first philosophy, the most general of the sciences, must also address the most fundamental principles — the common axioms — that are used in all reasoning. Thus, first philosophy must also concern itself with the principle of non-contradiction (PNC): the principle that “the same attribute cannot at the same time belong and not belong to the same subject and in the same respect”.
\footnote{See Aristotle \cite{Aristotle2}, Book IV, Part 3, 1005b19–20.}
\textbf{This, Aristotle says, is the most certain of all principles, and it is not just a hypothesis. It cannot, however, be proved, since it is employed, implicitly, in all proofs, no matter what the subject matter. It is a first principle, and hence is not derived from anything more basic.}
\end{quotation}

The LNC is one of Aristotle's three "Laws of Thought",
\footnote{See Gottlieb \cite{Gottlieb}, and Boole's \cite{Boole}, pp. 48-49, Proposition IV. See also Stillwell \cite{Stillwell}, p.99: "In fact, if \textit{p} + \textit{q} is taken to mean '\textit{p} or \textit{q} but not both,' then the algebraic rules of propositional logic become exactly the same as those of mod 2 arithmetic."
}
and is an axiom or theorem in classical logic 
\footnote{See Gabbay, \cite{Gabbay}, Chapter 2.6.}
(e.g. \textit{Principia Mathematica}),
and is an axiom or theorem in many non-classical logics (e.g. in intuitionism, but not in multi-valued logics).
\footnote{But see Priest et al. \cite{Priest2}: "dialetheism amounts to the claim that there are \textit{true contradictions}."} 
The LNC in sequent form 
\footnote{See Horn \cite{Horn},  Gottlieb \cite{Gottlieb}; Grishin \cite{Grishin}; and Smith \cite{Smith}, \S 11.} 
is: $\vdash \lnot (A \land \lnot A)$. Its verbal characterizations include “opposite assertions cannot both be true simultaneously”, and "no unambiguous statement can be both true and false".
\footnote{See Perzanowski \cite{Perzanowski} p.22, para.4: "The Principle of Non-Contradiction occurs in at least four versions: METAPHYSICAL — no object can, at the same time be and not be such-and-such; LOGICAL — no unambiguous statement can be both true and false; PSYCHOLOGICAL — nobody really and seriously has contradictory experiences, i.e., nobody really sees and does not see (hears and does not hear) simultaneously, etc.; ETHICAL — no one in his right mind would simultaneously demand (or perform) A and not-A."} 
\footnote{An example use of LNC in the context of the RH is found in Edwards \cite{Edwards}, chapter 9, p.202, citing Landau \cite{Landau}, which uses the LNC to prove the theorem that "\textit{if there are only a finite number of exceptions to the Riemann hypothesis, then $S(t)$ cannot be bounded below}".}

According to LNC, a function $f(s)$ of variable $s$ cannot be \textit{both} convergent \textit{and} divergent at \textit{any} value of $s$. Therefore, Riemann's $\zeta(s)$ violates the LNC, because it claims that $\zeta(s)$ is convergent at all values of $s\in \mathbb{C}$ in half-plane $\text{Re}(s)\le1$ (except $s=1$), while also the Dirichlet series "expression" of $\zeta(s)$ is \textit{proven to be divergent} at the same values of $s$. Therefore, in all logics that have LNC as an axiom or theorem, Riemann's $\zeta(s)$ must be false.

Furthermore, in logics that assume the principle of "explosion" (ECQ) (e.g. classical and intuitionistic logics), violation of LNC causes \textit{any} proposition to be "trivially true". 
\footnote{See, e.g. Kleene \cite{Kleene2}, p.101. according to which ECQ ($A, \neg A \vdash B$) is valid in both classical and intuitionistic logics.}
In contrast, paraconsistent bivalent logics reject "explosion" (ECQ), by rejecting the axioms (e.g. disjunctive syllogism and/or disjunction introduction) that lead to ECQ.
\footnote{See Mortansen \cite{Mortansen} and Priest et al. \cite{Priest3}.} In paraconsistent logics, unrelated propositions are no longer "trivially true", but propositions that are directly related to the contradictory proposition remain invalid. (In a paraconsistent bivalent logic, any proposition is false if it assumes that Riemann's analytic continuation of $\zeta(s)$ is true).

The Riemann Hypothesis (RH) states that "all the zeros of $\zeta(s)$ are on the critical line $\text{Re}(s)=0.5$." Because $\zeta(s)$ has no zeros, RH is a proposition with non-existent subjects ("vacuous subjects"). When RH is rephrased as "if $\zeta(s)$ has zeros, then all zeros are on the critical line $\text{Re}(s)=0.5$", the RH is both true and false in classical logic, according to material implication. This result, of being both true and false, violates the LNC.
\footnote{This result of RH being both true and false (a "paradox") is inconsistent with other results, such as Hasse's proof of the RH analogue for elliptic curves of genus 1 (see e.g. Milne \cite{Milne3}, p.3), and Deligne's proof of Weil's conjecture III (see e.g. Milne \cite{Milne3}, p.49). All of these alleged proofs include a violation of the LNC, caused by the analytic continuation of the Zeta function, and the consequently false determinations that the Zeta function has a pole and zeros, that its functional equation is valid, etc. }

Russell's \textit{On Denoting} (which is not a formal logic, but is still relevant to this situation) states that a sentence with a non-existent subject (e.g. the RH) can be interpreted as either a true statement or as a false one. If the RH is interpreted as "there exist zeros of $\zeta(s)$, and it is not the case that any of them are located off of the critical line $\text{Re}(s)=0.5$", then the RH is false, because $\zeta(s)$ has no zeros. However, the alternative interpretation is "it is not the case that there exist zeros of $\zeta(s)$ and any of them are located off of the critical line $\text{Re}(s)=0.5$". This version of RH is true, because it indeed is not the case that there exist zeros of $\zeta(s)$. So according to Russell, the ambiguity of RH means that it can be interpreted as either true or false (and thus is both).
\footnote{
Note that Russell's \textit{On Denoting} assumes that the LEM is true, so it differs from intuitionistic logic.
}

In contrast, some non-classical logics (e.g. multi-valued logics) and philosophical texts (Frege's \textit{Über Sinn und Bedeutung}, Strawson's \textit{On Referring}) reject the LEM, thereby enabling a third state in addition to "true" and "false". For example, Frege's \textit{Über Sinn und Bedeutung} holds that propositions with vacuous subjects (e.g. the RH) lack any truth-value, so they are neither true nor false. Strawson's reasoning in \textit{On Referring} states that questions with "vacuous subjects" (such as the RH) are "absurd" and therefore not asked, thereby inherently creating three truth-values (true, false, absurd), thereby rejecting the LEM.
\footnote{Note: the 160 year history of the RH should be sufficient evidence to refute this argument.}

\subsection{LNC is the Test for Consistency of an Axiomatic System} \label{Consistency}

Langer \cite{Langer} further defines the "axiomatic method" as follows (emphasis added in bold font):
\footnote{See Langer \cite{Langer}, pp.185-186.}

\begin{quotation}
All we ask of a postulate [axiom] is (1) that it shall belong to the system, i.e. be expressible entirely in the language of the system ["\textit{coherence}"]; (2) that it shall imply further propositions of the system ["\textit{contributiveness}"]; \textbf{(3) that it shall not \textit{contradict} any other accepted postulate, or any proposition implied by such another postulate ["\textit{consistency}"]}; and (4) that it itself shall not be implied by other accepted postulates, jointly or singly taken ["\textit{independence}"].
\end{quotation}
 
Langer \cite{Langer} also states that "Contradictory theorems cannot follow from consistent postulates."
\footnote{See Langer \cite{Langer}, p.202.}
Therefore, the LNC is the test for \textit{consistency} of a axiom system. According to Carnap \cite{Carnap}, the LNC is a "sentential formula" that is a tautology.
\footnote{See Carnap \cite{Carnap} p.26: "T8-1. The following formulas are tautologies and hence L-true", followed by two alternate expression of the LEM, and the LNC: (a) $p\lor \sim p$, (b)$\sim p\lor p$, and (c) $\sim (p.\sim p)$.} 
\footnote{See also Carnap \cite{Carnap} p.42: "For suppose it is not raining here now ... E.g. the modal sentence "it is impossible that it is raining and it is not raining" is true, whereas the sentence "it is impossible that it is raining" (produced therefrom by the indicated replacement) is false - for while it is not the case that it is raining here now, this case is nevertheless logically possible. Thus symbolic languages with modality symbols are generally not extensional.} Further according to Carnap (emphasis added in bold font)\cite{Carnap}:
\footnote{See Carnap \cite{Carnap} p.173.}

\begin{quotation}
An [Axiomatic System] AS is said to be \textit{inconsistent} provided that among its theorems is one of the form ${\frakfamily\fraklines\textfrak{S}}_i$ and another of the form $\sim{\frakfamily\fraklines\textfrak{S}}_i$. An AS is said to be \textit{consistent} provided that is not inconsistent. In view of T6-15,
\footnote{See Carnap \cite{Carnap} p.23: "The class comprising the sentential formulas ${\frakfamily\fraklines\textfrak{S}}_i$ and $\sim{\frakfamily\fraklines\textfrak{S}}_i$ L-implies every sentential formula; and likewise the conjunction ${\frakfamily\fraklines\textfrak{S}}_i.\sim{\frakfamily\fraklines\textfrak{S}}_i$" L-implies every sentential formula." This corresponds to "explosion" / \textit{ex contradictione (sequitur) quodlibet} (ECQ).} 
any sentence of the language is derivable from ${\frakfamily\fraklines\textfrak{S}}_i$ and $\sim{\frakfamily\fraklines\textfrak{S}}_i$ together; the theorems of an inconsistent AS therefore include all the sentences of the language L', and the AS in consequence is trivial and useless for practical purposes. \textbf{Consistency is thus an obvious requisite of any non-trivial AS.}
\end{quotation}

Tarski \cite{Tarski} affirms these comments (emphasis added in bold font):
\footnote{See Tarski \cite{Tarski} p.135.}

\begin{quotation}
A deductive theory is called CONSISTENT or NON-CONTRADICTORY if no two asserted statements of this theory contradict each other, or, in other words, if of any two contradictory statements (cf. \S 7) at least one cannot be proved. A theory is called COMPLETE, on the other hand, if of any two contradictory sentences formulated exclusively in the terms of the theory under consideration (and the theories preceding it) at least one sentence can be proved in this theory. Of a sentence which has the property that its negation can be proved in a given theory, it is usually said that it can be DISPROVED in that theory. \textbf{In this terminology we can say that a deductive theory is consistent if no sentence can be both proved and disproved in it[.]}
\end{quotation}

Also Langer \cite{Langer} reaffirms these comments (emphasis added in bold font. Italic font is in the original):
\footnote{See Langer \cite{Langer} p.135.}

\begin{quotation}
The same theorem may follow from more than one possible selection of premises ... But \textbf{\textit{contradictory theorems can never follow from consistent postulates}}. No matter how widely developed the system, how far removed a theorem may be from the original assumptions, they and they only are its ultimate premises; if two theorems in a system are incompatible, and there has been no error in the process of deduction, then the postulates, no matter how obvious and simple they appear, are inconsistent[.]
\end{quotation}

However, contradictory theorems ocassionally \textit{\textbf{do}} arise from consistent postulates. \textbf{These are "paradoxes"}. They arise, for example, from propositions that have "vacuous subjects". 
\footnote{See, e.g. Grattan-Guinness \cite{Grattan-Guinness}, p.338: "Questions of form should be distinguished from those concerning the existence assumptions that have to be abandoned in each case (citing Grattan-Guinness \cite{Grattan-Guinness2}). For example, there is no barber who shaves those and only those who do not shave themselves, thus there is no barber (seemingly Russell's reaction in (citing Russell, equivalent to \cite{Russell4}, p.101)); by contrast, eliminating Russell's paradoxical class affects set theory and logic quite fundamentally, as he was to find for several years to come."}
They also arise from contradictory self-referential statements (e.g. the Liar Paradox).
\footnote{
See, e.g. Grattan-Guinness \cite{Grattan-Guinness}, p.338: "Russell also did not much consider the logical forms of the paradoxes. In Cantor's and Burali-Forti's results [paradoxes of set theory], given the premise \textit{p that} there exists a greatest cardinal or ordinal respectively, opposing conclusions ($c$ and $\sim c$) are deduced about it:
\begin{equation}
p\supset c \text{ and } p \supset \sim c; \therefore \sim p.
\end{equation}
\textit{Reductio ad absurdum} proofs can have this logical structure, sometimes in the condensed form given by $c=p$:
\begin{equation}
p\supset \sim p; \therefore \sim p.
\end{equation}
(This is the version called 'reductio' in \textit{PM}, $\star 2\cdot01$, although without distinction of '$\supset$' from '$\therefore$' - or of \textit{reductio} from the method of indirect proof, which is effected by deducing contradictory consequences from $\sim p$). But with Russell's paradox, from the premise $r$ that his class exists, we deduce the following about the proposition $b$ that it belongs to itself:
\begin{equation}
r\supset.b \supset \sim b \text{ and } r \supset .\sim b \supset b; \therefore r \supset .b \equiv \sim b.
\end{equation}
The differences may be reconciled via \textit{reductio}, so no basic issue arises; in its terms, the paradoxes of the greatest numbers and of naming exemplify the first form while Russell's, the liar and Grelling's take the second."}

Furthermore, the above quotations stress the importance of the LNC, according to majority opinion. Therefore they are derived from a logic that assumes bivalence (most likely a classical logic), and thus ignore multi-valued logics (which tolerate contradictions). This is discussed in greater detail later in this paper. 

Morever, Hilbert's "formalist" program was "to justify classical mathematics by reducing it to a formal system whose consistency should be established by finitistic (hence constructive) means."
\footnote{See Moschovakis \cite{Moschovakis}, 2nd para.} 
\footnote{See also Brouwer \cite{Brouwer3}, p.88: "In the domain of finite sets in which the formalistic axioms have an interpretation perfectly clear to the Intuitionists, unreservedly agreed to by them, the two tendencies differ solely in their method, not in their results; this becomes quite different however in the domain of infinite or transfinite sets, where, mainly by the application of the axiom of inclusion, quoted above, the formalist introduces various concepts, entirely meaningless to the Intuitionist ..."}
At the beginning of the 20th century, Hilbert's "formalist" program "was the most powerful contemporary rival to L.E.J. Brouwer’s developing Intuitionism."
\footnote{See Moschovakis \cite{Moschovakis}, 2nd para.}
According to  Brouwer, "Hilbert was mistaken in claiming that consistency is all that is needed for mathematical existence",
\footnote{See Davis \cite{Davis}, p.95.}
\footnote{See also Brouwer \cite{Brouwer3}, p.90: "Although the formalists must admit contradictory results as mathematical if they want to be consistent, there is something disagreeable for them in a paradox like that of Burali-Forti because at the same time the progress of their arguments is guided by the \textit{principium contradictionis}, i.e., by the rejection of the simultaneous validity of two contradictory properties. For this reason the axiom of inclusion has been modified ..."}
and furthermore, 
 
\begin{quotation}
\textit{to exist} [Brouwer's italics] in mathematics means: to be constructed by intuition; and the question whether a certain language is consistent, is not only unimportant in itself, it is also not a test for mathematical existence.
\footnote{See Davis \cite{Davis}, p.95, fn. 19, citing Brouwer's dissertation \textit{On the Foundations of Mathematics}, in \cite{Brouwer}, p.96.}
\end{quotation}
 
Brouwer also "correctly predicted [Gödel's proof] that any attempt to prove the consistency of complete induction on the natural numbers would lead to a vicious circle."
\footnote{See Moschovakis \cite{Moschovakis}, 2nd para., citing Brouwer's 1912 essay \textit{Intuitionism and Formalism Brouwer}.}

\section{Aristotelian Logic - Axiomatic Method, and Three Laws of Thought}

If Riemann's analytic continuation of Dirichlet series $\zeta(s)$ to half-plane $\text{Re}(s)\le1$ is \textit{true}, then the two contradictory definitions of $\zeta(s)$ in that half-plane (divergent and convergent) violate all three of Aristotle's "Laws of Thought",
\footnote{The LOI, LEM, and LNC.}
for all values of $s$ in half-plane $\text{Re}(s)\le 1$ (except at $s=1$). This is because in Aristotelian classical logic, the connectives $\rightarrow$, $\land$, $\lor$, and $\neg$, are all definable in terms of the others.
\footnote{See Bezhanishvili et al. \cite{Bezhanishvili}  p.3: "Heyting proved that the axioms in Figure 1.1 are independent—none is derivable from the others—and stated that, in contrast to classical logic, in intuitionistic logic none of the connectives $\rightarrow$, $\land$, $\lor$, or $\neg$ is definable in terms of the others (as was proved in Wajsberg \cite{Wajsberg}, McKinsey \cite{McKinsey})."}

Aristotle's "Laws of Thought" are an ancient minimalistic axiomatic system, consisting of three axioms. It is the foundation of traditional logic. 
\footnote{See Fontainelle \cite{Fontainelle}, p.216: "[LEM] does not hold true for multi-valued logics (see page 229) and [LNC] does not hold true when we encounter a paradox (see page 218)."}
As stated by Russell \cite{Russell} at Chapter VII, "On Our Knowledge of General Principles": 
\footnote{See also Boole, \cite{Boole}, which discusses the LOI in Chapter II, pp.34-36, Para.12-13; the LNC in proposition IV, Chapter III, p.49; and the LEM in pp.8 and 99-100, and in proposition II, Chapter III, p.48}
\footnote{See also Brittanica \cite{Brittanica}, citing Dorbolo \cite{Dorbolo}}
\footnote{According to Priest \cite{Priest1} p.139, both LNC and LEM as defined in Aristotle's \textit{Metaphysics}, Book 4, "are not \textit{logical} principles for Aristotle, but \textit{metaphysical} principles, governing the nature of beings \textit{qua} beings. By the time one gets to Leibniz, however, the Laws have been absorbed into the logical canon."}
\footnote{See also Grattan-Guinness \cite{Grattan-Guinness}, p.148: "After stating the identity law as '$x \prec x$' for proposition $x$, Pierce stated  ... that the '\textit{principle of contradiction}' and of '\textit{excluded middle}' were written on p.177 respectively as '$x \prec \overline{x}$' and '$\overline{x} \prec x$'."}
 
\begin{quotation}
Aristotle's "Laws of Thought" are: 

(1) The law of identity [LOI]: 'Whatever is, is.'

(2) The law of contradiction [LNC]: 'Nothing can both be and not be.'

(3) The law of excluded middle [LEM]: 'Everything must either be or not be.'
\end{quotation}

All three axioms are inherited into classical logic, either as axioms or as theorems (e.g. the LNC and LEM are theorems in Russell's  \textit{Principia Mathematica},
\footnote{See Grattan-Guiness \cite{Grattan-Guinness}, p.390: "This theorem [$\star 2 \cdot 11$] was 'the law of excluded middle', a metalaw to us; others of this status included the laws of contradiction and of double negation ($\star 3 \cdot 24$ and $\star 4 \cdot 13$ respectively)."}
and LOI is also referred to as "Material Equivalence"
\footnote{See Lee \cite{Lee}, pp.193-194, 251} 
or "Leibniz's Law" 
\footnote{See e.g., Grattan-Guinness \cite{Grattan-Guinness}, p.447: "[Ramsey] called his new primitive notion a 'function in extension', symbolized '$\phi_e$'; under it and the interpretation of quantification the Leibnizian form of identity
\begin{equation}
x=y := .'(\phi_e).\phi_e x \equiv \phi_e y'
\end{equation}
was acceptable, for it covered all possible associations of proposition and individual and so would be a tautology if $x$ were identical with $y$ and a contradiction otherwise (citing Ramsey \cite{Ramsey}."}). In regards to predicate calculus,  Lemmon \cite{Lemmon} states that:
 
\begin{quotation}
Other theorems, corresponding at the predicate calculus level to the laws of non-contradiction (37),
\footnote{See Lemmon \cite{Lemmon}, p.50, wherein the LNC is proven based on \textit{Reductio ad Absurdum} (RAA) that is discussed in Lemmon \cite{Lemmon}, pp.26-27 and 39-40.}
identity (38),
\footnote{See Lemmon \cite{Lemmon}, p.51, wherein the LOI is proven based on the \textit{Rule of Conditional Proof} (CP) that is discussed in Lemmon \cite{Lemmon}, pp.14-18 and 39-40.}
and excluded middle (44), 
\footnote{See Lemmon \cite{Lemmon}, p.52, wherein the LEM is proven based on the rules (including the RAA) that are discussed in Lemmon \cite{Lemmon}, pp.39-40.}
all of whose proofs are easy, are:

129. $\vdash (x) \neg (Fx \& \neg Fx)$;

130. $\vdash (x) (Fx \rightarrow Fx)$;

131. $\vdash (x) (Fx \lor \neg Fx)$.

As in the propositional calculus, theorems here may be thought of as conveying logical truths, propositions true simply on logical grounds.
\end{quotation}

Some non-classical logics, such as intuitionistic logics (including minimal logic), have the LOI and LNC, but reject the LEM when a proposition has neither been proved or disproved. The LNC is also the test for simple consistency of a propositional calculus.
\footnote{See Carnap \cite{Carnap} p.173: "Consistency is thus an obvious requisite of any non-trivial [Axiomatic System]." See also Tarski et al. \cite{Tarski2}, p.28: "If [theorem] $T$ is inconsistent, two sentences $\Phi$ and $\neg \Phi$ are valid in $T$". See also Tarski et al.'s \cite{Tarski2} example in pp.46-47.}

The three "Laws of Thought" originate in Aristotle's \textit{De interpretatione}, the second of the six texts in Aristotle's \textit{Organon}.
\footnote{See Aristotle \cite{Aristotle} \textit{De interpretatione}.} 
However, the Law of Non-Contradiction (LNC), which is also discussed elsewhere in Aristotle's works, 
\footnote{See Gottlieb \cite{Gottlieb}, citing Aristotle's \textit{Metaphysics} IV (Gamma) 3–6, especially 4; \textit{De Interpretatione}; and \textit{Posterior Analytics} I, chapter 11.} 
is older. The LNC is discussed elsewhere, such as in Euclid (which is approximately contemporaneous with Aristotle),
\footnote{See Hardy \cite{Hardy3}, p.19 : "The proof [of the existence of an infinity of prime numbers] is by \textit{reductio ad absurdum}, and \textit{reductio ad absurdum}, which Euclid loved so much, is one of a mathematician's favourite weapons. It is a far finer gambit than any chess gambit: a chess player may offer the sacrifice of a pawn or even a piece, but a mathematician offers \textit{the game}."}
in Plato's Socratic dialogues
\footnote{In the Socratic dialogue \textit{Republic}, Plato Socrates states: "It is obvious that the same thing will not do or suffer opposites in the same respect in relation to the same thing and at the same time." See Priest \cite{Priest1} pp.137-138, citing Hamilton \cite{Hamilton}, p.436b.}
\footnote{In the Socratic dialogue \textit{Euthyphro}, Socrates uses the LNC in an argument. See Smith \cite{Smith2}, p.29: "Socrates next contends that if Euthyphro's definition of piety is right, then there must be objects that are contemporaneously pious and impious, since they are loved and hated by the gods at the same time. Euthyphro realizes the absurdity of the proposition and is forced to review his understanding of what it is to be pious." }
(all of which predate Aristotle), including the Socratic dialogue \textit{Parmenides},
\footnote{But Priest \cite{Priest1} interprets the Socratic dialogue \textit{Parmenides} as advocating against the LNC: "Even if all things come to partake of both [the form of like and the form of unlike], and by having a share of both are both like and unlike one another, what is there surprising in that? ... when things have a share in both or are shown to have both characteristics, I see nothing strange in that, Zeno, nor yet in a proof that all things are one by having a share in unity and at the same time many by sharing in plurality. But if anyone can prove that what is simple unity itself is many or that plurality itself is one, then shall I begin to be surprised." See Priest \cite{Priest1} p.138, citing Hamilton \cite{Hamilton}, p.129b,c.} 
\footnote{But Brownstein \cite{Brownstein}, pp.49-50, interprets the same section of \textit{Parmenides} as agreeing with the LNC. Brownstein assumes that \textit{a} is a red circle, \textit{b} is a red square, and \textit{c} is a green circle.
\begin{quotation}
Thus \textit{a} and \textit{b} are qualitatively similar to one another [in color] ... but dissimilar to \textit{c}. ... Thus \textit{a} and \textit{c} are similar to each other [in shape] while both are not similar to \textit{b}. We might describe this situation as one in which objects \textit{a}, \textit{b}, and \textit{c} are both alike and unlike ... Plato makes it clear that he does not regard the kind of situation I have described as an absurdity at all.
\end{quotation}
However, we \textit{cannot} describe this situation as one in which objects \textit{a}, \textit{b}, and \textit{c} are both circles and non-circles, or both red and non-red.}
and is attributed to Parmenides and other "pre-Socratics".
\footnote{See Cohen \cite{Cohen2}, pp.75: "It is Parmenides, (one of the pre-Socratic philosophers in the 5th century BCE) who is credited with originally setting out ... the law of noncontradiction, put also as 'Never will this prevail, that what is not is,' by Plato in \textit{The Sophist}."} 

Bertrand Russell argued that these three "laws" are either axioms or theorems of logic:
\footnote{See Russell \cite{Russell}, Chapter VII: "On Our Knowledge of General Principles".}

\begin{quotation}
'[A]nything implied by a true proposition is true' ... is one of a certain number of self-evident logical principles. Some at least of these principles must be granted before any argument or proof becomes possible. When some of them have been granted, others can be proved, though these others, so long as they are simple, are just as obvious as the principles taken for granted. For no very good reason, three of these principles have been singled out by tradition under the name of 'Laws of Thought.' 
\end{quotation}

Lee \cite{Lee} stated that these "are so general and intuitive that their general forms are accepted as laws of logic."
\footnote{See Lee \cite{Lee}, pp.3-4.} 
Minto \cite{Minto} added that "[i]t is even said that all the doctrines of Deductive or Syllogistic Logic may be educed from them."
\footnote{See also Minto \cite{Minto}, p.29}
However, none of these commentators envisioned the existence of non-classical logics, despite the fact that Aristotle himself discussed the existence of future contingent propositions in \textit{De interpretatione} \S 9, in \textit{Organon}.
\footnote{See Aristotle \cite{Aristotle}.}

(Note: In the "classical logic" of Whitehead and Russell's \textit{Principia Mathematica} ("\textit{PM}"), Aristotle's three "Laws of Thought" are theorems derived from other axioms.
\footnote{See Langer \cite{Langer}, p.305, which states that the LNC is proved in Th. 3.24 of \textit{Principia Mathematica}, where the two famous authors state that "in spite of its fame, we have found few occasions for its use."} 
\footnote{See also Andrews \cite{Andrews2}, p.54, which states that all three of the Laws of Thoughts are theorems in \textit{Principia Mathematica}: The LNC in Th. *3.24, and the LEM in Th. *2.11. and the Principle of Identity in Th. *2.08.} 
\footnote{See also Whitehead and Russell's \cite{Whitehead2} discussion of equivalence and Th. *4.01 on p.115 ("It is obvious that two propositions are equivalent when, and only when, both are true or both are false.").}
\footnote{See also Whitehead and Russell's \cite{Whitehead2} discussion of the 'Law of Identity' on pp. 22-23, 39, and 92-93. In addition, see the discussion of Th. *2.08 on p.99: "I.e. any proposition implies itself. This is called the 'principle of identity' and referred to as 'Id.' It is not the same as the 'law of identity' ('$x$ is identical with $x$'), but the law of identity is inferred from it (cf. *13.15)."})

\subsection{The Law of Identity (LOI)} \label{LOI}

The Law of Identity (LOI) is also known as the "Law of Equivalence", and as "Leibniz's Law". This law is also referred to as the "Identity of Indiscernibles".
\footnote{See Forrest, \cite{Forrest}, which formulates it as: "[I]f, for every property $F$, object $x$ has $F$ if and only if object $y$ has $F$, then $x$ is identical to $y$. Or in the notation of symbolic logic:
$\forall(Fx \leftrightarrow Fy) → x=y$."} 
It is the first of Aristotle's three "Laws of Thought", 
\footnote{See, e.g. Russell \cite{Russell}, at Chapter VII: "On Our Knowledge of General Principles".} 
and is an axiom of classical and intuitionistic propositional logics. 

The Law of Identity (LOI) states that a proposition ($P$) "is the same with itself and different from another". This can be written as $P \equiv P$. In the notation of Whitehead and Russell's \textit{Principia Mathematica}
\footnote{See Langer \cite{Langer}, p.307, sequent (*4.2).}, 
the corresponding propositional logic sequent is:  $\vdash .p \equiv p$.
\footnote{See also Russell \cite{Russell}, Chapter VII: "The law of identity: 'Whatever is, is.'" In the context of the RH, if $\zeta(s)$ is both convergent and divergent at any value of $s$, then it both 'is' and 'is not' divergent there, violating the LOI.} 
Taski's version of "Leibniz's Law" is "\textit{$x=x$ if, and only if, $x$ has every property which $x$ has.}".
\footnote{See also Tarski, \cite{Tarski}, p.56: "Leibniz's Law" can be simplified to "\textit{$x=x$ if, and only if, $x$ has every property which $x$ has.}" See also Tarski, \cite{Tarski}, p.57: "\textit{$y=x$ if, and only if, $y$ has every property which $x$ has, and $x$ has every property which $y$ has.}" This is clearly not the case with the divergent Dirichlet series $\zeta(s)$ and the convergent Riemann $\zeta(s)$ in half-plane $\text{Re}(s)\le1$.}
\footnote{See also Sruton \cite{Scruton}, pp.144-146: "But \textit{what} is identity? Philosophers agree on the following four characteristics: ... (ii) Identity is reflexive: everything is identical with itself: (x)(x = x)".}

Another definition of the LOI, in the context of logical discourse, is that the definition of a proposition must be consistent throughout a logical discourse (e.g., the proof of a mathematical theorem). Changing the definition of a proposition in the course of a logical discourse is "equivocation". Aristotle states that "[t]he identity of subject and of predicate must not be 'equivocal'." 
\footnote{See Aristotle \cite{Aristotle} \textit{De interpretatione} \S 6, in \textit{Organon}.}

In the context of the RH, if the Dirichlet series $\zeta(s)$ and Riemann's $\zeta(s)$ are both \textit{true}, then the LOI is violated in half-plane $\text{Re}(s)\le1$ (except at $s=1$), because the two different definitions of $\zeta(s)$ produce \textit{two different values} of $\zeta(s)$ at each value of $s$. In that half-plane, a divergent $\zeta(s)$ is \textit{not equivalent} to a convergent $\zeta(s)$. So if both definitions are true, $\zeta(s)$ \textit{is not equivalent to itself}.

Therefore, $\zeta(s)$ cannot have a plurality of definitions that produce more than one values of $\zeta(s)$ at the same value of $s$. This would be "equivocation", and would mean that $\zeta(s)$ is different from itself ($\zeta(s) \ne \zeta(s)$), thereby violating the LOI.

\subsection{The Law of the Excluded Middle (LEM)} \label{LEM}

The Law of the Excluded Middle (LEM) is another of the three Aristotelian "Laws of Thought". It states that every proposition is either true or false, and thus cannot be both (hence the "excluded middle"). Another interpretation of the LEM that \textit{only one} of a proposition $p$, and its negation $\neg p$, is true ($p \lor \neg p$). 
\footnote{See Aristotle \cite{Aristotle} \textit{De interpretatione} \S 9, in \textit{Organon}.} 

According to Aristotle:
\footnote{See Aristotle \cite{Aristotle} \textit{De interpretatione} \S 9, in \textit{Organon}.}

\begin{quotation}
In the case of that which is or which has taken place, propositions, whether positive or negative, must be true or false. Again, in the case of a pair of contradictories, either when the subject is universal and the propositions are of a universal character, or when it is individual, as has been said, one of the two must be true and the other false[.]
\end{quotation}

The sequent of the LEM is written as: $\forall P \vdash (P \lor \neg P)$. Counter-intuitively, the truth table of the logical disjunction "$\lor$" is that of the Boolean "Inclusive OR", not that of the Boolean "Exclusive OR (XOR)". 
\footnote{See Aloni \cite{Aloni}; and Horn \cite{Horn}, \S 2: "LEM and LNC".} 
(So only in logics that have \textit{both} the LEM \textit{and} the LNC is the middle indeed excluded). 

The LEM is rejected by some non-classical logics, such as intuitionistic logics, \footnote{See Moschovakis \cite{Moschovakis}.} by multi-valued logics (e.g. 3VL), and also by the informal logics described in Frege's \textit{Über Sinn und Bedeutung}, Strawson's \textit{On Referring}, and Russell's \textit{On Denoting}.

In the context of the RH, if both the Dirichlet series $\zeta(s)$ and Riemann's $\zeta(s)$ are true, the LEM is violated throughout half-plane $\text{Re}(s)\le1$, (except at $s=1$"), because then \textit{both} the proposition $p$ ("$\zeta(s)$ is divergent") and its negation $\neg p$ ("$\zeta(s)$ is convergent") are true. 
Under one interpretation, the disjunction ("or") in the LEM is non-exclusive. 
\footnote{See Aloni \cite{Aloni}, 1st para.: "In logic, disjunction is a binary connective ($\lor$) classically interpreted as a truth function the output of which is true if at least one of the input sentences (disjuncts) is true, and false otherwise."}
So under this interpretation, the LEM merely states that propositions $p$ and $\neg p$ cannot \textit{both} be false. However, even with this stricter interpretation of LEM, Riemann's $\zeta(s)$ still violates the LEM in classical logic, dueto its Law of Double Negation Elimination ($\neg (\neg p)$ = p). 

If we assume both definitions of $\zeta(s)$ are true ($p$ = divergent, and $\neg p$ = convergent), and then negate them both, then both $q = \neg (p)$ = convergent, and $\neg q = \neg (\neg p)$ = divergent.  Also, both $q$ and $\neg q$ are false, because they are negations of true statements. This result \textit{is undeniably} a violation of the LEM. In classical logic, due to its Law of Double Negation Elimination ($\neg (\neg p)$ = p), we violate the LEM with the original two propositions! 
\footnote{Because $\neg q = \neg (\neg p)$ and $\neg (\neg p) = p$, therefore $\neg q = p$. So both $q$ and $\neg q$ being false is the same as $p$ and $\neg p$ being false.} 

\subsection{The Law of Non-Contradiction (LNC)} \label{LNC}

The Law of the Non-Contradiction (LNC) is the the third axiom of Aristotle's "Laws of Thought". The LNC states that a proposition ($P$) and its negation ($\neg P$) cannot both be true simultaneously. One expression of this law 
\footnote{See Horn \cite{Horn},  Gottlieb \cite{Gottlieb}; Grishin \cite{Grishin}; and Smith \cite{Smith}, \S 11.} 
is the sequent: $\forall P \vdash \neg(P \land \neg P)$.

Another expression of LNC is that "no unambiguous statement can be both true and false."
\footnote{See Perzanowski \cite{Perzanowski} p.22, para.4: "The Principle of Non-Contradiction occurs in at least four versions: METAPHYSICAL — no object can, at the same time be and not be such-and-such; LOGICAL — no unambiguous statement can be both true and false; PSYCHOLOGICAL — nobody really and seriously has contradictory experiences, i.e., nobody really sees and does not see (hears and does not hear) simultaneously, etc.; ETHICAL — no one in his right mind would simultaneously demand (or perform) A and not-A."} 
Yet another version is that one of a proposition ($P$), or its negation ($\lnot P$), is true. 
\footnote{See Langer \cite{Langer}, pp.262-283, and 300: "Any proposition is either true or false".  This version is true only in logics that assume the LEM. It is a major issue for the intuitionists, and fails in multi-valued logics.} According to Aristotle:

\begin{quotation}
A simple proposition is a statement, with meaning, as to the presence of something in a subject or its absence, in the present, past, or future, according to the divisions of time.
\footnote{See Aristotle \cite{Aristotle} \textit{De interpretatione} \S 5, in \textit{Organon}.} 

An affirmation is a positive assertion of something about something, a denial a negative assertion ... Those positive and negative propositions are said to be contradictory which have the same subject and predicate.
\footnote{See Aristotle \cite{Aristotle} \textit{De interpretatione} \S 6, in \textit{Organon}.}

We see that in a pair of this sort both propositions cannot be true[.]
\footnote{See Aristotle \cite{Aristotle} \textit{De interpretatione} \S 7, in \textit{Organon}.}
\end{quotation}

In the context of the RH, the LNC is violated if the Dirichlet series $\zeta(s)$ and Riemann's $\zeta(s)$ are both true, because then $\zeta(s)$ has two \textit{contradictory} values (divergence and convergence) at all values of $s$ in half-plane $\text{Re}(s)\le1$ (except at $s=1$).
\footnote{See Carnap \cite{Carnap} p.173: "Consistency is thus an obvious requisite of any non-trivial [Axiomatic System]." See also Tarski et al. \cite{Tarski2}, p.28: "If [theorem] $T$ is inconsistent, two sentences $\Phi$ and $\neg \Phi$ are valid in $T$". See also Tarski et al.'s \cite{Tarski2} example in pp.46-47.}

\subsection{Aristotle's Laws of Thought, Applied to the Zeta Function}

When applied to the Zeta function $\zeta(s)$, the LOI holds that $\zeta(s)$ cannot have two different values at any value of $s$, 
\footnote{Note: Also according to the formal definition of a function, $\zeta(s)$ cannot have two different values at any value of $s$.}
because this would mean that the proposition "$\zeta(s)$" is not equal to itself  ($P\not\equiv P$).

The LNC is more specific. It states that a proposition $P$ and its contradiction $\neg P$ cannot \textit{both} be true simultaneously. Using the function $\zeta(s)$ as an example, $\zeta(s)$ cannot be both convergent and divergent at the same value of $s$, because this would mean that proposition $P$ and its negation $\neg P$ were both true. 

So given that the Dirichlet series $\zeta(s)$ is \textit{proven} to be divergent throughout half-plane $\text{Re}(s)\le1$, the LOI and the LNC hold that Riemann's $\zeta(s)$ cannot be valid at any value of $s$ in half-plane $\text{Re}(s)\le1$ (except at $s=1$), nor can any other analytic continuation of $\zeta(s)$.
\footnote{See Carnap \cite{Carnap}, p.18: "A sentential formula is said to be \textit{L-false} (or logically false, or contradictory) in case its range is the null range, i.e. it is false for every value-assignment. Every L-false sentence is evidently false; moreover, its falsity resides entirely in the sense of the sentence and is independent of the facts."} 

So, according to all classical and intuitionistic propositional logics that have LOI and LNC as axioms or theorems, $\zeta(s)$ is defined exclusively by the Dirichlet series (which has no zeros). This means that the zeros of the Riemann Hypothesis (RH) \textit{do not exist}. The non-existent zeros are "vacuous subjects" of a proposition, like "the present King of France" in Bertrand Russell's famous proposition: "The present King of France is bald".
\footnote{See Russell \cite{Russell2}, pp.483-485 and 490.}

\subsection{LNC and the Two Contradictory Zeta Functions }

Riemann's version of $\zeta(s)$ violates the LNC.

In logic, the law of identity is the first of the three classical laws of thought. It states that "each thing is the same with itself and different from another". ... In logical discourse, violations of the Law of Identity (LOI) result in the informal logical fallacy known as equivocation.

If analytic continuation of $\zeta(s)$ is true, that the alternative version of $\zeta(s)$ is convergent for all $s \in \mathbb{C}$, $s\ne1$, then $\zeta(s)$ is \textit{both} convergent \textit{and} divergent throughout half-plane $\text{Re}(s)\le1$, where Riemann's $\zeta(s)$ and the Dirichlet series $\zeta(s)$ disagree. The sole exception is the pole at $s=1$, where both the Dirichlet series $\zeta(s)$ and Riemann's $\zeta(s)$ agree on divergence. 

In other words, analytic continuation of $\zeta(s)$ claims that throughout half-plane $\text{Re}(s)\le1$, for all $s \in \mathbb{C}$ (except $s=1$), $s\ne1$, both a proposition ($P$) and its negation ($\neg P$) are simultaneously true. So this claim contradicts the LNC,  which states that a proposition ($P$) and its negation ($\lnot P$) cannot both be true simultaneously ($\lnot (A \land \lnot A)$). Thus the LNC and Riemann's $\zeta(s)$ cannot both be true. 

The Law of Non-Contradiction (LNC) is "derivable in classical as well as in intuitionistic constructive propositional calculus",
\footnote{See Grishin\cite{Grishin}.} 
so Riemann's $\zeta(s)$ violates the LNC in \textit{both} the classical and the intuitionistic schools of propositional logic.
So in both of these logics, the LNC and proof of Dirichlet series $\zeta(s)$ divergence in half-plane $\text{Re}(s)\le1$ together are sufficient to falsify Riemann's version of $\zeta(s)$.

\subsection{ECQ is a Medieval Addition to Aristotelian Logic}

In Aristotelian logic, if the LNC is violated, the result is \textit{ex contradictione (sequitur) quodlibet} ("ECQ"), which is also called the "principle of explosion". This is the law that \textit{any} proposition can be proven from a contradiction. So due to ECQ, any argument containing a contradiction is "trivially true". 

In the context of the RH, this means that analytic continuation of $\zeta(s)$ violates the LNC and triggers ECQ, because the Dirichlet series $\zeta(s)$ is proven to be divergent there. If assumed to be true, this so-called analytic continuation of $\zeta(s)$ triggers ECQ ("explosion"). 
\footnote{See e.g. Gelbart et al. \cite{Gelbart}, Abstract: "we describe the two major methods for proving the analytic continuation and functional equations of $L$-functions: the method of integral representations, and the method of Fourier expansions of Eisenstein series."}
\footnote{See also Gelbart et al. \cite{Gelbart}, p.78, which states: 
\begin{quotation}
To analytically continue $\zeta(s)$, basically 'the constant term' is enough: reading through the spectral proof of the analytic continuation of $\phi(s)$ for $E(z,s)$, one demonstrates that $\xi(s)$ is holomorphic everywhere, save for simple poles at s = 0 and 1.
\end{quotation}
However, analytic continuation of $\zeta(s)$ violates the LNC and thus is false, and Dirichlet series $\zeta(s)$ has neither poles nor zeros.} 
\footnote{Moreover, $L$-functions are generalizations of the Riemann $\zeta(s)$ function (whose analytic continuation violates LNC). So, analytically-continued $L$-functions violate LNC, as do the arguments that assume that analytically-continued $L$-functions are true, e.g. those described in Gelbart et al. \cite{Gelbart}, p.65: "The Dirichlet $L$-functions $L(s, \chi)$ satisfy the properties \textbf{E}, \textbf{BV}, and \textbf{FE} analogous to those of $\zeta(s)$ (which corresponds to the trivial character)", citing Davenport \cite{Davenport}. See also Gelbart et al. \cite{Gelbart}, pp.60-61, for the definitions of properties \textbf{E}ntirety (\textbf{E}), \textbf{V}ertical strips (\textbf{BV}), and \textbf{F}unctional \textbf{E}quation (\textbf{FE}).}

\section{Classical Logics} 

\subsection{Definition}

Whitehead and Russell's \textit{Principia Mathematica} is referred to as "\textit{the}" classical logic.
\footnote{See F. E. Andrews \cite{Andrews2}, p.54, footnote 3: "In this century the logic of \textit{Principia Mathematica} [henceforth PM] has so succeeded that it is now called "Classical logic"".} 
\footnote{See also Priest \cite{Priest4}, p.xvii, "Around the turn of the twentieth century, a major revolution occurred in logic. Mathematical techniques of a quite novel kind were applied to the subject, and a new theory of what is logically correct was developed by Gottlob Frege, Bertrand Russell and others. This theory has now come to be called \textit{the} 'classical logic'. The name is rather inappropriate, since the logic has only a somewhat tenuous connection with logic as it was taught and understood in Ancient Greece or the Roman Empire. But it is classical in another sense of that term, namely standard."} Other examples of classical logic include George Boole's algebraic reformulation of Aristotelian logic,
\footnote{See Boole \cite{Boole}, especially Propositions III and IV on pp. 48-49, that correspond to the LEM and LNC, respectively.} 
and the second-order logic found in Gottlob Frege's \textit{Begriffsschrift} (when applied to "judgable content").
\footnote{See Lotter \cite{Lotter}, \S 3a: "Frege's early semantics is based on the notion of a conceptual content, that is, it is based on that part of meaning that is relevant for logical inferences. The class of conceptual contents in turn is divided up into judgable and non-judgable ones, whereby the former are logically composed of and can be decomposed into the latter. What Frege may have had in mind – although he does not put it exactly this way – with his distinction between judgable and non-judgable contents is the following consideration: a judgable content is such that we can reasonably either affirm or deny it"}
\footnote{Note: Riemann \cite{riemann1859number} was published in 1859. Riemann died in 1866. Sigwart's work was published in 1873, Frege \cite{Frege2} in 1892, Russell \cite{Russell2} in 1905, and Russell's \textit{Principia Mathematica} in 1910-1913. So Riemann had no knowledge of any of these before his death. In contrast, Boole \cite{Boole} discusses the LOI in Chapter II, pp.34-36, Para.12-13; the LNC in proposition IV, Chapter III, p.49; and the LEM in pp.8 and 99-100, and in proposition II, Chapter III, p.48. Boole \cite{Boole} was published in 1854, a few years before Riemann's 1859 paper, but Riemann does not appear to have been aware of it or its implications.}

"Classical logics" are logics that assume the following as axioms or theorems:
\footnote{See Wikipedia \cite{Classical}, citing Gabbay, \cite{Gabbay}, Chapter 2.6. See also Lee \cite{Lee}, p.251.}

\begin{tabular}{ |p{4cm}||p{3cm}|p{3cm}|p{3cm}|  }
 \hline
 \multicolumn{4}{|c|}{Axioms of Classical Logic} \\
 \hline
 Name & Synonym & Sequent & Description\\
 \hline  \hline
Law of Non-Contradiction (LNC) &  & $\lnot (p \land \lnot p)$ & \\
\hline
Principle of Explosion  & \textit{Ex Contradictione Quodlibet} (ECQ) & $\forall p,  \forall q$: 
$(p \land \lnot p)\vdash q$ & \\
 \hline
 Law of the Excluded Middle (LEM) & & $p \lor \lnot p$ &  \\
 \hline  \hline
Double Negation (DN) & Double Negative Elimination & $p \equiv \lnot \lnot p$ &  \\
 \hline
Monotonicity of Entailment & Weakening & $p \vdash q$ &  Adding presumption $a$ results in $p, a \vdash q$\\
 \hline
Idempotency of Entailment & Contraction & $p, p, a \vdash q$ &  Deleting one of presumptions $p$ results in $p, a \vdash q$\\
  \hline
Commutativity (Com) of Conjunction &  & $(p \land q)\equiv (q \land p)$ & \\
  \hline
De Morgan's Duality (DeM) &  &  $\lnot(p \land q) \equiv (\lnot p \lor \lnot q)$ & Every logical operator is dual to another\\
  \hline
DeM continued &  &  $\lnot(p \lor q) \equiv (\lnot p \land \lnot q)$& \\
  \hline  \hline
\end{tabular}

Also, most semantics of classical logic are \textit{bivalent}, meaning all of the possible denotations of propositions can be categorised as either true or false.
\footnote{See Wikipedia \cite{Classical}, citing Gabbay, \cite{Gabbay}, Chapter 2.6. See also Lee \cite{Lee}, p.251.}
Any higher-order logic that is based on a "classical logic" inherits all of these properties, \textit{in addition to} the three "Laws of Thought".
\footnote{See Sakharov \cite{First-Order-Logic}: "The set of axiom schemata of first-order predicate calculus is comprised of the axiom schemata of propositional calculus together with the two following axiom schemata."} 
\footnote{See also Andrews \cite{Andrews} p.201: "So far we have been concerned with first-order logic, and its subsystem propositional calculus, which we might regard as zeroth-order logic."} 
\footnote{See also Kleene \cite{Kleene}, p.74: "The predicate calculus includes the propositional calculus."}

In classical propositional logics, the three "Laws of Thought" can be either axioms or theorems. For example, Kleene lists all three of the "Laws of Thought" as axioms of a classical propositional calculus. \footnote{See Kleene \cite{Kleene}, p.8: "Now we make one further assumption about the atoms, which is characteristic of classical logical logic. We assume that each atom (or the proposition it expresses) is either \textit{true} or \textit{false} but not both."} \footnote{See also Kleene \cite{Kleene}, p.16, formulas *1, *50, and *51.} In contrast, in \textit{Principia Mathematica}, the three "Laws of Thought" are theorems.
\footnote{See Langer \cite{Langer}, p.305, which states that the LNC is proved in Th. 3.24 of \textit{Principia Mathematica}. Whitehead and Russell: "[I]n spite of its fame, we have found few occasions for its use."} 
\footnote{See also Andrews \cite{Andrews2}, p.54, which states that all three of the Laws of Thoughts are theorems in \textit{Principia Mathematica}: The LNC in Th. 3.24, and the LEM in Th. 2.11. and the Principle of Identity in Th. 2.08. (Russell states that the 'Law of Identity' is inferred later in \textit{PM} from the Principle of Identity). See \textit{Principia Mathematica} to *56, Cambridge, 1967, pp. 99, 101, 111.}

\subsection{Relationship Between Math and Logic} \label{math_logic}

Russell hoped to prove that symbolic logic (and more specifically, his version of "classical" logic) is "practically identical" to mathematics.
\footnote{See Russell's \cite{Russell3} definition of mathematics in p.157, para.106: "This definition brought Mathematics into very close relation to Logic, and made it practically identical with Symbolic Logic."} 
\footnote{See also Scruton \cite{Scruton}, p.77: "... as Russell believed, that mathematics is, in the last analysis, merely logic in another guise."}
But after the publications of Gödel's incompleteness theorems, Russell's "logicism" project had to be abandoned,
\footnote{See Scruton, \cite{Scruton}, p.395: "The final blow to the logicist programme was struck by Gödel, in his famous meta-mathamatical proof that there can be no proof of the completeness of arithmetic which permits a proof of its consistency, and vice versa."}
\footnote{See also Wikipedia \cite{LiarParadox} (citing Crossley et al.  \cite{Crossley}, pp. 52–53): "Roughly speaking, in proving the first incompleteness theorem, Gödel used a modified version of the liar paradox, replacing 'this sentence is false' with 'this sentence is not provable', called the 'Gödel sentence G'."}
and Hilbert's "formalism" project had to be abandoned too.
\footnote{See Scruton \cite{Scruton}, p.395: "It follows too that we cannot treat mathematics as Hilbert wished, merely as strings of provable formulae: the theory of 'formalism' is false."}
Gödel proved that Russell's classical logic is incomplete, by showing that it cannot decide paradoxes such as the Liar's paradox (which is undecidable in a bivalent logic). 
\footnote{See Grattan-Guinness \cite{Grattan-Guinness}, p.512: "Both logicism and formalism now had to be set aside in their current forms, although $PM$ still provided a main source for many basic notions in mathematical logic. However in assuming bivalency, the theorem did not affect intuitionism ... Further, it had no major effect on mathematicians; apart from their general uninterest in foundations, it used a far more formal notion of proof than even their most 'rigorous' practitioners entertained, so that it would not have seemed to bear upon their concerns."}

However, Russell's concept of "logic" was limited to his version of "classical" logic. Over the course of the past century, a wide variety of non-classical logics have been developed, that differ from classical logic.
\footnote{See e.g., Priest's \cite{Priest4} book on non-classical logics.}
So if mathematics is "merely logic in another guise", 
\footnote{Scruton's phrasing of Russell's argument, in \cite{Scruton}, p.77.}
then \textit{which} logic corresponds (or logics correspond) to mathematics? 

Russell assumed that classical logic would be the logic equivalent to all mathematics. This was disproved by Gödel's Incompleteness Theorems. But perhaps the entire body of mathematics is equivalent to some other non-classical logic, for example a 3VL (such as Priest's "Logic of Paradox" $LP$) which rejects LNC and thus the requirement for consistency? Or perhaps the entire body of mathematics is equivalent to the entire inconsistent body of logic? In regards to the last question, Brouwer \cite{Brouwer3} stated:
\footnote{See Brouwer \cite{Brouwer3}, p.84, citing Mannoury \cite{Mannoury}.}
 
\begin{quotation}
To the philosopher or to the anthropologist, but not to the mathematician, belongs the task of investigating why certain systems of symbolic logic rather than others may be effectively projected upon nature. Not to the mathematician, but to the psychologist, belongs the task of explaining why we believe in certain systems of symbolic logic and not in others, in particular why we are averse to the so-called contradictory systems in which the negative as well as the positive of certain propositions are valid.
\end{quotation}

Brouwer's intuitionism "differs from [Russell's] logicism by treating logic as a part of mathematics rather than as the foundation of mathematics".
\footnote{See Moschovakis \cite{Moschovakis}, 2nd para.}

Brouwer's intuitionism also differs from Russell's logicism by treating language as "having nothing to do with mathematics", 
\footnote{See Vafeiadou et al. \cite{Vafeiadou}, p.2, citing Brouwer \cite{Brouwer}, p.79:
"Hilbert's formalist program was doomed to failure because 'language ... is a means ... for the communication of mathematics but ... has nothing to do with mathematics' and is not essential for it." }
whereas logic has always treated mathematical propositions as being a subset of all logical propositions, wherein propositions are linguistic constructs.

Boole showed \footnote{See Boole \cite{Boole}.} that Aristotle's "Laws of Thought" can be represented by an algebra (which is a specific subset of mathematics). By extension, each non-classical logic can have its own corresponding mathematical representation, and consequently, the entire body of logic can be represented by a \textit{subset} of mathematics. The following was Brouwer's opinion: "Far from mathematics being logic (as Frege and Russell had maintained), logic itself is derived from mathematics."
\footnote{See Davis \cite{Davis}, p.95.}
\footnote{See also Curry \cite{Curry}. p.265: "Recent foundational studies (recursive arithmetic, combinatory logic including the theories of lambda conversion, Post's formalized syntax, etc.) show that important theories can be constructed without the aid of any logical calculus, and that these are sufficient for portions of mathematics; so that logic is founded on mathematics, as the intuitionists have long held, rather than the reverse."}

However, historically logic did not originate from mathematics. It originated from philosophy, as method for regulating all arguments, not only mathematical ones. 

If Brouwer is correct, and logic is derived from mathematics, and also mathematics is inconsistent, then we are left with Timon's argument regarding the impossibility of proof. If the many conflicting logics (classical, intuitionistic, multi-valued logic) are derived from the contradictory body of mathematics, are there any "self-evident" general principles? If no, then deductive proof is impossible, because everything will have to be proved by means of something else, and all argument will be either circular or an endless chain hanging from nothing.

Finally, in contrast to Russell, Hilbert, and Brouwer, Wittgenstein defined philosophy (which since Aristotle is defined as including logic) as "all those primitive propositions which are assumed as true without proof by the various sciences".
\footnote{See Wikipedia \cite{Wittgenstein}, citing Klagge et al. \cite{Klagge} p.332, citing Nedo et al. \cite{Nedo} p.89.
} 
According to this argument, the body of logic does not "correspond" to mathematics, and is not derived from mathematics. Instead, the body of \textit{philosophy} (which includes logic) is the foundation of mathematics. The propositions of logic are underlying (and often unstated) assumptions of mathematics. 

According to this view, the laws of logic cannot be determined by mathematical considerations. Instead, they are (and must be) determined according to philosophical considerations. Given that even in ancient Greece there were rival schools of philosophy (and even of logic: e.g. the differences between Aristotle and the Skeptics), it is follows that differing schools of philosophy give birth to different logics, which in turn give birth to different schools of thought in each of the sciences (that are offspring og the different schools of philosophy). 

Also implicit in this argument is that each version of these "primitive propositions" reflects a philosophical worldview ("Weltanschauung"). Especially if Timon is correct in that these "primitive propositions" are not agreed upon, and are impossible to prove. In the context of logics, this would mean that one logic would be selected instead of others because its "absurdities" (or "paradoxes") are held to be less problematic than the "absurdities" (or "paradoxes") of other logics. Such is the case currently with the popularity of the "classical logic" of Russell and Whitehead's \textit{Principia Mathematica}, despite the problems arising from its material implication and ECQ.

However, there are several problems with Wittgenstein's approach. First, it assumes that science is based on deductive reasoning, which starts with philosophy and ends with science. But for the most part, this is not the case. Instead, most areas of science are primarily based on inductive reasoning. Mathematics is the exception, because it does involve a great deal of deductive reasoning. But even mathematics requires other types of reasoning (e.g. inductive, abductive, and analogical). What type of reasoning led to Riemann to propose his famous hypothesis? Surely not deductive reasoning.

Moreover, if Wittgenstein's definition of philosophy is true, \textit{which} subset(s) of the entire body of logic form(s) the foundation of math? Today, there are many logics (classical and non-classical) that contradict one another by least one axiom. 
\footnote{See Priest \cite{Priest4}.}
Also, certain "primitive propositions" are impermissible in certain logics (e.g. "paradoxes" in classical logic, due to LNC), but are permissible in other logics (e.g. paradoxes in MVLs, due to the absence of LNC). So which "primitive propositions" are included, for example, in the foundations of mathematics? Is LNC itself included in the foundations of mathematics?

More specifically, regarding the Riemann Hypothesis (RH), which logic is assumed to be in its foundation? This paper shows that in logics that have the LNC and ECQ as axioms (e.g. classical and intuitionistic logics) 
\footnote{See, e.g. Kleene \cite{Kleene2}, p.101. according to which both classical and intuitionistic logics have ECQ ($A, \neg A \vdash B$) as a theorem.}
the Riemann Zeta Function $\zeta(s)$ violates the LNC, so ECQ renders "trivially true" any proof that assumes Riemann's $\zeta(s)$ is true. So  $\zeta(s)$ is defined by its Dirichlet series, which has no zeros, which means that the Riemann Hypothesis (RH) is directed to an empty set, so both the RH and its negation are both "vacuously true". So RH is an unresolvable paradox in these logics.

In contrast, Priest's "Logic of Paradox" ($LP$), which is Kleene's three-valued logic (3VL) with the third truth-value assigned to paradoxes.
\footnote{See e.g. the following articles on "Logic of Paradox": Priest \cite{Priest5}, Priest \cite{Priest7}, and Hazen et al. \cite{Hazen}.} 
$LP$ enables the RH to be used in logical argument (thanks to $LP$'s rejection of the LEM and the LNC).

This result is consistent with Wittgenstein's argument, because by selecting a foundation logic for the derivation of Riemann's $\zeta(s)$, and for the use of the RH, we are selecting foundation propositions that are assumed to be true.

Moreover, the use of $LP$ as the underlying logic of the RH shows that paradoxes (such as RH) are not "a triviality unworthy of serious consideration",
\footnote{See Priest, \cite{Priest5}, p.219}
or a source of catastrophes (according to ECQ), but instead are an important element in logic and in mathematics. 

\subsection{LNC and Bivalence are Assumed by Gödel's and Tarski's Theorems}

\subsubsection{Gödel's First Incompleteness Theorem}

This use of $LP$ as the underlying logic of RH also renders Gödel's first incompleteness theorem irrelevant. 
\footnote{See Kripke \cite{Kripke}, p.714: "The proof by Gödel and Tarski that a language cannot contain its own semantics applied only to languages without truth gaps.)"}
Gödel's first incompleteness theorem is:
 
\begin{quotation}
\textbf{First Incompleteness Theorem}: "Any consistent formal system F within which a certain amount of elementary arithmetic can be carried out is incomplete; i.e., there are statements of the language of F which can neither be proved nor disproved in F." 
\footnote{Wikipedia \cite{IncompletenessTheorems}, citing Raatikainen \cite{Raatikainen}.}
\end{quotation}

Gödel's first incompleteness theorem has been called a "restatement of the Liar paradox", 
\footnote{See also Wikipedia \cite{LiarParadox} (citing Crossley et al.  \cite{Crossley}, pp. 52–53): "Roughly speaking, in proving the first incompleteness theorem, Gödel used a modified version of the liar paradox".}
and of course the Liar Paradox is neither true nor false (or it is both). This is problematic in \textit{classical logic}, due t LNC and ECQ. In contrast, in $LP$, a third truth-value is assigned to the Liar Paradox (and to all other paradoxes).

So one interpretation of Gödel's first incompleteness theorem is that it is merely a tautology: that paradoxes exist, and that classical logic cannot cope with paradoxes (due to the LNC and ECQ, and the lack of a third truth-value).  In $LP$, Gödel's first incompleteness theorem can be interpreted as another tautology: that there exist propositions that have a third truth-value (neither true nor false).

\subsubsection{Gödel's Second Incompleteness Theorem}

The use of $LP$ as the underlying logic also provides a new interpretation of Gödel's second incompleteness theorem. Gödel's second incompleteness theorem is:
 
\begin{quotation}
\textbf{Second Incompleteness Theorem}: "Assume F is a consistent formalized system which contains elementary arithmetic. Then $ F\not \vdash {\text{Cons}}(F)$." 
\footnote{Wikipedia \cite{IncompletenessTheorems}, citing Raatikainen \cite{Raatikainen}.}
\end{quotation}

In Priest's "Logic of Paradox" ($LP$), Gödel's second incompleteness theorem can be interpreted as a tautology: the canonical consistency statement $\text{Cons}(LP)$ is not provable in $LP$, because $LP$ rejects the LNC, and tolerates inconsistency (i.e. statements with the third truth-value).

Moreover, if $LP$ is indeed the foundational logic underlying the RH problem, does applying the axioms of $LP$ to solve the RH correspond to "adding new rules from 'outside' of number theory in order to solve RH"? No, because the axioms of $LP$ are \textit{inherited} into the axioms of the RH problem. 

\subsubsection{Tarski's Undefinability Theorem}

Tarski's undefinability theorem states that arithmetical truth cannot be defined in arithmetic, and more generally that truth in any sufficiently strong formal system cannot be defined within the system.
\footnote{See Wikipedia, \cite{UndefinabilityTheorem}.}
But the theorem does not prevent truth in that system from being defined in a stronger system.
\footnote{See Wikipedia, \cite{UndefinabilityTheorem}: "For example, the set of (codes for) formulas of first-order Peano arithmetic that are true in $N$ is definable by a formula in second order arithmetic. Similarly, the set of true formulas of the standard model of second order arithmetic (or $n$-th order arithmetic for any $n$) can be defined by a formula in first-order Zermelo–Fraenkel set theory (both ZF and ZFC)."}
The results of this paper, which state that RH (a conjecture in number theory) is false in intutitionistic logic, but is a paradox that triggers ECQ in classical logic, and a pararadox that does \textit{not} cause ECQ in certain 3VLs, is entirely consistent with Tarski's undefinability theorem. The truth-value of the RH is defined by the logical context in which it resides.

As disclosed in Tarski's \textit{The Concept of Truth in Formalized Languages}, 
\footnote{See Tarski, \cite{Tarski3}, p.158, 162.} Tarski uses the Liar paradox in the proof of his Undefinability Theorem, just like Gödel used it in his first incompleteness theorem. The LNC is "Theorem 1" of Tarski's theorem (and the LNC is presented without proof). 
\footnote{See Tarski, \cite{Tarski3}, p.197.}
Instead, it is described as "an almost immediate consequence of [Definitions] 22 and 23."
\footnote{See Tarski, \cite{Tarski3}, p.193 and 195, respectively.}
\footnote{The fact that LNC is a theorem in Tarski's model, rather than an axiom, is not particularly important.}

Therefore, Tarski's Undefinability Theorem is inapplicable in a logic without the LNC, such as a 3VL with truth-value gaps.
\footnote{See Kripke \cite{Kripke}, p.714: "The proof by Gödel and Tarski that a language cannot contain its own semantics applied only to languages without truth gaps.)"}
McGee \cite{McGee} states this, and goes further:
 
\begin{quotation}
Tarski’s analysis leaves open the prospect that we can develop a fully satisfactory theory of truth for a substantial fragment of English; also the prospect that we can develop a theory of truth for English as a whole which, while not fully satisfying our intuitions, is none the less useful and illuminating. Both prospects have been substantially advanced by Saul Kripke’s \cite{Kripke} \textit{Outline of a Theory of Truth}, which exploits the idea that there are truth-value gaps.
\footnote{See McGee \cite{McGee}, which cites Kripke \cite{Kripke}. Kripke's \cite{Kripke}, p.700 states: "One appropriate scheme for handing connectives is Kleene's strong three-valued logic". Footnote 18 on Kripke's \cite{Kripke}, p.700 cites Kleene's \cite{Kleene2} (1952 ed.) description of 3VL in pp.332-340. In the footnote, Kripke states: "'Undefined' is not an \textit{extra} truth-value". So Kripke's use of Kleene's 3VL is similar to Frege's "truth-value gaps", and different from Priest's use of Kleene's 3VL for "truth-value gluts" in \textit{LP}.}
\footnote{See also Kripke \cite{Kripke}, p.711: "So far we have assumed that truth gaps are to handled according to the methods of Kleene. It is by no means necessary to do so. Just about any scheme for handling truth-value gaps is usable, provided that the basic property of the monotonicity of $\phi$ is preserved; that is, provided that extending the interpretation of $T(x)$ never changes the truth-value of any sentence of $L$, but at most gives truth-values to previously undefined cases."}
\end{quotation}

Tarski's reliance on the LNC is reiterated in \textit{Some Observations on the Concepts of $\omega$-Consistency and $\omega$-Completeness} \cite{Tarski4}, where Tarski expressly constructs a symbolical language which, "[i]n spite of its great simplicity ... suffices for the expression of every idea which can be formulated in [Whitehead and Russell's] \textit{Principia Mathematica}."
\footnote{See Tarski \cite{Tarski4}, p.279 and footnote 3.}

Kremer \cite{Kremer} adds:
\begin{quotation}
We had to wait until the work of Kripke \cite{Kripke} and of Martin \& Woodruff \cite{Martin} for a systematic formal proposal of a semantics for languages with their own truth predicates. The basic thought is simple: take the offending sentences, such as [the liar paradox], to be \textit{neither true nor false}. Kripke, in particular, shows how to implement this thought for a wide variety of languages, in effect employing a semantics with three values, \textit{true, false and neither}.
\footnote{Kremer's footnote: "Kripke prefers to treat neither not as a third truth value but as the absence of a truth value." [Author's supplemental footnote: Note that this is Frege's interpretation as well].}
It is safe to say that Kripkean approaches have replaced Tarskian pessimism as the new orthodoxy concerning languages with their own truth predicates.
\end{quotation}

\subsection{The Variety of Non-Classical Logics}

Each non-classical logic is non-classical because it rejects at least one of the axioms of classical logic. 
\footnote{See Wikipedia \cite{Classical}: "Classical logic (or standard logic) is an intensively studied and widely used class of formal logics. Each logical system in this class shares characteristic properties", citing Gabbay, \cite{Gabbay}, Chapter 2.6.}
\footnote{See Priest \cite{Priest4}. See also Sadegh-Zadeh \cite{Sadegh-Zadeh}, p.1030: "Consequently, a large number of such non-classical logics have developed. ... Each of them effectively dismantles the classical logic in a particular way."}
For example, intuitionistic logics reject the Law of the Excluded Middle (LEM), Double Negation (DN), and part of De Morgan's laws.
\footnote{See Wikipedia \cite{Classical}, citing Gabbay, \cite{Gabbay}, Chapter 2.6.}
\footnote{See also Bezhanishvili et al. \cite{Bezhanishvili}  p.4: "From [\textit{intuitionistic propositional calculus}] IPC one obtains a system equivalent to the classical propositional calculus (CPC) used in Principia
by adding any of the following axioms:

$p \lor ¬p$ (excluded middle);

$\neg\neg¬p \rightarrow p$ (double negation elimination);

$((p \rightarrow q) \rightarrow p) \rightarrow p$ (Peirce’s law)."}
\footnote{See also Bezhanishvili et al. \cite{Bezhanishvili}  p.4, fn. 3: "According to Mints \cite{Mints}, p. 701: “Russell anticipated intuitionistic logic by clearly distinguishing propositional principles implying the law of the excluded middle from remaining valid principles. In fact, he states what was later called Peirce’s law."}

Multi-valued logics reject bivalence, allowing for additional truth-values (not just "true" and "false"). Examples of multi-valued logics include three-valued logics (3VL), 
\footnote{Id. Initially developed by Jan Łukasiewicz. Another 3VL, with a slightly different truth table, was developed by Kleene. See the Wikipedia entry on 3VL \cite{3VL}.} 
and infinitely-valued logics ("fuzzy logic") 
\footnote{Id. "Fuzzy logic" permits truth-values to be any Real number between 0 and 1.}

Paraconsistent logics (e.g., relevance logic) reject the Principle of Explosion (ECQ).
\footnote{Id.}
Relevance logic, linear logic, and non-monotonic logic reject monotonicity of entailment;
\footnote{Id.}
Non-reflexive logic (i.e. "Schrödinger logic") rejects or restricts the law of identity.
\footnote{Id., citing da Costa et al. \cite{daCosta}.}

The non-classical logics that reject the LEM include intuitionistic logics (e.g., as formalized by Kleene), multi-valued logics (e.g., 3VL), and Frege's \textit{Begriffsschrift} (when applied to "non-judgable content")
\footnote{See Lotter \cite{Lotter}, \S 3a}.
Higher-order logics that are based on these non-classical logics also reject the LEM.
\footnote{See Sakharov \cite{First-Order-Logic}: "The set of axiom schemata of first-order predicate calculus is comprised of the axiom schemata of propositional calculus together with the two following axiom schemata".} 
\footnote{See also Andrews \cite{Andrews} p.201: "So far we have been concerned with first-order logic, and its subsystem propositional calculus, which we might regard as zeroth-order logic."} 
\footnote{See also Kleene \cite{Kleene}, p.74: "The predicate calculus includes the propositional calculus".}

\section{Paraconsistency and Dialetheism} 

\subsection{Dialetheism Rejects the LNC (for Paradoxes)}

Aristotle introduced the LNC as “the most certain of all principles" ("\textit{firmissimum omnium principiorum}", according to the Medieval theologians). 
\footnote{See Priest et al. \cite{Priest2}, citing Book 
$\Gamma$ of Aristotle's \textit{Metaphysics} (1005b24).
}
"The LNC has been an (often unstated) assumption, felt to be so fundamental to rationality that some claim it \textit{cannot} be defended."
\footnote{See Priest et al. \cite{Priest2}, citing Lewis \cite{Lewis}.
}

Yet "[f]rom the very dawn of Greek thought ... these principles [of LNC] have been contested, first by some rhetoricians and sophists, later on by certain metaphysicists, and recently even by several logicians and mathematicians."
\footnote{See Perzawoski \cite{Perzanowski}, p.22, para.5.} Dialetheism is this view that rejects the LNC, by holding that there exist propositions that are simultaneously true and false (i.e. paradoxes / antinomies).
\footnote{See Priest et al. \cite{Priest2}.}

"As a challenge to the LNC, therefore, dialetheism assails what most philosophers take to be unassailable common sense, calling into question the rules for what can be called into question".
\footnote{See Priest et al. \cite{Priest2}, citing Woods \cite{Woods2003}, Woods \cite{Woods2005}, and Dutilh-Novaes \cite{Dutilh}.}
One of the first logicians to question the status of the logical version of the LNC was Jan Łukasiewicz, father of the Polish school of logic.
\footnote{See Perzawoski 
\cite{Perzanowski}, p.23, para.7.}
In Łukasiewicz's article \textit{On the Principle of Consistency in Aristotle} \cite{Lukasiewicz} \cite{Lukasiewicz2} \cite{Lukasiewicz3},
"Łukasiewicz endorsed only the ethical version of the principle of non-contradiction".
\footnote{See Perzawoski \cite{Perzanowski}, p.23, para.7.} 

Dialetheism argues that \textit{some} propositions are true, \textit{some} are false, and \textit{some} are paradoxes that have a third truth-value. Therefore, dialetheism does \textit{not} reject the LNC for \textit{all} propositions. 
\footnote{But see Beziau \cite{Beziau}. His "trivial dialetheism" argument is based on the false assumption that dialetheism rejects the LNC for \textit{all} propositions, which would indeed reduce 3VL to a single valued logic (i.e. a triviality equivalent to ECQ).}

\subsection{Paraconsistent Logics Accept LNC But Reject ECQ}

A paraconsistent logic rejects "explosion" (ECQ). 
\footnote{See da Costa et al. \cite{daCosta2}, p.1: "It is natural then to put the question whether it is possible to develop a logic in which contradictions can be mastered, in which there are inoffensive or, at least, not dangerous contradictions. The creation of paraconsistent logic by the first author of the present paper (da Costa), more than thirty years ago, brought an affirmative answer to this question. We shall retrace here the history of this invention that has contributed to the subversion of the usual conception of logic."
}
Paraconsistency must be distinguished from dialetheism.
\footnote{But see Priest et al. \cite{Priest3}, citing Asmus \cite{Asmus}.}
"In the literature, especially in the part of it that contains objections to paraconsistent logic, there has been some tendency to confuse paraconsistency with dialetheism (the philosophy that contradictions exist)."
\footnote{Priest et al. \cite{Priest3}}

Paraconsistent logic (logic that rejects ECQ) does \textit{not} entail dialetheism. "Paraconsistency is a property of a consequence relation, whereas dialetheism is a view about truth  ... The fact that one can construct a model where a contradiction holds but not every sentence of the language holds (or where this is the case at some world) does not mean that the contradiction is true \textit{per se}. " 
\footnote{Priest et al. \cite{Priest3}}

The following quotes provide the rationale for paraconsistent logic: 
\footnote{
See McKubre-Jordens \cite{McKubre-Jordens}.
}
\begin{quotation}
Suppose I have proved that the Russell set is and is not a member of itself. Why should it follow from this that there is a donkey braying loudly in my bedroom? ...

The question of relevance (just what has a donkey to do with set theory?) is one that has plagued classical logic for a long time, and is one that makes classical logic a hard pill to swallow to first-time students of logic, who are often told that 'this is the way it is' in logic. Fortunately for those students, paraconsistency provides an alternative.
\end{quotation}

Stanisław Jaśkowski (Łukasiewicz's pupil) produced a paraconsistent logic which "accommodates" paradoxes, and allows for their investigation, without "explosion" (ECQ).
\footnote{See Jaśkowski \cite{Jaskowski}, p.1: "Examples of convincing reasonings which nevertheless yield two contradictory conclusions were the reason why others sometimes disagreed with the Stagirite’s [Aristotle's] firm stand. That was why Aristotle’s opinion was not in the least universally shared in antiquity. His opponents included Heraclitus of Ephesus, Antisthenes the Cynic, and others (cf. Łukasiewicz \cite{Lukasiewicz} \cite{Lukasiewicz2} \cite{Lukasiewicz3}, p. 1). In the early l9th century Heraclitus’ idea was taken up by Hegel, who opposed to classical logic a new logic, termed by him dialectics, in which co-existence of two contradictory statements is possible."}
\footnote{See also Perzanowski 
\cite{Perzanowski2}, p.1: "Any educated person knows, or at least should know, that most cases of incoherences, impossibilities and — in a theoretical framework — paradoxes are rather suspicious members of a domain", and also p.1, fn. 1: "With exceptions of Hegel, Hegelians, etc." See also Perzanowski's  further unflattering comments regarding "Inconsistency believers" in \cite{Perzanowski2}, p.19.}
\footnote{See also Perzanowski 
\cite{Perzanowski}, p.23, para.8.}
"Jaśkowski’s point of departure was a discourse, the situation of a discussion. When one asks: \textit{Is it the case that A?}, and does not know the answer, one often considers both possibilities at once. Likewise, when defending \textit{A}, one respects, at least during a honest discussion, an opponent who claims \textit{not-A}. Which logic applies here?"
\footnote{See Perzawoski 
\cite{Perzanowski}, p.23, para.10.}

"Firstly, [Jaśkowski] created a discursive calculus \textbf{D2}, which fulfilled all the formal criteria we tend to impose on interesting paraconsistent logics. Secondly, his construction in its deep structure enables us to consider inconsistencies [paradoxes] occurring in a theory \textbf{\textit{T}} as contingent statements in a related modal theory \textbf{M(\textit{T})} playing the role of its metatheory. Thirdly, it often allows for the consistent examination of a given inconsistency [paradox]. Sometimes even for the understanding of its mechanism and sources"
\footnote{See Perzanowski 
\cite{Perzanowski}, p.24, para.13: }

According to paraconsistent logics (because they reject ECQ), theorems that assume that Riemann's $\zeta(s)$ is true (e.g. the RH) do not result in "trivial truth". 

One notable subset of paraconsistent logics is that of relevance logics. In relevance logics, "a conditional with a contradictory antecedent that does not share any propositional or predicate letters with the consequent cannot be true (or derivable)."
\footnote{See Wikipedia \cite{Relevance}, citing Routley et al. \cite{Routley} and Mares \cite{Mares}.}

Therefore, according to relevance logics (because they assume the LNC), Riemann's $\zeta(s)$ violates the LNC and does not trigger ECQ for propositions that do not recite, and are unrelated to, the Riemann $\zeta(s)$.

\section{Intuitionistic Logics} 

\subsection{Intuitionistic Logics Reject the LEM (in Regards to Proof)} \label{Intu}

Classical logics assume both the LEM and the LEC, os in these logics, exactly one of ($P$) and ($\neg P$) can be true. 
\footnote{See Plisko \cite{Plisko}; and Stanford \cite{Stanford}.}
This use of LEM together with LNC enables the technique of proving that ($P$) is true, by instead proving that ($\neg P$) is false. This technique is called "\textit{proof by contradiction}". According to classical logic, which has the LEM as an axiom, proof by contradiction is a valid form of proof.

However, in logics that reject LEM (e.g. intuitionistic logics and multi-valued logics), proof by contradiction is \textit{not} a valid form of proof.
\footnote{See Bauer \cite{Bauer}, p.482, \S 1.2: "Proof by contradiction,
or \textit{reductio ad absurdum} in Latin, is the reasoning principle:
\begin{quotation}
\textit{If a proposition $P$ is not false, then it is true.}
\end{quotation}
In symbolic form it states that $\neg \neg P \Rightarrow P$ for all propositions $P$, and is equivalent to excluded middle."  (Note that $\neg \neg P \Rightarrow P$ in classical logic is Double Negation Elimination, which like the LEM is rejected by intuitionist logic). }

Some non-classical logics reject the LEM, and thus also reject proof by contradiction. The intuitionistic school of logic, founded by Brouwer, and formalized by Heyting,
\footnote{See Haack, \cite{Haack}, pp.216-220, citing Brouwer \cite{Brouwer2} and Heyting \cite{Heyting}: "Because he regarded mathematics as essentially mental, and hence thought of mathematical and, \textit{a fortiori}, logical formalism as relatively unimportant, Brouwer didn't give a formal system of the logical principles which are intuitionistically valid. However, intuitionistic logic was formailzed by Heyting, who gives these axioms ...", wherein axiom (10) is ECQ: $\neg p \rightarrow (p \rightarrow q)$. Also, axiom (11) is $((p \rightarrow q) \& (p \rightarrow \neg q) \rightarrow \neg p)$.}
\footnote{So in intuitionistic logic, the problem of "vacuous subjects" and other similar paradoxes do not exist, because the propositions that create them are held to be false. In intuitionistic logic, RH is false, because the zeros are proven to not exist.}
is one of the non-classical schools of logic that reject the LEM in certain instances. According to Moschovakis \cite{Moschovakis}:
 
\begin{quotation}
Intuitionistic logic can be succinctly described as classical logic without the Aristotelian law of excluded middle (LEM) ($A \lor \neg A$) or the classical law of double negation elimination ($\neg \neg A \to A$), but with the law of contradiction  ($A \to B) \to ((A \to \neg B) \to \neg A$) and \textit{ex falso quodlibet} ($\neg A \to (A \to B)$)
\end{quotation}
 
Kleene \cite{Kleene2}, p.120, *51, Remark 1 agrees, proving that "either of $\lnot \lnot A \supset A$ [Principle of Double Negation] or $A \lor \lnot A$ [LEM] can be chosen as the one non-intuitionistic postulate of the classical system."
\footnote{See Kleene \cite{Kleene2}, p.82, Postulate 8 of the "Postulates for the propositional calculus" and the comment regarding "$\circ$" on p.82, the discussion surrounding Postulate $8^I$ on p.101, and Remark 1 on p.120.} 
\footnote{See also Haack \cite{Haack}, p.218: "Heyting's logic lacks some classical theorems; notably, neither $'p \lor \neg p'$, nor $'\neg \neg p \rightarrow p'$, are theorems. However, the double negation of all classical theorems are valid in intuitionistic logic."}
\footnote{But see also Bezhanishvili et al. \cite{Bezhanishvili}  p.4: "From [\textit{intuitionistic propositional calculus}] IPC one obtains a system equivalent to the classical propositional calculus (CPC) used in Principia
by adding any of the following axioms:

$p \lor ¬p$ (excluded middle);

$\neg\neg¬p \rightarrow p$ (double negation elimination);

$((p \rightarrow q) \rightarrow p) \rightarrow p$ (Peirce’s law)."}
(Kleene \cite{Kleene2} also states that the LNC is valid in intuitionistic logic).
\footnote{See Kleene \cite{Kleene2}, p.119, law *50. $\vdash \lnot (A \land \lnot A)$, which is not marked with "$\circ$". See also Kleene \cite{Kleene2}, p., p.101, discussing that ECQ is valid in intuitionistic logics.}
\footnote{See also Haack \cite{Haack}, p.218, axiom (10), which is ECQ: "$\neg p \rightarrow (p \rightarrow q)$". Also, axiom (10) shows that Heyting's intuitionism has the LNC, because there is no ECQ without the LNC.}
\footnote{But see Haack \cite{Haack}, p.218: "Heyting's is not the only, although it is the best entrenched, system of intuitionistic logic: in fact, Johansson's logic [citing \cite{Johansson}], which lacks the tenth axiom [ECQ], has, arguably, a better claim properly to represent the logical principles which are acceptable by intuitionist standards."}

However, the Brouwer–Heyting–Kolmogorov (BHK) interpretation of intuitionistic logic \textit{does} assume the LEM, in the following circumstance: (1) either a proof of proposition ($P$) exists, or (2) an impossibility proof exists for ($P$). 
\footnote{See the discussion of the Brouwer-Heyting-Kolmogorov interpretation in Iemhoff \cite{Iemhoff}, \S 3.1.} 
According to Iemhoff \cite{Iemhoff}:

\begin{quotation}
The BHK-interpretation is not a formal definition because the notion of construction is not defined and therefore open to different interpretations. Nevertheless, already on this informal level one is forced to reject one of the logical principles ever-present in classical logic: the principle of the excluded middle $(A\lor \neg A)$. According to the BHK-interpretation[,] this statement holds intuitionistically if the creating subject knows a proof of A[,] or a proof that A cannot be proved. In the case that neither for A nor for its negation a proof is known, the statement $(A\lor \neg A)$ does not hold. 
\end{quotation}
Further according to Iemhoff \cite{Iemhoff}:

\begin{quotation}
Indeed, there are propositions, such as the Riemann hypothesis, for which there exists currently neither a proof of the statement nor of its negation. Since knowing the negation of a statement in intuitionism means that one can prove that the statement is not true, this implies that both $A$ and $\lnot A$ do not hold intuitionistically, at least not at this moment.  
\end{quotation}

As for the relationship between the LEM and the Riemann hypothesis, the situation is more interesting than as described in the quote above. Chapter \ref{RH} of this paper discusses the relationship between the LEM and the Riemann hypothesis in greater detail.

\subsection{Minimal Logic Rejects Both LEM and ECQ}

One variant of intuitionistic logic is minimal logic. Minimal logic rejects not only LEM, but also ECQ ($\bot \vdash B$). 
\footnote{See Bezhanishvili et al. \cite{Bezhanishvili}, pp.3-4: "An alternative tradition to the formalization of intuitionistic logic, starting with Kolmogorov \cite{Kolmogorov}, leads to a weaker logical calculus, now known as \textit{minimal calculus} [Johansson, \cite{Johansson}]. The distinguishing feature of the minimal calculus is that the formula $\neg p \rightarrow (p \rightarrow q)$, corresponding to the principle \textit{ex falso quodlibet}, is not a theorem. Though the historical debate over the intuitionistic acceptability of \textit{ex falso quodlibet} is interesting, here we focus only on Heyting’s formalization of intuitionistic propositional logic as IPC."}
\footnote{See also Bezhanishvili et al. \cite{Bezhanishvili}, p.3, fn.2: "Kolmogorov’s \cite{Kolmogorov} propositional calculus is in fact equivalent to the implication-negation fragment of minimal calculus (see Plisko \cite{Plisko2})."}
(However, minimal logic does derive a special case of ECQ ($\bot \vdash \neg B$)). Adding ECQ to minimal logic results in intuitionistic logic, and adding the Law of the Excluded Middle (LEM), Double Negation (DN), or Pierce's Law to intuitionistic logic results in classical logic.
\footnote{See Wikipedia \cite{Minimal}, citing Johansson \cite{Johansson} and Troelstra et al. \cite{Troelstra}, p.37.
}
\footnote{See also Bezhanishvili et al. \cite{Bezhanishvili}  p.4: "From [\textit{intuitionistic propositional calculus}] IPC one obtains a system equivalent to the classical propositional calculus (CPC) used in Principia
by adding any of the following axioms:

$p \lor ¬p$ (excluded middle);

$\neg\neg¬p \rightarrow p$ (double negation elimination);

$((p \rightarrow q) \rightarrow p) \rightarrow p$ (Peirce’s law)."}

\section{The Derivation of Riemann's Zeta Function is Not Valid in Logics with LNC}  \label{Riemann-Invalid_2}

\subsection{As Predicted by LNC, the Derivation of Riemann's Zeta Contains Contradictions, and Thus is Invalid}

Moreover, as predicted by the LNC's holding that Riemann's analytic continuation of $\zeta(s)$ is false, the derivation of Riemann's $\zeta(s)$ is invalid in logics with LNC. The derivation uses Cauchy's integral theorem, but contradicts the theorem's prerequisites.
\footnote{See Whittaker et al. \cite{Whittaker}, top of p.87: "If there are two paths $z_0AZ$ amd $z_0BZ$ from $z_0$ to $Z$, and if $f(z)$ is a function of $z$ analytic at all points on these curves and throughout the domain encircled by these two paths, then $\int_{z_0}^Zf(z)\,dz$ has the same value of integration, whether the path of integration is $z_0AZ$ or $z_0BZ$."} 

Riemann used Cauchy's integral theorem to find the limit of the Hankel contour as the Hankel contour approaches the branch cut of $f(s)=\log(-s)$ for $s \in \mathbb{C}$. 

But by definition, all of the points on the branch cut $f(s)=\log(-s)$ have no value at non-negative Real values of $s$. Because they have no value, the function is also non-holomorphic at these points on half-axis $s\ge0$.
\footnote{However, for $s \in \mathbb{C}$, there exists a definition for the branch cut of $f(s)=\log(-s)$ that assigns to it the values of $f(s)=\log(|s|)$ (and remains undefined at $s=0$). This definition contradicts the definition of logarithms of Real numbers. See Encyclopedia of Math \cite{EoM}: "The single-valued branch of this function defined by $\ln (z) = \ln |z| + i \arg (z)$, where $\arg (z)$ is the principal value of the argument of the complex number $z$, $-\pi < \arg(z) \le \pi$, is called the principal value of the logarithmic function."}
\footnote{See Whittaker et al. \cite{Whittaker}, p.244, which states that "by §5.2 corollary 1, the path of integration may be deformed (without affecting the value of the integral) into the path of integration which starts at $\rho$, proceeds along the Real axis to $\lambda$, describes a circle of radius $\lambda$ counter-clockwise round the origin and returns to $\rho$ along the Real axis". The cited "§5.2 corollary 1" appears at the top of Whittaker et al.'s p.87, and is the path equivalence corollary of Cauchy's integral theorem, discussion of which begins on Whittaker et al.'s p.85.}

Moreover, the Hankel contour is either open or closed. In both cases, the Hankel contour violates prerequisites of Cauchy's integral theorem. If it is open (at $s = +\infty$), then it violates the prerequisite of Cauchy's integral theorem that there be \textit{two different paths} connecting two points (in other words, that the contour be closed). If, on the other hand, it is closed (for example at at $x=\infty$),
\footnote{See Whittaker et al. \cite{Whittaker}, p.245: "We shall write $\int_\infty^{(0+)}$ for $\int_C$, meaning thereby that the path of integration starts at 'infinity' on the Real axis, encircles the origin in the positive direction, and returns to the starting point." }
then it encircles the non-holomorphic points of the branch cut, which contradicts another prerequisite of Cauchy's integral theorem (that all points within the contour be holomorphic). 

Therefore, regardless of whether the Hankel contour is interpreted as open or closed, it contradicts prerequisites of Cauchy's integral theorem.
\footnote{See Edwards \cite{Edwards}, pp.10-11.}
\footnote{See Whittaker et al. \cite{Whittaker}, pp.85-87, 244-45 and 266.} This is described in greater detail in Chapter \ref{Riemann-Invalid_2}.

\subsection{The Derivation of Riemann's Zeta Function, Part 1} \label{Deriv_1}

Riemann's version of the Zeta function $\zeta(s)$, that he claims is an alternative "expression" of the Dirichlet series $\zeta(s)$ that "remains valid for all $s$", is derived as follows.
\footnote{See Riemann \cite{riemann1859number} pp.1-2; Edwards \cite{Edwards}, pp.9-11; and Whittaker et al. \cite{Whittaker}, pp.265-266.}

First, Riemann begins with Euler's factorial function (written here in Gauss's notation, as used in Edwards \cite{Edwards}):
\footnote{See Edwards, \cite{Edwards}, p.8, footnote, discussing Legendre's notation.}
\begin{equation}
\prod(s) = \int_{0}^{\infty}e^{-x}x^{s}\,dx 
\end{equation}
The above equation is valid for $s>-1$. So, therefore for $s>0$,
\begin{equation}
\prod(s-1) =  \int_{0}^{\infty}e^{-x}x^{s-1}\,dx 
\end{equation}
Substitution of $nx$ for $x$ in Euler's integral expression
for $\prod(s-1)$ results in:
:
\footnote{See Edwards, \cite{Edwards}, p.9}
\begin{equation}
\prod(s-1) =  \int_{0}^{\infty}e^{-nx}(nx)^{s-1}\,dx
\end{equation}
Extracting the $n^{s-1}$ term from the integral (because $n^{s-1}$ is independent of $x$) results in:
\begin{equation}
\prod(s-1) = n^{s-1} \int_{0}^{\infty}e^{-nx}x^{s-1}\,dx
\end{equation}
Rearranging terms results in:
\begin{equation}
\int_{0}^{\infty}e^{-nx}x^{s-1}\,dx =  \frac{\prod(s-1)}{ n^{s-1}}
\end{equation}
Only if we assume that $n^{s} \approx n^{s-1}$ do we obtain the result used in Riemann's paper:
\footnote{See Edwards, \cite{Edwards}, p.9, and Riemann \cite{riemann1859number}, p.1.}
\begin{equation}
\int_{0}^{\infty}e^{-nx}x^{s-1}\,dx = \frac{\prod(s-1)}{n^s}
\end{equation}
wherein  $(s>0, n=1,2,3,\ldots)$. This error in Riemann's analytic continuation of the Dirichlet series $\zeta(s)$ is minor compared to what is described next.

\subsection{The Derivation of Riemann's Zeta Function, Part 2} \label{Deriv_2}

Next, Riemann takes the last equation of the preceding section,
\footnote{See Riemann \cite{riemann1859number}, p.1.}
\begin{equation}
\int_{0}^{\infty}e^{-nx}x^{s-1}\,dx = \frac{\prod(s-1)}{n^s}
\end{equation}
On the left side of the equation, Riemann uses the equation
\footnote{See Edwards \cite{Edwards}, p.9, footnote, citing Abel and Chebyshev.}
$\sum_{n=1}^{\infty} r^{-n} = (r-1)^{-1}$ to replace $e^{-nx}$ in the integral  with $1/(e^{x}-1)$. On the right side of the equation, Riemann sums the term $1/n^{s}$, from $n = 1$ to $\infty$, thereby obtaining:
\begin{equation}
\int_{0}^{\infty} \frac{x^{s-1}}{e^x-1}\,dx = \prod(s-1) \cdot\sum_{n=1}^{\infty} \frac{1}{n^s}
\end{equation}
By definition, $\zeta(s) = \sum n^{-s}$, so the above equation can be rewritten as:
\begin{equation} \label{Eq1}
\int_{0}^{\infty} \frac{x^{s-1}}{e^x-1}\,dx = \prod(s-1) \cdot \zeta(s)
\end{equation}
Next, Riemann considers the following integral:
\begin{equation}
\int_{+\infty}^{+\infty} \frac{(-x)^{s}}{(e^{x}-1)} \cdot \frac{dx}{x}
\end{equation}
Edwards \cite{Edwards} states: 
\footnote{See Edwards \cite{Edwards}, p.10.}

\begin{quotation}
The limits of integration are intended to indicate a path of integration which begins at $+\infty$ , moves to the left down the positive Real axis, circles the origin once once in the positive (counterclockwise) direction, and returns up the positive Real axis to $+\infty$. The definition of $(-x)^s$ is $(-x)^s = \exp[s\cdot \log(-x)]$, where the definition of $\log(-s)$ conforms to the usual definition of $\log(z)$ for $z$ not on the negative Real axis as the branch which is Real for positive Real $z$; thus $(-x)^s$ is not defined on the positive Real axis and, strictly speaking, the path of integration must be taken to be slightly above the Real axis as it descends from $+\infty$ to $0$ and slightly below the Real axis as it goes from $0$ back to $+\infty$.
\end{quotation}

This is the Hankel contour. 
\footnote{See Whittaker et al. \cite{Whittaker}, pp.244-45 and 266.} 
When written in three terms, with the first term a slight distance above the Real axis as it descends from $+\infty$ to $\delta$, the middle term representing the circle with radius $\delta$ around the origin, and the third term a slight distance below the Real axis as it goes from $\delta$ back to $+\infty$, it is:
\footnote{See Edwards \cite{Edwards}, p.10.}
\begin{equation} \label{Hankel}
\int_{+\infty}^{\delta} \frac{(-x)^{s}}{(e^{x}-1)}\cdot \frac{dx}{x} + \int_{|z|=\delta} \frac{(-x)^{s}}{(e^{x}-1)}\cdot \frac{dx}{x} + \int_{\delta}^{+\infty} \frac{(-x)^{s}}{(e^{x}-1)}\cdot \frac{dx}{x}
\end{equation}
In regards to the middle of these three terms (the circle), Edwards \cite{Edwards} states: 
\footnote{See Edwards \cite{Edwards}, p.10.}

\begin{quotation}
[T]he middle term is $2\pi i$ times the average value of $(-x)^s\cdot (e^{x}-1)^{-1}$ on the circle $|x|=\delta$ [because on this circle $i \cdot d \theta = (dx/x)$]. Thus the middle term approaches zero as $\delta \to 0$ provided $s>1$ [because $x(e^{x}-1)^{-1}$ is nonsingular near $x=0$]. The other two terms can then be combined to give[:]
\end{quotation}

\begin{equation}
\int_{+\infty}^{+\infty} \frac{(-x)^{s}}{e^{x}-1} \cdot \frac{dx}{x} = \lim_{\delta \to 0} \Big[ \int_{+\infty}^{\delta}  \frac {\exp[s(\log x - i\pi)]}{(e^{x}-1)}\cdot \frac{dx}{x} + \int_{\delta}^{+\infty} \frac{\exp[s(\log x + i\pi)]}{(e^{x}-1)}\cdot \frac{dx}{x} \Big]
\end{equation}
resulting in 
\begin{equation}
\int_{+\infty}^{+\infty} \frac{(-x)^{s}}{e^{x}-1} \cdot \frac{dx}{x} = (e^{i\pi s} - e^{-i\pi s})\cdot \int_{0}^{\infty} \frac{x^{s-1}\,dx}{e^{x}-1}
\end{equation}
Since $(e^{i\pi s} - e^{-i\pi s}) = 2i\sin(\pi s)$, this can be rewritten as
\begin{equation}  
\int_{+\infty}^{+\infty} \frac{(-x)^{s}}{e^{x}-1} \cdot \frac{dx}{x} = 2i\sin(\pi s)\cdot \int_{0}^{\infty} \frac{x^{s-1}\,dx}{e^{x}-1}
\end{equation}
Rearranging the terms results in:
\begin{equation}  \label{Eq2}
\int_{0}^{\infty} \frac{x^{s-1}\,dx}{e^{x}-1} = \frac{1}{2i\sin(\pi s)} \cdot \int_{+\infty}^{+\infty} \frac{(-x)^{s}}{e^{x}-1} \cdot \frac{dx}{x}
\end{equation}
The left sides of Equations \ref{Eq1} and \ref{Eq2} are identical, so Riemann equates the right sides of Equations \ref{Eq1} and \ref{Eq2}, resulting in Equation \ref{Eq3}:
\begin{equation}  \label{Eq3}
\int_{+\infty}^{+\infty} \frac{(-x)^{s}}{e^{x}-1} \cdot \frac{dx}{x} = 2i\sin(\pi s)\cdot \prod(s-1) \cdot \zeta(s)
\end{equation}
Then, both sides of the equation are multiplied by $\prod(-s)\cdot s/ 2\pi is$, resulting in
\begin{equation}  
\frac{\prod(-s)\cdot s}{2\pi is} \cdot \int_{+\infty}^{+\infty} \frac{(-x)^{s}}{e^{x}-1} \cdot \frac{dx}{x} =  \frac{\prod(-s)\cdot s}{2\pi is} \cdot 2i\sin(\pi s)\cdot \prod(s-1) \cdot \zeta(s)
\end{equation}
The $s$ terms on the left side cancel out, as do the $2i$ terms on the right side, so
\begin{equation} \label{Eq4}
\frac{\prod(-s)}{2\pi i} \cdot \int_{+\infty}^{+\infty} \frac{(-x)^{s}}{e^{x}-1} \cdot \frac{dx}{x} =  \frac{\prod(-s)\cdot \prod(s-1) \cdot s}{\pi s} \cdot \sin(\pi s)\cdot \zeta(s)
\end{equation}
Next, the identity 
\footnote{See Edwards, \cite{Edwards}, p.8, Eq.5, citing "any book which deals with [the] factorial function or the '$\Gamma$-function', for example Edwards \cite{Edwards2}, pp.421-425." }
$\prod(s) = s\cdot \prod(s-1)$ is substituted into Eq. \ref{Eq4}, resulting in
\begin{equation} \label{Eq5}
\frac{\prod(-s)}{2\pi i} \cdot \int_{+\infty}^{+\infty} \frac{(-x)^{s}}{e^{x}-1} \cdot \frac{dx}{x} =  \frac{\prod(-s)\cdot \prod(s)}{\pi s} \cdot \sin(\pi s)\cdot \zeta(s)
\end{equation}
Finally, the identity
\footnote{See Edwards \cite{Edwards}, p.8, Eq. 6.}
$\sin(\pi s) = \pi s\cdot \Big[\prod(-s)\prod(s)\Big]^{-1}$ is substituted into the right side of Eq. \ref{Eq5}, resulting in
\begin{equation} 
\zeta(s) = \frac{\prod(-s)}{2\pi i} \cdot \int_{+\infty}^{+\infty} \frac{(-x)^{s}}{e^{x}-1} \cdot \frac{dx}{x}
\end{equation}
This is the Riemann Zeta Function.
\footnote{See Edwards \cite{Edwards}, pp.10-11. especially Eq.3.}

\subsection{The Hankel Contour}

In regards to Equation \ref{Hankel} above:
\footnote{See Edwards \cite{Edwards}, pp.10-11. See also Whittaker et al. \cite{Whittaker}, p.244.} 

\begin{equation} 
\label{Hankel2}
\begin{aligned}
\int_{+\infty}^{+\infty} \frac{(-x)^{s}}{(e^{x}-1)} \cdot \frac{dx}{x} = & \int_{+\infty}^{\delta} \frac{(-x)^{s}}{(e^{x}-1)}\cdot \frac{dx}{x} \\
& + \int_{|z|=\delta} \frac{(-x)^{s}}{(e^{x}-1)}\cdot \frac{dx}{x} \\
& + \int_{\delta}^{+\infty} \frac{(-x)^{s}}{(e^{x}-1)}\cdot \frac{dx}{x}
\end{aligned}
\end{equation}

Edwards \cite{Edwards} states:
\footnote{See Edwards \cite{Edwards}, p.10.}

\begin{quotation}
[T]hus $(-x)^s$ is not defined on the positive Real axis and, strictly speaking, the path of integration must be taken to be slightly above the Real axis as it descends from $+\infty$ to $0$ and slightly below the Real axis as it goes from $0$ back to $+\infty$.
\end{quotation}
 
This is the Hankel contour. Riemann copied this solution directly from Hankel's derivation of the Gamma function $\Gamma(s)$. 
\footnote{See Whittaker et al. \cite{Whittaker}, pp.244-245 and 266.}  Riemann uses the Hankel contour in Equation \ref{Eq3}. But what is the basis for equating the branch cut of $f(x)=\log(-x)$ to the limit of the Hankel Contour as the Hankel contour approaches the branch cut? Remember that, by definition, all points on the branch cut have no defined value. Equating the branch cut to the limit of the Hankel contour is a \textit{de facto} assignment of values to points that, by definition, have no value. Remember that in for $x \in \mathbb{R}$, the exponential function $y=\exp{x}$ has \textit{no values} of $x$ which result in $y$ being a non-positive number.

Riemann \cite{riemann1859number} and Edwards \cite{Edwards} fail to provide any reason, much less a mathematically valid reason, for equating the "\textbf{strictly speaking}" interpretation of the "first contour" on the left side of Eq. \ref{Hankel2} (points that, by definition, have no value, and thus are also non-holomorphic: points of the branch cut of $f(x)=\log(-x)$), to the "\textbf{non-strictly speaking}" interpretation of the "first contour" (the Hankel contour that is "slightly above the Real axis as it descends from $+\infty$ to $0$ and slightly below the Real axis as it goes from $0$ back to $+\infty$", on the right side of Eq. \ref{Hankel2}). 

Unlike Riemann \cite{riemann1859number} and Edwards \cite{Edwards}, Whittaker et al. \cite{Whittaker} does provide a reason: \textit{\textbf{the path equivalence corollary of Cauchy's integral theorem}} is given as the basis for equating the Hankel contour to the branch cut.
\footnote{See Whittaker et al. \cite{Whittaker}, p.244, which states that "by §5.2 corollary 1, the path of integration may be deformed (without affecting the value of the integral) into the path of integration which starts at $\rho$, proceeds along the Real axis to $\lambda$, describes a circle of radius $\lambda$ counter-clockwise round the origin and returns to $\rho$ along the Real axis". The cited "§5.2 corollary 1" appears at the top of Whittaker et al.'s p.87, and is the path equivalence corollary of Cauchy's integral theorem, discussion of which begins on Whittaker et al.'s p.85.} 
\textbf{However, this basis is not mathematically valid.} Both the Hankel contour and the branch cut \textit{\textbf{contradict}} the prerequisites of the Cauchy integral theorem, 
\footnote{See Whittaker et al. \cite{Whittaker}, p.85: "If $f(z)$ is a function of $z$, analytic at all points ... inside a contour $C$, then $\int_{(C)} f(z)\,dz = 0$.". The integrated function must be analytic (holomorphic) at all points inside the contour of integration.} 
and therefore \textit{\textbf{also contradict}} the prerequisites of its corollary. 
\footnote{See Whittaker et al. \cite{Whittaker}, top of p.87: "If there are two paths $z_0AZ$ amd $z_0BZ$ from $z_0$ to $Z$, and if $f(z)$ is a function of $z$ analytic at all points on these curves and throughout the domain encircled by these two paths, then $\int_{z_0}^Zf(z)\,dz$ has the same value of integration, whether the path of integration is $z_0AZ$ or $z_0BZ$."} 
\textit{\textbf{These contradictions invalidate the derivation of Riemann's version of $\zeta(s)$}} in logics with LNC.

\subsection{Cauchy's Integral Theorem and Its Path Equivalence Corollary}

\subsubsection{Cauchy's Integral Theorem}

Cauchy's integral theorem states that if function $f(z)$ of complex variable $z$ is holomorphic at all points on a simple closed curve ("contour") $C$, and if $f(z)$ is holomorphic at all points inside the contour, then the contour integral of $f(z)$ is equal to zero: 
\footnote{See Whittaker et al. \cite{Whittaker}, p.85.}
\begin{equation} \label{eq:2.1}
\int_{(C)} f(z)\cdot dz = 0
\end{equation}

\subsubsection{Path Equivalence Corollary}

The path equivalence corollary of Cauchy's integral theorem 
\footnote{See Whittaker et al. \cite{Whittaker}, p.87, Corollary 1.} 
states that: 

(1) if there exist four points $z_0$, $Z$, $A$, and $B$ on the Cartesian plane representing the complex domain, and the two points $z_0$ and $Z$ are connected by two distinct paths $z_0AZ$ and $z_0BZ$ (one path going through $A$, the other path goin through $B$), and 

(2) if function $f(z)$ of complex variable $z$ is holomorphic at all points on these two distinct paths $z_0AZ$ and $z_0BZ$, and $f(z)$ is holomorphic at all points enclosed by these two paths,  

(3) then any line integral connecting the two points $z_0$ and $Z$ inside this region (bounded by $z_0AZ$ and $z_0BZ$) has the same value, regardless of whether the path of integration is $z_0AZ$, or $z_0BZ$, or any other path disposed between $z_0AZ$ and $z_0BZ$. 

\subsection{Prerequisites of Cauchy Integral Theorem are Contradicted}

Riemann used Cauchy's integral theorem to find the limit of the Hankel contour as the Hankel contour approaches the branch cut of $f(x)=\log(-x)$ at $x \in \mathbb{C}$. But  by definition, $\log(-x)$ has no value (and thus is non-holomorphic) on half-axis $x\in \mathbb{R}, x\ge0$. The Hankel contour is either open, or closed, at $x = +\infty$ (the latter enclosing non-holomorphic points). In both cases, the Hankel contour violates prerequisites of Cauchy's integral theorem.

If the Hankel contour is open, the Cauchy integral theorem (which only applies to closed contours) cannot be used. In the alternative, if the Hankel contour is indeed closed at $+\infty$ on the branch cut,
\footnote{See Whittaker et al. \cite{Whittaker}, p.245: "We shall write $\int_\infty^{(0+)}$ for $\int_C$, meaning thereby that the path of integration starts at 'infinity' on the Real axis, encircles the origin in the positive direction, and returns to the starting point." }
then the Hankel contour still contradicts the requirements of the Cauchy integral theorem. This is because the closed Hankel contour encloses the entire branch cut of $f(z)$, and the branch cut consists entirely of non-holomorphic points. Also, there would be a non-holomorphic point on the Hankel contour itself, at the point where it intersects the branch cut at $+\infty$ on the Real axis. These reasons disqualify the use of the Cauchy integral theorem with the Hankel contour.

For these reasons it is not valid to use the Cauchy integral theorem's path equivalence corollary to find the limit of the Hankel contour as the Hankel contour approaches the branch cut of $f(x)=\log(-x)$ at $x \in \mathbb{C}$. So the derivation of Riemann's $\zeta(s)$ violates the LNC. For the same reasons, Hankel's derivation of the Gamma function $\Gamma(s)$ violates the LNC.  
\footnote{See Whittaker et al. \cite{Whittaker}, pp.244-246, \S 12.22 "Hankel's expression of $\Gamma(z)$ as a contour integral", citing Hankel \cite{Hankel}, p.7.}

\subsection{Strictly Speaking, the Points on the Hankel Contour Have No Defined Value}

Further in regards to the Hankel contour of Equation \ref{Hankel}: 
\footnote{See also Whittaker et al. \cite{Whittaker}, p.266.} 

\begin{equation} 
\begin{aligned}
\int_{+\infty}^{+\infty} \frac{(-x)^{s}}{(e^{x}-1)} \cdot \frac{dx}{x} = & \int_{+\infty}^{\delta} \frac{(-x)^{s}}{(e^{x}-1)}\cdot \frac{dx}{x} \\
& + \int_{|z|=\delta} \frac{(-x)^{s}}{(e^{x}-1)}\cdot \frac{dx}{x} \\
& + \int_{\delta}^{+\infty} \frac{(-x)^{s}}{(e^{x}-1)}\cdot \frac{dx}{x}
\end{aligned}
\end{equation}

Edwards \cite{Edwards} states (emphasis added in bold font): 
\footnote{See Edwards \cite{Edwards}, p.10.}

\begin{quotation}
The limits of integration are intended to indicate a path of integration which begins at $+\infty$ , moves to the left down the positive Real axis, circles the origin once once in the positive (counterclockwise) direction, and returns up the positive Real axis to $+\infty$. The definition of $(-x)^s$ is $(-x)^s = \exp[s\cdot \log(-x)]$, where the definition of $\log(-x)$ conforms to the usual definition of $\log(z)$ for $z$ not on the negative Real axis as the branch which is Real for positive Real $z$; \textbf{thus $(-x)^s$ is not defined on the positive Real axis and, strictly speaking, the path of integration must be taken to be slightly above the Real axis as it descends from $+\infty$ to $0$ and slightly below the Real axis as it goes from $0$ back to $+\infty$.}
\end{quotation}

So as $\delta \to 0$  (provided $s>1$), the middle (circular) term disappears, and the hankel contour is represented by the two "linear" terms, as follows:
\begin{equation} 
\int_{+\infty}^{+\infty} \frac{(-x)^{s}}{(e^{x}-1)} \cdot \frac{dx}{x} = \int_{\infty}^{+0} \frac{(-x)^{s}}{(e^{x}-1)}\cdot \frac{dx}{x} + \int_{+0}^{\infty} \frac{(-x)^{s}}{(e^{x}-1)}\cdot \frac{dx}{x}
\end{equation}

However, strictly speaking, the notation misleadingly indicates that both of the two terms on the right side of the equation are located directly on the branch cut of $\log(-x)$ (where by definition these two terms have no value). 

Any assignment of a value to the right side of this equation (as written) contradicts the definition of the logarithm $f(x)=\log(-x)$, which \textit{by definition} has no value for all non-negative Real values of $x$. 

\section{Multi-Valued Logics (Including Three-Valued Logics)} 

\subsection{Multi-Valued Logics Reject the LEM} \label{MVL}

Multi-valued logics reject the LEM, because they are not bivalent. They have at least one truth-value in addition to the two bivalent truth-values ("true" and "false"). For example, Frege's \textit{Über Sinn und Bedeutung} ("On Sense and Denotation") "claimed that an utterance of a sentence containing a non-referring singular is neither true nor false."
\footnote{See Marques, \cite{Marques} p.70, and Frege \cite{Frege2}.} 
More specifically, Frege states the following:

\begin{quotation}
The sentence ‘Scylla has six heads’ is not true, but the sentence ‘Scylla does not have six heads’ is not true either; for it to be true the proper name ‘Scylla’ would have to designate
something.
\footnote{See Marques, \cite{Marques} p.71, citing Frege \cite{Frege}, p.127. "Scylla" refers to the creature from Greek mythology.}
\end{quotation}
 
Therefore, according to Frege's logic, a proposition can have no truth-value, which means that a proposition has three possible states: true, false, or neither. As Marques \cite{Marques} states (emphasis added):
\footnote{See Marques \cite{Marques}, p.71}
 
\begin{quotation}
This gives expression to two natural ideas: i) a sentence such as ‘Scylla does not have six heads’ is the negation of ‘Scylla has six heads’; and ii) ‘Scylla has six heads’ is false if and only if its negation is true (that is, if ‘Scylla does not have six heads’ is true). \textit{When a sentence has no truth-value, the result of embedding the sentence, for instance under the scope of negation, also can have no truth-value.} 
\end{quotation}
 
See also Milne \cite{Milne2}, who states: "Frege holds that any sentence containing a bearerless name in a direct/non-oblique context is neither true nor false."
\footnote{See Milne \cite{Milne2}, p.473, citing Frege's "Logic". Milne's \cite{Milne2} p.474 reproduces Smiley's truth tables for Frege's three-valued logic (Citing Smiley \cite{Smiley}, pp.125-35.). The "third value" in these truth tables is an absence of any truth-value ("a truth-value gap"). Milne's \cite{Milne2} p.474: "Beware! The bar [symbol] is not a third truth-value, it signifies the absence of a truth-value. Where both \textbf{A} and \textbf{B} have truth-values, the connectives behave classically." }
 
\begin{quotation}
Frege holds that \textit{any} sentence containing a bearerless name in a direct/non-oblique context is neither true nor false. ... He terms the thought expressed by such a sentence 'fictitious' and a 'mock thought' ('Logic', p.130); they are such exactly and only in that they fail to be about actually existing objects. In particular he says 'Scylla has six heads' is not true, and 'Scylla does not have six heads' is not true. Lack of a bearer for a singular term spreads lack of truth-value pervasively to logically complex sentences.
\end{quotation}

\subsection{Priest's Three-Valued Logic Rejects LNC and ECQ}

One of the theorems in classical logics, and most non-classical logics (but not multi-valued logics!) is the "Principle of Explosion". In Latin: \textit{Ex Contradictione (Sequitur) Quodlibet} (ECQ): "from contradiction, anything (follows)".

According to this theorem, the result of a contradiction (a violation of LNC) is that any statement whatsoever can be proven. In other words, "a false proposition implies any proposition". 
\footnote{See Langer \cite{Langer}, p.284.} 
So a single contradiction in a theorem results in an "explosion" of false theorems that incorrectly assume the original contraction to be true. 

\subsection{Three-Valued Logics (3VLs) Bypass Aristote's Three Laws of Thought}

Kleene describes Brouwer's intuitionism as follows:  
\begin{quotation}
In 1908, Brouwer, in a paper entitled 'The untrustworthiness of the principles of logic', challenged the belief that the rules of the classical logic, which have come down to us essentially from Aristotle (384--322 B.C.) have an absolute validity, independent of the subject matter to which they are applied.
\footnote{See Kleene \cite{Kleene2}, p. 46.}
\end{quotation}

Brouwer's intuitionism is skeptical of Aristotelian logic, for reasons similar to those of Timon of Phlius, the Pyrrhonist philosopher. 
\footnote{Therefore Heyting's formalist version of intuitionism might be a misinterpretation of Brouwer's core argument.}
However, it must be said that Aristotle himself "expressed reservations about bivalence",
\footnote{See Haack \cite{Haack}, p.204, citing  Aristotle's \cite{Aristotle} \textit{De interpretatione} \S 9, in \textit{Organon}, which includes the famous statement: "A sea-fight must either take place to-morrow or not, but it is not necessary that it should take place to-morrow, neither is it necessary that it should not take place, yet it is necessary that it either should or should not take place to-morrow." The propositions about the future ("it will take place tomorrow" and "it will not take place tomorrow") are best understood according to probability theory (or 3VL, or fuzzy logic), not according to bivalent logic, which cannot assign truth-values to these two propositions. Aristotle correctly states: "One may indeed be more likely to be true than the other, but it cannot be either actually true or actually false." }
a core assumption of Aristotelian logic.
\footnote{See also Pelletier \cite{Pelletier2}, and Pelletier et al. \cite{Pelletier3}}

Regarding three-valued logics (3VLs), which adopt trivalence instead of bivalence, they satisfy both Aristotle's and Brouwer's rationales for a third truth-value. Regarding Priest's three-valued logic ("$LP$"), it satisfies yet another reason for a third truth-value: the need to work with paradoxes such as the liar paradox. $LP$, like all 3VLs, rejects the LEM, but it is unique in also rejecting the LNC and the LOI. This completely rejects Aristote's three laws of thought, satisfying Brouwer's above-cited argument.  
\footnote{Note also that Heyting's intuitionism, which completely (and mechanistically) removes the LEM and Double Negation from classical logic, is not really an accurate representation of Brouwer's above-cited argument for opposing the LEM in regards to the truth-value of propositions that have not been proven or disproven.}

Also, this rejection of LNC and LEM exposes a fundamental flaw in the axiomatic method advocated by Aristotle, according to which, "all deduction had to start, like Euclid, from general principles regarded as self-evident."
\footnote{See Russell \cite{Russell5}, p.234.}
\footnote{See also Mendell \cite{Mendell}, \S 2: "Aristotle's discussions on the best format for a deductive science in the \textit{Posterior Analytics} reflect the practice of contemporary mathematics as taught and practiced in Plato's Academy, discussions there about the nature of mathematical sciences, and Aristotle's own discoveries in logic. Aristotle has two separate concerns. One evolves from his argument that there must be first, unprovable principles for any science, in order to avoid both circularity and infinite regresses."}
\footnote{See also Mendell \cite{Mendell}, \S 2: "Aristotle distinguishes (\textit{Posterior Analytics} i.2) Two sorts of starting points for demonstration, \textbf{axioms} and \textbf{posits}. An \textbf{axiom} (\textit{axiôma}) is a statement worthy of acceptance and is needed prior to learning anything. Aristotle's list here includes the most general principles such as non-contradiction and excluded middle, and principles more specific to mathematicals, e.g., when equals taken from equals the remainders are equal."}
\footnote{See also Lemmon \cite{Lemmon}, pp.173-174: "The main burden of traditional logic is to distinguish, of the 256 possible patterns, which are \textit{valid} and which are \textit{invalid}. Two quite separate approaches are used, which yield the same result ... The second method, which is Aristotle's own, is to accept as valid certain 'self-evident' patterns in the first figure and then, using principles as (1)-(9), to \textit{deduce} the valid patters of the remaining figures. This method is traditionally known as \textit{reduction to the first figure}, and is said to take two forms, direct and indirect reduction. Roughly speaking, in indirect reduction the valid pattern is deduced by [\textit{Reductio ad absurdum}] RAA ..."} 
Rejection of LNC and LEM validates the argument presented by Timon of Phlius, the Pyrrhonist philosopher, that there are no "self-evident" general principles. 

Haack \cite{Haack} provides another famous example of axioms that were not "self-evident" after further consideration:
 
\begin{quotation}
Frege confidently supposed that the principles of his logical system were self-evident, until Russell showed that they were inconsistent! 
\footnote{See Haack \cite{Haack}, p.153.}
\end{quotation}
\begin{quotation}
Frege's response to the discovery of [Russell's Paradox] was to concede that he'd never really thought that the relevant axiom was \textit{quite} as self-evident as the others - a comment which may well induce a healthy skepticism about the concept of self-evidence. 
\footnote{See Haack \cite{Haack}, p.10.}
\end{quotation}

Which leads to the question: What does it mean to claim that some proposition is self-evident? According to Haack \cite{Haack} (who agrees with Timon of Phlius):
 
\begin{quotation}
Presumably, something to the effect that it is obviously true. But once it has been put like this, the difficulty with the concept of self-evidence cannot be disguised. The fact that a proposition is obvious is, sadly, no guarantee that it's true. (It is pertinent that different people, and different ages, find different and even incompatible propositions - that some men are naturally slaves, that all men are equal... - 'obvious'.) 
\footnote{See Haack \cite{Haack}, pp.235-236.}
\end{quotation}
In addition, in regards to the specific example of "self-evident" axioms of Frege's logicisim: 
\begin{quotation}
Whether one says that Frege's inconsistent axioms only \textit{seemed} self-evident, but couldn't really have been, or that they \textit{were} self-evident but unfortunately weren't true, self-evidence must fail to supply an epistemological guarantee; because either (on the latter assumption) a proposition may be self-evident but false, or else (on the former assumption) though if a proposition is self-evident then it is, indeed, true, one has no certain way to tell when a proposition is really self-evident. 
\footnote{See Haack \cite{Haack}, p.236.}
\end{quotation}

\subsection{In 3VLs, Material Implication from Paradoxes Does Not Result in ECQ }

In contrast to classical logic, the  three-valued logics (3VLs) created by Łukasiewicz,
\footnote{See Łukasiewicz \cite{Lukasiewicz4}, which presents the first 3VL, that was "later on criticized by Suszko \cite{Suszko}" but "later used by Asenjo, da Costa and D’Ottaviano and Priest, to develop paraconsistent systems of logic." (See Béziau \cite{Beziau2}, p.25, last para.). A predecessor of paraconsistent and multi-valued logics was Vasiliev. See da Costa et al. \cite{daCosta2} and Bazhanov \cite{Bazhanov}.}
Bochvar,
\footnote{See Bochvar \cite{Bochvar}, and see also Urquhart \cite{Urquhart2}, pp.252-253, \S 1.6: "The work of the Russian logician Bochvar \cite{Bochvar} represents a new philosophical motivation for many-valued logic; its use as a means of avoiding the logical paradoxes. His system introduces the intermediate value \textit{I} in addition to the classical values \textit{T} and \textit{F}. His idea is to avoid logical paradoxes such as Russell's and Grelling's by declaring the crucial sentences involving them to be meaningless (having the value \textit{I})."}
Kleene, 
\footnote{See Kleene, \cite{Kleene3}.}
and Frege (in his \textit{Über Sinn und Bedeutung}), 
\footnote{See Frege, \cite{Frege2}. It has a truth-value gap instead of a third truth-value.}
all reject LEM, as does intuitionism (in a more limited manner),
\footnote{Specifically, in regards to proof and absence of proof.}
thereby allowing propositions (e.g. the RH) to be "\textit{neither} true nor false". 
\footnote{But see Woleński \cite{Wolenski}, \S 3.3, which states that 3VL rejects the LNC (emphasis added): "Sentences about future contingent states of affairs are natural candidates for having the third value ($1/2$). For example, the sentence “I will visit Warszawa next year”, is neither true nor false, it is merely possible and has the value $1/2$. Its negation has the same value. This idea led to three-valued logic ... \textit{This means that the laws of contradiction and excluded middle do not hold in three-valued logic}."}

Further in contrast, Priest's "Logic of Paradox" ($LP$) expressly rejects LNC by allowing for paradoxes, thereby allowing propositions (e.g. the RH) to be \textit{both} true and false. 
\footnote{See Priest \cite{Priest6}. These four possibilities (true, false, both, neither) form the "catuskoti" (or "tetralemma") of early Buddhist logic, which rejects the LNC. }
Moreover, it is proven that in Frege's logic that the state of "\textit{neither} true \textit{nor} false" implies the state of "\textit{both} true and false", and vice versa. 
\footnote{See Milne's \cite{Milne2}, p.475 (citing Heidelberger \cite{Heidelberger}): "Putting that all together we get,
\begin{quotation}
It's not true that $P$ and it's not false that $P$ only if it's both true that $P$ and false that $P$.
\end{quotation}
In short, everywhere we think there’s a truth-value gap, there’s also a‘glut’! (And vice versa!)".}
\footnote{See also Priest \cite{Priest6}, p.27, which recites an explanation different from Heidelberger's \cite{Heidelberger}: "Notably, assuming De Morgan’s laws, ... $\lnot(A \lor \lnot A)$ is equivalent to ... $A \land \lnot A$".}

Also, according to the truth tables of Łukasiewicz, Kleene's and Priest's 3VLs, when $p$ has the 3rd truth-value, the material implication "if $p$, then $q$" is \textit{not} always true. This is a rejection of ECQ, because certain contradictions (those with the 3rd truth-value) do \textit{not} imply "trivial truth" for all other propositions.

\subsection{Gödel's Incompleteness Theorems and Tarski's Undefinability Theorem Assume Bivalence and LNC, and Thus are Irrelevant in 3VL}

Further regarding Łukasiewicz's 3VL, due to the third truth-value, both the LEM
and the LNC fail.
\footnote{See Haack \cite{Haack2}, p.5.} 
Priest's $LP$ goes further, by expressly defining the third truth-value as "\textit{both} true and false". So Gödel's second incompleteness theorem and Tarski's undefinability theorem (both of which assume bivalence, and the LNC) are rendered irrelevant in 3VL. Due to the absence of LNC, inconsistency is permitted in 3VL, and due to the absence of ECQ, inconsistency does not lead to triviality. It is the axiom LNC in classical logic (the foundational logic of math) that forces mathematics to be either consistent or trivial.

\subsection{Three Logical Frameworks for Dealing With Paradoxes}

According to Perzanowski \cite{Perzanowski2}, there are "at least" three logical responses to inconsistencies:
\footnote{See Perzanowski \cite{Perzanowski2}, p.11, para.19. Note that the fourth combination (rejecting LNC but accepting ECQ) is not possible, because a violation of LNC is a prerequisite for ECQ.}
\footnote{Superficially, Perzanowski's three responses to inconsistencies appear to be related to Lakatos's \cite{Lakatos} three methods ("monster-barring",  "monster-adjustment", and exception handling) to respond to counter-examples to mathematical theorems. However, Lakatos's three methods are relevant only when the LNC is accepted as an axiom, and therefore do not map to Perzanowski's three responses to inconsistencies. But note that Zermelo–Fraenkel set theory (both ZF and ZFC) has the LNC as an axiom, and is both a "monster-barring" and an inconsistency "enemy" foundation for set theory.}

\textbf{(1) Inconsistency "enemies":} This is the approach of logics (e.g. classical and intuitionistic logics) that accept both LNC and ECQ. Contradictions are not permitted to exist. The existence of a contradiction is a sign of logical "disease" (according to Tarski), so a single inconsistency trivializes everything. Therefore, every inconsistency must be discovered and quarantined. 

In regards to the RH, both RH and its negation $\neg$RH are paradoxes, and have either a truth-value glut, or a truth-value gap. These results are impermissible in classical logic, due to the LNC and LEM. So in classical logic, all theorems that assume that RH is true are invalid ("trivially true" due to ECQ) and unsound (due to the false assumption that AC of $\zeta(s)$ is true).

\textbf{(2) Paradox "believers":} This is the approach of logics (e.g. 3VL) that hold that paradoxes exist, and must be accounted for. Therefore, these logics assign the third truth-value to paradoxes, thereby bypassing both LNC and ECQ. 
\footnote{See Perzanowski \cite{Perzanowski2}, p.11, par.19, footnote 8: "The position has rather a long tradition, starting with the Sophists, Nicolas of Cusa, Hegel and Hegelians of several types (including the dialectic philosophers). In our time the position is defended by several Australian philosophers, including the late Richard Routley (later Sylvan), Chris Mortensen, and, under the name of \textit{dialethism}, by Graham Priest."}

Therefore, in such logics, the fact that RH and its negation $\neg$RH are paradoxes is not a catastrophe. In the case of 3VL, RH and its negation $\neg$RH are given the third truth-value. The truth tables of three-valued logic are applied accordingly.

\textbf{(3) Inconsistency "investigators":} This is the approach of paraconsistent logics that do not accept paradoxes, but do not want to trivialize the entire system due to the discovery of a paradox. These logics accept LNC but reject ECQ. They believe that “[i]n formal logic, a contradiction is the signal of a defeat; but in the evolution of real knowledge it marks the first step in progress towards victory." 
\footnote{See Whitehead \cite{Whitehead1}, Ch.11, p.187.}

Therefore, in such logics. the fact that RH and its negation $\neg$RH are paradoxes is \textit{not} fatal (does not necessarily cause ECQ). In the case of relevance logic, for example, RH and $\neg$RH only cause ECQ for propositions that are directly relevant to the RH and $\neg$RH. The truth tables of relevance logic are applied accordingly.

\subsubsection{In Classical Logic, RH Violates the LNC and Triggers ECQ}

In classical logic, because $\zeta(s)$ has no zeros, by material implication both RH ("\textit{all zeros are on the critical line}") and $\neg$RH ("\textit{not all zeros are on the critical line}") are "vacuously true". RH can be rephrased as "for all $s$, if $\zeta(s)=0$, then $\text{Re}(s)=0.5$", and $\neg$RH can rephrased as "\textit{not} for all $s$, if $\zeta(s)=0$, then $\text{Re}(s)=0.5$". Both are true, because according to material implication in classical logic, a false proposition (in this case, $\zeta(s)=0$) implies anything.

Also, RH and $\neg$RH are both \textit{false} by conjunction, when for example, $\neg$RH is rephrased as "there exists $\zeta(s)=0$ and $\text{Re}(s)\ne0.5$", and RH is rephrased as "\textit{not} (there exists $\zeta(s)=0$ and $\text{Re}(s)\ne0.5$)". However, the negations of these last two propositions are \textit{true}. 
Therefore, both RH and $\neg$RH are true, so RH is a semantic paradox, and therefore is also a contradiction.
\footnote{See Haack \cite{Haack}, pp.137-138: "it is possible to classify the paradoxes in two distinct groups, those which essentially involve set-theoretical concepts, such as '$\in$' and 'ordinal number', and those which essentially involve semantic concepts, such as 'false', 'false of ...', and 'definable'." }
\footnote{See also Bolander \cite{Bolander}, which adds a third group of paradoxes: "Epistemic paradoxes". These are similar to semantic paradoxes, except that "the central concept involved is knowledge rather than truth".}
\footnote{In contrast, see See Haack \cite{Haack}, p.138: "Russell himself, however, didn't think of the paradoxes as falling into two distinct groups, \textit{because he thought that they all as the result of one fallacy}, from violations of the 'vicious circle principle'." (Emphasis in the original).}

In classical logic, any conjecture that assumes the truth of a contradiction (\textit{such as the RH}) is, due to LNC and ECQ, "trivially true".

\subsubsection{In 3VL, There is Neither LNC Nor ECQ}

Three-value logics (3VL) avoid some of the paradoxes of classical logic, such as the paradoxes of implication. They do so by adding a third truth-value. As stated in Haack \cite{Haack}: "The proponent of a 3-valued logic ... seems to claim that there are valid arguments/logical truths of classical logic[,] the informal analogues of which aren't valid/logically true, so that classical logic is \textit{actually incorrect}",
\footnote{See Haack \cite{Haack}, p.222}
and "This explains in a more precise way the idea
\footnote{Citing Haack \cite{Haack}, ch.9 \S 3.}
... that deviant logics pose a more serious challenge than extended logics to classical logic."

In regards to the laws of classical logic, Haack is correct. All of 3VLs discussed in this paper (Frege's, Łukasiewicz's, Post's,  Bochvar's, Kleene's, and Priest's version thereof) bypass the LEM of classical logic.
\footnote{See, e.g., Haack \cite{Haack2}, pp.4-7, and Hazen et al. \cite{Hazen2}, and Hazen et al. \cite{Hazen3}.}
Priest's and Bochvar's 3VLs go further and assign the 3rd truth-value to paradoxes, thereby rejecting the LNC. Haack argues that also the other 3VLs (Łukasiewicz's, Post's, and Kleene's non-Priest version) reject the LNC. 
\footnote{See, e.g., Haack \cite{Haack2}, pp.4-7.} 

However, when comparing the truth tables of different 3VLs (Frege's, Łukasiewicz's, Post's,  Bochvar's, Kleene's, and Priest's $LP$) to those of classical logic, Haack is \textit{incorrect}. The truth tables of classical logic are included, \textit{in their entirety}, within the truth tables of all of these 3VLs. The 3VLs "extend" the truth tables of classical logic to a third truth-value. So Haack's classification of "deviant logics" and "extended logics" is misleading, because the truth tables of the so-called "deviant logics" are "extended" versions of classical logic's truth tables, and do not contradict any value in classical logic's truth tables.

In Frege's 3VL, the third truth-value is "\textit{neither} true nor false". But it is proven in Frege's logic that \textit{neither} implies \textit{both}, and vice versa.
\footnote{See Milne's \cite{Milne2} p.475 (citing Heidelberger \cite{Heidelberger}): "In short, everywhere we think there’s a truth-value gap, there’s also a‘glut’! (And vice versa !)".}
\footnote{See also Priest \cite{Priest6}, p.27, which recites an explanation for this phenomenon that is different from Heidelberger's \cite{Heidelberger}: "Notably, assuming De Morgan’s laws, ... $\lnot(A \lor \lnot A)$ is equivalent to ... $A \land \lnot A$".}
\footnote{See also Bolander \cite{Bolander}, \S 3.2.2 "Extensions and Alternatives to Kripke’s Theory of Truth", which states (emphasis in the original): "The choice is between \textit{truth-value gaps} and \textit{truth-value gluts}: A truth-value gap is a statement with no truth-value, neither true or false (like \textit{undefined} in Kleene’s strong three-valued logic), and a truth-value glut is a statement with several truth-values, e.g. both true and false (like in the paraconsistent logic LP). There are also arguments in favour of allowing both gaps and gluts, e.g. by letting the set of truth-values form of a bilattice [citing Fitting \cite{Fitting} and Odintsov et al. \cite{Odintsov}]. The simplest non-trivial bilattice has exactly four values, which in the context of truth-values are interpreted as: \textit{true}, \textit{false}, $\bot$  (neither true nor false), and $\top$  (both true and false). For a more extensive discussion of Kripke’s theory, its successors and rivals, see the entry on the liar paradox [citing Beall et al. \cite{Beall}]."}
According to the truth tables of Kleene's, Łukasiewicz's, and Priest's 3VLs, the proposition "if RH, then $p$" is \textit{true} by material implication, if $p$ is true. This is consistent with ECQ, because RH has the third truth-value.
However, in these same 3VLs, the proposition "if RH, then $p$" is \textit{false} by material implication, if $p$ is false. This is \textit{inconsistent} with ECQ, because according to ECQ, the result should be trivial truth.

In regards to the proposition "if RH, then $p$", there is a difference of opinion between Łukasiewicz and Kleene/Priest regarding the value of material implication when $p$ has the third truth-value "neither/both". In Łukasiewicz's 3VL, the proposition "if RH, then $p$" is true when $p$ has the third truth-value "neither/both". In Kleene/Priest, the proposition has the third truth-value when $p$ has the third truth-value. 
\footnote{See Kleene \cite{Kleene2}, p.335.}
\footnote{See also Haack \cite{Haack}, pp.206-208.}
\footnote{See Wikipedia \cite{3VL}: "The Łukasiewicz Ł3 has the same tables for AND, OR, and NOT as the Kleene logic given above, but differs in its definition of implication in that "unknown implies unknown" is \textit{true}. This section follows the presentation from [Malinowski \cite{Malinowski}]."}
Both interpretations are \textit{servicable for reasoning purposes} since these rules at least have the property that they will do not lead us from an assumption having a truth-value of "true", or a truth-value glut (that includes the truth-value of "true") ... to a false conclusion.
\footnote{Paraphrasing the reasoning for accepting material implication in classical logic. See Lemmon \cite{Lemmon}, p.60, citing Chapter 2, \S 4, pp.75-82.} 
In other words, both Łukasiewicz's and Kleene/Priest's 3VL material implication provide "truth preservation".

Priest's "Logic of Paradox" ($LP$) is Kleene's 3VL, and therefore assigns the third truth-value (\textit{both} true and false) to paradoxes. In $LP$, the material implication "if RH, then $p$" has the third truth-value (the same truth-value as $p$).  

But Priest's version of Kleene's three-valued logic (3VL) (which Priest calls "Logic of Paradox" $LP$), paradoxes such as the RH are assigned the third-truth-value (\textit{both} true and false).
\footnote{But Haack \cite{Haack}, p. 211, and ch.8, \S 2, argues that "this kind of approach to the paradoxes is apt to from the frying pan - the Liar paradox - to the fire - the Strengthened Liar ('this sentence is either false or paradoxical', true if false or paradoxical, false or paradoxical if true)." The counter-argument to Haack is Priest's concept of a "truth-value glut". If we assume the existence of a third-truth-value, and therby bypass LEM and LNC, then both Haack's argument and also Priest's can be true. Otherwise we have yet another paradox.}
According to $LP$'s material implication, the truth-value of material implication "if RH, then $p$" is the same as the truth-value of $p$. 

So in classical logic, a paradox (e.g. RH) implies ECQ, and thus implies trivial truth. But in $LP$, a paradox does not imply ECQ. Instead, a paradox can imply non-trivial truth, falsity, or the third truth-value. Also, both Gödel's second incompleteness theorem and Tarski's undefinability theorem are irrelevant in $LP$, because $LP$ rejects the LNC, but both of these theorems assume the LNC.

\subsubsection{In Intuitionistic Logic, RH is False }

Kleene states the following in regards to intuitionistic logic:
 
\begin{quotation}
An existence statement \textit{there exists a natural number $n$ having the property $P$}, or briefly \textit{there exists an n such that $P(n)$}, has its intuitionistic meaning as a partial communication (or abstract) of a statement giving a particular example of a natural number $n$ which has the property $P$, or at least giving a method by which in principle one could find such an example.
\footnote{See Kleene \cite{Kleene2}, p. 49.}
\end{quotation}
 
So, for example, prior to Wiles's proof of Fermat's last theorem ("FLT"), Intuitionists would reject any non-constructive existence proof (which is acceptable in classical logic), such as: "If FLT is true, then the number 5013 has the property $P(n)$, and if FLT is false, then the number 10 has the property $P(n)$."  
\footnote{See e.g., Kleene \cite{Kleene2}, p. 50.}

Kleene's implementation of intuitionistic logic is based on that of Hilbert and Ackerman, Hilbert and Bernays, Gentzen, etc.,
\footnote{See Kleene \cite{Kleene2}, pp. 69.}
and is identical to classical logic, but without both the Principle of Double Negation and the LEM.
\footnote{See Kleene \cite{Kleene2}, p.120, *51, Remark 1 states that "either of $\lnot \lnot A \supset A$ [Principle of Double Negation] or $A \lor \lnot A$ [LEM] can be chosen as the one non-intuitionistic postulate of the classical system."}
However, Kleene's implementation of intuitionistic logic, which eliminates the LEM \textit{completely},is wrong. Intuitionists \textit{do accept} the LEM, \textit{but only if} there is a constructive existence proof, or disproof. Therefore, a more accurate implementation of intuitionistic logic is Kleene's 3VL. Kleene states the following regarding his own 3VL:
 
\begin{quotation}
We further conclude from the introductory discussion that, for the definitions of partial recursive operations, $t, f, u$ must be susceptible of another meaning besides (i) 'true', 'false', 'undefined', namely (ii) 'true', 'false', 'unknown (or value immaterial)'. Here 'unknown' is a category into which we can regard any proposition as falling, whose value we either do not know or choose for the moment to disregard; and it does not then exclude the other two possibilities 'true' and 'false'."
\footnote{See Kleene \cite{Kleene2}, p.335.}
\end{quotation}

In other words, this interpretation of 3VL implements what the intuitionists argued: that in the absence of a constructive proof or disproof, a proposition has an 'unknown' truth-value. The LEM becomes relevant after a classical truth-value is obtained, therefore of a constructive proof (or disproof).

\section{The 3rd Truth-Value, Truth-Value Gluts, and Truth-Value Gaps}

\subsection{Truth-Value Glut: RH is Both True and False}

\subsubsection{Classical Logic:  Russell's "On Denoting"}

Russell's \textit{On Denoting} \cite{Russell2} (which like axiomatic set theory, has the LEM as an axiom) holds that a proposition with a vacuous subject (e.g. the Riemann hypothesis) is ambiguous, because it can be interpreted in two ways. Therefore, depending on how such a statement (e.g. the RH) is interpreted, it can be either true or false. (In its ambiguous state, it has both meanings).

In contrast, the negation ("the present King of France is not bald") can be interpreted as the conjunction of the following three propositions:
\footnote{See Batty \cite{Batty}, "1. Russell Recap".}

i. There is at least one King of France. $\exists x(Kx)$

ii. There is at most one King of France. $(x)(y)(Kx \land Ky \to x = y)$

iii. Whatever is King of France is not bald. $(x)(Kx \to \lnot Bx)$

When these three propositions are conjoined, we get: "There is one and only one present King of France and he is not bald." In standard logical notation, this first sentence is:
\footnote{See Russell \cite{Russell2}, p.490, and Jacquette \cite{Jacquette}, pp.5-6.}
\[\exists x \Big(Kx \land (\forall y)\Big((Ky\to x=y \Big) \land \lnot Bx\Big)\]
This sentence is \textit{false}, because it quantifies over a non-existent entity. ("There is one and only one present King of France" is false). 

A second interpretation of the sentence is: "It is not the case that that there exists a present King of France and he is bald". The second interpretation is \textit{true}, because it is indeed not the case that that there exists a present King of France. In standard logical notation, this second sentence is:
\footnote{See Russell \cite{Russell2}, p.490, and Jacquette \cite{Jacquette}, pp.5-6.}
\[\lnot \exists x \Big(Kx \land (\forall y)\Big((Ky\to x=y \Big) \land Bx\Big)\] 

If the RH is interpreted according to Russell's first interpretation, as "there exist zeros of $\zeta(s)$ and they are not located off of the critical line $\text{Re}(s)=0.5$", then the RH is \textit{false}, because it quantifies over non-existent entities (the non-existent zeros of $\zeta(s)$). 

However, if the RH is interpreted according to Russell's second interpretation, as "it is not the case that there exist zeros of $\zeta(s)$ and they are located off of the critical line $\text{Re}(s)=0.5$", then it is \textit{true}, because indeed it is \textit{not} the case that there exist zeros of $\zeta(s)$. (Note: Both axiomatic set theory and Russell's \textit{On Denoting} assume that the LEM is true).

Moreover, if we apply Russell's first interpretation to the RH, and to its negation $\neg$RH ("not all zeros of $\zeta(s)$ are on the critical line $\text{Re}(s)=0.5$"), then paradoxically \textit{both} are \textit{false}.

\subsection{Truth-Value Gap: RH is Neither True Nor False}

In the alternative, some logics that reject the LEM hold the Riemann hypothesis to be neither true nor false, because in these logics, some propositions are not assigned a (classical) truth-value. 

In those systems that embrace truth-value gaps (Strawson, Frege) or non-classically-valued systems (Łukasiewicz, Bochvar, Kleene), some sentences or statements are not assigned a (classical) truth-value. However, in the specific case of Strawson's \textit{On Referring}, its reasoning is inapplicable to the Riemann hypothesis, for reasons that will be discussed later in this paper.

\subsubsection{Intuitionistic Logic}

Brouwer presented his theory of intuitionism, a philosophy of the foundations of mathematics, in \textit{Intuitionism and Formalism} (1913).
\footnote{See Brouwer \cite{Brouwer3}.}
As Davis \cite{Davis} explains:

\begin{quotation}
For Brouwer, some propositions can neither be said to be true or to be false; these are propositions for which no method is currently known by means of which this can be decided one way or the other. Hilbert's original proof of Gordon's conjecture used the law of the excluded middle in the way mathematicians usually do: he showed that denying the conjecture would lead to a contradiction. To Brouwer such a proof was unacceptable.
\footnote{See Davis \cite{Davis}, p.95.}
\end{quotation}

This summary is repeated by Iemhoff \cite{Iemhoff}:

\begin{quotation}
According to the BHK-interpretation[,] this statement [LEM] holds intuitionistically if the creating subject knows a proof of A[,] or a proof that A cannot be proved. In the case that neither for A nor for its negation a proof is known, the statement $(A\lor \neg A)$ does not hold. 
\end{quotation}

Brouwer did not object to the LNC, and thus the LNC is included in intuitionistic logic. The LNC, in combination with the proof that the Dirichlet series $\zeta(s)$ is divergent in the half-plane $\text{Re}(s)\le1$ means that $\zeta(s)$ has no zeros. In light of these facts, according to intuitionism, the LEM \textit{does} hold for both of the propositions "$\zeta(s)$ has zeros" and "$\zeta(s)$ has no zeros". 

However, because $\zeta(s)$ has no zeros, the RH is directed to "vacuous subjects". Therefore, no proof is possible for either RH or for its negation $\neg$RH. So according to intuitionism, the LEM \textit{does not} hold for either the RH or for its negation $\neg$RH. And therefore, according to intuitionism, the RH has no truth-value. (It is a "truth-value gap").

The RH being a paradox provides a stronger argument than Brouwer's against the LEM: the LEM cannot be used to hold that theorems are either true or false, because some theorems are paradoxes (and thus require a 3rd truth -value).

\subsubsection{Russell's Argument}

Russell's \textit{On Denoting} (1905) preceded L.E.J. Brouwer's intuitionism by a few years, and presents ideas that are shared with intuitionism in regards to "vacuous subjects": (1) requiring a proof of existence in order to use LEM, and (2) in the absence of proof of existence, abandoning the LEM and assigning a "truth-value gap" to the proposition.

According to Russell's "Theory of Descriptions", a proposition with a vacuous subject "'C has the property $\phi$' is false for all values of $\phi$".
\footnote{See Russell \cite{Russell2} p.490.}
So, according to Russell, "the present King of France is bald" is "certainly false", and "the present King of France is not bald" is also false if it means "There is an entity which is now King of France and is not bald", but true if it means "It is false that there is now an entity which is now King of France and is not bald".
\footnote{See Russell \cite{Russell2} p.490. See also Pelletier et al. \cite{Pelletier}, and Haack \cite{Haack2}, p.15.}

However, it is important to note that according to Venn's "Modern" Square of Opposition, these results of Russell's "Theory of Descriptions" creates a paradox (which violates the LNC, one of the theorems in Whitehead and Russell's classical logic). According to Venn's "Modern" Square of Opposition, if both "the present King of France is bald" is false, and "the present King of France is not bald" is also false, then both of the propositions "all Kings of France are bald" and "all Kings of France are not bald" are true.
\footnote{See Parsons \cite{Parsons}.}

In addition, Russell unknowingly also presents an \textit{alternative} argument in favor of abandoning the LEM, but in which a proposition with a vacuous subject (e.g. "the present King of France is bald") is assigned a third truth-value (a "truth-value gap") instead of the truth-value of "false":
\footnote{See Russell \cite{Russell2} p.485.}
 
\begin{quotation}
By the law of the excluded middle [LEM], either 'A is B' or 'A is not B' must be true. Hence either "the present King of France is bald" or "the present King of France is not bald" must be true. Yet if we enumerated the things that are bald, and then the things that are not bald, we should not find the present King of France in either list.   
\end{quotation}
 
(Russell fails to mention the obvious conclusion: that because we do not find the present King of France in either list, it means that both propositions are neither true nor false).

\subsubsection{Frege's Argument}

Speranza et al.'s \cite{Speranza} quotation of Christoph Sigwart presents the essence of Frege's argument regarding truth-value gaps: 
\footnote{See Speranza et al. \cite{Speranza}, p.148.}

\begin{quotation}
For Strawson, as for his intellectual predecessor Frege [1892], the notion of presupposition has semantic status as a necessary condition on true or false assertion ... In fact, the earliest pragmatic treatments of the failure of existential presupposition predate Frege's analysis by two decades. Here is Christoph Sigwart [1873] on the problem of vacuous subjects:

"As a rule, the judgement A is not B presupposes the existence of A in all cases when it would be presupposed in the judgement A is B ... 'Socrates is not ill' presupposes in the first place the existence of Socrates, because only on the presupposition [Voraussetuzung] of his existence can there be any question of his being ill."
(Sigwart [1873/1895: 122], ...)
\end{quotation}

\subsubsection{Strawson's "On Referring"}

Accordingly, Aristotle's and Russell's logics (which assume the LEM) hold the RH to be false, but Frege's and Strawson's logics hold that the RH cannot be used to make a true or false assertion (thereby rejecting the LEM).
\footnote{See Horn \cite{Horn}, "4. Gaps and Gluts: LNC and Its Discontents".}  More specifically, according to Horn:

\begin{quotation}
In those systems that do embrace truth-value gaps (Strawson, arguably Frege) or non-classically-valued systems (Łukasiewicz, Bochvar, Kleene), some sentences or statements are not assigned a (classical) truth-value; in Strawson's famous dictum, the question of the truth-value of “The king of France is wise”, in a world in which France is a republic, simply fails to arise. The negative form of such vacuous statements, e.g. “The king of France is not wise”, is similarly neither true nor false. This amounts to a rejection of LEM, as noted by Russell [in "On Denoting"].
\end{quotation}

In contrast to Russell's \textit{On Denoting} \cite{Russell2}, Strawson's \textit{On Referring} \cite{Strawson} states that a statement with a vacuous subject (a subject term that has no referent, e.g. "the present King of France") is \textit{not} false. Instead, it is "absurd" and therefore not asked. So, it is neither true nor false (and thus belongs in a third category, whose existence is a rejection of LEM). Strawson provides the following example:

\begin{quotation}
A literal-minded and childless man asked whether all his children are asleep will certainly not answer "Yes" on the ground that he has none; but nor will he answer "No" on this ground. Since he has no children,
the question does not arise.
\end{quotation}

However, Strawson assumes that the potential questioner knows that the question has a vacuous subject. The 160 year history of the RH shows that this is not always the case. 

In the context of the Riemann Hypothesis, Strawson's argument is clearly wrong. Over the course of the past 160 years, many mathematicians \textit{have} asked if all of the zeros of $\zeta(s)$ indeed are on the critical line $\text{Re}(s)=0.5$. The question \textit{has} arisen, because in contrast to Strawson's examples ("the present King of France", the children of a man well-known to be childless), it has not been common knowledge that Riemann's $\zeta(s)$ violates LNC (or that $\zeta(s)$ thus has no zeros). Instead, Riemann's $\zeta(s)$ was widely assumed \textit{to indeed have zeros}. So an axiom of Strawson's logic (common knowledge that subject of the question is vacuous) is clearly false in the context of the Riemann Hypothesis.

\subsection{Comparison of Truth-Value Gluts to Truth-Value Gaps}

\subsubsection{Comparison of Truth Tables}

The Three-Valued Logic Truth Tables shown below are those of Frege, Kleene, Bochvar, and Łukasiewicz. There are others, but a full discussion is beyond the scope of this paper.
\footnote{See e.g. Ciucci et al. \cite{Ciucci}.}

Remember that the "third value" in the Frege truth tables is the \textit{absence} of any truth-value ("a truth-value gap"). As Milne \cite{Milne2} states: "Beware! The bar [ - ] is not a third truth-value, it signifies the absence of a truth-value. Where both [variables] have truth-values, the connectives behave classically. "
\footnote{See Milne \cite{Milne2}, p.473, citing Frege's "Logic". Milne's \cite{Milne2} p.474 reproduces Smiley's truth tables for Frege's three-valued logic (Citing Smiley \cite{Smiley}, pp.125-35.).} 
\footnote{See also Haack \cite{Haack}, p.212.}

The "Kleene" and "Łukasiewicz" tables are "essentially those of Kleene’s and Łukasiewicz’s three valued logics", respectively. 
\footnote{See Priest et al. \cite{Priest3}, \S 3.6.}
\footnote{See also Haack \cite{Haack}, pp.206-208.}
The "Bochvar" tables are those of yet another 3VL, which (unlike Kleene's and Łukasiewicz's) was originally intended as a solution to semantic paradoxes.
\footnote{See Haack \cite{Haack}, pp.206-208, citing Bochvar \cite{Bochvar}.}
Bochvar adds an "assertion operator" (presented here as "T"), which means something like "It is true that:". The "external connectives are defined as follows: $\neg A = \neg TA$, $A \& B = TA \& TB$, $A \lor B = TA \lor TB$, $A \rightarrow B = TA \rightarrow TB$.

The values in all of the truth tables presented here are: $\CIRCLE$ =  \textit{True (only)}, $\Circle$ = \textit{False (only)}, and $\XBox$ = \textit{Both (True and False)}. The conditional $\rightarrow$, follows Kleene’s three valued logic,
\footnote{See Priest et al. \cite{Priest3}, \S 3.6.}
\footnote{See Kleene \cite{Kleene2}, p.335.}
\footnote{See also Wikipedia \cite{3VL}: "The Łukasiewicz Ł3 has the same tables for AND, OR, and NOT as the Kleene logic given above, but differs in its definition of implication in that 'unknown implies unknown' is \textbf{true}", citing Malinowski \cite{Malinowski}}
and material equivalence $\leftrightarrow$, is defined as "means the same as".

\begin{table}
\begin{tabular}[t]{c|c}
$\lnot$ & \\ \hline
$\CIRCLE$ & $\Circle$ \\ 
- & -\\ 
$\Circle$ & $\CIRCLE$ \\    
\end{tabular}
\hfill
\begin{tabular}[t]{c|c|c|c}
$\land$ & $\CIRCLE$ & - & $\Circle$ \\ \hline
$\CIRCLE$ & $\CIRCLE$ & - & $\Circle$ \\ 
- & - & - & - \\ 
$\Circle$ & $\Circle$ & - & $\Circle$ \\    
\end{tabular}
\hfill
\begin{tabular}[t]{c|c|c|c}
$\lor$ & $\CIRCLE$ & - & $\Circle$ \\ \hline
$\CIRCLE$ & $\CIRCLE$ & - & $\CIRCLE$ \\ 
- & - & - & - \\ 
$\Circle$ & $\CIRCLE$ & - & $\Circle$ \\    
\end{tabular}
\hfill
\begin{tabular}[t]{c|c|c|c}
$\rightarrow$ & $\CIRCLE$ & - & $\Circle$ \\ \hline
$\CIRCLE$ & $\CIRCLE$ & - & $\Circle$ \\ 
- & - & - & - \\ 
$\Circle$ & $\CIRCLE$ & - & $\CIRCLE$ \\    
\end{tabular}
\hfill
\begin{tabular}[t]{c|c|c|c}
$\leftrightarrow$ & $\CIRCLE$ & - & $\Circle$ \\ \hline
$\CIRCLE$ & $\CIRCLE$ & - & $\Circle$ \\ 
- & - & - & - \\ 
$\Circle$ & $\Circle$ & - & $\CIRCLE$ \\    
\end{tabular}
\caption{Frege's Truth Tables (Truth-Value Gaps)}

\begin{tabular}[t]{c|c}
$\lnot$ & \\ \hline
$\CIRCLE$ & $\Circle$ \\ 
$\XBox$ & $\XBox$\\ 
$\Circle$ & $\CIRCLE$ \\    
\end{tabular}
\hfill
\begin{tabular}[t]{c|c|c|c}
$\land$ & $\CIRCLE$ & $\XBox$ & $\Circle$ \\ \hline
$\CIRCLE$ & $\CIRCLE$ & $\XBox$ & $\Circle$ \\ 
$\XBox$ & $\XBox$ & $\XBox$ & $\Circle$ \\ 
$\Circle$ & $\Circle$ & $\Circle$ & $\Circle$ \\    
\end{tabular}
\hfill
\begin{tabular}[t]{c|c|c|c}
$\lor$ & $\CIRCLE$ & $\XBox$ & $\Circle$ \\ \hline
$\CIRCLE$ & $\CIRCLE$ & $\CIRCLE$ & $\CIRCLE$ \\ 
$\XBox$ & $\CIRCLE$ & $\XBox$ & $\XBox$ \\ 
$\Circle$ & $\CIRCLE$ & $\XBox$ & $\Circle$ \\    
\end{tabular}
\hfill
\begin{tabular}[t]{c|c|c|c}
$\rightarrow$ & $\CIRCLE$ & $\XBox$ & $\Circle$ \\ \hline
$\CIRCLE$ & $\CIRCLE$ & $\XBox$ & $\Circle$ \\ 
$\XBox$ & $\CIRCLE$ & $\XBox$ & $\XBox$ \\ 
$\Circle$ & $\CIRCLE$ & $\CIRCLE$ & $\CIRCLE$ \\    
\end{tabular}
\hfill
\begin{tabular}[t]{c|c|c|c}
$\leftrightarrow$ & $\CIRCLE$ & $\XBox$ & $\Circle$ \\ \hline
$\CIRCLE$ & $\CIRCLE$ & $\Circle$ & $\Circle$ \\ 
$\XBox$ & $\Circle$ & $\CIRCLE$ & $\Circle$ \\ 
$\Circle$ & $\Circle$ & $\Circle$ & $\CIRCLE$ \\
\end{tabular}
\caption{Kleene's Truth Tables (Truth-Value Gluts)}

\begin{tabular}[t]{c|c}
$\lnot$ & \\ \hline
$\CIRCLE$ & $\Circle$ \\ 
$\XBox$ & $\XBox$\\ 
$\Circle$ & $\CIRCLE$ \\    
\end{tabular}
\hfill
\begin{tabular}[t]{c|c|c|c}
$\land$ & $\CIRCLE$ & $\XBox$ & $\Circle$ \\ \hline
$\CIRCLE$ & $\CIRCLE$ & $\XBox$ & $\Circle$ \\ 
$\XBox$ & $\XBox$ & $\XBox$ & $\XBox$ \\ 
$\Circle$ & $\Circle$ & $\XBox$ & $\Circle$ \\    
\end{tabular}
\hfill
\begin{tabular}[t]{c|c|c|c}
$\lor$ & $\CIRCLE$ & $\XBox$ & $\Circle$ \\ \hline
$\CIRCLE$ & $\CIRCLE$ & $\XBox$ & $\CIRCLE$ \\ 
$\XBox$ & $\XBox$ & $\XBox$ & $\XBox$ \\ 
$\Circle$ & $\CIRCLE$ & $\XBox$ & $\Circle$ \\    
\end{tabular}
\hfill
\begin{tabular}[t]{c|c|c|c}
$\rightarrow$ & $\CIRCLE$ & $\XBox$ & $\Circle$ \\ \hline
$\CIRCLE$ & $\CIRCLE$ & $\XBox$ & $\Circle$ \\ 
$\XBox$ & $\XBox$ & $\XBox$ & $\XBox$ \\ 
$\Circle$ & $\CIRCLE$ & $\XBox$ & $\CIRCLE$ \\    
\end{tabular}
\hfill
\begin{tabular}[t]{c|c}
$T$ &  \\ \hline
$\CIRCLE$ & $\CIRCLE$  \\ 
$\XBox$ & $\Circle$  \\ 
$\Circle$ & $\Circle$  \\
\end{tabular}
\caption{Bochvar's Truth Tables (Truth-Value Gluts)}

\begin{tabular}[t]{c|c|c|c}
$\rightarrow$ & $\CIRCLE$ & $\XBox$ & $\Circle$ \\ \hline
$\CIRCLE$ & $\CIRCLE$ & $\XBox$ & $\Circle$ \\ 
$\XBox$ & $\CIRCLE$ & $\CIRCLE$ & $\XBox$ \\ 
$\Circle$ & $\CIRCLE$ & $\CIRCLE$ & $\CIRCLE$ \\    
\end{tabular}
\caption{"Łukasiewicz" Material Implication  (Other Operators are Same as Kleene's)}
\end{table}

\subsubsection{Every Truth-Value Gap Implies a Glut}

The following natural deduction rules in classical logic fail in Frege's Truth-value gap logic:
v-introduction, $\rightarrow$-introduction (conditional proof), \textit{reductio ad absurdum}, \textit{ex falso quodlibet} (ECQ), the law of the excluded middle (LEM)
\footnote{See Milne's \cite{Milne2} p.474.}

However, enough of classical logic remains valid to prove the following:
\footnote{See Milne's \cite{Milne2} p.475 (citing Heidelberger \cite{Heidelberger}): "In short, everywhere we think there’s a truth-value gap, there’s also a‘glut’! (And vice versa !)".}
 
\begin{quotation}
It's not true that $P$ and it's not false that $P$ only if it's both true that $P$ and false that $P$.
\end{quotation}
 
So in Frege's logic, whenever there is a truth-value gap, there is also a truth-value glut (and vice versa). 
\footnote{See also Priest \cite{Priest6}, p.27, which recites an explanation for this phenomenon that is different from Heidelberger's \cite{Heidelberger}: "Notably, assuming De Morgan’s laws, ... $\lnot(A \lor \lnot A)$ is equivalent to ... $A \land \lnot A$".}
This is a paradox, because there is a contradiction here.
\footnote{See Heis \cite{Heis}: "Frege, of course, would resolve this paradox by prescribing that a logically perfected language have no bearerless names. Milne \cite{Milne2} advocates instead adopting a semantic (as opposed to Frege's functional) theory of negation. He rejects Frege's solution because it precludes a plausible semantics for ordinary language, and because the set-theoretic paradoxes show that even a scientific language such as Frege's own needs to allow for the possibility of singular terms (like "the extension of $x\notin x$") that are nevertheless bearerless."}
\footnote{See also Scruton \cite{Scruton}, p.63: "Frege argued that there are just two 'truth-values' as he called them: the true and the false. He therefore suggested that a sentence will refer to one or other of two things: truth (the true) or falsehood (the false)."}
\footnote{See also Scruton \cite{Scruton}, p.72: "Just as 'the golden mountain' lacks a reference, therefore, the sentence 'the golden mountain is hidden' lacks a truth-value."}
But according to another interpretation, this is a \textit{not} a paradox, because the difference lies in the definition of tautologies.
\footnote{See Wikipedia \cite{3VL}, citing Look \cite{Look}: "In these truth tables, the \textit{unknown} state can be thought of as neither \textit{true} nor \textit{false} in Kleene logic, or thought of as both \textit{true} and \textit{false} in Priest logic. The difference lies in the definition of tautologies. Where Kleene logic's only designated truth-value is $T$, Priest logic's designated truth-values are both $T$ and $U$. In Kleene logic, the knowledge of whether any particular \textit{unknown} state secretly represents \textit{true} or \textit{false} at any moment in time is not available. However, certain logical operations can yield an unambiguous result, even if they involve at least one \textit{unknown} operand."} 
\footnote{See Wikipedia \cite{3VL}, citing Look \cite{Look}: "Kleene logic has no tautologies (valid formulas) because whenever all of the atomic components of a well-formed formula are assigned the value Unknown, the formula itself must also have the value Unknown."} 
\footnote{See Wikipedia \cite{3VL}, citing Look \cite{Look}: "[Priest's] Logic of Paradox ($LP$) has the same truth tables as Kleene logic, but it has two designated truth-values instead of one; these are: \textit{True} and \textit{Both} (the analogue of \textit{Unknown}), so that $LP$ does have tautologies but it has fewer valid inference rules."}

However, if this is indeed a paradox, then we should always apply logic based on truth-value gluts (e.g. Kleene's three-valued logic) instead of logic based on truth-value gaps (e.g. Frege's logic), because the former is "truth preserving".
\footnote{See Priest et al. \cite{Priest3}, \S 3.6: "Let t [true] and b [both] be the designated values. These are the values that are preserved in valid inferences. If we define a consequence relation in terms of preservation of these designated values, then we have the paraconsistent logic $LP$. In $LP$. ECQ is invalid", citing Priest \cite{Priest5}.}
\footnote{See also MacFarlane \cite{MacFarlane}, p.14: "The idea is to keep the classical idea that validity is truth preservation, but give up the classical assumption that the same sentence cannot be both true and false."}

Perhaps the most interesting result in the "Kleene" three-valued truth tables is that of material implication, $(A\rightarrow B)$. In classical logic, the material implication $A\to B$ is equivalent to $\neg A\lor B$ (this can be seen in the "Frege" truth tables).  So it is true if $A$ is false, regardless of whether $B$ is true or false.

In a three-valued logic, the material implication $A\to B$ remains equivalent to $(\neg A\lor B)$. So if $A$ is "both true and false", then the material implication \textit{is not false}, regardless of the value of $B$. This can be seen in the "Kleene" (But Not "Frege") truth tables.

So if the RH has the third truth-value ("both true and false"), then in classical logic all theorems that assume RH is true cannot be proven true (or proven false). But in three-valued logic, material implication holds that they might be proven true.

\section{State Table of $\zeta(s)$ and Truth Tables of RH}

\subsection{The State Table of the Zeta Function}

It is said that a picture is worth a thousand words, so we begin our discussion by filling out the state table of $\zeta(s)$, as a function of: the truth/falsity of the Law of Non-Contradiction (LNC), and the truth/falsity of the analytic continuation of $\zeta(s)$.  

\begin{table}[h] 
\begin{tabular}{cc|c|c|}
 & & \multicolumn{2}{c}{\{LNC\}} \\
 &  & True & False  \\ 
\hline
 & True & Divergent \& Convergent & Divergent \& Convergent \\ 
\{Analytic & & (Paradox)* & (Paradox)$\dagger$ \\ 
\cline{2-4}
Contin. of $\zeta(s)$\} & False & Divergent & Divergent  \\ 
 & & &  \\  
\hline
\end{tabular}
\label{tab:TT_Zeta}
\textbf{\caption{State Table of $\zeta(s)$ in Half-Plane $\text{Re}(s)\le1$}} 
.

(* = Violates LNC. In logics that have ECQ as a theorem, this triggers ECQ.)

($\dagger$ = In certain Multi-Valued Logics, paradoxes are assigned a 3rd truth-value.)
\end{table}

The Dirichlet series $\zeta(s)$ is \textit{proven} to be divergent throughout the half-plane $\text{Re}(s)\le1$. So if the analytic continuation ("AC") of $\zeta(s)$ to half-plane $\text{Re}(s)\le1$ is true, then the function $\zeta(s)$ is a paradox in that half-plane, because it is \textit{both} convergent \textit{and} divergent at every value of $s$ in that half-plane (except at the pole at $s=1$). \footnote{In contrast, Weierstrass's analytic continuation lacks such a direct contradiction. See Chapter  \ref{Weierstrass} of this paper for more details.}

If the analytic continuation (AC) of $\zeta(s)$ is true, then the state-value of $\zeta(s)$ in half-plane $\text{Re}(s)\le1$ is "paradox". In a logic with both LNC and ECQ, this violation of LNC triggers ECQ. In contrast, in a certain 3VLs and 4VLs, the LNC can be bypassed, by assigning a third truth-value to paradoxes.

But if the AC of $\zeta(s)$ is false,
\footnote{Due to, for example, Riemann's analytic continuation of $\zeta(s)$ being invalid. See Chapters \ref{Riemann-Invalid_1} and \ref{Riemann-Invalid_2}.}
then $\zeta(s)$ is exclusively defined by its Dirichlet series, which is divergent throughout the half-plane $\text{Re}(s)\le1$, and has no zeros and no poles, both in logics where the LNC holds, and in logics where the LNC fails. 

\subsection{The Truth Table of the Riemann Hypothesis}

As shown in the preceding state table, the state of $\zeta(s)$ depends on whether the analytic continuation of $\zeta(s)$ is true or false. If the analytic continuation is true, $\zeta(s)$ is a paradox in half-plane $\text{Re}(s)\le1$.  If false, $\zeta(s)$ is exclusively defined by the Dirichlet series, and thus is divergent in said half-plane.

Next we discuss the truth table of the Riemann hypothesis, as a function of the state of $\zeta(s)$ in half-plane $\text{Re}(s)\le1$, and of three different classes of logic (classical, intuitionistic, and 3VLs that assign a 3rd truth-value to paradoxes).  

\begin{table}[h] \label{TT_RH}
\begin{tabular}{c|c|c|c|} \\
 & \multicolumn{3}{c}{\{Logics\}} \\
 & Classical & Intutionistic & 3VLs That Assign   \\ 
 & & & a 3rd Truth-Value   \\  
  & & & to Paradoxes  \\  
\hline
 Paradox $\zeta(s)$ & Trivially True & Trivially True & 3rd Truth-Value$\dagger$ \\ 
 (Convergent \& Divergent) & (due to ECQ) & (due to ECQ) & (due to $\zeta(s)$) \\ 
 \hline
 Dirichlet Series $\zeta(s)$ & Paradox* \& ECQ & False** & 3rd Truth-Value $\ddagger$ \\ 
  & (due to no Zeros) & (due to no Zeros) & (due to no Zeros) \\ 
 \hline
\end{tabular}
 \label{tab:TT_RH}
\textbf{\caption{Truth Table of the Riemann Hypothesis (RH)}}

.

(* = Both RH and anti-RH ("All zeros are off the critical line") are true, due to "vacuous zeros". This violates LNC and triggers ECQ.)

(** = RH's zeros are proven to be unconstructable.)

($\dagger$ = In Bochvar's 3VL, material implication has the 3rd truth-value if it is from a paradox to any other proposition. In Priest's $LP$, there is no material implication for paradoxes.)

($\ddagger$ = In all 3VLs, the material implication of a false 1st proposition to any 2nd proposition (and to the negation of the 2nd proposition) is true, resulting in a paradox.)
\end{table}

\subsubsection{If Analytic Continuation of $\zeta(s)$ is False}

If analytic continuation (AC) of $\zeta(s)$ is false, 
\footnote{We show in this paper that Riemann's alleged proof of the analytic continuation of $\zeta(s)$ is false. There are other alleged proofs, so unfortunately this result is not dispositive by itself. 
See e.g. Titchmarsh et al. \cite{Titchmarsh}, \S 2.1 to \S 2.10, pp.13-27, which lists seven such proofs.}
then $\zeta(s)$ is exclusively defined by its Dirichlet series, which is divergent throughout half-plane $\text{Re}(s)\le1$, and has neither zeros nor poles. In this scenario, the RH refers to "vacuous zeros" that do not exist. 

In classical logic, both material implication and Venn's "Modern" Square of Opposition hold that in the case of "vacuous zeros", both the RH and its negation the anti-RH ("All zeros of $\zeta(s)$ are off the critical line") are true, which means that RH is a paradox that violates the LNC, and triggers ECQ. 

In intuitionistic logic, in the case of "vacuous zeros", the RH is false, because Dirichlet series $\zeta(s)$ is proven to have no zeros. So the zeros of RH are proven to be unconstructable. 

In contrast, in the 3VLs discussed in this paper, in the case of "vacuous zeros", both the RH and its negation the anti-RH ("All zeros of $\zeta(s)$ are off the critical line") are true, which means that RH is a paradox that is assigned the third truth-value.

\subsubsection{If Analytic Continuation of $\zeta(s)$ is True and LNC is True}

If the analytic continuation of $\zeta(s)$ is true, it creates the paradox of $\zeta(s)$ being both convergent and divergent throughout half-plane $\text{Re}(s)\le1$. All paradoxes violate the LNC, so the analytic continuation of $\zeta(s)$ and the LNC cannot\textit{both} hold true simulataneously.
\footnote{Intuitionism rejects this use of the LEM. See Brouwer \cite{Brouwer4}, p.23: "The axiom of the \textit{solvability of all problems} as formulated by Hilbert in 1900 \cite{Hilbert1901} is equivalent to the logical Principle of the Excluded Middle; therefore, since there are no sufficient grounds for this axiom and since logic is based on mathematics - and not vice versa - the use of the Principle of the Excluded Middle is \textit{not permissible} as part of a mathematical proof", and p.27, fn.4: "However, in his more recent lecture \textit{Axiomatic Thinking} \cite{Hilbert1917}, (p.412), Hilbert qualifies the question of the solvability of all mathematical problems by calling it a question still to be solved." }
In any logic that has both LNC and ECQ (e.g. classical and intuitionistic logics), this violation of the LNC triggers ECQ, which in turn renders any other proposition "trivially true". Here, it is the Riemann hypothesis which is rendered "trivially true" by ECQ. In fact, even in the stricter "relevance logics", which require that the antecedent and consequent of an implication to be "relevantly" related,
\footnote{See Wikipedia \cite{Relevance}, citing Routley et al. \cite{Routley} and Mares \cite{Mares}: "Relevance logic aims to capture aspects of implication that are ignored by the 'material implication' operator in classical truth-functional logic, namely the notion of relevance between antecedent and conditional of a true implication."}
the RH is "trivially true" due to ECQ, because the RH is directly related to the function $\zeta(s)$. 

\subsubsection{If Analytic Continuation of $\zeta(s)$ is True and LNC is False}

The upper right-most entry of RH's truth table is where the analytic continuation of $\zeta(s)$ is true, and the LNC is false. In this scenario, neither classical nor intuitionistic logic be used, because both logics have LNC. 

What is needed is a logic that permits paradoxes, such as the example 3VLs discussed in this paper: Bochvar's 3VL and Priest's $LP$. However, even in these two 3VLs, the RH is an unprovable paradox. In Bochvar's 3VL, material implication has the 3rd truth-value if it is from a paradox to any other proposition. 

So in Bochvar's 3VL, if ($\zeta(s)=0$) has the 3rd truth-value, then its material implication to any other proposition has the 3rd truth -value ("paradox"). So the result is always the 3rd truth -value ("paradox"). 

In Priest's $LP$, there is no material implication for paradoxes. More specifically, if ($\zeta(s)=0$) has the 3rd truth-value, then its material implication to a true proposition is a "quasi-valid" truth - but only if there are no paradoxical statements involved. Which is the case here, because the AC of $\zeta(s)$ renders $\zeta(s)=0$ a paradox.
\footnote{See Priest, \cite{Priest5}, p.235, \S IV.8: "The proposal is that we allow ourselves quasi-valid inferences even though they are not generally valid. We do know that quasi-valid inferences are truth preserving provided that there are no paradoxical sentences involved (see Section IV.1). Hence, if we were certain that we were not dealing with paradoxical sentences, we could use quasi-valid rules with a clear conscience. "}
But here there is a paradoxical statement involved. So again we are stuck with a truth-value of "paradox".

In contrast, if ($\zeta(s)=0$) is a false proposition, then it materially implies anything, which incldes both ($\text{Re}(s) = 1/2$) and ($\text{Re}(s) \ne 1/2$). Again, a paradox. See the two tables immediately below, Table \ref{tab:TT_RH_Mat} of RH, and Table \ref{tab:TT_Anti_RH} of Anti-RH ("All zeros of $\zeta(s)$ are off the critical line.").

\begin{table}[h] 
\begin{tabular}{cc|c|c|}
 & & \multicolumn{2}{c}{$\text{Re}(s)=1/2$}  \\
 &  & True & False  \\ 
\hline
 & Paradox (Convergent & No Implication (Priest's $LP$) & No Implication (Priest's $LP$) \\ 
$\zeta(s) = 0$ & to Zero, \& Divergent) & Paradox (Bochvar's 3VL) & Paradox (Bochvar's 3VL) \\ 
\cline{2-4}
 & False (Convergent to  & True & True  \\ 
 & Not Zero, \& Divergent) & &  \\  
\hline
\end{tabular}
\label{tab:TT_RH_Mat}
\textbf{\caption{RH as Material Implication, if AC of $\zeta(s)$ is True and LNC is False}} 

\end{table}

\begin{table}[h] 
\begin{tabular}{cc|c|c|}
 & & \multicolumn{2}{c}{$\text{Re}(s)\ne1/2$}  \\
 &  & True & False  \\ 
\hline
 & Paradox (Convergent  & No Implication (Priest's $LP$) & No Implication (Priest's $LP$) \\ 
$\zeta(s) = 0$ & to Zero, \& Divergent) & Paradox (Bochvar's 3VL) &  Paradox (Bochvar's 3VL) \\ 
\cline{2-4}
 & False (Convergent to  & True & True  \\ 
 & Not Zero, \& Divergent) & &  \\  
\hline
\end{tabular}
\label{tab:TT_Anti_RH}
\textbf{\caption{Anti-RH as Material Implication, if AC of $\zeta(s)$ is True and LNC is False}} 

\end{table}

\section{Some Implications in Mathematics} 

The falsity of analytic continuation of $\zeta(s)$, such that $\zeta(s)$ is exclusively defined by Dirichlet series $\zeta(s)$, has far-reaching implications. Some of these implications are discussed below.

\subsection{Prime Number Theorem}

Borwein et al. \cite{Borwein} states: 
"The proof of the prime number theorem relies on showing that $\zeta(s)$ has no zeros of the form $1 + it$ for $t \in \mathbb{R}$", 
\footnote{See Borwein et al. \cite{Borwein}, p.61}
and also states:
\begin{quotation}
In fact, this statement is equivalent to the prime number theorem, namely 
\begin{equation}
\pi(x) \sim \frac{x}{\log x}, x \rightarrow \infty
\end{equation}
(a problem that required a century of mathematics to solve).
\footnote{See Borwein et al. \cite{Borwein}, p.16.}
\end{quotation}
and further states:
\begin{quotation}
[W]e present part of de la Vallée Poussin’s proof of the prime number theorem (see Section 12.4); in particular, we prove that $\zeta(1 + it) \ne 0$ for $t \in \mathbb{R}$.
\footnote{See Borwein et al. \cite{Borwein}, p.9.}
\end{quotation}
Edwards \cite{Edwards} concurs:
\begin{quotation}
Since $\text{Re} \rho \le 1$ for all $\rho$ (by the Euler product formula - see Section 1.9), this amounts to proving that there are no roots $\rho$ [of Riemann's $\zeta(s)$] on the line $\text{Re}( s) = 1$. Thus, given von Mangoldt's 1894 formula for $\psi(x)$, the proof of the prime number theorem can be reduced to proving that there are no roots $\rho$ on the line $\text{Re}(s) = 1$ and to proving that the above limit can be evaluated termwise.
\footnote{See Edwards \cite{Edwards}, p.68.}
\end{quotation}
Edwards \cite{Edwards} also states that:
\begin{quotation}
Hadamard's proof that there are no roots $\rho$ on $\text{Re}(s) = 1$ is given in Section 4.2. 
De la Vallée Poussin admitted that Hadamard's proof was the simpler of the two, and although simpler proofs have since been found (see Section 5.2), Hadamard's is perhaps still the most straightforward and natural proof of this fact.
\footnote{See Edwards \cite{Edwards}, p.69.}
\end{quotation}

Borwein et al. \cite{Borwein} concludes with: "Thus [the prime number theorem] follows from the truth of the Riemann hypothesis." \footnote{Id.} Unfortunately, Borwein is wrong. There is no such relationship. (It appears that Borwein arrives at this conclusion because de la Vallée Poussin’s proof assumes that Riemann's analytic continuation of $\zeta(s)$ is true.)

When $\zeta(s)$ is defined as Riemann's $\zeta(s)$, the proof that $\zeta(1 + it) \ne 0$ for all $t \in \mathbb{R}$ is  "\textit{nontrivial}" (according to Borwein).
\footnote{See Borwein et al. \cite{Borwein}, p.16: "However, the proof that the zero-free region includes the
vertical line $\text{Re}(s) = 1$ (i.e., $\zeta(1 + it) \ne 0$ for all $t \in \mathbb{R}$) is already nontrivial."}
But when $\zeta(s)$ is defined as the Dirichlet series $\zeta(s)$, the proof of this theorem is \textit{trivial}: The Dirichlet series of $\zeta(s)$ has no zeros, so $\zeta(1 + it) \ne 0$ for all $t \in \mathbb{R}$.

\subsection{Analogues of the RH}

There exist analogues of the RH that (allegedly) have been proven to be true. These proofs need to be revisited, due to the invalidity of Riemann's AC of $\zeta(s)$. These analogues are invalid due to violating the LNC (for the same reasons that Riemann's AC of $\zeta(s)$ violates the LNC), and they are also unsound, due to falsely assuming that Riemann's AC of $\zeta(s)$ is true. See, for example: 
\begin{enumerate}
    \item Hasse's proof of the RH for elliptic curves of genus 1, 
\footnote{See Milne \cite{Milne3}, p.3.}
    \item Weil's proof of the RH for elliptic curves of arbitrary genus $g$,
\footnote{See also Jannsen \cite{Jannsen}, pp.4-5: "More generally one can show the following result which goes back to E. Artin and F.K. Schmidt: for a smooth projective (geometrically irreducible) curve $X$ of genus $g$ over $\mathbb{F}_q$ one has:
\begin{equation}
Z(X, T) = \frac{P(T)}{(1 − T)(1 − qT)}    
\end{equation}
where $P(T)$ is a polynomial of degree $2g$ in $\mathbb{Z}[T]$, with constant coefficient $1$. Furthermore Hasse (for $g = 1$, as well as for elliptic curves) and Weil (for arbitrary $g$) proved that the zeros of $P(q^{−s})$ lie on the line $\text{Re}(s) = 1/2$. Applied to $\zeta(X, s) = Z(X, q^{−s})$ this proves the analogoue (conjectured by Artin) of the Riemann hypothesis in the case of function fields."}
 and 
    \item Deligne's proof of the Weil conjecture III (which is the function field analogue of the Grand Riemann Hypothesis).
\footnote{See Milne \cite{Milne3}, p.49.} 
\end{enumerate}
All of these alleged proofs include a violation of the LNC, caused by the analytic continuation of the Zeta function, and the consequently false determinations that the Zeta function has a pole and zeros, that its functional equation is valid, etc.

For example, the Weil-conjecture expressly assumes that analytic continuation of $\zeta(s)$ is valid:
\footnote{See Jannsen \cite{Jannsen}, p.5.}
 
\begin{quotation}
\textbf{Weil-conjecture} (proved by Deligne in 1973): Let $X$ be a geometric irreducible smooth projective variety$\mathbb{F}_q$. Define
\begin{equation}
Z(X, T) = \exp(\sum_{n=1}^{\infty} |X(\mathbb{F}_{q^n})|\frac{T^n}{n})\in \mathbb{Q}[[T]] 
\end{equation}
Then the following holds

\textbf{I}: $Z(X, T)$ is rational, i.e., in $\mathbb{Q}(T)$.
(In particular, this implies the existence of a meromorphic continuation of the zeta-function
$\zeta(X, s) = Z(X, q^{−s})$, for which the series initially only converges for $\text{Re}(s) >> 0$).
\end{quotation}
 
This "meromorphic continuation" of the zeta-function $\zeta(X, s)$ violates the LNC for the same reason that the "meromorphic continuation" of $\zeta(s)$ violates the LNC: The series $\zeta(X, s)$ is "initially" convergent  \textit{only} for $\text{Re}(s) >> 0$ (and thus "initially" must be divergent for all other values of $s$). The series $\zeta(X, s)$ cannot be both divergent and convergent for $\text{Re}(s) << 0$

Therefore, the zeros of the "meromorphic continuation" of $\zeta(X, s)$ do not exist, and the $\zeta(X, s)$ analogue of the RH is false in intuitionistic logic, a paradox that triggers ECQ in classical logic, and has a third truth-value in a 3VL. The author conjectures that the same applies to all other allegedly proven analogues of the RH.

\subsection{L-Functions, the Modularity Theorem, and the Hasse-Weil Theorem}

\subsubsection{L-Functions}

Katz et al. \cite{Katz} states:
\footnote{See Katz et al. \cite{Katz}, pp.3-4.}
\begin{quotation}
The Riemann Zeta Function is but the first of a zoo of zeta and $L$-functions for which we can ask similar questions. There are the Dirichlet $L$-functions $L(s,\chi)$ defined as follows: $q \ge 1$ is an integer, $\chi : (\mathbb{Z}/q\mathbb{Z})^{*} \rightarrow C^{*}$ a (primitive) character and we extend $\chi$ to $\mathbb{Z}$ by making it periodic, and $\chi(m) = 0$ if $(m, q) \ne 1$. Then
\begin{equation}
L(s, \chi) = \sum_{n=1}^{\infty} \chi(n)n^{−s} = \prod_{p} (1 − \chi(p)p^{−s})^{-1}.    
\end{equation}
\end{quotation}
Dirichlet $L$-functions are generalizations of Riemann's $\zeta(s)$:
\footnote{See Ash et al. \cite{Ash}, p.200: "Dirichlet's $L$-functions can be thought of as a generalization of the Riemann zeta-function $\zeta(s)$. In the next section, we will describe a monster generalization of $\zeta(s)$ called the Hasse-Weil zeta-function."}
\begin{quotation}
By analytic continuation, [the Dirichlet $L$-series, $L(s, \chi) = \sum \chi(s)/n^s$] can be extended to a meromorphic function on the whole complex plane, and is then called a Dirichlet $L$-function and also denoted $L(s, \chi)$.
\footnote{See Wikipedia \cite{Dirichlet_L-function}.}
\end{quotation}
Also, note that:
\begin{quotation}
Just as the Riemann zeta function is conjectured to obey the Riemann hypothesis, so the Dirichlet $L$-functions are conjectured to obey the generalized Riemann hypothesis.
\footnote{Id.}
\end{quotation}

In logics that have LNC and ECQ, the analytic continuation of $\zeta(s)$ violates the LNC, and triggers ECQ. Riemann's $\zeta(s)$ is merely one example of a Dirichlet $L$-function. Therefore, generalizations of Riemann's $\zeta(s)$ (such as $L$-functions) are unsound, because they falsely assume that Riemann's analytic continuation of $\zeta(s)$ is valid.

\subsubsection{Modularity Theorem}

Sutherland \cite{Sutherland} concisely describes the Modularity theorem (previously called the Taniyama-Shimura conjecture), as follows:
\begin{quotation}
 Every elliptic curve $E$/$\mathbb{Q}$ is modular.
\footnote{See Sutherland \cite{Sutherland}, p.13, Theorem 25.33.}
\end{quotation}
Weisstein \cite{Weisstein2} provides a more detailed description of the Modularity theorem:
\begin{quotation}
In effect, the conjecture says that every rational elliptic curve is a modular form in disguise. Or, more formally, the conjecture suggests that, for every elliptic curve $y^2=Ax^3+Bx^2+Cx+D$ over the rationals, there exist nonconstant modular functions $f(z)$ and $g(z)$ of the same level $N$ such that
\begin{equation}
[f(z)]^2=A[g(z)]^2+Cg(z)+D.  
\end{equation}
Equivalently, for every elliptic curve, there is a modular form with the same Dirichlet $L$-series.
\footnote{See Weisstein \cite{Weisstein2}.}
\footnote{See also Frey \cite{Frey}, \S5.2, p.19: "Theorem 5.1 - \textit{Tanayama's and the Hasse-Weil conjecture is equivalent with the existence of a non-trivial map $\phi: X_0(N_E)\to E$ defined over $\mathbb{Q}$}. We call an elliptic curve $E$ over $\mathbb{Q}$ \textbf{modular} if a map $\phi$ like in the theorem exists. With this notation we can reformulate Taniyama‘s conjecture:
Conjecture 4 (Taniyama-Shimura-Weil) — Every elliptic curve
defined over $\mathbb{Q}$ is modular."}
\end{quotation}
The above-cited quote from Weisstein \cite{Weisstein2} expressly refers to "Dirichlet $L$-\textit{series}" (not "Dirichlet $L$-\textit{functions}").
\footnote{See also Sutherland \cite{Sutherland}, \S 25.8: "Although we defined the $L$-function of an elliptic curve using an Euler product, we can always expand this product to obtain a Dirichlet series".} 
As long as the $L$-series are not analytically continued to become $L$-functions, they do not violate the LNC. Bruin \cite{Bruin} states the following regarding the relationship between the Modularity theorem and analytic continuation of $L$-functions of elliptic curves:
\begin{quotation}
The modularity theorem implies that $L$-functions of elliptic curves over $\mathbb{Q}$ admit an analytic continuation to all of $\mathbb{C}$. This is not at all obvious and there is no known direct way to prove it. 
\end{quotation}
However, according to material implication, a true proposition cannot materially imply a false proposition. So if the Modularity theorem is indeed true, it \textit{cannot imply} that $L$-functions of elliptic curves over $\mathbb{Q}$ admit an analytic continuation to all of $\mathbb{C}$.
(In fact, the existence of $L$-functions in general remains an unproven conjecture). 
\footnote{See Wikipedia \cite{L-function}: "It is this (conjectural) meromorphic continuation to the complex plane which is called an $L$-function."}
\footnote{See e.g. the statement in Bombieri \cite{Bombieri}, p.5: "Not a single example of validity or failure of a Riemann hypothesis for an $L$-function is known up to this date." }

\subsubsection{Hasse-Weil Theorem (a Corollary of the Modularity Theorem)}

Wiles \cite{Wiles} states the following in regards to analytic continuation of $L$-functions, and the Hasse-Weil conjecture:
\begin{quotation}
Then we can define the incomplete $L$-series of C (incomplete because we omit the Euler factors for primes $p|2\Delta$) by
\begin{equation}
L(C, s) := \prod_{p|2\Delta} (1 − a_{p}p^{−s} + p^{1−2s})^{−1}
\end{equation}
We view this as a function of the complex variable $s$ and this Euler product is then known to converge for $\text{Re}(s) > 3/2$. A conjecture going back to Hasse (see the commentary on 1952(d) in [Weil \cite{Weil}]) predicted that $L(C, s)$ should have a holomorphic continuation as a function of s to the whole complex plane. This has now been proved [citing Wiles \cite{Wiles2}, Taylor et al. \cite{Wiles3}, and Breuill et al. \cite{Breuil}.]
\footnote{See Wiles \cite{Wiles}, p.2. This result violates the LNC, for the same reasons that analytic continuation of $\zeta(s)$ violates the LNC. } 
\end{quotation}
In fact, the first sentences of Wiles \cite{Wiles2} state the following:
\begin{quotation}
An elliptic curve over $\mathbb{Q}$ is said to be modular if it has a finite covering by a modular curve of the form $X_0(N)$. Any such elliptic curve has the property
that its Hasse-Weil zeta function has an analytic continuation and satisfies a functional equation of the standard type. 
\end{quotation}
But this cannot be true. The analytic continuation of the Hasse-Weil zeta function violates the LNC, and the "functional equation of the standard type" is not valid.

Sutherland \cite{Sutherland} also discloses the relationship between the Modularity theorem and the Hasse-Weil conjecture:
\begin{quotation}
When $E$ is modular, the $L$-function of $E$ is necessarily the $L$-function of a modular form, and this implies that $L_{E}(s)$ has an analytic continuation and satisfies a functional equation, since this holds for the $L$-function of a modular form ...
\footnote{See Sutherland \cite{Sutherland}, p.9: "Theorem 25.25 (Hecke). Let $f \in S_k(\Gamma_0(N))$. The $L$-function $L_f(s)$ extends analytically to a holomorphic function on $\mathbb{C}$, and the normalized $L$-function $\overline{L}_f(s) = N^{s/2}(2\pi)^{−s}\Gamma(s)L_f(s)$
satisfies the functional equation
$\overline{L}_f(s) = \pm \overline{L}_f(k − s)$."}
Prior to the proof of the Modularity theorem, this was an open question known as the Hasse-Weil conjecture; we record it here as a corollary to the Modularity Theorem.
\footnote{See also Sutherland \cite{Sutherland}, \S 25.9}
\end{quotation}

The Hasse-Weil theorem is a "corollary of the Modularity theorem" (according to Sutherland \cite{Sutherland}).
\footnote{See also Frey \cite{Frey}, p.17: "Conjecture 2 (Hasse-Weil) — $L_{E}(s)$ has an analytic continuation to
$\mathbb{C}$ satisfying the following functional equation..."; and p.19: "Taniyama stated the following conjecture: Conjecture 3. — Assume that the Hasse-Weil conjecture is true for the $L$-series $L_{E}(s) = \sum^{\infty}_{n=1} b_n n^{−s}$. Then $f_{E}(z) := \sum^{\infty}_{n=1} b_n e^{2\pi inz}$ is a cusp form."}
When an elliptic curve $E$ is modular, the Modularity theorem implies that the Hasse-Weil conjecture is true.
\footnote{See also Frey \cite{Frey}, p.20, Theorem 5.3.}
Moreover, Wiles's proof of Fermat's last theorem  assumes that the properties of modular elliptic curves (including the Hasse-Weil conjecture) are true.
\footnote{See Frey \cite{Frey}, pp.20-22, Theorem 5.3, \S6, and \S7.}

But the Hasse-Weil theorem is unsound, because it falsely assumes that the analytic continuation used to create Dirichlet $L$-functions is valid. Also, material implication in classical in intuitionistic logics holds that a true proposition cannot imply a false proposition. 

Therefore, if the Modularity theorem is true, it cannot materially imply a false Hasse-Weil theorem. According material implication, if the Hasse-Weil theorem is false, the Modularity theorem must be false as well. So the Modularity theorem must be false, and its progeny (e.g. Wiles's "proof" of Fermat's last theorem) must also be false. These "theorems" should never have been called theorems, because they are built upon unproven conjectures (i.e. Dirichlet $L$-functions).

\subsection{The Birch and Swinnerton-Dyer Conjecture} \label{BSD}

As discussed above, analytic continuation of the Dirichlet series $\zeta(s)$ to half-plane $\text{Re}(s)\le1$ violates the LNC, because the  Dirichlet series $\zeta(s)$ is proven to be divergent throughout half-plane $\text{Re}(s)\le1$. Riemann's $\zeta(s)$  violates not only the LNC in said half-plane, but also the Law of Identity (LOI) and the definition of a "function" in set theory  (due to the one-to-two relationship of domain to range). 

The Dirichlet series exclusively defines $\zeta(s)$, so at $s=1$, $\zeta(s)$ is the "harmonic series", which is proven to be divergent by the "Integral test for convergence".
\footnote{See, e.g., Guichard et al.'s \cite{Guichard1}, discussion of the \textit{Integral test for convergence}, at Theorems 13.3.3 and 13.3.4 and their proofs.}
This confirms that $\zeta(1)\ne0$. Also, the invalidity of analytic continuation of $\zeta(s)$ (in logics with LNC), and thus of analytic continuation of $L$-functions, disposes of the Landau-Siegel zero, "which no one believes exists". 
\footnote{See Conrey \cite{Conrey}, p.351: "The ineffectivity comes about from the assumption that some $L$-function actually has a real zero near $1$. Such a hypothetical zero of some $L$-function, which no one believes exists, is called a Landau-Siegel zero."}

According to Clay Mathematics Institute \cite{Clay}, this resolves the Birch and Swinnerton-Dyer (BSD) Conjecture in favor of finiteness: 
\begin{quotation}
[T]his amazing conjecture asserts that if $\zeta(1)$ is equal to 0, then there are an infinite number of rational points (solutions), and conversely, if $\zeta(1)$ is not equal to 0, then there is only a finite number of such points.
\end{quotation}
However, the BSD conjecture is unsound, because it falsely assumes that the analytic continuation used to create Dirichlet $L$-functions is valid.

\subsubsection{Hasse–Weil Zeta Function - 1st Example} \label{HWZ}

Further in regards to the BSD Conjecture, one example of the Hasse–Weil zeta function is for a nonsingular plane projective curve $C$, given by a homogeneous equation $F(x, y, z)=0$ with integer coefficients of degree $d$. 
\footnote{See Ash et al. \cite{Ash}, p.201.} 
\begin{quotation}
Let's continue our example with $C=P^{1}$, a projective line. To get the Hasse-Weil zeta function [we solve:]
\footnote{See Ash et al. \cite{Ash}, p.204.} 
\begin{equation} \label{eq:11.4}
Z(P^1, s) = \prod_p (1 - \frac{1}{p^s})^{-1} \cdot (1 - \frac{p}{p^s})^{-1}
\end{equation}
\end{quotation}
The Euler product of the Riemann zeta-function is:
\footnote{See Ash et al. \cite{Ash}, p.175, Eq. 11.13.}
\begin{equation} \label{eq:11.5}
\zeta(s) = \prod_p (1 - \frac{1}{p^s})^{-1}
\end{equation}
When the Euler product (Eq. \ref{eq:11.5}) is substituted into the Hasse-Weil zeta function (Eq. \ref{eq:11.4}), the result is:
\footnote{See Ash et al. \cite{Ash}, p.204.} 
\begin{equation} 
Z(P^1, s) = \zeta(s)\cdot  \zeta(s-1)
\end{equation} 

Given that given that analytic continuation of $\zeta(s)$ violates the LNC, $\zeta(s)$ is defined exclusively by the Dirichlet series $\zeta(s)$, which has no zeros. Neither $\zeta(s)$ nor $\zeta(s-1)$ can equal zero. Therefore, the Hasse-Weil zeta function $Z(P^1, s)$, which is the product of $\zeta(s)$ and $\zeta(s-1)$, is non-zero for all $s \in \mathbb{C}$.

\subsubsection{Hasse–Weil Zeta Function - 2nd Example}

Another version of the Hasse–Weil zeta function 
holds that the zeta function $Z_{E, \mathbf{Q}}(s)$ of elliptic curve $E$ over rational number field $\mathbb{Q}$ of conductor $N$ is:
\footnote{See Wikipedia \cite{HasseWeil}, citing Silverman \cite{Silverman} \S C.16, and Serre \cite{Serre}.} 
\begin{equation} Z_{E, \mathbf{Q}}(s) = \frac{\zeta(s)\cdot  \zeta(s-1)}{L(E,s)}
\end{equation} 

Again, given that analytic continuation of $\zeta(s)$ violates the LNC, $\zeta(s)$ is defined exclusively by the Dirichlet series $\zeta(s)$, which has no zeros. So neither $\zeta(s)$ nor $\zeta(s-1)$ can equal zero, and therefore their product, which is this version of the Hasse-Weil zeta function $Z_{E, \mathbf{Q}}(s)$, is non-zero for all values of $s \in \mathbb{C}$. 

Moreover, rearranging the terms of this Hasse–Weil zeta function produces: 
\begin{equation} 
L(E,s) = \frac{\zeta(s)\cdot \zeta(s-1)}{Z_{E, \mathbf{Q}}(s)}
\end{equation} 

Since neither $\zeta(s)$ nor $\zeta(s-1)$ can equal zero, $L(E,s) \ne 0$ for all $s \in \mathbb{C}$. So at $s=1$, the function $L(E, 1) \ne 0$. Given this result, all modular elliptic curves $E$ have rank 0, and thus are finite. Thus resolving the BSD conjecture to finiteness.
\footnote{See Wiles \cite{Wiles} (citing Kolyvagin \cite{Kolyvagin}): "Kolyvagin showed in 1990 that for modular elliptic curves, if $L(C, 1) \ne 0$ then $r = 0$ and if $L(C, 1) = 0$ but $L'(C, 1) \ne 0$ then $r = 1$".} 

However, $L$-functions are generalizations of Riemann's $\zeta(s)$, whose analytic continuation violates the LNC. Likewise, $L$-functions are divergent throughout a half-plane, and their "analytic continuation" to this half-plane is invalid, because it violates the LNC. This result is \textit{inconsistent} with Wiles's \cite{Wiles} proof of analytic continuation of $L(C, s)$. That proof is invalid, because the analytic continuation of $L$-functions violate the LNC.
\footnote{See Wiles \cite{Wiles}, p.2: "A conjecture going back to Hasse ... predicted that $L(C, s)$ should have a holomorphic continuation as a function of $s$ to the whole complex plane. This has now been proved", citing, \textit{inter alia}, Wiles's \cite{Wiles2} proof of Fermat's last theorem.} 

\subsection{Finiteness of the Tate–Shafarevich Group and the Brauer Group}

According to Totaro \cite{Totaro} and \cite{Totaro2}, the resolution of the BSD Conjecture also resolves equivalent conjectures. Totaro \cite{Totaro2} lists a few:
\footnote{See Totaro \cite{Totaro2}, page 578.}
\begin{quotation}
To spell out the relations between the Tate conjecture and finiteness problems,
let $X$ be a smooth projective surface over a finite field $k$, and let $f$ be a morphism with connected fibers from $X$ onto a smooth projective curve $C$. Assume that the generic fiber $F$ of $f$, which is a curve over the function field $k(C)$, is smooth over $k(C)$. Let $J$ be the Jacobian of $F$; thus $J$ is an abelian variety over the global field $k(C)$. Then the following are equivalent:
\footnote{Citing Ulmer \cite{Ulmer}, Proposition 5.1.2 and Theorem 6.3.1.}
\begin{itemize}
\item the Tate conjecture holds for divisors on $X$;
\item the Brauer group of $X$ is finite;
\item the Tate–Shafarevich group of $J$ is finite;
\item the Birch–Swinnerton-Dyer conjecture holds for $J$.
\end{itemize}
\end{quotation}

As discussed in the preceding section, the BSD conjecture is unsound, because it falsely assumes that the analytic continuation of $\zeta(s)$ is true. Moreover, as discussed in the following section(s), the Tate conjecture is \textit{unsound}, due to same false assumption. Therefore, the equivalence between the BSD conjecture and the Tate conjecture, as described by Totaro, is correct. 

Moreover, Wiles's \cite{Wiles} official Clay Foundation description of the BSD conjecture states the following: 

\begin{quotation}
There is an analogous conjecture for elliptic curves over function fields. It has been proved in this case by Artin and Tate \cite{Tate} that the $L$-series has a zero of order at least $r$, but the conjecture itself remains unproved. In the function field case it is now known to be equivalent to the finiteness of the Tate–Shafarevich group. 
\footnote{Wiles's \cite{Wiles} p.2, citing Tate \cite{Tate}, and citing Milne's \cite{Milne} Corollary 9.7.}
\end{quotation}
If these conjectures are indeed equivalent to the BSD conjecture, \footnote{As stated by Totaro  \cite{Totaro2}, page 578.}
then they too are unsound. The analytic continuation of $L$-series that takes place in the cited Tate \cite{Tate} reference violates the LNC.
\footnote{See Tate \cite{Tate}, p.416: "It is generally conjectured that $L_S$ has an analytic continuation throughout
the $s$-plane. This general conjecture, which in principle underlies those of
Birch and Swinnerton-Dyer, has been verified in some special cases, notably for
$A$ of $C.M.$-type (Weil, Deuring, Shimura), in which case $L_S$ can be identified
as a product of Hecke $L$-series, and for some elliptic curves related to modular function fields, when $L_S$ can be related to modular forms (Eichler, Shimura)."}
The other finiteness conjectures are "inspired" by the BSD conjecture.
\footnote{See Tate \cite{Tate}, p.426: "Inspired by the work of Birch and Swinnerton-Dyer, in the way explained
below, Mike Artin and I conjecture ... The Brauer group $Br(X)$ is finite "}

\subsection{The Tate Conjecture}

\subsubsection{The Tate conjecture, Argument 1}

Regarding the Tate conjecture, Totaro \cite{Totaro} states:
 
\begin{quotation}
Tate and Milne proved the equivalence of two problems, the Tate conjecture for elliptic surfaces over finite fields and the Birch-Swinnerton-Dyer conjecture for elliptic curves over global fields of positive characteristic. Both problems remain open. See for example Ulmer’s notes \cite{Ulmer2} on elliptic curves over function fields. 
\end{quotation}
 
Ulmer \cite{Ulmer2}, p.6, \S 3, discloses the following regarding Zeta functions over a finite field:
 
\begin{quotation}
Let $\chi$ be a variety over the finite field $\mathbb{F}_q$. 

It follows that $\zeta(\chi, s)$ has a meromorphic continuation to the whole $s$ plane, with poles on the lines $\text{Re}(s) \in \{0, \ldots, \text{dim} \chi \}$ and zeroes on the lines $\text{Re}(s) \in \{1/2, \ldots, \text{dim} −1/2\}$. This is the analogue of the Riemann hypothesis for $\zeta(\chi, s)$.

... Thus $\zeta(C, s)$ has simple poles for $s \in \frac{2 \pi i}{\log q}\mathbb{Z}$ and $s \in 1 + \frac{2 \pi i}{\log q}\mathbb{Z}$ and its zeroes lie on the line $\text{Re}(s) = 1/2$.  
\end{quotation}
 
This "meromorphic continuation" of $\zeta(\chi, s)$ is analogous to that of the Dirichlet series $\zeta(s)$ in the original Riemann hypothesis. The simple poles are analogous to the simple pole of Riemann's $\zeta(s)$, and of course the "zeroes on the lines $\text{Re}(s) \in \{1/2, \ldots, \text{dim} −1/2\}$" are analogous to the RH's zeros on the line $\text{Re}(s) = 1/2$. 

Ulmer \cite{Ulmer2}, pp.31-32, then discloses Tate's first and second conjectures, as follows (emphasis added):
\begin{quotation}
\textbf{Conjecture 9.2} $(T_2(\chi))$. We have
\begin{equation}
\text{Rank NS}(\chi) = − \text{ord}_{s=1} \zeta(\chi, s)
\end{equation}

\textbf{Note that by the Riemann hypothesis, the poles of $\zeta(\chi, s)$ at $s = 1$ come from $P_2(\chi, q^{−s})$.} More precisely, using the cohomological formula (4.1) of Lecture $0$ for $P_2$, we have that the order of pole of $\zeta(\chi, s)$ at $s = 1$ is equal to the multiplicity of $q$ as an eigenvalue of $Fr_q$ on $H^2(\overline{\chi}, \mathbb{Q}_ℓ)$.

Thus we have a string of inequalities:
\begin{equation}
\text{Rank NS}(\chi) \le \text{dim}_{\mathbb{Q}_ℓ} H^2(\overline{\chi}, \mathbb{Q}_ℓ)^{Fr_q = q} \le − \text{ord}_{s=1} \zeta(\chi, s)   \end{equation}

[Tate's first conjecture] $T_1(\chi)$ is that the first inequality is an equality and [Tate's second] conjecture $T_2(\chi)$ is that the leftmost and rightmost integers are equal. It follows trivially that $T_2(\chi)$
implies $T_1(\chi)$. Tate proved the reverse implication.
\end{quotation}

Prior to the "meromorphic continuation" discussed in Ulmer \cite{Ulmer2}, either $\zeta(\chi, s)$ is divergent at the values of $s$ covered by said "continuation", or $\zeta(\chi, s)$ has no value at these "pre-continuation" values of $s$. In both scenarios, the "meromorphic continuation" results in two conflicting definitions for certain values of $s$, thereby violating LNC and triggering ECQ.

Therefore, as discussed regarding the original Riemann Hypothesis pertaining to $\zeta(s)$, "the meromorphic continuation" of $\zeta(s)$ violates the LNC, and is invalid in logics with LNC, so $\zeta(\chi, s)$ has neither zeros nor poles. 

Moreover, Tate's conjectures are unsound, because they falsely assume that Riemann's "meromorphic continuation" of $\zeta(s)$ is valid.

In intuitionistic logic, the proof that the poles and zeros of Riemann's $\zeta(s)$ are non-existent is sufficient to render the Tate conjecture false. In contrast, in classical logic, a material implication with a "vacuous subject" (such as a proposition regarding a non-existent pole) is both true and false, resulting in an undecidable paradox that violates LNC and triggers ECQ. In certain 3VLs (e.g. Priest's $LP$), such a paradox does not violate LNC, and does not cause ECQ. 

\subsubsection{The Tate conjecture, Argument 2}

Milne \cite{Milne4} states the following (emphasis added):
\footnote{See Milne \cite{Milne4}, p.3.} 
\begin{quotation}
THEOREM 1.4. Let $X$ be a variety over $\mathbb{F}$ of dimension $d$, and let $r \in N$. The following statements are equivalent:

(a) $T^r(X,l)$ and $E^r(X,l)$ are true for a single $l$.

(b) $T^r(X,l)$, $S^r(X,l)$, and $T^{d-r}(X,l)$ are true for a single $l$.

(c) $T^r(X,l)$, $E^r(X,l)$, $S^r(X,l)$, $E^{d-r}(X,l)$, and $T^{d-r}(X,l)$ are true for all $l$, and the $\mathbb{Q}$-subspace $A_l^r(X)$ of $T_l^r(X)$  generated by the algebraic classes is a Q-structure on $T_l^r(X)$, i.e. $A_l^r(X) \bigotimes_\mathbb{Q} \mathbb{Q}_l \simeq  T_l^r(X)$

(d) \textbf{the order of the pole of the zeta function $Z(X,t)$ at $t=q^{-r}$ is equal to the rank of the group of numerical equivalence classes of algebraic cycles of codimension $r$.}

\end{quotation}
 
However, the original Riemann $\zeta(s)$ is not valid in logics with LNC in the half-plane of the analytic continuation. The resulting exclusive definition of $\zeta(s)$, the Dirichlet series $\zeta(s)$, is convergent in one half-plane, divergent in the other half-plane, \textbf{and has neither zeros nor poles}.
\footnote{See Rowland et al. \cite{Rowland}: "The word 'pole' is used prominently in a number of very different branches of mathematics. Perhaps the most important and widespread usage is to denote a singularity of a complex function."}
\footnote{See also Wikipedia \cite{Poles}, "Definitions": "The characterization of zeros and poles implies that zeros and poles are isolated, that is, every zero or pole has a neighbourhood that does not contain any other zero and pole.".}
This applies to the generalizations of $\zeta(s)$, as well.  Therefore, clause (d) of Milne's \cite{Milne4} Theorem 1.4, which Milne calls  "the full Tate conjecture"
\footnote{See Milne \cite{Milne4}, p.3, discussion of Theorem 1.4, last line.}
applies to \textbf{a pole that does not exist}.

In intuitionistic logic, the proof that the pole is non-existent is sufficient to render the Tate conjecture false. In contrast, in classical logic, a material implication with a "vacuous subject" (such as a proposition regarding a non-existent pole) is both true and false, resulting in an undecidable paradox that violates LNC and triggers ECQ. In certain 3VLs, such a paradox does not violate LNC, and does not cause ECQ. 

The above analysis also applies to clauses (a), (b) and (c) of Milne's \cite{Milne4} Theorem 1.4.

\subsection{The Hodge conjecture}

Several references expressly state that the Tate conjecture is equivalent to the Hodge conjecture in the case of abelian varieties of $CM$-type. 

As discussed above, the Tate conjecture is unsound due to its false assumption that the analytic continuation of the Zeta function $Z(X,t)$ is true, and that consequently to $\zeta(s)$ has a pole. This falsity of the Tate conjecture results in the falsity of the Hodge conjecture too, because the two are equivalent in the case of abelian varieties of $CM$-type. The invalidity of the Hodge conjecture in this one specific case is sufficient to invalidate it in general.

Gordon \cite{Gordon} states at 
page 364, \S 11.2:
\begin{quotation}
The main result of Pohlmann \cite{Pohlmann} is that for abelian varieties of $CM$-type, the Hodge and Tate conjectures are equivalent. Then that the validity of the Tate conjecture for an abelian variety $A$ implies the validity of the Hodge conjecture for $A$ has been proved by Piatetskii-Shapiro \cite{Piatetskii-Shapiro}, Deligne (unpublished) and Deligne \cite{Deligne}. Borovoi \cite{Borovoi} extends the result of Piatetskii-Shapiro \cite{Piatetskii-Shapiro}, and Borovoi \cite{Borovoi2} contains a weaker version of the main theorem of Deligne \cite{Deligne}, from which Tate implies Hodge for abelian varieties follows as a corollary.
\end{quotation}

\begin{table}[h]
\begin{tabular}{lll} \\
\textbf{Year} & \textbf{Author} & \textbf{Topic} \\
1968 & Pohlmann \cite{Pohlmann} & Hodge if and only if Tate for $CM$-type \\
1971 & Piatetskii-Shapiro \cite{Piatetskii-Shapiro} & Tate implies Hodge \\
1974 & Borovoi \cite{Borovoi} & Tate implies Hodge \\
1977 & Serre \cite{Serre2} & Connections between Hodge and Tate conj. \\
1982 & Deligne \cite{Deligne} & Absolute Hodge cycles, Tate implies Hodge \\
\end{tabular}
\caption{\textbf{Chronological listing of work on the Hodge conjecture for abelian varieties.} (See Gordon \cite{Gordon}, p.366.)}
\end{table}

The cited Deligne \cite{Deligne} reference discloses the following:
\begin{quotation}
COROLLARY 6.2. Let $A$ be an abelian variety over $\mathbb{C}$. If Tate’s conjecture is true for $A$, then so also is the Hodge conjecture.
\footnote{See Deligne \cite{Deligne}, p.43.}
\end{quotation}

\begin{quotation}
REMARK 6.3. The last result was first proved independently by Piatetskii-Shapiro \cite{Piatetskii-Shapiro} and Deligne (unpublished) by an argument similar to that which concluded the proof of the main theorem. \textbf{(Corollary 6.2 is easy to prove for abelian varieties of $CM$-type; in fact, Pohlmann \cite{Pohlmann} shows that the two conjectures are equivalent in that case.)}  
\footnote{See Deligne \cite{Deligne}, p.43.}
\end{quotation}

Shioda \cite{Shioda} provides more details (emphasis added in bold:
\begin{quotation}
\textit{Abelian varieties of $CM$ type} ([Pohlmann \cite{Pohlmann}, \S 2]). In this case, Pohlmann gave a combinatorial description of the Hodge ring $\mathcal{B}^{*}(A)$ in terms of the action of the $CM$ field on the complex cohomology $H^{*}(A, \mathbf{C})$, \textbf{and proved the equivalence of the Hodge Conjecture and the Tate Conjecture for this type of abelian varieties.} There is given an explicit example (due to Mumford) of a 4-dimensional abelian variety of $CM$ type such
that $\mathcal{B}^{2}(A)\ne \mathcal{D}^{2}(A)$, for which Hodge (A, 2) is still unknown. 

According to Mumford \cite{Mumford}, an abelian variety $A$ is of $CM$ type in the extended sense (i.e. isogenous to a product of abelian varieties of $CM$ type in the usual sense) if and only if its Hodge group $\text{Hg}(A)$ is an algebraic torus. We have $\text{dim Hg}(A) \le \text{dim} A$, and $A$ is called \textit{non-degenerate} if equality holds (Kubota \cite{Kubota}, Ribet \cite{Ribet}). For an abelian variety $A$ of $CM$ type, the two conditions (i) $A$ is non-degenerate and (ii) $\mathcal{B}^{2}(A) = \mathcal{D}^{2}(A)$ seem closely related. A recent result of Ribet and Lenstra (private communication in May 1981) shows that (i) and (ii) are indeed equivalent if $A$ is an abelian variety with the $CM$ field which is an abelian extension of $\mathbf{Q}$.  Hazama \cite{Hazama} shows that if $A$ is simple, then (i) implies (ii) in general. 
\footnote{See Shioda \cite{Shioda}, page 60.}
\end{quotation}

Moreover, Beauville \cite{Beauville} discloses that:
\begin{quotation}
For most abelian varieties, the Hodge conjecture holds for trivial reasons: the algebra of Hodge classes is generated in degree one. 
\footnote{See Beauville \cite{Beauville}, p. 12, Corollary 5.5: If the algebra $Hdg^{*}(X)$ is generated by $Hdg^{1}(X)$, the Hodge conjecture holds for $X$.}
This is the case in particular:

* for a general abelian variety [Mattuck \cite{Mattuck}];

* for a product of elliptic curves [Tate \cite{Tate2}];

* for a simple
\footnote{A complex torus $T$ is simple if the only complex subtori it contains are $(0)$ and $T$.} 
abelian variety of dimension $p$, where $p$ is a prime number [Tankeev \cite{Tankeev}].
\footnote{See Beauville \cite{Beauville}, p.14.}
\end{quotation}

Note that Beauville's \cite{Beauville} statement that "for most abelian varieties, the Hodge conjecture holds for trivial reasons", was not originally intended to refer to "trivial truth" (as per ECQ). Yet ironically, the Hodge conjecture is indeed "trivially true" as per ECQ. The analytic continuation of the Zeta function violates the LNC, and triggers ECQ. This renders unsound any conjecture that assumes that the analytic continuation of the Zeta function is true, and consequently that the Zeta function has poles and zeros. Tate \cite{Tate2}, whose title is "Algebraic cycles and poles of zeta functions", does precisely this. Therefore, the Tate conjecture is unsound, due to false assumptions. 

Given that the Tate and Hodge conjectures are equivalent for the "trivial case" discussed in Tate \cite{Tate2}, and also for "abelian varieties of $CM$ type", the Hodge conjecture is unsound, because it is equivalent to the unsound Tate conjecture in these instances.

So the result in classical and intuitionistic logics is that the Tate and Hodge conjectures violate LNC and trigger ECQ, even if this can only be proven for the specific instances of "a product of elliptic curves" and "abelian varieties of $CM$ type". The unsoundness of the Hodge conjecture in these specific instances is sufficient to invalidate it in all other instances (in logics with LNC).

\subsection{Other Number Theory Conjectures}

\subsubsection{The Generalized Riemann Hypothesis (GRH), Extended Riemann Hypothesis (ERH), and Grand Lindelöf Hypothesis (GLH) }

According to the Wikipedia entry on the Generalized Riemann Hypothesis (GRH) \cite{GRH}:
 
\begin{quotation}
Various geometrical and arithmetical objects can be described by so-called global $L$-functions, which are formally similar to the Riemann zeta-function. One can then ask the same question about the zeros of these $L$-functions, yielding various generalizations of the Riemann hypothesis.  ...

Global $L$-functions can be associated to elliptic curves, number fields (in which case they are called Dedekind zeta-functions), Maass forms, and Dirichlet characters (in which case they are called Dirichlet $L$-functions). 

When the Riemann hypothesis is formulated for Dedekind zeta-functions, it is known as the extended Riemann hypothesis (ERH) 
\footnote{See also Chandrasekharan \cite{Chandrasekharan}, p.4: "If, on the other hand, one \textit{assumes} the 'extended Riemann hypothesis', that not only the Riemann zeta-function but all the $L$-functions, modulo $q$, of Dirichlet, have all their zeros in the critical strip on the critical line, one would get ..."}
and when it is formulated for Dirichlet $L$-functions, it is known as the generalized Riemann hypothesis (GRH).
\footnote{See discussions of the GRH in Iwaniec et al. \cite{Iwaniec} and Sarnak \cite{Sarnak}.}
\end{quotation}

There exist additional hypothesis derived from the RH, such as the Lindelöf Hypothesis,
\footnote{See the Wikipedia entry on the Lindelöf Hypothesis \cite{Lindelof}.} 
and the Grand Lindelöf Hypothesis (GLH), which is a generalization of the Lindelöf hypothesis.
\footnote{See discussions of the GLH in Iwaniec et al. \cite{Iwaniec} and Sarnak \cite{Sarnak}.}

All of these hypotheses are generalizations of the RH, and like the RH, they too falsely assume the truth of analytic continuation, and assume the existence of non-existent zeros. Their truth-values correspond to those of the RH, according to the logic applied (classical, intuitionistic, 3VL, etc.). 

\subsubsection{The Bloch-Kato Conjecture}

According to Boston \cite{Boston}, "the Bloch-Kato conjecture [is] a vast generalization of the Birch and Swinnerton-Dyer conjecture":
\footnote{See Boston \cite{Boston}, page cxvii.}
 
\begin{quotation}
In this way we can restate the desired inequality in terms of the order of a Selmer group being bounded by a special value of an $L$-function, and we have a case of the Bloch-Kato conjecture, a vast generalization of the Birch and Swinnerton-Dyer conjecture. 
\end{quotation}
Furthermore, according to Bellaïche \cite{Bellaiche}:
 
\begin{quotation}
In the case where $V = V_p(E)$, the Bloch-Kato conjecture is closely related to the Birch and Swinnerton-Dyer conjecture, so all results about the Birch and Swinnerton-Dyer conjecture give a result for the Bloch-Kato conjecture. For example, the combination of results of Gross-Zagier and Kolyvagin shows that for if $\text{ord}_{s=0} L\big(V_p(E), s\big) \le 1$, the Bloch-Kato conjecture is known for $V = V_p(E)$.
\footnote{See Bellaïche \cite{Bellaiche}, p. 50.}
\end{quotation}
 
The Bloch-Kato conjecture falsely assumes that the analytic continuation of $L$-functions is valid, that they have zeros, etc. So the Bloch-Kato conjecture is unsound in logics that have the LNC. 

Bellaïche \cite{Bellaiche} states that there is no direct relation between the Grand Riemann Hypothesis and the Bloch-Kato conjectures for special values (the "Tamagawa number conjecture"):

\begin{quotation}
However, be aware that there is no direct relation between the Grand Riemann Hypothesis, which is interested in the zeros of $L(V, s)$ on $(w + 1)/2 < \text{Re}(s) < w/2 + 1$ and the Bloch-Kato conjecture, which is concerned by the zeros of
$L(V, s)$ at integers.
\footnote{See Bellaïche \cite{Bellaiche}, p.39.}
\end{quotation}
 
Bellaïche overlooks the fact that $L$-functions are generalizations of Riemann's $\zeta(s)$. In logics that have the LNC, the analytic continuation of Riemann's $\zeta(s)$ is invalid, as is the "meromorphic continuation" of $L$-functions. 
\footnote{See Bellaïche \cite{Bellaiche}, p. 38, Conjecture 3.1: "Then the function $L(V, s)$ admits a meromorphic continuation on all the complex plane." See also Bellaïche \cite{Bellaiche}, p. 44: "We assume that the $L$-function $L(V, s)$ has a meromorphic continuation to the entire plane, in accordance to Conjecture 3.1."}
Bellaïche refers to analytical continuation as a "mysterious process",
\footnote{See Bellaïche \cite{Bellaiche}, p. 44.}
when in fact it is an invalid process in any logic that has the LNC. As with the RH and the BSD conjecture, the Bloch-Kato conjecture refers to non-existent zeros. Both the Grand Riemann Hypothesis and the Bloch-Kato conjecture are unsound.

\subsection{P vs. NP}
\label{P vs. NP}

Venn's "Modern" Square of Opposition (see Figure \ref{fig:venn_square2} below) resolves the $P$ vs. $NP$ question, by showing that $P \ne NP$. The author has not found any reference that applies this technique of logic to solve this specific problem.   

\begin{figure} [ht]
    \centering
    \includegraphics[scale=1.5]{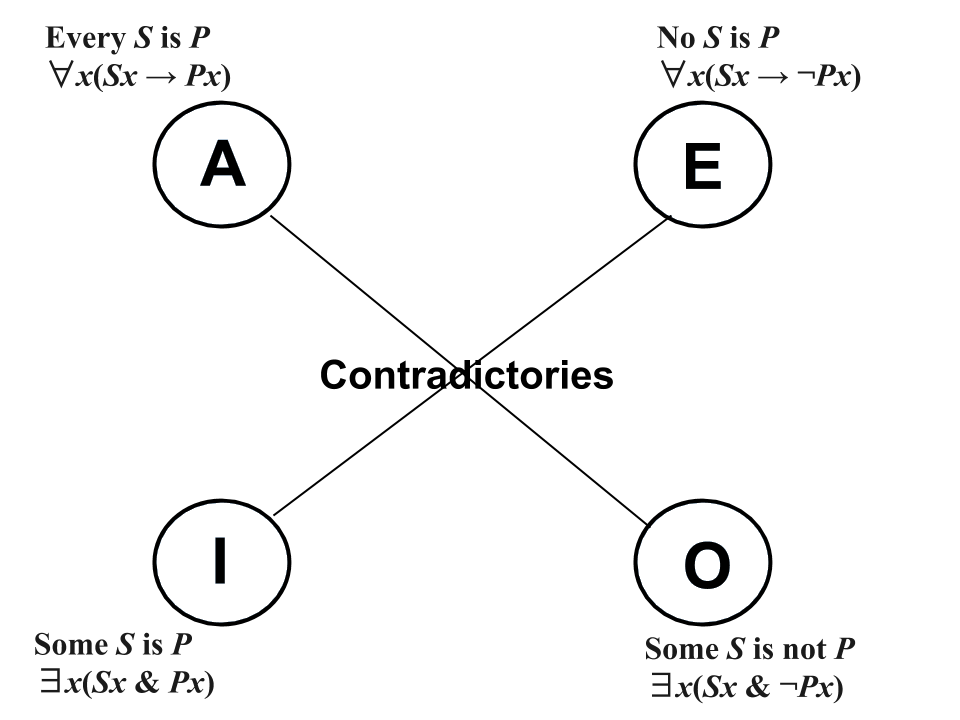}  
    \caption{Venn's "Modern" Square of Opposition}
    \label{fig:venn_square2}
\end{figure}

First, we assume classical logic as the foundational logic. The definition of "equivalence" in Whitehead and Russell's \textit{Principia Mathematica} \cite{Whitehead2} is: "two propositions are equivalent when, and only when, both are true or both are false" 
\footnote{See Whitehead and Russell's \cite{Whitehead2}, p.115, discussion of equivalence and Th. *4.01.}
\textit{PM} is a bivalent logic, and thus has no third truth-value such as "indeterminate" or "paradox". 

Next, we accept the conventional definitions of "$P$" as "the set of problems solvable in polynomial time", and "$NP$" as "the set of problems verifiable in polynomial time". Furthermore, "$P$-complete" is defined as the set of problems \textit{proven} to be in $P$, and "$NP$-complete" is defined as the set of problems \textit{proven} to be in $NP$.

So based on the definition of "equivalence" in \textit{Principia Mathematica} \cite{Whitehead2}, $P=NP$ only if both of the following propositions $A_1$ and $A_2$ are true:

\begin{quotation}
Proposition $A_1$: $P  \rightarrow NP$ 

Proposition $A_2$: $NP  \rightarrow P$ 
\end{quotation}
Moreover, in classical logic $P=NP$ is the same as $\neg P = \neg NP$, and $\neg P = \neg NP$ is the same as both of the following propositions $A_3$ and $A_4$ being true: 
\begin{quotation}
Proposition $A_3$: $\neg P  \rightarrow \neg NP$ 

Proposition $A_4$: $\neg NP  \rightarrow \neg P$
\end{quotation}

We analyze these four propositions ($A_1$, $A_2$, $A_3$, and $A_4$), and show that neither of the pairs ($A_1$ and $A_2$) or ($A_3$ and $A_4$) consists of two true propositions. Therefore, $P \ne NP$.

\subsubsection{Proposition A1: "Every P is NP"}

This is proposition $A_1$ and its related propositions, according to Venn's "Modern" Square of Contradiction:
\footnote{Compare Aristotle's and Venn's "Modern" Square of Contradiction (SoC) at Parsons \cite{Parsons}.}

\begin{quotation}
Proposition $A_1$: Every $P$ is $NP$: $\forall x(Px \rightarrow NPx)$

Proposition $E_1$: No $P$ is $NP$: $\forall  x(Px \rightarrow \neg NPx)$

Proposition $I_1$: Some $P$ is $NP$: $\exists x(Px \& NPx)$

Proposition $O_1$: Some $P$ is not $NP$:  $\exists x (Px \& \neg NPx)$  
\end{quotation}

Proposition $I_1$: There exists a problem that is solvable in polynomial time ($P$) AND is verifiable in polynomial time ($NP$). This proposition is TRUE (e.g. any $P$-complete problem. One example is the Circuit Value Problem (CVP). In fact, all $P$ problems are $NP$).

Proposition $O_1$: There exists a problem that is solvable in polynomial time ($P$) AND is not verifiable in polynomial time ($\neg$NP). This proposition is FALSE. (All $P$ problems are $NP$. All problems are verifiable in polynomial time, due to the existence of polynomial time sorting algorithms).

Therefore, because "A" contradicts "O", and "E" contradicts "I", we can determine the following:

Proposition $A_1$: For all problems, IF a problem is solvable in polynomial time ($P$) THEN it is verifiable in polynomial time ($NP$). This proposition is TRUE. (The confirmation is that all $P$ problems are $NP$).  

Proposition $E_1$: For all problems, IF a problem is solvable in polynomial time ($P$) THEN it is NOT verifiable in polynomial time ($\neg$NP). This proposition is FALSE. (The confirmation is that all $P$ problems are $NP$. All problems are verifiable in polynomial time, due to the existence of polynomial time sorting algorithms).

So $A_1$ is TRUE.

\subsubsection{Proposition A2: "Every NP is P"}

For proposition $A_2$, which is a material implication in the opposite direction of proposition $A_1$, Venn's "Modern" Square of Opposition is:

\begin{quotation}
Proposition $A_2$: Every $NP$ is $P$: $\forall x(NPx \rightarrow Px)$

Proposition $E_2$: No $NP$ is $P$: $\forall  x(NPx \rightarrow \neg Px)$

Proposition $I_2$: Some $NP$ is $P$: $\exists x(NPx \& Px)$

Proposition $O_2$: Some $NP$ is not $P$:  $\exists x (NPx \& \neg Px)$  
\end{quotation}

Proposition $I_2$: There exists a problem that is verifiable in polynomial time ($NP$) AND is solvable in polynomial time ($P$). This proposition is TRUE (e.g. a $P$-complete problem. All $P$ problems are $NP$).

Proposition $O_2$: There exists a problem that is verifiable in polynomial time ($NP$) AND is not solvable in polynomial time ($\neg$P). This proposition is TRUE (e.g. an $NP$-complete problem, such as Travelling Salesman Problem, which grows exponentially).

Therefore, because "A" contradicts "O", and "E" contradicts "I", we can determine the following:

Proposition $A_2$: For all problems, IF a problem is verifiable in polynomial time ($NP$) THEN it is solvable in polynomial time ($P$). This proposition is FALSE.  

Proposition $E_2$: For all problems, IF a problem is verifiable in polynomial time ($NP$) THEN it is NOT solvable in polynomial time ($\neg$P). This proposition is FALSE.

So both of the propositions $A_2$ and $E_2$ are FALSE. Some problems are verifiable in polynomial time ($NP$) and solvable in polynomial time ($P$). Other problems are verifiable in polynomial time ($NP$) BUT NOT solvable in polynomial time ($\neg P$).

This result, that proposition $A_2$ is FALSE, is sufficient to show that $P\ne NP$. 

\subsubsection{Proposition A3: "Every 'not P' is 'not NP'"}

For the sake of completeness, we also we also evaluate the pair of propositions $A_3$ and $A_4$, in order to show that $\neg P\ne \neg NP$. According to Venn's "Modern" Square of Contradiction, this is Proposition $A_3$ and its related propositions:

\begin{quotation}
Proposition $A_3$: Every $\neg P$ is $\neg NP$: $\forall x(\neg Px \rightarrow \neg NPx)$

Proposition $E_3$: No $\neg P$ is $\neg NP$: $\forall  x(\neg Px \rightarrow NPx)$

Proposition $I_3$: Some $\neg P$ is $\neg NP$: $\exists x(\neg Px \& \neg NPx)$

Proposition $O_3$: Some $\neg P$ is $NP$:  $\exists x (\neg Px \&  NPx)$  
\end{quotation}

Proposition $I_3$: There exists a problem that is NOT solvable in polynomial time ($\neg P$) AND is NOT verifiable in polynomial time ($\neg NP$). This proposition is FALSE (e.g. All problems are verifiable in polynomial time, due to the existence of polynomial time sorting algorithms.)

Proposition $O_3$: There exists a problem that is NOT solvable in polynomial time ($\neg P$) AND is verifiable in polynomial time ($NP$). This proposition is TRUE (e.g. Due to sorting algorithms being in polynomial time)

Therefore, because "A" contradicts "O", and "E" contradicts "I", we can determine the following:

Proposition $A_3$: For all problems, IF a problem is NOT solvable in polynomial time ($\neg P$) THEN it is NOT verifiable in polynomial time ($\neg NP$). This proposition is FALSE (All problems are verifiable in polynomial time, due to the existence of polynomial time sorting algorithms). 

Proposition $E_3$: For all problems, IF a problem is NOT solvable in polynomial time ($\neg P$) THEN it is verifiable in polynomial time ($\neg$NP). This proposition is TRUE (confirmed, due to sorting algorithms being in polynomial time).

So $A_3$ is FALSE. This is sufficient to show that $\neg P \ne \neg NP$.

\subsubsection{Proposition A4: "Every 'not NP' is 'not P'"}

For proposition $A_4$, which is a proposition in the opposite direction of proposition $A_3$, Venn's "Modern" Square of Opposition is:

\begin{quotation}
Proposition $A_4$: Every $\neg NP$ is $\neg P$: $\forall x(\neg NPx \rightarrow \neg Px)$

Proposition $E_4$: No $\neg NP$ is $\neg P$: $\forall  x(\neg NPx \rightarrow  Px)$

Proposition $I_4$: Some $\neg NP$ is $\neg P$: $\exists x(\neg NPx \& \neg Px)$

Proposition $O_4$: Some $\neg NP$ is $P$:  $\exists x (\neg NPx \&  Px)$  
\end{quotation}

Proposition $I_4$: There exists a problem that is NOT verifiable in polynomial time ($\neg NP$) AND is NOT solvable in polynomial time ($\neg P$). This proposition is FALSE. (All problems are verifiable in polynomial time, due to the existence of polynomial time sorting algorithms).

Proposition $O_4$: There exists a problem that is NOT verifiable in polynomial time ($\neg NP$) AND is solvable in polynomial time ($\neg$P). This proposition is FALSE. (All problems are verifiable in polynomial time, due to the existence of polynomial time sorting algorithms).

Therefore, because "A" contradicts "O", and "E" contradicts "I", we can determine the following:

Proposition $A_4$: For all problems, IF a problem is NOT verifiable in polynomial time ($\neg NP$) THEN it is NOT solvable in polynomial time ($\neg P$) . This proposition is TRUE, because "O" is FALSE, and also due to material implication, which is always true if the antecedent is a "vacuous subject" (as it is in this proposition). 

Proposition $E_4$: For all problems, IF a problem is NOT verifiable in polynomial time ($\neg NP$) THEN it is solvable in polynomial time ($\neg$P). This proposition is TRUE, because "I" is FALSE, and also due to material implication, which is always true if the antecedent is a "vacuous subject" (as it is in this proposition).

So $A_4$ and $E_4$ together form a PARADOX. 

\subsubsection{$P \ne NP$, Because A2 is False and A4 is a Paradox}

In classical logic, in order for $P=NP$, both $A_1$ and $A_2$ must be true, or both $A_3$ and $A_4$ must be true. 
\footnote{See Whitehead and Russell's \cite{Whitehead2} discussion of equivalence and Th. *4.01 on p.115 ("It is obvious that two propositions are equivalent when, and only when, both are true or both are false.").}
But that is not the case. Out of the pair $A_1$ and $A_2$, $A_2$ is FALSE, so only $A_1$ is TRUE. Out of the pair $A_3$ and $A_4$, $A_4$ is a PARADOX, so only $A_3$ is exclusively TRUE. So $P \ne NP$.

\section{Some Implications in Physics} \label{Phys}

The invalidity of analytic continuation of $\zeta(s)$ in logics with LNC means that in physics, "Zeta Function Regularization" violates the LNC and triggers ECQ, thereby rendering "trivially true" every physics model that uses it. 

This paper points out a few articles in the physics literature where this "regularization" (it actually is a "contradiction") is used, in models pertaining to Yang-Mills theory, the Casimir Effect, Quantum Electrodynamics (QED), Chromodynamics (QCD), Supersymmetry (SUSY), Quantum Field Theory (QFT), and Bosonic String Theory.

\subsection{Riemann Zeta Function Regularization}

Physicists have a procedure they call "Riemann zeta function regularization", that replaces the Dirichlet series $\zeta(s)$ with Riemann's $\zeta(s)$, whenever the former produces divergent values. This "regularization" introduces a contradiction whenever it is used, thus violating the LNC, and rendering the relevant mathematical proof "trivially true" in any logic with both LNC and ECQ. 
\footnote{See Bell \cite{Bell}, p.33, citing Dirac \cite{Dirac}: "[Dirac] divided the difficulties of quantum mechanics into two classes, those of the first class and those of the second. The second-class difficulties were essentially the infinities of relativistic quantum field theory. Dirac was very disturbed by these,
and was not impressed by the 'renormalisation' procedures by which they are circumvented. Dirac tried hard to eliminate these second-class difficulties, and urged others to do likewise."}
\footnote{See Dirac \cite{Dirac}: "I am inclined to suspect that the renormalization theory is something that will not survive in the future, and that the remarkable agreement between its results and experiment should be looked on as a fluke."}

In addition, all physics arguments (e.g. two-dimensional Yang-Mills theory) that falsely assume that $\zeta(s)$ is convergent for values of $s$ in half-plane $\text{Re}(s)\le1$, \footnote{See, e.g. Witten \cite{witten1991}: Eq. 2.20, Eq.2.32, and Eq.3.22, etc.}
even without explicit reference to "Riemann zeta function regularization", are unsound, and thus "trivially true" due to LNC and ECQ.

Hawking \cite{Hawking} describes the use of Riemann Zeta function regularization as:

\begin{quotation}
... a technique for obtaining finite values to path integrals for fields (including the gravitational field) on a curved spacetime background or, equivalently, for evaluating the determinants of differential operators such as the four-dimensional Laplacian or D'Alembertian. 
\footnote{See Hawking \cite{Hawking}, p.133, \S 1 Introduction, 1st para.}
\end{quotation}
 
According to Dittrich \cite{Dittrich}: 

\begin{quotation}
[In] many local relativistic quantum field theory models of elementary particles, ... Riemann’s results are of utmost importance for handling infinities with the aid of his zeta-function regularization.
\footnote{See Dittrich \cite{Dittrich}, p.3.}
\end{quotation}

Moreover, according to Bilal et al. \cite{Bilal}:
\begin{quotation}
We emphasize the close relationship between zeta function methods and arbitrary spectral cutoff regularizations in curved spacetime. This yields, on the one hand, a physically sound and mathematically rigorous justification of the standard zeta function regularization at one loop and, on the other hand, a natural generalization of this method to higher loops. In particular, to any Feynman diagram is associated a generalized meromorphic zeta function.
\footnote{See Bilal et al. \cite{Bilal}, Abstract.}
\end{quotation}
This despite the following:
\begin{quotation}
In spite of its power and elegance, the zeta function approach suffers from two important drawbacks. The first drawback, shared with dimensional regularization, is the absence of any obvious reason for why precisely it works. Even though replacing sums like $\sum_{n>0} n$ by $\zeta_{R}(−1) = −1/12$ is a perfectly well-defined procedure in the mathematical sense, it is abstract and unphysical.
\footnote{Note: It is not "well-defined procedure in the mathematical sense", and in fact is illogical.}
It is clear that the analytic continuation subtracts the divergence, as required, but it is very unclear how it does so explicitly and why the remaining finite part is the actual correct physical value.
\footnote{See Bilal et al. \cite{Bilal}, 4th page.}
\end{quotation}

\subsection{Yang-Mills Theory} 
\label{Yang}

Witten \cite{witten1991} describes two-dimensional quantum Yang-Mills Theory (YMT) from three different "points of view": 

\begin{enumerate}
\item Standard physical methods, 
\item Relating YMT to the large \textit{k} limit of three-dimensional Chern-Simons theory, and two-dimensional conformal field theory, and 
\item Relating the weak coupling limit of YMT to the theory of Reidmeister-Ray-Singer Torsion.
\end{enumerate}

The abstract of Witten \cite{witten1991} states that the results obtained from these three points of view are in agreement, and "give formulas for the volumes of the moduli spaces of representations of fundamental groups of two dimensional surfaces." However, each of these three points of view use Riemann's version of $\zeta(s)$, which is invalid for values of $s$ in half-plane $\text{Re}(s)\le 1$. So all three "points of view" of 2D YMT are "trivially true" in logics with LNC and ECQ. \footnote{See Witten \cite{witten1991}, p.154, description of Eq. 1.2: "... and $\Sigma$ [is] a Riemann surface of genus $g$, one finds $Vol(\textit{M})=2\cdot (2\pi^2)^{1-g} \cdot \zeta(2g-2)$, where $\zeta(s)$ is the Riemann zeta function".  But Riemann's $\zeta(s)$ violates LNC, so $Vol(\textit{M})$ is divergent at g=0 and g=1.} 
\footnote{Witten's \cite{witten1991}, p.154, description of Eq. 1.2 also refers (in regards to Eq. 3.18) to the Hurwitz zeta function and Dirichlet $L$-functions. These generalizations of Riemann's $\zeta(s)$ inherit Riemann's $\zeta(s)$ falsity in half-plane $\text{Re}(s)\le 1$.}
\footnote{Witten's \cite{witten1991}, p.174 (last para.) states: "We will formulate this in a way that exhibits the relation to IRF models - which also appear, after a much more difficult analysis, in computing Wilson line expectation values in three dimensional Chern-Simons theory [25]. For convenience, we will consider first the case that $\Sigma$ has genus zero." However, Eq.1.2 with $g=0$ produces a divergent $Vol(\textit{M})$.}
\footnote{Witten's \cite{witten1991}, p.159, description of Eq. 2.20: "With an explicit choice (such as zeta function regularization) for defining the determinants that appear in evaluating the left and right-hand sides of (2.20), an \textit{a priori} computation of \textit{$\Delta v$} can be given." The so-called "Zeta function regularization" replaces the Dirichlet series $\zeta(s)$ with the false Riemann $\zeta(s)$, for values of $s$ in half-plane $\text{Re}(s)\le 1$.} 
\footnote{See also Witten's \cite{witten1991}, p.161, description of Eq. 2.28: "We will ensure this by using the zeta function definition of determinants [3]". But the Dirichlet series $\zeta(s)$ is divergent at $s=0$.}
\footnote{See also Witten's \cite{witten1991}, p.178, Eq. 3.8, which includes Riemann's $\zeta(s)$. Eq. 3.8 is divergent at $g=0$ and $g=1$.} 
\footnote{See also Witten's \cite{witten1991}, p.180, Eq.3.22: "The Hurwitz zeta function ... is then continued holomorphically throughout the complex $z$ plane, except for a pole at $z=1$". This is false.}
\footnote{See also Witten's \cite{witten1991}, p.201, Eq. 4.95: "with $\zeta(s)$ the Riemann zeta function". This is divergent if $\text{Re}(s)\le 1$.}

Moreover, Aguilera-Damia et al. \cite{Aguilera-Damia} applies Zeta-function regularization to $N=4$ super-Yang-Mills theory:
\begin{quotation}
Using $\zeta$-function regularization, we study the one-loop effective action of fundamental strings in
$AdS_5$ × $S^5$ dual to the latitude $\frac{1}{4}$-BPS Wilson loop in $N=4$ super-Yang-Mills theory. To avoid certain ambiguities inherent to string theory on curved backgrounds we subtract the effective action of the holographic $\frac{1}{2}$-BPS Wilson loop. We find agreement with the expected field theory result at first order in the small latitude angle expansion but discrepancies at higher order.
\footnote{See Aguilera-Damia et al. \cite{Aguilera-Damia}, abstract.}
\end{quotation}
 
So because of zeta-function regularization, $N=4$ super-Yang-Mills theory is "trivially true" in logics with LNC and ECQ.

\subsection{Casimir Effect, QED, and QCD}

Dittrich \cite{Dittrich} states that "Riemann Zeta Function Regularization" is used to derive the Casimir effect. 
\footnote{See Dittrich \cite{Dittrich}, pp.30-34.}  
Tong \cite{Tong} confirms that this is the case for "Casimir Energy".
\footnote{See Tong \cite{Tong}, pp.38-40. Tong's discussion on Casimir Energy begins on p.38 with the following quote attributed to Ramanujan, in a letter to G.H.Hardy: "“I told him that the sum of an infinite no. of terms of the series: $1 + 2 + 3 + 4 + ... =  1/12$ under my theory. If I tell you this you will at once point out to me the lunatic asylum as my goal." Ramanujan was aware that this equation violates the rules of arithmetic.}

Dittrich \cite{Dittrich} also states that: "The same procedure finds application in QED and QCD." 
\footnote{See Dittrich \cite{Dittrich}, p.34.} 
If true, then the Casimir effect, Quantum Electrodynamics (QED), and Quantum Chromodynamics (QCD) are all "trivially true" in logics with LNC and  ECQ.

\subsection{Supersymmetry (SUSY)}

According to Elizalde \cite{Elizalde}, Supersymmetry (SUSY) incorporates Riemann Zeta Function Regularization:

\begin{quotation}
Regularization and renormalization procedures are essential issues in contemporary physics — without which it would simply not exist, at least in the form known today (2000). They are also essential in supersymmetry calculations. Among the different methods, zeta-function regularization — which is obtained by analytic continuation in the complex plane of the zeta-function of the relevant physical operator in each case — might well be the most beautiful of all. Use of this method yields, for instance, the vacuum energy corresponding to a quantum physical system (with  constraints of any kind, in principle).
\footnote{See Elizalde \cite{Elizalde}, 1st para. It appears that "\textit{with}  constraints" should be "\textit{without} constraints".}
\end{quotation}
 
Therefore, due to the use of zeta-function regularization, Supersymmetry (SUSY) is "trivially true" in logics with LNC and ECQ. 

\subsection{Quantum Field Theory (QFT)}

According to Elizalde \cite{Elizalde}, Riemann Zeta Function Regularization is also used in Quantum Field Theory (QFT):

\begin{quotation}
These mathematically simple-looking relations involve very deep physical concepts (no wonder that understanding them took several decades in the recent history of quantum field theory, QFT). The zeta-function method is unchallenged at the one-loop level, where it is rigorously defined and where many calculations of QFT reduce basically (from a mathematical point of view) to the computation of determinants of elliptic pseudo-differential operators ...
\footnote{See Elizalde \cite{Elizalde}, 2nd para.}
\end{quotation}
 
Penrose \cite{Penrose} goes further, saying that:
 
\begin{quotation}
Whatever philosophical position is taken on this issue, renormalization is an essential feature of modern QFT. Indeed, as things stand, there is no accepted way of obtaining finite answers without such an 'infinite rescaling' procedure applied not necessarily only to charge, or mass, but to other quantities also. Theories in which this kind of procedure works are called \textit{renormalizable}. In a renormalizable QFT, it is possible to collect together all the divergent parts of the Feynman graphs into a finite number of 'parcels' which can be 'scaled away' by renormalization, any remaining divergent expressions being deemed to cancel out with each other \footnote{See Penrose \cite{Penrose}, \S26.9, p.678.}
\end{quotation}

However, if QFT uses "zeta-function normalization", then it is "trivially true" in logics with LNC and ECQ. Moreover, it is clear that physicists care neither about this specific logical problem (violation of LNC and ECQ), nor about the more general problem of the logical foundations of mathematics. Here is Penrose \cite{Penrose} again:

\begin{quotation}
It is a common standpoint, among particle physicists, to take renormalization as a selection principle for proposed theories. Accordingly any non-renormalization theory would be automatically rejected as inappropriate to Nature. 
\footnote{See Penrose \cite{Penrose}, \S26.9, p.678.}
\end{quotation}
So any particle physics theory without a glaring logical contradiction is "inappropriate to Nature"?
\footnote{Note that Niels Bohr held that violation of LNC is a core principle of quantum physics. He chose the motto "Contraria Sunt Complementa" ("Opposites are Complementary") for his coat of arms, when inducted into the Danish Order of the Elephant in 1947. See Wikipedia \cite{Bohr}, citing Wheeler \cite{Wheeler}.}
Penrose also states:
\begin{quotation}
Many (and perhaps even most) physicists would take the view that the framework of QFT is 'here to stay', and that the blame for any inconsistencies (these being usually from infinities coming from divergent integrals, or from divergent sums, or both) lies in the particular scheme to which QFT is being applied, rather than in the framework of QFT itself.
\footnote{See Penrose \cite{Penrose}, \S26.1, p.656.}
\end{quotation}
So in summary: according to many physicists, any particle physics theory without the contradiction inherent in "renormalization" is automatically rejected, but the blame for any inconsistencies in accepted theories does not lie in the framework of QFT itself. This is madness. (Especially because vacuous subjects generate paradoxes. It is entirely possible that certain hypothesized particles do not exist, and hence generate paradoxes).    
\subsection{Bosonic String Theory}

There are several examples in Bosonic string theory of the use of Riemann's $\zeta(s)$, and the functional equation of a relationship between $\zeta(s)$ and $\zeta(1-s)$.

The He et al. \cite{He} reference links Riemann's $\zeta(s)$ to expressions of the Veneziano amplitude
\footnote{See Wikipedia \cite{Veneziano1}, citing Veneziano \cite{Veneziano2}. See also Turco et al. \cite{Turco} (unpublished). }
that describe the scattering of four bosonic open strings with tachyonic masses. This is based on work by Freund et al. \cite{Freund}, whose abstract states: 
\begin{quotation}
We show that the Veneziano and Virasoro-Shapiro four-particle scattering amplitudes can be factored in terms of an infinite product of non-archimedean string amplitudes. This factorization is equivalent to the functional equation for the Riemann zeta function.
\end{quotation}

Toppan \cite{Toppan} provides a description of the heat-kernel method and of generalized Riemann's zeta-functions associated to elliptic operators. (These generalized zeta-functions violate the LNC in half of their respective domains, just as the original Riemann $\zeta(s)$ does). Toppan \cite{Toppan} then defines their role in defining one-loop partition functions for Euclidean Field Theories. 

Toppan \cite{Toppan} then applies these results to the Polyakov functional quantization of the closed bosonic string, to
derive its critical dimensionality of $D=26$. Núñez \cite{Nunez} confirms the use of Zeta function regularization in obtaining the "trivially true" dimensionality of   $D=26$ .
\footnote{See Núñez \cite{Nunez}, Eq.105, bottom of p.17 to top of p.18.}
Therefore, Bosonic String theory is "trivially true" in logics with LNC and ECQ.

Moreover, in Bosonic string theory, the mass of states in lightcone gauge is:
\footnote{See PhysicsOverflow \cite{PhysicsOverflow}, citing Tong \cite{Tong}, Eq. at top of p.39.}
\begin{equation} \label{M^2}
M^{2} = \frac{4}{\alpha'}\Big[\sum_{n=1}^{\infty}\alpha^{i}_{−n}\alpha^{i}_{n} + \frac{D-2}{2}\Big(\sum_{n=1}^{\infty}n \Big)\Big]  
\end{equation}
The Dirichlet series $\zeta(s)$ is divergent at $s=-1$, but Riemann's $\zeta(s)$ at $s=-1$ is:
\footnote{See PhysicsOverflow \cite {PhysicsOverflow}, citing Tong \cite{Tong}, Eq. at middle of p.39.}
\begin{equation}
\zeta(−1) = −1/12
\end{equation}
If the value of Riemann's $\zeta(s)$ at $s=-1$ is substituted for the divergent Dirichlet series $\zeta(s)$ at $s=-1$ (thereby violating the LNC, and triggering ECQ), the mass of states is:
\footnote{See PhysicsOverflow \cite {PhysicsOverflow}, citing Tong \cite{Tong}, Eq.2.26, p.39.}
\begin{equation}
M^{2} = \frac{4}{\alpha'}\Big( N − \frac{(D-2)}{24} \Big)
\end{equation}
At the ground state $N=0$, the formula simplifies to:
\footnote{See PhysicsOverflow \cite {PhysicsOverflow}, citing Tong \cite{Tong}, Eq.2.27, p.40.}
\begin{equation}
M^{2} = \frac{-(D−2)}{6\cdot \alpha'}
\end{equation}
which corresponds to a particle with an imaginary mass, known as a tachyon. Moreover, at the first excited state ($N=1$), the Equation \ref{M^2} is massless ($M^{2}=0$) at $D=26$. 

These results are "trivially true" in logics with LNC and ECQ. If $\zeta(s)$ is defined by the Dirichlet series, then the mass of states in lightcone gauge (Equation \ref{M^2}) is as follows:
At $D=2$,
\begin{equation} 
M^{2} = \frac{4\cdot N}{\alpha'} 
\end{equation}
At all other values of $D$, the value of $M^{2}$ is divergent. Moreover, Equation \ref{M^2} is massless (i.e. $M^2=0$) only if both $D=2$ and $N=0$, or if $D=2$ and $\alpha'$ is infinitesimal. Moreover, at $D\ne 2$ and $\alpha'$ is infinitesimal, the value of $M^{2}$ is divergent.

\subsection{Riemann's Zeta Function and the Failure of LOI in Quantum Physics}

As stated in a previous section of the present paper, the convergent Riemann's $\zeta(s)$) in half-plane $\text{Re}(s)\le 1$ (except at $s=1$), where the Dirichlet series $\zeta(s)$ is proven divergent, is a violation of the Law of Identity (LOI) for the function $\zeta(s)$. In this scenario, $\zeta(s)$ is not equal to itself in half-plane $\text{Re}(s)\le 1$ (except at $s=1$).

Given that quantum physics extensively uses Riemann's $\zeta(s)$ in "zeta-function regularization", it is not surprising to see published articles that state that LOI fails in quantum physics (emphasis added):
\footnote{See French, \cite{French}.}
 
\begin{quotation}
However, it has also been argued that quantum physics is in fact compatible with a metaphysics of individual objects, but that such objects are indistinguishable in a sense \textit{which leads to the violation of Leibniz’s famous Principle of the Identity of Indiscernibles}. This last claim has recently been contested in a way that has reinvigorated the debate over the impact of the theory.
\end{quotation}

This leads to the questions: what remains of quantum physics if "renormalization" (e.g. use of Riemann's $\zeta(s)$) is no longer permitted? Will LOI hold true in whatever remains? An additional question: will whatever remains be able to explain experimental results?

\subsection{3VL in Physics}

\subsubsection{Schrödinger's Cat}

Classical logic is the assumed logical foundation of the "Schrödinger Cat" illustration of the Copenhagen interpretation of quantum mechanics, as evidenced by concerns about its violation of the LNC. As Baggott \cite{Baggott} states:
\begin{quotation}
On the surface, it really seems as though we ought to be able to resolve this paradox with ease. But we can't. There is obviously no evidence for peculiar superposition states of live-and-dead things or of 'classical' macroscopic objects of any description.
\footnote{See Baggott \cite{Baggott}, pp.133-134.}
\end{quotation}
In classical logic, the contradictory statements of "the cat is alive" and "the cat is dead" would violate the LNC, and trigger ECQ, if both were true simultaneously.
\footnote{See Griffiths \cite{Griffiths}, \S 10.1 "Schrödinger's Cat", citing Schrödinger \cite{Schroedinger}. See also da Costa et al. \cite{daCosta}, abstract, which argues  that "Schrödinger logics" (Non-reflexive logics) are "logical systems in which [Leibniz's] principle of identity is not true in general." See also Penrose \cite{Penrose}, \S29, especially \S29.7 - \S29.9, pp.804-812, which discusses the "paradox of Schrödinger's Cat", but only in the context of a 2VL.}

A more appropriate logical foundation for the "Schrödinger Cat" scenario is a 3VL. For example, Łukasiewicz's 3VL  has a third truth-value of "unknown", which is relevant for the state where it is unknown if the cat is alive or dead.
\footnote{Arguably, so do the "future contingents" discussed in Aristotle's \textit{De interpretatione} \S 9, in \textit{Organon}.}
Moreover, in a 3VL such as Łukasiewicz', a logical proposition having the 3rd state does not result in the entire model being "trivially true" due to LNC and ECQ. It is due to LNC and ECQ that "Schrödinger's Cat" is usually discussed in the context of probability theory (which has no truth-values) rather than logic (which does). 
\footnote{See Fronhöfer \cite{Fronhoefer}, p.2.}
\footnote{Also see Baggott \cite{Baggott}, pp.131-135, which discusses the paradox of "Schrödinger's Cat", but fails to consider non-classical logics.}
3VLs provide a truth-functional way of addressing the paradox, without "trivial truth", and without resorting to probability theory. 

Other logics that are appropriate for the "Schrödinger Cat" scenario are intuitionistic logic (that rejects LEM) and its variants, such as minimal logic (that rejects both LEM and ECQ). \footnote{See Wikipedia \cite{Minimal}, citing Johansson \cite{Johansson} and Troelstra et al. \cite{Troelstra}, p.37.}
Intuitionistic logic is applicable here because while the chamber containing the cat is sealed, outside observers cannot prove either that the cat is dead, or that it is alive. In other words, we have no "proof" for the cat being alive or dead. These logics acknowledge that there exist instances when neither proposition $A$ nor its negation $\neg A$ can be proven, which in classical logic would violate the LEM. As with Aristotle's "future contingents", and probability theory, intuitionistic logic acknowledges that there exist conditions of uncertainty, due to the limits of human knowledge. So in intuitionistic logic, we must acknowledge the limits of our knowledge, and concede that we do not know the status of the cat. This is an Epistemological issue.

\subsubsection{Particle/Wave Duality}

Another example of contradiction in physics is the particle/wave duality. In classical physics, which pertains to "large scale" phenomena, particles and waves are mutually exclusive categories. So in classical logic, "quantum scale" assumptions such as the dual nature of matter (and light) are paradoxes that, due to LNC and ECQ, would cause the "trivial truth" of classical physics theories. (Note that particle/wave duality also violates the LEM). 

As Penrose \cite{Penrose} states:
\begin{quotation}
These kinds of consideration led to the conclusion that an ordinary particle displays wavelike behavior, this having a universal relationship to the particle's rest-mass as determined by the Planck and de Broglie formulae. But, in the previous two decades, a converse to this had already been established, demonstrating that entities previously thought of as purely wavelike - basically Maxwell's oscillating electric and magnetic fields as the constituents of light
\footnote{Citing Penrose \cite{Penrose}, \S19.2} - had also to be viewed as having a \textit{particulate} nature, again consistent with Planck and de Broglie formulae. The most convincing evidence for this was in the \textit{photoelectric effect} ...
\footnote{Penrose \cite{Penrose},\S21.4, p.501. See also \S21.5, pp.505-507, and \S21.7, pp.511-515. }
\end{quotation}

Given that particle/wave duality is observed at the "quantum scale" but not at the "large scale", the logical foundation of the math used at the large scale can be classical logic. But at the "quantum scale", phenomena such as particle/wave duality must be described in a logic that rejects the LNC and ECQ (i.e. a non-classical logic). 

As with the paradox of "Schrödinger's Cat", a 3VL with a third truth-value is a good candidate. One obvious candidate is Priest's 3VL, which has a 3rd truth-value corresponding to "truth-value gluts". This is appropriate because light is both particle and wave, simultaneously. In such a 3VL, the dual nature of light is assigned to a third truth-value, instead of to a contradiction that causes the entire model to be "trivially true" due to LNC and ECQ. 

\subsubsection{Galilean Relativity and Special Relativity}

Cohen \cite{Cohen2} discloses the following issue raised by Immanuel Kant:
\begin{quotation}
Kant's early work is characterized by an attempt to identify internal contradictions in abstract metaphysical theories derived from pure logic. For example, Kant is concerned that although in logic either A or not-A is true, in reality, something can be both A and not-A. A physical object like a table on a train, for instance, can be both in motion and motionless since it depends on the position of the observer.
\footnote{See Cohen \cite{Cohen2}, pp.76, and 238-239, citing Kant \cite{Kant}, pp.203-242.}
\end{quotation}
However, Galileo preceeded Kant by over a century in raising this issue:
\begin{quotation}
Galilean invariance or Galilean relativity states that the laws of motion are the same in all inertial frames. Galileo Galilei first described this principle in 1632 in his \textit{Dialogue Concerning the Two Chief World Systems} using the example of a ship travelling at constant velocity, without rocking, on a smooth sea; any observer below the deck would not be able to tell whether the ship was moving or stationary.
\footnote{See Wikipedia \cite{Galilean_inv}. See also Penrose \cite{Penrose}, \S17.2 "Spacetime for Galilean relativity", pp.385-387.}
\end{quotation}

Kant's example has two independent "frames of reference", with each "frame of reference" having its own observer who is unaware of the other. Kant's example presents several problems regarding the LNC:

(1) The below-deck observer would determine that an object (such as a table) fixed to the ship is motionless. Physical experiment would confirm this result. But an  observer outside the ship would determine that the ship (and thus the table attached to it) is in motion relative to some other point. The table is both in motion and motionless - but not to the same observer. Each of the answers is subjectively true to its respective observer. To resolve the dilemma, either the two observers need to communicate with one another. or a third observer is needed to objectively determine that only the outside observer is correct.

(2) Also there is a possibility that the outside observer is in agreement with the below-deck observer, and both are wrong. For example, if the outside observer is on a spaceship travelling parallel to the below-deck observer's spaceship, in a featureless area of outer space, both observers will determine that the table is motionless (the wrong answer). Again a third observer is needed, with access to additional information (e.g. a reference point), in order to determine that the first two observers are wrong.

Therefore, the question is not "whether or not the table is in motion", but rather "whether or not the table is in motion in relation to point $x$ in space". So if no observer can observe "point x", and all observers are in an inertial state,
\footnote{Both constant speed and stillness are "inertial states". }
this necessitates a 3rd truth-value (e.g. "indeterminate") for the question "whether or not the table is in motion in relation to point $x$." This 3rd truth-value renders LNC and ECQ irrelevant. 

As with with the paradox of "Schrödinger's Cat", the 3rd truth-value in Galilean Relativity is necessitated by the observer's lack of critical information, not by some other characteristic of reality.

Norton \cite{Norton} adds the following, in regards to "relativity of simultaneity" in Einstein's special theory of relativity:
\footnote{See also Wikipedia \cite{Relativity}.}
\begin{quotation}
The relativity of simultaneity adds to the repertoire of quantities that are relative and not absolute. There is no absolute fact to whether a spaceship is moving uniformly or is at rest. It can only be said to be at rest relative to another body. There is no absolute fact as to whether a rod is a foot long or a process lasts for one minute. They can only true with respect an observer with a definite state of motion. To this list we add that there is no absolute fact to whether two spatially separated events are simultaneous; or whether two spatially separated clocks are synchronous. These can only be true relative to an observer with a definite state of motion.
\end{quotation}
So the truth value of statements pertaining to simultaneity also should be assigned the 3rd truth value (unless a specific frame of reference is specified). 

\subsubsection{Popper, Bohr, Einstein, and Bell}

Moreover, some propositions with the 3rd truth-value fail Popper's "falsifiability" test for scientific conjectures, because they are paradoxes that are both true and false (or neither).
\footnote{See Wikipedia \cite{Popper}: "To say that a given statement (e.g., the statement of a law of some scientific theory)—call it "T"—is "falsifiable" does not mean that "T" is false. Rather, it means that, \textbf{if} "T" is false, \textbf{then} (in principle), "T" could be \textbf{shown} to be false, by observation or by experiment. Popper's account of the logical asymmetry between verification and falsifiability lies at the heart of his philosophy of science."}
One example of this is the Riemann Hypothesis. 

Moreover, if Niels Bohr is correct regarding contradiction being an inherent characteristic of quantum physics, \footnote{See Wikipedia \cite{Bohr}, citing Wheeler \cite{Wheeler}.}
then its underlying logic must be able to cope with paradoxes, and thus must be non-classical. The classical logic that underlies mathematics (and thus classical physics too) is unable to cope with paradoxes, due to LNC and ECQ. If there is to be a unification of classical and quantum physics, it can only happen if the foundational logic is a non-classical logic (that accepts the paradoxes of quantum physics).

It is also noted that the paradoxes of both (a) value of the Schrödinger wave function prior to "collapse" (according to the Copenhagen interpretation), and (b) whether events are simultaneous in Einstein's special theory of relativity, are due to limits of what can be known. Another example is Heisenberg's uncertainly principle. In all of these cases, observers are barred from knowing the truth-value of a proposition.

Einstein "wanted things out there to have properties, whether or not they were measured".
\footnote{See Mermin \cite{Mermin}, p.38, citing Pais \cite{Pais}: "We often discussed his notions on objective reality. I recall that during one walk Einstein suddenly stopped, turned to me and asked whether I really believed that the moon exists only when I look at it." The best response to this question is to invoke David Hume's arguments regarding the “problem of induction”. See Henderson \cite{Henderson}.}
But there are limits to observer knowledge, even after measurement, as Heisenberg's uncertainty principle clearly shows.
This in turn raises another issue: What truth-value do we assign to unknowable propositions? Classical logic does not have an answer for this. So a non-classical logic (such as a 3VL) must be used instead.

Another famous philosophical question raised by quantum physics is: "If a tree falls in the forest, and there’s nobody around to hear, does it make a sound?"
\footnote{See Baggott \cite{Baggott2}, and Wikipedia \cite{Tree}.}
The physicist John Bell asked it in the following form: 
\begin{quotation}
What exactly qualifies some physical systems to play the role of 'measurer'? Was the wavefunction of the world waiting to jump for thousands of years until a single-celled living creature appeared? Or did it have to wait a little longer, for some better qualified system ... with a PhD?
\footnote{See Baggott \cite{Baggott}, p.134, citing Bell \cite{Bell}, p.34.}
\end{quotation}
As Karl Popper would have gladly explained, Bell's questions are "unfalsifiable" philosophical questions that fall outside of the purview of science. 
\footnote{See, for example, debates pertaining to "consciousness" in Van Gulick \cite{VanGulick}.}

\section{Conclusion}

Analytic continuation of $\zeta(s)$ violates the LNC, because it contradicts the proven divergence of the Dirichlet series $\zeta(s)$ in the half-plane $\text{Re}(s)\le1$. According to Aristotle's LOI, LEM, and LNC, any "analytic continuation" of $\zeta(s)$ to the half-plane $\text{Re}(s)\le1$ is false.  
Therefore, in logics that include LNC and ECQ, the falsity of analytic continuation of $\zeta(s)$ renders "trivially true" all arguments that falsely assume the truth of the "analytic continuation" of $\zeta(s)$. 

Moreover, because the analytic continuation of $\zeta(s)$ is false, the Dirichlet series exclusively defines $\zeta(s)$, and therefore $\zeta(s)$ has no zeros. Thus, both the Riemann Hypothesis (RH) and anti-RH ("All zeros of $\zeta(s)$" are off the critical line") are true propositions, due to their "vacuous subjects": the non-existent zeros of Dirichlet series $\zeta(s)$. This paradoxical result violates the LNC. So in classical and intuitionistic logics, ECQ renders "trivially true" all "proofs" that assume RH is true. In 3VLs that assign the 3rd truth -value to paradoxes (e.g. Buchavar's 3VL and Priest's "$LP$), the RH has the 3rd truth-value. 

The result that RH is a paradox causes all conjectures that assume it is true to be "trivially true" in logics with LNC and ECQ. This result, and also the result in 3VLs that RH has the third truth-value, is inconsistent with "proofs" of analogues of the RH, which claim to prove that the analogues of RH are "exclusively true" (not paradoxes). See e.g. (1) Hasse's proof of the RH for elliptic curves of genus 1, 
\footnote{See Milne \cite{Milne3}, p.3.}
(2) Deligne's proof of the Weil conjecture III,
\footnote{See Milne \cite{Milne3}, p.49.}
and (3) Weil's proof of the RH for elliptic curves of arbitrary genus $g$. \footnote{See also Jannsen \cite{Jannsen}, pp.4-5.}

All of these alleged proofs include a violation of the LNC, caused by the analytic continuation of an analogue of $\zeta(s)$, and the consequently false determinations that: the analogue of $\zeta(s)$ has a pole and zeros, that its functional equation is valid, etc.

Moreover, as stated in Chapter \ref{Consistency}, Langer's \cite{Langer} statement that "[c]ontradictory theorems cannot follow from consistent postulates"
\footnote{See Langer \cite{Langer}, p.202.}
is wrong. Contradictory theorems \textit{do} follow from consistent postulates, if the theorems are directed to "vacuous subjects", or if the postulates result in self-reference. Therefore, MacFarlane's \cite{MacFarlane} quote citing Meyer \cite{Meyer} on this topic requires clarification: 
 
\begin{quotation}
There’s no good reason to assume that mathematics must be consistent. If math is about a supersensible realm of objects, why should we assume they’re like ordinary empirical objects with respect to consistency? But if math is a free human creation, why can’t it be inconsistent?
\begin{quotation}
... for certain purposes an inconsistent system might be more useful, more beautiful, and even—at the furthest metaphysical limits—as the case may be, more accurate. 
\footnote{MacFarlane \cite{MacFarlane}, p.1, citing Meyer \cite{Meyer}, p.814.
}
\end{quotation}
Classical logic \textit{forces} math to be consistent or trivial.
\footnote{See MacFarlane \cite{MacFarlane}, p.1.}
\end{quotation}
 
However, it is not \textit{classical logic} per se that forces math to be consistent or trivial, rather it is ECQ and its prerequisite LNC that do so. Any logic that has LNC and ECQ (e.g. Intuitionistic logic) would, if assumed to be the foundation logic of math, force math to be consistent or trivially true. Moreover, Intuitionistic logic would do so in a more restrictively than classical logic, due to an insistence on constructive proof.  

In regards to a 3VL or 4VL as a possible foundation logic instead of classical logic, Hazen et al. \cite{Hazen3} (citing Dunn \cite{Dunn}), states that if one tries to formulate a second-order logic of a 3VL or 4VL, the resultant system collapses to its classical counterpart. 
\footnote{See Hazen et al. \cite{Hazen3}, p.507: "We are not sure what general morals to draw from all this. An obvious one to draw from the negative results of Sects. 6 and 7 is that many non-classical logics do not have well-behaved Second Order versions: something already shown, in a different way and for different logics in [Dunn \cite{Dunn}]."}
\footnote{See Dunn \cite{Dunn}, p.261: "In Dunn \cite{Dunn2} it was shown (among other things) that if one tries to formulate second-order quantum logic with a certain minimal principle of extensionality, one is doomed to failure in the sense that the resultant system collapses to its classical counterpart. It was remarked in In Dunn \cite{Dunn2} that this result is generalizable to a large class of non-classical logics, and this is the point of the present paper."}
Moreover, another Hazen et al. article (\cite{Hazen}) states: "it will be extremely difficult to appeal to [Priest's] second-order $LP$ for the purposes that its proponents advocate, until some deep, intricate, and hitherto unarticulated metaphysical advances are made." 
\footnote{See Hazen et al. \cite{Hazen}, abstract: "The logic of paradox, $LP$, is a first-order, three-valued logic that has been advocated by Graham Priest as an appropriate way to represent the possibility of acceptable contradictory statements. Second-order $LP$ is that logic augmented with quantification over predicates. As with classical second-order logic, there are different ways to give the semantic interpretation of sentences of the logic. The different ways give rise to different logical advantages and disadvantages, and we canvass several of these, concluding that it will be extremely difficult to appeal to second-order $LP$ for the purposes that its proponents advocate, until some deep, intricate, and hitherto unarticulated metaphysical advances are made."}

Note also that the most important "cost" of having classical logic as the foundational logic of math is that LNC and ECQ force math to be \textit{incomplete}, as formally proven by Gödel in his first incompleteness theorem (that utilizes the Liar paradox). If a multi-valued logic (e.g. 3VLs and 4VLs) were able to be the foundational logic of mathematics, and thus could assign a 3rd or 4th truth-value to paradoxes, this would make mathematics \textit{complete} (at the cost of being \textit{inconsistent}).

On an unrelated note: In his \textit{The History of Modern Philosophy}, Bertrand Russell states the following:
\footnote{See Russell \cite{Russell5}, p.202.}
\begin{quotation}
Throughout modern times, practically every advance in science, in logic, or in philosophy has had to be made in the teeth of the opposition from Aristotle's disciples.
\end{quotation}
Ironically, Russell's \textit{Principia Mathematica} includes Aristotle's "Laws of Thought" (LOI, LEM, and LNC) as theorems, which \textit{de facto} makes Russell a disciple of Aristotle. So Russell's statement can be interpreted as the liar's paradox. Also, Łukasiewicz's 3VL was derived from Aristotle's future contingents, was an advance in logic, and was not made "in the teeth of opposition from Aristotle's disciples."

On another unrelated note, the RH has been described as "[e]legant, crisp, falsifiable, and far-reaching" and "the epitome of what a conjecture should be".
\footnote{See Sarnak \cite{Sarnak}, first page; and Iwaniec et al. \cite{Iwaniec}, p.712.} 
In fact, the RH is a paradox, and thus \textit{not} falsifiable. So according to Karl Popper's philosophy of science, RH is not a "scientific question". This highlights an implicit assumption of Karl Popper's philosophy of science: paradoxes do not exist. 

Finally, we note that in the twenty years since the initial announcement of the Millenium Problems, \textit{none} of the official descriptions of the problems have ever listed "paradox" as a possible answer, nor has the mathematical community argued that it should be listed as a possible answer. This demonstrates that the mathematical community has still not internalized the results of Gödel's famous work.

\section*{Acknowledgments}

The author thanks the reviewer(s) and the editorial staff. The research for this paper did not receive any funding from any funding agency in the public, commercial, or not-for-profit sectors. All research was performed in the author's off-duty time.

\singlespacing
\bibliographystyle{acm}
\bibliography{sample.bib}

\begin{thebibliography}{100}

\bibitem{Abramowitz}
{\sc Abramowitz, M., and Stegun, I.~A.}, Eds.
\newblock {\em Handbook of mathematical functions}, 9th~ed.
\newblock National Bureau of Standards, 1970.
\newblock \url{http://people.math.sfu.ca/~cbm/aands/abramowitz_and_stegun.pdf}.

\bibitem{Aguilera-Damia}
{\sc Aguilera-Damia, J., Faraggi, A., Zayas, L.~P., Rathee, V., and Silva, G.}
\newblock Zeta-function regularization of holographic wilson loops.
\newblock {\em Physical Review D 98}, 4 (2018), 046011.
\newblock
  \url{https://journals.aps.org/prd/abstract/10.1103/PhysRevD.98.046011}.

\bibitem{Albert}
{\sc Albert, H.}
\newblock {\em Traktat über kritische Vernunft}.
\newblock Mohr, 1968.

\bibitem{Aloni}
{\sc Aloni, M.}
\newblock Disjunction.
\newblock In {\em The Stanford Encyclopedia of Philosophy}, E.~N. Zalta, Ed.,
  winter 2016~ed. Metaphysics Research Lab, Stanford University, 2016.
\newblock
  \url{https://plato.stanford.edu/archives/win2016/entries/disjunction/}.

\bibitem{Andrews2}
{\sc Andrews, F.~E.}
\newblock The principle of excluded middle then and now: Aristotle and
  principia mathematica.
\newblock {\em Animus}, 1 (1996), 53--66.
\newblock \url{https://core.ac.uk/download/pdf/11700623.pdf}.

\bibitem{Andrews}
{\sc Andrews, P.~B.}
\newblock {\em An Introduction to Mathematical Logic and Type Theory: To Truth
  Through Proof}, 2nd~ed.
\newblock Kluwer Academic Publishers, 2002.
\newblock
  \url{https://books.google.com/books?id=nV4zAsWAvT0C&pg=PA201#v=onepage&q&f=false}.

\bibitem{EoM}
{\sc Anonymous}.
\newblock Logarithmic function.
\newblock In {\em Encyclopedia of Mathematics}, {Last} modified on 7 {February}
  2011~ed. Springer-Verlag, 2011.
\newblock
  \url{http://www.encyclopediaofmath.org/index.php?title=Logarithmic_function&oldid=16249}.

\bibitem{Stanford}
{\sc Anonymous}.
\newblock {Glances} {Ahead}: {More} to {Think} {About}. {IV}. {The} {Law} of
  {Excluded} {Middle}.
\newblock In {\em An {Introduction} to {Philosophy}}. Stanford University
  Philosophy Dept., 2014.
\newblock
  \url{https://web.archive.org/web/20140731134904/http://web.stanford.edu/~bobonich/glances%20ahead/IV.excluded.middle.html
  }.

\bibitem{StackExchange2}
{\sc Anonymous}.
\newblock Equivalence of $a \rightarrow b$ and $\lnot a \lor b$, 2016.
\newblock
  \url{https://web.archive.org/web/20161229224743/https://math.stackexchange.com/questions/243949/equivalence-of-a-rightarrow-b-and-lnot-a-vee-b}.

\bibitem{Square}
{\sc Anonymous}.
\newblock Square of opposition.
\newblock In {\em The Internet Encyclopedia of Philosophy (IEP)}. Internet
  Encyclopedia of Philosophy, 2018.
\newblock
  \url{https://web.archive.org/web/20181224033855/https://www.iep.utm.edu/sqr-opp/}.

\bibitem{Apostol}
{\sc Apostol, T.~M.}
\newblock An elementary view of euler's summation formula.
\newblock {\em The American Mathematical Monthly 106}, 5 (May 1999), 409--418.
\newblock \url{https://www.jstor.org/stable/2589145}.

\bibitem{Aristotle2}
{\sc Aristotle}.
\newblock Metaphysics, 1994-2000.
\newblock \url{http://classics.mit.edu/Aristotle/metaphysics.mb.txt}, Last
  updated 2000.

\bibitem{Aristotle}
{\sc Aristotle}.
\newblock Organon: On interpretation, 2015.
\newblock \url{https://ebooks.adelaide.edu.au/a/aristotle/interpretation/},
  Last updated Wednesday, July 15, 2015.

\bibitem{Ash}
{\sc Ash, A., and Gross, R.}
\newblock {\em Elliptic Tales}.
\newblock Princeton University Press, 2012.

\bibitem{Asmus}
{\sc Asmus, C.}
\newblock Paraconsistency on the rocks of dialetheism.
\newblock {\em Logique et Analyse 55}, 217 (2012), 3–21.

\bibitem{Baggott2}
{\sc Baggott, J.}
\newblock Quantum theory: If a tree falls in the forest...
\newblock In {\em OUP Blog}. Oxford University Press, February 2011.
\newblock
  \url{https://web.archive.org/web/20110219141459/https://blog.oup.com/2011/02/quantum/}.

\bibitem{Baggott}
{\sc Baggott, J.}
\newblock {\em Farewell to Reality: How Modern Physics has Betrayed the Search
  for Scientific Truth}.
\newblock Pegasus Books, 2013.

\bibitem{Batty}
{\sc Batty, C.}
\newblock Strawson, “on referring”, Sept. 2009.
\newblock
  \url{http://www.uky.edu/~cebatt2/teaching/PHI_515_Sep2009/PHI_515-Strawson.pdf}.

\bibitem{Bauer}
{\sc Bauer, A.}
\newblock Five stages of accepting constructive mathematics.
\newblock {\em Bulletin of the American Mathematical Society 54}, 3 (July
  2017), 481–498.
\newblock \url{http://dx.doi.org/10.1090/bull/1556}.

\bibitem{Bazhanov}
{\sc Bazhanov, V.}
\newblock The dawn of paraconsistency: Russia's logical thought in the turn of
  xx century.
\newblock {\em Manuscrito 34}, 1 (2011), 89--98.
\newblock \url{http://www.scielo.br/pdf/man/v34n1/a04v34n1.pdf}.

\bibitem{Beall}
{\sc Beall, J., Glanzberg, M., and Ripley, D.}
\newblock Liar paradox.
\newblock In {\em The Stanford Encyclopedia of Philosophy}, E.~N. Zalta, Ed.,
  fall 2017~ed. Metaphysics Research Lab, Stanford University, 2017.
\newblock
  \url{https://plato.stanford.edu/archives/fall2017/entries/liar-paradox/}.

\bibitem{Beauville}
{\sc Beauville, A.}
\newblock The hodge conjecture, 2019.
\newblock
  \url{https://web.archive.org/web/20190221225036/https://pdfs.semanticscholar.org/3505/0485dfc7fbc7268d6a3de41e9b103d109c72.pdf}.

\bibitem{Beeson}
{\sc Beeson, M.~J.}
\newblock {\em Foundations of Constructive Mathematics}.
\newblock Springer-Verlag, 1985.
\newblock \url{https://link.springer.com/book/10.1007/978-3-642-68952-9}.

\bibitem{Bell}
{\sc Bell, J.}
\newblock Against 'measurement'.
\newblock {\em Physics World\/} (1990), 33--40.

\bibitem{Bellaiche}
{\sc Bellaïche, J.}
\newblock An introduction to the conjecture of bloch and kato, 2009.
\newblock
  \url{https://web.archive.org/web/20190201063004/http://virtualmath1.stanford.edu/~conrad/BSDseminar/refs/BKintro.pdf}.

\bibitem{Belnap}
{\sc Belnap, N.}
\newblock A useful four-valued logic.
\newblock In {\em Modern Uses of Multiple-Valued Logic}, J.~M. Dunn and
  G.~Epstein, Eds., vol.~2. D. Reidel Publishing Company, 1977, pp.~8--37.
\newblock \url{https://doi.org/10.1007/978-94-010-1161-7_2}.

\bibitem{Bezhanishvili}
{\sc Bezhanishvili, G., and Holliday, W.~H.}
\newblock A semantic hierarchy for intuitionistic logic, 2018.
\newblock [Online; accessed 17-March-2019].
  \url{https://escholarship.org/uc/item/2vp2x4rx}.

\bibitem{Beziau}
{\sc B\'eziau, J.-Y.}
\newblock Trivial dialetheism and the logic of paradox.
\newblock {\em Logic and Logical Philosophy 25\/} (2016), 51–56.
\newblock
  \url{https://pdfs.semanticscholar.org/f71e/241d0f89c7f162952fe00ec6c095210582f9.pdf}.

\bibitem{Beziau2}
{\sc B\'eziau, J.-Y.}
\newblock The lvov-warsaw school: A true mythology.
\newblock In {\em The Lvov-Warsaw School. Past and Present}, Ángel Garrido and
  U.~Wybraniec-Skardowska, Eds. Birkhäuser, 2018, pp.~779--815.
\newblock \url{http://www.jyb-logic.org/papirs/jyb-lws1.pdf}.

\bibitem{Bilal}
{\sc Bilal, A., and Ferrari, F.}
\newblock Multi-loop zeta function regularization and spectral cutoff in curved
  spacetime.
\newblock {\em Nucl.Phys. B. 877\/} (2013), 956--1027.
\newblock \url{https://arxiv.org/abs/1307.1689}.

\bibitem{Birdsong}
{\sc Birdsong, S.}
\newblock Summary of convergence tests.
\newblock
  \url{http://math2.uncc.edu/~sjbirdso/calc%20II-sum08/series/convergence%20tests.pdf
  }. [Online; accessed 9-July-2017].

\bibitem{Bochvar}
{\sc Bochvar, D.}
\newblock On a three-valued logical calculus and its application to the
  analysis of contradictories.
\newblock {\em Matematiceskij sbornik [Rec. Math. N.S.] 4(46)}, 2 (1938),
  287–308.
\newblock \url{http://mi.mathnet.ru/eng/msb/v46/i2/p287}.

\bibitem{Boeckh}
{\sc Boeckh, A.}
\newblock {\em Gesammelte kleine Schriften}, vol.~VI.
\newblock Druck und Verlag von B. G. Teubner, 1866.
\newblock \url{https://archive.org/details/augustboeckhsge05bratgoog/page/n9}.

\bibitem{Bolander}
{\sc Bolander, T.}
\newblock Self-reference.
\newblock In {\em The Stanford Encyclopedia of Philosophy}, E.~N. Zalta, Ed.,
  fall 2017~ed. Metaphysics Research Lab, Stanford University, 2017.
\newblock
  \url{https://plato.stanford.edu/archives/fall2017/entries/self-reference/}.

\bibitem{Bombieri}
{\sc Bombieri, E.}
\newblock Problems of the millennium: the riemann hypothesis.
\newblock
  \url{http://www.claymath.org/sites/default/files/official_problem_description.pdf}.
  [Online; accessed 12-March-2019].

\bibitem{Boole}
{\sc Boole, G.}
\newblock {\em An Investigation of the Laws of Thought}.
\newblock Macmillan and Co., 1854.
\newblock \url{https://books.google.com/books?id=DqwAAAAAcAAJ}.

\bibitem{Borovik}
{\sc Borovik, A.~V.}
\newblock {\em Mathematics under the Microscope}.
\newblock The University of Manchester, 2007.
\newblock \url{http://eprints.ma.man.ac.uk/844/1/covered/MIMS_ep2007_112.pdf}.

\bibitem{Borovoi}
{\sc Borovoi, M.~V.}
\newblock On the action of the galois group on rational cohomology classes of
  type (p,p) of abelian varieties.
\newblock {\em Matematiceskij Sbornik (Recueil Math\'ematique de la Soci\'et\'e
  Math\'ematique de Moscou) 94}, 4 (1974), 649–652,656.
\newblock (English translation in Math. USSR Sbornik 23 (1974), 613–616).

\bibitem{Borovoi2}
{\sc Borovoi, M.~V.}
\newblock The shimura-deligne schemes mc(g,h) and the rational cohomology
  classes of type (p,p) of abelian varieties.
\newblock {\em Problems in Group Theory and Homological Algebra\/} (1977),
  3–53.

\bibitem{Boston}
{\sc Boston, N.}
\newblock The proof of fermat's last theorem, Spring 2003.
\newblock
  \url{https://web.archive.org/web/20040325174327/https://www.math.wisc.edu/~boston/869.pdf}.

\bibitem{Bridges}
{\sc Bridges, D., and Palmgren, E.}
\newblock Constructive mathematics.
\newblock In {\em The Stanford Encyclopedia of Philosophy}, E.~N. Zalta, Ed.,
  summer 2018~ed. Metaphysics Research Lab, Stanford University, 2018.
\newblock
  \url{https://plato.stanford.edu/archives/sum2018/entries/mathematics-constructive/}.

\bibitem{Bromwich}
{\sc Bromwich, T. J.~I.}
\newblock The relation between the convergence of series and of integrals.
\newblock {\em Proc. Lond. Math. Soc. 6\/} (1908), 327--338.
\newblock \url{https://archive.org/stream/proceedingslond12socigoog}.

\bibitem{Bromwich2}
{\sc Bromwich, T. J.~I., and MacRobert, T.~M.}
\newblock {\em An Introduction to the Theory of Infinite Series}, 3rd~ed.
\newblock American Mathematical Society, 1991.

\bibitem{Brouwer3}
{\sc Brouwer, L.}
\newblock Intuitionism and formalism.
\newblock {\em Bull. Amer. Math. Soc. 20}, 2 (1913), 81--96.
\newblock \url{https://projecteuclid.org/euclid.bams/1183422499}.

\bibitem{Brouwer2}
{\sc Brouwer, L.}
\newblock Historical background, principles and methods of intuitionism.
\newblock {\em South African Journal of Science 49}, 3-4 (Oct 1952), 139.
\newblock \url{https://hdl.handle.net/10520/AJA00382353_3009}.

\bibitem{Brouwer}
{\sc Brouwer, L.}
\newblock On the foundations of mathematics.
\newblock In {\em Collected Works}, A.~Heyting, Ed., vol.~I. North Holland
  Publishing Co., 1975.

\bibitem{Brouwer4}
{\sc Brouwer, L. E.~J.}
\newblock Intuitionist set theory.
\newblock In {\em From Brouwer to Hilbert}, P.~Mancosu, Ed. Oxford University
  Press, 1998, ch.~1, pp.~23--27.

\bibitem{Brownstein}
{\sc Brownstein, D.}
\newblock Aspects of the problem of universals.
\newblock In {\em University of Kansas Humanistic Studies}, vol.~44. University
  of Kansas Publications, 1973, pp.~49--50.
\newblock
  \url{https://kuscholarworks.ku.edu/bitstream/handle/1808/6511/humseries.44.Aspects_of_Universals.pdf?sequence=1}.

\bibitem{Bruin}
{\sc Bruin, P.}
\newblock What is. . . an l-function?
\newblock
  \url{https://web.archive.org/web/20190518200001/https://www.math.leidenuniv.nl/~pbruin/L-functions.pdf}.
\newblock Zurich Graduate Colloquium, 30 October 2012.

\bibitem{Cahen}
{\sc Cahen, E.}
\newblock Sur la fonction $\zeta(s)$ de riemann et sur des fonctions analogues.
\newblock {\em Annales de l'\'ecole Normale 11}, 3 (1894), 75--164.
\newblock \url{https://archive.org/details/surlafonctionzet00caheuoft}.

\bibitem{Calinger}
{\sc Calinger, R.~S.}
\newblock {\em Leonhard Euler: Mathematical Genius in the Enlightenment}.
\newblock Princeton University Press, 2015.
\newblock \url{https://books.google.com/books?id=jSv3BgAAQBAJ&dq}.

\bibitem{Cameron}
{\sc Cameron, P.}
\newblock Peter cameron's blog: Quotes, 10 2011.
\newblock
  \url{https://web.archive.org/web/20111007124611/https://cameroncounts.wordpress.com/quotes/}.

\bibitem{Cardano}
{\sc Cardano, G.}
\newblock {\em The Book of My Life}.
\newblock E. P. Dutton \& Co., Inc., 1930.
\newblock
  \url{https://archive.org/details/GirolamoCardanoTheBookOfMyLifeDeVitaPropriaLiber}.

\bibitem{Morgan}
{\sc Carlson, J.}, Ed.
\newblock {\em 100 Years of Topology: Work Stimulated by Poincar\'e’s
  Approach to Classifying Manifolds\/} (June 2010), vol.~19 of {\em Clay
  Mathematics Proceedings}, Institut Henri Poincar\'e, Clay Mathematics
  Institute.
\newblock http://www.claymath.org/library/proceedings/cmip19.pdf.

\bibitem{Carnap}
{\sc Carnap, R.}
\newblock {\em Introduction to Symbolic Logic and its Applications}.
\newblock Dover Publications, 1958.

\bibitem{Chandrasekharan}
{\sc Chandrasekharan, K.}
\newblock The work of enrico bombieri.
\newblock In {\em Proceedings of the 1974 International Congress of
  Mathematicians}, R.~D. I.~James, Ed. Canadian Mathematical Congress, 1975,
  pp.~3--10.
\newblock
  \url{https://www.mathunion.org/fileadmin/ICM/Proceedings/ICM1974.1/ICM1974.1.ocr.pdf}.

\bibitem{Breuil}
{\sc Christophe~Breuil, Brian~Conrad, F.~D., and Taylor, R.}
\newblock On the modularity of elliptic curves over q: wild 3-adic exercises.
\newblock {\em J. Amer. Math. Soc. 14}, 4 (10 2001), 843–939.
\newblock \url{https://www.jstor.org/stable/827119}.

\bibitem{Church}
{\sc Church, A.}
\newblock 'logic', in the encyclopaedia britannica, xiv edition, chicago 1959.
\newblock {\em The Journal of Symbolic Logic 23\/} (1959), 22–29.

\bibitem{Ciucci}
{\sc Ciucci, D., and Dubois, D.}
\newblock A map of dependencies among three-valued logics.
\newblock {\em Information Sciences 250\/} (2013), 162--177.
\newblock \url{https://hal.archives-ouvertes.fr/hal-01122868}.

\bibitem{Clark}
{\sc Clark, P.~L.}
\newblock Dirichlet series, 2010.
\newblock
  \url{https://web.archive.org/web/20100611225105/http://math.uga.edu/~pete/4400dirichlet.pdf}.

\bibitem{Clay}
{\sc {Clay Mathematics Institute}}.
\newblock {The Birch and Swinnerton-Dyer Conjecture}.
\newblock In {\em Millenium Problems}. Clay Mathematics Institute, 2018.
\newblock
  \url{https://web.archive.org/web/20180111100858/http://www.claymath.org/millennium-problems/birch-and-swinnerton-dyer-conjecture}.
  [Online; accessed 15-January-2018].

\bibitem{Cohen2}
{\sc Cohen, M.}
\newblock {\em The Philosophy Bible: The Definitive Guide to the Last 3,000
  Years of Thought}, 1st~ed.
\newblock Firefly Books Ltd., 2016.

\bibitem{Cohen}
{\sc Cohen, S.~M.}
\newblock Aristotle's metaphysics.
\newblock In {\em The Stanford Encyclopedia of Philosophy}, E.~N. Zalta, Ed.,
  winter 2016~ed. Metaphysics Research Lab, Stanford University, 2016.
\newblock
  \url{https://plato.stanford.edu/archives/win2016/entries/aristotle-metaphysics/}.

\bibitem{Coleman}
{\sc Coleman, M.}
\newblock Background analytic continuation, 2017.
\newblock
  \url{https://web.archive.org/web/20171118051044/http://www.maths.manchester.ac.uk/~mdc/MATH41022/Background/Background%20Analytic%20Continuation.pdf}.

\bibitem{Conrad}
{\sc Conrad, K.}
\newblock Dirichlet series, 2019.
\newblock
  \url{https://web.archive.org/web/20190325042037/https://kconrad.math.uconn.edu/math5121s18/handouts/dirichletseries.pdf}.

\bibitem{Conrey}
{\sc Conrey, J.~B.}
\newblock The riemann hypothesis.
\newblock {\em Notices of the American Mathematical Society 50}, 3 (March
  2003), 341--353.
\newblock \url{https://www.ams.org/notices/200303/fea-conrey-web.pdf}.

\bibitem{Copi}
{\sc Copi, I.~M., and Cohen, C.}
\newblock {\em Introduction to Logic}, 12th~ed.
\newblock Pearson Education, 2005.

\bibitem{Courant}
{\sc Courant, R., and Robbins, H.}
\newblock {\em What is Mathematics?}
\newblock Oxford University Press, 1996.

\bibitem{Crossley}
{\sc Crossley, J., and Ash, C.}
\newblock {\em What is mathematical logic?}
\newblock Oxford University Press, 1972.

\bibitem{Curry}
{\sc Curry, H.~B.}
\newblock Review: Grundzüge der theoretischen logik (3rd edition).
\newblock {\em Bull. Amer. Math. Soc. 59}, 3 (1953), 263–267.
\newblock
  \url{http://www.ams.org/journals/bull/1953-59-03/S0002-9904-1953-09701-4/S0002-9904-1953-09701-4.pdf}.

\bibitem{daCosta2}
{\sc da~Costa, N., B\'eziau, J.-Y., and Bueno, O.}
\newblock Paraconsistent logic in a historical perspective.
\newblock {\em Logique et Analyse 38}, 150-152 (1995), 111--125.
\newblock Special Issue Dedicated to the Memory of Leo Apostel.
  \url{https://www.jstor.org/stable/44084538}.

\bibitem{daCosta}
{\sc da~Costa, N., and Krause, D.}
\newblock Schrödinger logics.
\newblock {\em Studia Logica 53}, 4 (December 1994), 533–550.
\newblock \url{https://link.springer.com/article/10.1007%2FBF01057649}.

\bibitem{Davenport}
{\sc Davenport, H.}
\newblock {\em Multiplicative Number Theory}, 3rd~ed., vol.~74 of {\em Graduate
  Texts in Mathematics}.
\newblock Springer-Verlag, 2000.
\newblock Revised and with a preface by Hugh L. Montgomery. MR 2001f:11001.

\bibitem{Davis}
{\sc Davis, M.}
\newblock {\em Engines of Logic}.
\newblock W. W. Norton \& Company, Inc., 2000.

\bibitem{Davis2}
{\sc Davis, W.~A., and the Georgetown Logic~Group}.
\newblock {\em An Introduction to Logic}.
\newblock Kunos Press, 2007.

\bibitem{Decker}
{\sc Decker, H.}
\newblock Historical and computational aspects of paraconsistency in view of
  the logic foundation of databases.
\newblock In {\em Semantics in Databases: Second International Workshop}, K.-D.
  S. B.~T. Leopoldo~Bertossi, Gyula O.H.~Katona, Ed., vol.~2582 of {\em Lecture
  Notes in Computer Science (LNCS)}. Springer-Verlag, January 2001, pp.~63--81.
\newblock \url{https://link.springer.com/chapter/10.1007/3-540-36596-6_4}.

\bibitem{Deligne}
{\sc Deligne, P., and Milne, J.}
\newblock Hodge cycles on abelian varieties.
\newblock In {\em Hodge cycles, motives, and Shimura varieties}, J.~O. A. S.
  K.-y. Deligne, P.;~Milne, Ed., vol.~900 of {\em Lecture Notes in Mathematics
  (LNM)}. Springer-Verlag, 1982, pp.~9--100.
\newblock \url{https://www.jmilne.org/math/Documents/Deligne82.pdf}.

\bibitem{p-series}
{\sc Department~of Mathematics, O. S.~U.}
\newblock p-series, 1996.
\newblock
  \url{https://web.archive.org/web/20050414083205/https://math.oregonstate.edu/home/programs/undergrad/CalculusQuestStudyGuides/SandS/SeriesTests/p-series.html}.
  [Online; accessed 9-July-2017].

\bibitem{Derbyshire}
{\sc Derbyshire, J.}
\newblock {\em Prime Obsession: Bernhard Riemann and The Greatest Unsolved
  Problem in Mathematics}.
\newblock Joseph Henry Press, 2003.

\bibitem{Dirac}
{\sc Dirac, P. A.~M.}
\newblock The evolution of the physicist’s picture of nature.
\newblock {\em Scientific American 208\/} (May 1963), 45--53.
\newblock
  \url{https://blogs.scientificamerican.com/guest-blog/the-evolution-of-the-physicists-picture-of-nature/}.

\bibitem{Dittrich}
{\sc Dittrich, W.}
\newblock On riemann’s paper, “on the number of primes less than a given
  magnitude”.
\newblock \url{https://arxiv.org/pdf/1609.02301.pdf}. Published 2 Aug 2017.

\bibitem{Dorbolo}
{\sc Dorbolo, J.}
\newblock {Aristotle: Laws of Thought}.
\newblock In {\em InterQuest: Introduction to Philosophy (PHL201)}. Oregon
  State University, 2002.
\newblock
  \url{http://oregonstate.edu/instruct/phl201/modules/Philosophers/Aristotle/aristotle_laws_of_thought.html}.

\bibitem{Dummett}
{\sc Dummett, M.}
\newblock {\em Elements of Intuitionism}, 2nd~ed., vol.~39 of {\em Oxford Logic
  Guides}.
\newblock Clarendon Press, 2000.

\bibitem{Dunn2}
{\sc Dunn, J.~M.}
\newblock Quantum mathematics.
\newblock In {\em PSA: Proceedings of the Biennial Meeting of the Philosophy of
  Science Association\/} (1980), P.~Asquith and R.~Giere, Eds., The Philosophy
  of Science Association, pp.~512--531.
\newblock
  \url{https://www.journals.uchicago.edu/doi/10.1086/psaprocbienmeetp.1980.2.192608}.

\bibitem{Dunn}
{\sc Dunn, J.~M.}
\newblock The impossibility of certain higher-order non-classical logics with
  extensionality.
\newblock In {\em Philosophical Analysis: A Defense by Example}, D.~F. Austin,
  Ed., vol.~39 of {\em Philosophical Studies Series}. Kluwer Academic
  Publishers, 1988, p.~261–279.
\newblock \url{https://rd.springer.com/chapter/10.1007/978-94-009-2909-8_16}.

\bibitem{Dutilh}
{\sc Dutilh-Novaes, C.}
\newblock Contradiction: the real challenge for paraconsistent logic.
\newblock In {\em Handbook of Paraconsistency}, J.~B\'eziau, W.~Carnielli, and
  D.~Gabbay, Eds. College Publications, 2008.

\bibitem{Duering}
{\sc Düring, I.}
\newblock Aristotle in the ancient biographical tradition.
\newblock {\em Acta Universitatis Gothoburgensis 63\/} (1957).
\newblock
  \url{https://books.google.com/books/about/Aristotle_in_the_ancient_biographical_tr.html?id=ThwXAQAAMAAJ}.

\bibitem{Brittanica}
{\sc {Editors of Encyclopedia Brittanica}}.
\newblock {Laws of Thought}.
\newblock In {\em Encyclopaedia Britannica}. Encyclopaedia Britannica, Inc.,
  2016.
\newblock
  \url{https://web.archive.org/web/20170320072527/https://www.britannica.com/topic/laws-of-thought}.

\bibitem{Edwards2}
{\sc Edwards, H.~M.}
\newblock {\em Advanced Calculus}.
\newblock Houghton, 1969.

\bibitem{Edwards}
{\sc Edwards, H.~M.}
\newblock {\em Riemann's Zeta Function}.
\newblock Dover Publications, 2001.

\bibitem{Elizalde}
{\sc Elizalde, E.}
\newblock Zeta-function method for regularization.
\newblock In {\em Encyclopedia of Mathematics}, {Last} modified on 29
  {December} 2015~ed. Springer-Verlag, 2015.
\newblock
  \url{https://www.encyclopediaofmath.org/index.php/Zeta-function_method_for_regularization}.

\bibitem{Euler3}
{\sc Euler, L.}
\newblock Methodus uniuersalis seriemm convergentium summas quam proxime
  inueniendi.
\newblock {\em Commentarii academiae scientarum Petropolitanae 8\/} (1736),
  8--9.
\newblock Comment No. 46 from Enestroemiani’s Index.
  \url{http://eulerarchive.maa.org/pages/E046.html}.

\bibitem{Euler4}
{\sc Euler, L.}
\newblock Methodus universalis series summandi ulterius promota.
\newblock {\em Commentarii academiae scientarum Petropolitanae 8\/} (1736),
  147--158.
\newblock Comment No. 55 from Enestroemiani’s Index.
  \url{http://eulerarchive.maa.org/pages/E055.html}.

\bibitem{Euler2}
{\sc Euler, L.}
\newblock Variae observationes circa series infinitas.
\newblock {\em Commentarii academiae scientarum Petropolitanae 9}, 1744 (1737),
  160–188.
\newblock Comment No. 72 from Enestroemiani’s Index.
  \url{http://eulerarchive.maa.org/pages/E072.html}.

\bibitem{Kleene3}
{\sc Fitting, M.}
\newblock Kleene’s three valued logics and their children.
\newblock {\em Fundamenta Informaticae 20}, 1,2,3 (1994), 113--131.
\newblock
  \url{https://pdfs.semanticscholar.org/c224/2a131ef0cb0b65f46638a1110a84165a3608.pdf}.

\bibitem{Fitting}
{\sc Fitting, M.}
\newblock Bilattices are nice things.
\newblock In {\em Self-reference}, V.~F.~H. Thomas~Bolander and S.~A. Pedersen,
  Eds. University of Chicago Press, 2006.
\newblock
  \url{https://pdfs.semanticscholar.org/705d/47eb02020ae3ef514bb1cdd4a9dcd6eac14a.pdf}.

\bibitem{Fontainelle}
{\sc Fontainelle, E.}
\newblock Logic.
\newblock In {\em Trivium: The Classical Liberal Arts of Grammar, Logic, and
  Rhetoric}, J.~Martinueau, Ed. Bloomsbury USA, 2016.

\bibitem{Forrest}
{\sc Forrest, P.}
\newblock The identity of indiscernibles.
\newblock In {\em The Stanford Encyclopedia of Philosophy}, E.~N. Zalta, Ed.,
  winter 2016~ed. Metaphysics Research Lab, Stanford University, 2016.
\newblock
  \url{https://plato.stanford.edu/archives/win2016/entries/identity-indiscernible/}.

\bibitem{Frege2}
{\sc Frege, G.}
\newblock On sense and nominatum (ueber sinn und bedeutung).
\newblock In {\em Readings in Philosophical Analysis}, H.~Feigl and W.~Sellars,
  Eds. Appleton-Century-Crofts, Inc., 1949.
\newblock \url{http://www.uh.edu/~garson/SenseandReference.pdf}.

\bibitem{Frege}
{\sc Frege, G.}
\newblock Logic.
\newblock In {\em Posthumous Writings}, H.~Hermes, F.~Kambartel, and
  F.~Kaulbach, Eds. The University of Chicago Press, 1979, pp.~126--151.

\bibitem{French}
{\sc French, S.}
\newblock Identity and individuality in quantum theory.
\newblock In {\em The Stanford Encyclopedia of Philosophy}, E.~N. Zalta, Ed.,
  fall 2015~ed. Metaphysics Research Lab, Stanford University, 2015.
\newblock \url{https://plato.stanford.edu/archives/fall2015/entries/qt-idind/}.

\bibitem{Freund}
{\sc Freund, P. G.~O., and Witten, E.}
\newblock Adelic string amplitudes.
\newblock {\em Phys. Lett. B.}, 199 (1987), 191.
\newblock \url{https://doi.org/10.1016/0370-2693(87)91357-8}.

\bibitem{Frey}
{\sc Frey, G.}
\newblock The way to the proof of fermat's last theorem.
\newblock {\em Annales de la Facult\'e des Sciences de Toulouse XVIII\/}
  (2009), 5--23.
\newblock \url{http://afst.cedram.org/item?id=AFST_2009_6_18_S2_5_0}.

\bibitem{Fronhoefer}
{\sc Fronhöfer, B.}
\newblock {Introduction to Many-Valued Logics}, 2011.
\newblock
  \url{https://web.archive.org/web/20131225052706/http://www.wv.inf.tu-dresden.de/Teaching/SS-2011/mvl/mval.HANDOUT2.pdf}.
  Archived on Dec. 25, 2013.

\bibitem{Gabbay}
{\sc Gabbay, D.~M.}
\newblock Classical vs non-classical logic.
\newblock In {\em Handbook of Logic in Artificial Intelligence and Logic
  Programming}, D.~Gabbay, C.~Hogger, and J.~Robinson, Eds., vol.~2. Oxford
  University Press, 1994, ch.~2.6.

\bibitem{Gardner2}
{\sc Gardner, M.}
\newblock {\em Martin Gardner's Sixth Book of Mathematical Games from
  Scientific American}.
\newblock Univ of Chicago Press, 1971.

\bibitem{Gardner}
{\sc Gardner, M.}
\newblock {\em Aha! Gotcha}.
\newblock W. H. Freeman \& Co., 1982.

\bibitem{Gelbart}
{\sc Gelbart, S.~S., and Miller, S.~D.}
\newblock Riemann's zeta function and beyond.
\newblock {\em Bull. Amer. Math. Soc. 41}, 1 (2004), 59--112.
\newblock
  \url{http://www.ams.org/journals/bull/2004-41-01/S0273-0979-03-00995-9/home.html}.

\bibitem{Gordon}
{\sc Gordon, B.~B.}
\newblock A survey of the hodge conjecture for abelian varieties.
\newblock \url{ https://arxiv.org/abs/alg-geom/9709030}. Published 26 Sep 1997.
\newblock Appendix B in "A Survey of the Hodge Conjecture" by James D. Lewis,
  2nd ed. Print ISBN: 978-1-4704-2852-5.

\bibitem{Gottlieb}
{\sc Gottlieb, P.}
\newblock Aristotle on {Non-contradiction}.
\newblock In {\em The Stanford Encyclopedia of Philosophy}, E.~N. Zalta, Ed.,
  summer 2015~ed. Metaphysics Research Lab, Stanford University, 2015.
\newblock
  \url{https://plato.stanford.edu/archives/sum2015/entries/aristotle-noncontradiction/}.

\bibitem{Gourdon}
{\sc Gourdon, X., and Sebah, P.}
\newblock Numerical evaluation of the riemann zeta-function, 2003.
\newblock
  \url{http://numbers.computation.free.fr/Constants/Miscellaneous/zetaevaluations.pdf}.

\bibitem{Grattan-Guinness2}
{\sc Grattan-Guinness, I.}
\newblock Structural similarity or structuralism? comments on priest's analysis
  of the paradoxes of self-reference.
\newblock {\em Mind 107}, 428 (Oct 1998), 823--834.
\newblock \url{https://www.jstor.org/stable/2659786}.

\bibitem{Grattan-Guinness}
{\sc Grattan-Guinness, I.}
\newblock {\em The Search for Mathematical Roots, 1870–1940: Logics, Set
  Theories, and the Foundations of Mathematics from Cantor through Russell to
  Gödel}.
\newblock Princeton University Press, 2000.

\bibitem{Griffiths}
{\sc Griffiths, R.~B.}
\newblock The consistent histories approach to quantum mechanics.
\newblock In {\em The Stanford Encyclopedia of Philosophy}, E.~N. Zalta, Ed.,
  spring 2017~ed. Metaphysics Research Lab, Stanford University, 2017.
\newblock
  \url{https://plato.stanford.edu/archives/spr2017/entries/qm-consistent-histories/}.

\bibitem{Grishin}
{\sc Grishin, V.}
\newblock Contradiction, law of.
\newblock In {\em Encyclopedia of Mathematics}, {Last} modified on 17 {March}
  2014~ed. Springer-Verlag, 2014.
\newblock
  \url{https://www.encyclopediaofmath.org/index.php/Contradiction,_law_of}.

\bibitem{Guichard1}
{\sc Guichard, D., and Koblitz, N.}
\newblock The integral test.
\newblock In {\em Calculus: Late transcendentals}. Dept. of Mathematics,
  Whitman College, 2017.
\newblock
  \url{https://www.whitman.edu/mathematics/calculus_late_online/section13.03.html}.
  [Online; accessed 15-August-2017].

\bibitem{Gomez-Torrente}
{\sc Gómez-Torrente, M.}
\newblock Logical truth.
\newblock In {\em The Stanford Encyclopedia of Philosophy}, E.~N. Zalta, Ed.,
  spring 2019~ed. Metaphysics Research Lab, Stanford University, 2019.
\newblock
  \url{https://plato.stanford.edu/archives/spr2019/entries/logical-truth/}.

\bibitem{Godel}
{\sc Gödel, K.}
\newblock Über formal unentscheidbare sätze der principia mathematica und
  verwandter systeme i.
\newblock {\em Monatshefte für Mathematik und Physik 38}, 1 (1931), 173--198.
\newblock \url{http://www.w-k-essler.de/pdfs/goedel.pdf}.

\bibitem{Godel2}
{\sc Gödel, K.}
\newblock {\em On Formally Undecidable Propositions of Principia Mathematica
  and Related Systems}.
\newblock Dover Publications, Inc., 1992.

\bibitem{Haack}
{\sc Haack, S.}
\newblock {\em Philosophy of Logics}.
\newblock Cambridge University Press, 1978.

\bibitem{Haack3}
{\sc Haack, S.}
\newblock {\em Deviant Logic, Fuzzy Logic}.
\newblock University of Chicago Press, 1996.

\bibitem{Haack2}
{\sc Haack, S.}
\newblock Deviant logics.
\newblock In {\em Encyclopedia of Language and Linguistics}, R.~E. Asher, Ed.,
  revised 1997~ed., vol.~2. Pergamon Press, 1997, pp.~256--263.
\newblock \url{https://www.academia.edu/24533445/Deviant_Logic_1994_1997_}.

\bibitem{Hairer}
{\sc Hairer, E., and Wanner, G.}
\newblock {\em Analysis by Its History}.
\newblock Springer-Verlag, 1996.
\newblock \url{https://www.springer.com/us/book/9780387945514}.

\bibitem{Hankel}
{\sc Hankel, H.}
\newblock Die euler'schen integrale bei unbeschränkter variabilität des
  argumentes.
\newblock {\em Zeitschrift für Math, und Phys. 9\/} (1864), 1--21.
\newblock
  \url{https://books.google.com/books?id=M3daAAAAcAAJ&pg=PP1#v=onepage&q&f=false}.

\bibitem{Hardy3}
{\sc Hardy, G.}
\newblock {\em A Mathematician's Apology}.
\newblock Cambridge University Press. Republished by University of Alberta
  Mathematical Sciences Society., 1940.
\newblock \url{https://archive.org/details/AMathematiciansApology-G.h.Hardy}.

\bibitem{Hardy}
{\sc Hardy, G.}
\newblock {\em Divergent Series}.
\newblock Oxford University Press, 1949.
\newblock \url{https://archive.org/details/DivergentSeries}.

\bibitem{Hardy2}
{\sc Hardy, G., and Riesz, M.}
\newblock {\em The General Theory of Dirichlet's Series}.
\newblock Cambridge University Press. Reprinted by Cornell University Library
  Digital Collections., 1915.
\newblock \url{https://archive.org/details/cu31924060184441/page/n5}.

\bibitem{Havil}
{\sc Havil, J.}
\newblock {\em Gamma: Exploring Euler's Constant}.
\newblock Princeton University Press, 2003.

\bibitem{Hawking}
{\sc Hawking, S.~W.}
\newblock Zeta function regularization of path integrals in curved spacetime.
\newblock {\em Comm. Math. Phys. 55}, 2 (1977), 133--148.
\newblock \url{https://projecteuclid.org/euclid.cmp/1103900982}.

\bibitem{Hazama}
{\sc Hazama, F.}
\newblock Algebraic cycles on nonsimple abelian varieties.
\newblock {\em Duke Math. J. 58\/} (1980), 31--37.
\newblock \url{https://projecteuclid.org/euclid.dmj/1077307371}.

\bibitem{Hazen3}
{\sc Hazen, A., and Pelletier, F.}
\newblock Pecularities of some three- and four-valued second order logics.
\newblock {\em Logica Universalis 12}, 3-4 (10 2018).
\newblock \url{https://doi.org/10.1007/s11787-018-0214-7}.

\bibitem{Hazen2}
{\sc Hazen, A.~P., and Pelletier, F.~J.}
\newblock K3, Ł3, lp, rm3, a3, fde, m: How to make many-valued logics work for
  you.
\newblock \url{https://arxiv.org/abs/1711.05816}, 11 2017.
\newblock \url{https://www.uni-log.org/prize/Canada.pdf}.

\bibitem{Hazen}
{\sc Hazen, A.~P., and Pelletier, F.~J.}
\newblock Second-order logic of paradox.
\newblock {\em Notre Dame Journal of Formal Logic 59}, 4 (2018), 547--558.
\newblock \url{https://projecteuclid.org/euclid.ndjfl/1536653099}.

\bibitem{He}
{\sc He, Y.-H., Jejjala, V., and Minic, D.}
\newblock From veneziano to riemann: A string theory statement of the riemann
  hypothesis, 2015.
\newblock \url{https://arxiv.org/abs/1501.01975}.

\bibitem{Heidelberger}
{\sc Heidelberger, H.}
\newblock The indispensability of truth.
\newblock {\em American Philosophical Quarterly 5\/} (1968), 212--217.
\newblock \url{https://www.jstor.org/stable/20009275}.

\bibitem{Heis}
{\sc Heis, J.}
\newblock Book review - {The Cambridge Companion to Frege}, {Michael Potter and
  Tom Ricketts} (eds.),.
\newblock {\em Notre Dame Philosophical Reviews}, 2011.11.34 (2011).
\newblock \url{https://ndpr.nd.edu/news/the-cambridge-companion-to-frege/}.

\bibitem{Henderson}
{\sc Henderson, L.}
\newblock The problem of induction.
\newblock In {\em The Stanford Encyclopedia of Philosophy}, E.~N. Zalta, Ed.,
  spring 2019~ed. Metaphysics Research Lab, Stanford University, 2019.

\bibitem{Heyting2}
{\sc Heyting, A.}
\newblock Die intuitionistische grundlegung der mathematik.
\newblock {\em Erkenntnis 2}, 1 (1931), 106–115.
\newblock \url{https://www.jstor.org/stable/20011630}.

\bibitem{Heyting}
{\sc Heyting, A.}
\newblock {\em Intuitionism: An Introduction}, 2nd revised~ed.
\newblock North Holland Publishing Co., 1966.

\bibitem{Heyting3}
{\sc Heyting, A.}
\newblock The formal rules of intuitionistic logic.
\newblock In {\em From Brouwer to Hilbert}, P.~Mancosu, Ed. Oxford University
  Press, 1998, ch.~24, pp.~311--327.

\bibitem{Hilbert1901}
{\sc Hilbert, D.}
\newblock Mathematische probleme.
\newblock {\em Archiv der Mathematik und Physik\/} (1901), 44--63, 213–237.
\newblock \url{https://babel.hathitrust.org/cgi/pt?id=njp.32101075453553}.

\bibitem{Hilbert1917}
{\sc Hilbert, D.}
\newblock Axiomatisches denken.
\newblock {\em Mathematische Annalen 78\/} (1917), 405--415.
\newblock \url{http://eudml.org/doc/158776}.

\bibitem{Hildebrand}
{\sc Hildebrand, A.}
\newblock Introduction to analytic number theory, 2013.
\newblock
  \url{https://web.archive.org/web/20190326050948/https://faculty.math.illinois.edu/~hildebr/ant/main.pdf}.

\bibitem{Horn}
{\sc Horn, L.~R.}
\newblock Contradiction.
\newblock In {\em The Stanford Encyclopedia of Philosophy}, E.~N. Zalta, Ed.,
  spring 2014~ed. Metaphysics Research Lab, Stanford University, 2014.
\newblock
  \url{https://plato.stanford.edu/archives/spr2014/entries/contradiction/}.

\bibitem{Hurley}
{\sc Hurley, P.}
\newblock {\em A Concise Introduction to Logic}, 4th~ed.
\newblock Wadsworth Publishing, 1991.

\bibitem{Iemhoff}
{\sc Iemhoff, R.}
\newblock Intuitionism in the {Philosophy} of {Mathematics}.
\newblock In {\em The Stanford Encyclopedia of Philosophy}, E.~N. Zalta, Ed.,
  winter 2016~ed. Metaphysics Research Lab, Stanford University, 2016.
\newblock
  \url{https://plato.stanford.edu/archives/win2016/entries/intuitionism/}.

\bibitem{Ivic}
{\sc Ivić, A.}
\newblock {\em The Riemann Zeta-Function: The Theory of the Riemann
  Zeta-Function with Applications}.
\newblock John Wiley \& Sons, Inc., 1985.
\newblock \url{https://www.springer.com/us/book/9780387945514}.

\bibitem{Iwaniec}
{\sc Iwaniec, H., and Sarnak, P.}
\newblock Perspectives on the analytic theory of {L}-functions.
\newblock {\em GAFA, Geometric and Functional Analysis\/} (2000).
\newblock
  \url{http://web.math.princeton.edu/sarnak/Perspectives%20on%20the%20Analytic%20Theory%20of%20L-functions.pdf}.

\bibitem{Jacoby}
{\sc Jacoby, F.}
\newblock {\em Die Fragmente der griechischen Historiker}, vol.~244.
\newblock Weidmann, 1926-1957.
\newblock \url{https://referenceworks.brillonline.com/cluster/Jacoby%20Online
  }.

\bibitem{Jacquette}
{\sc Jacquette, D.}
\newblock Introduction: Logic, philosophy, and philosophical logic.
\newblock In {\em A Companion to Philosophical Logic}, D.~Jacquette, Ed. John
  Wiley \& Sons, 2006, pp.~5--6.
\newblock
  \url{https://books.google.com/books?id=pzf7_sT58PUC&printsec=frontcover#v=onepage&q&f=false}.

\bibitem{Jannsen}
{\sc Jannsen, U.}
\newblock Deligne’s proof of the weil-conjecture, 2015/16.
\newblock
  \url{http://www.mathematik.uni-regensburg.de/Jannsen/home/Weil-gesamt-eng.pdf},.

\bibitem{Jaskowski}
{\sc Jaśkowski, S.}
\newblock A propositional calculus for inconsistent deductive systems.
\newblock {\em Logic and Logical Philosophy 7\/} (1999), 35--56.
\newblock
  \url{http://http://apcz.umk.pl/czasopisma/index.php/LLP/article/viewFile/2616/2595}.

\bibitem{Jensen}
{\sc Jensen, J. L. W.~V.}
\newblock Om rÆkkers konvergens.
\newblock {\em Tidsskrift for mathematik 2}, 5 (1884), 63--72.
\newblock \url{https://www.jstor.org/stable/24540057}.

\bibitem{Johansson}
{\sc Johansson, I.}
\newblock Der minimalkalkül, ein reduzierter intuitionistischer formalismus.
\newblock {\em Compositio Mathematica 4}, 1 (1936), 119–136.
\newblock \url{http://www.numdam.org/numdam-bin/item?id=CM_1937__4__119_0}.

\bibitem{Kant}
{\sc Kant, I.}
\newblock Attempt to introduce the concept of negative magnitudes into
  philosophy (1763).
\newblock In {\em Theoretical Philosophy, 1755-1770}, D.~Walford and
  R.~Meerbote, Eds., 1st paperback~ed. Cambridge University Press, 2003.
\newblock \url{https://books.google.com/books?id=oooN5uQ9O6cC}.

\bibitem{Katz}
{\sc Katz, N.~M., and Sarnak, P.}
\newblock Zeros of zeta functions and symmetry.
\newblock {\em Bulletin of the American Mathematical Society 36}, 1 (February
  1999), 1--26.
\newblock
  \url{http://www.ams.org/journals/bull/1999-36-01/S0273-0979-99-00766-1/home.html}.

\bibitem{Klagge}
{\sc Klagge, J.~C., and Nordmann, A.}
\newblock {\em Ludwig Wittgenstein: Public and Private Occasions}.
\newblock Rowman \& Littlefield, 2003.

\bibitem{Kleene}
{\sc Kleene, S.~C.}
\newblock {\em Mathematical Logic}.
\newblock Dover Publications, 2002, 1967.

\bibitem{Kleene2}
{\sc Kleene, S.~C.}
\newblock {\em Introduction to Metamathematics}.
\newblock Ishi Press International, 2009, 1952.

\bibitem{Kline}
{\sc Kline, M.}
\newblock {\em Mathematics for the Nonmathematician}, 1985~ed.
\newblock Dover Publications, Inc., 1967.

\bibitem{Kolmogorov2}
{\sc Kolmogorov, A.~N.}
\newblock O principe tertium non datur.
\newblock {\em Matematiceskij Sbornik (Recueil Math\'ematique de la Soci\'et\'e
  Math\'ematique de Moscou) 32}, 4 (1925), 646--667.
\newblock (English translation in van Heijenoort (1967)).

\bibitem{Kolmogorov}
{\sc Kolmogorov, A.~N.}
\newblock On the principle of the excluded middle.
\newblock In {\em From Frege to Gödel: A Source Book in Mathematical Logic,
  1897-1931}, J.~van Heijenoort, Ed. Harvard University Press, 1967,
  pp.~414--437.

\bibitem{Kolyvagin}
{\sc Kolyvagin, V.}
\newblock Finiteness of $e(\mathbb{Q})$ and $sha(e,\mathbb{Q})$ for a class of
  weil curves.
\newblock {\em Math. USSR Izv. 32\/} (1989), 523–541.

\bibitem{Kremer}
{\sc Kremer, P.}
\newblock The revision theory of truth.
\newblock In {\em The Stanford Encyclopedia of Philosophy}, E.~N. Zalta, Ed.,
  winter 2016~ed. Metaphysics Research Lab, Stanford University, 2016.
\newblock
  \url{https://plato.stanford.edu/archives/win2016/entries/truth-revision/}.

\bibitem{Kripke}
{\sc Kripke, S.~A.}
\newblock Outline of a theory of truth.
\newblock {\em Journal of Philosophy 72\/} (1975), 690--716.
\newblock \url{https://www.jstor.org/stable/2024634}.

\bibitem{Kubota}
{\sc Kubota, T.}
\newblock On the field extension by complex multiplication.
\newblock {\em Trans. Amer. Math. Soc. 118\/} (1965), 113--122.
\newblock
  \url{https://www.ams.org/journals/tran/1965-118-00/S0002-9947-1965-0190144-8/S0002-9947-1965-0190144-8.pdf}.

\bibitem{Kuznetsov}
{\sc Kuznetsov, A.~V.}
\newblock On superintuitionistic logics.
\newblock In {\em Proceedings of the 1974 International Congress of
  Mathematicians}, R.~D. I.~James, Ed. Canadian Mathematical Congress, 1975,
  pp.~243--250.
\newblock
  \url{https://www.mathunion.org/fileadmin/ICM/Proceedings/ICM1974.1/ICM1974.1.ocr.pdf}.

\bibitem{Lakatos}
{\sc Lakatos, I.}
\newblock {\em Proofs and Refutations: The Logic of Mathematical Discovery},
  reissue~ed.
\newblock Cambridge University Press, 2015.

\bibitem{Landau}
{\sc Landau, E.}
\newblock Zur theorie der riemannschen zetafunktion.
\newblock {\em Vierteljahrsschr. Naturf. Ges. Zurich 56\/} (1911), 125--148.
\newblock \url{http://www.ngzh.ch/archiv/1911_56/56_1-2/56_8.pdf}.

\bibitem{Lande}
{\sc Lande, N.~P.}
\newblock {\em Classical Logic and Its Rabbit-Holes: A First Course}.
\newblock Hackett Publishing Co., Inc., 2013.

\bibitem{Langer}
{\sc Langer, S.~K.}
\newblock {\em An Introduction to Symbolic Logic}, 3rd revised~ed.
\newblock Dover Publications, 1967.

\bibitem{Lee}
{\sc Lee, S.-F.}
\newblock {\em Logic: A Complete Introduction}.
\newblock Hodder \& Stoughton, 2017.

\bibitem{Lemmon}
{\sc Lemmon, E.}
\newblock {\em Beginning Logic}.
\newblock Hackett Publishing Co., Inc., 1978.

\bibitem{Lewis}
{\sc Lewis, D.}
\newblock Letters to beall and priest.
\newblock In {\em The Law of Non-Contradiction. New Philosophical Essays},
  G.~Priest, J.~Beall, and B.~Armour-Garb, Eds. Oxford University Press, 2004,
  p.~176–177.

\bibitem{Look}
{\sc Look, B.~C.}
\newblock Symbolic logic ii, lecture 7: Beyond propositional logic.
\newblock
  \url{https://web.archive.org/web/20111109200756/http://www.uky.edu/~look/Phi520-Lecture7.pdf}.

\bibitem{Lotter}
{\sc Lotter, D.}
\newblock Gottlob frege: Language.
\newblock In {\em Internet Encyclopedia of Philosophy}, J.~Fieser and
  B.~Dowden, Eds. Internet Encyclopedia of Philosophy, 2018.
\newblock
  \url{https://web.archive.org/web/20180714153349/https://www.iep.utm.edu/freg-lan/}.

\bibitem{MacFarlane}
{\sc MacFarlane, J.}
\newblock Relevance logic, April 2016.
\newblock
  \url{https://web.archive.org/web/20171209105735/http://johnmacfarlane.net/142/relevance.pdf}.

\bibitem{Maclaurin}
{\sc Maclaurin, C.}
\newblock {\em A Treatise of Fluxions}.
\newblock William Baynes and William Davis, 1742.

\bibitem{Malinowski}
{\sc Malinowski, G.}
\newblock {Many-valued Logic} and its {Philosophy}.
\newblock In {\em Handbook of the History of Logic: The Many Valued and
  Nonmonotonic Turn in Logic}, D.~Gabbay and J.~Woods, Eds., vol.~8. Elsevier,
  2007, pp.~13--94.
\newblock
  \url{https://books.google.com/books?id=3TNj1ZkP3qEC&source=gbs_navlinks_s}.

\bibitem{Manin}
{\sc Manin, Y.~I., and Grell, B.}
\newblock Reviewed work: Provable and unprovable.
\newblock {\em Studia Logica 41}, 1 (1982), 84--86.
\newblock \url{https://www.jstor.org/stable/20015041}.

\bibitem{Mannoury}
{\sc Mannoury, G.}
\newblock Methodologisches und philosophisches zur elementarmathematik.
\newblock {\em Vierteljahrschrift für wissenschaftliche Philosophie und
  Soziologie 35\/} (1911), 263--265.
\newblock
  \url{https://www.jstor.org/stable/20114258?seq=1#page_scan_tab_contents}.

\bibitem{Mares}
{\sc Mares, E.}
\newblock Relevance logic.
\newblock In {\em The Stanford Encyclopedia of Philosophy}, E.~N. Zalta, Ed.,
  spring 2014~ed. Metaphysics Research Lab, Stanford University, 2014.
\newblock
  \url{https://plato.stanford.edu/archives/spr2014/entries/logic-relevance/}.

\bibitem{Marques}
{\sc Marques, T.}
\newblock {\em Bivalence and the Challenge of Truth-Value Gaps}.
\newblock PhD thesis, Univ. of Stirling, November 2003.
\newblock
  \url{https://www.researchgate.net/profile/Teresa_Marques5/publication/34243597_Bivalence_and_the_challenge_of_truth-value_gaps/links/54e790480cf27a6de10a8185/Bivalence-and-the-challenge-of-truth-value-gaps.pdf}.

\bibitem{Martin}
{\sc Martin, R.~L., and Woodruff, P.~W.}
\newblock On representing ‘true-in-l’ in l.
\newblock {\em Philosophia 5\/} (1975), 217–221.
\newblock \url{https://link.springer.com/article/10.1007%2FBF02379018 }.

\bibitem{Mattuck}
{\sc Mattuck, A.}
\newblock Cycles on abelian varieties.
\newblock {\em Proc. Amer. Math. Soc. 9\/} (1958), 88--98.
\newblock
  \url{https://www.ams.org/journals/proc/1958-009-01/S0002-9939-1958-0098752-1/S0002-9939-1958-0098752-1.pdf}.

\bibitem{McGee}
{\sc McGee, V.}
\newblock Semantic paradoxes and theories of truth.
\newblock In {\em Routledge Encyclopedia of Philosophy}. Taylor and Francis,
  1998.
\newblock
  \url{https://www.rep.routledge.com/articles/thematic/semantic-paradoxes-and-theories-of-truth/v-1}.

\bibitem{McKinsey}
{\sc McKinsey, J. C.~C.}
\newblock Proof of the independence of the primitive symbols of heyting’s
  calculus of propositions.
\newblock {\em The Journal of Symbolic Logic 6}, 4 (1939), 155–158.
\newblock \url{https://www.jstor.org/stable/2268715}.

\bibitem{McKubre-Jordens}
{\sc McKubre-Jordens, M.}
\newblock This is not a carrot: Paraconsistent mathematics.
\newblock {\em +Plus Magazine\/} (August 24, 2011).
\newblock \url{https://plus.maths.org/content/not-carrot}.

\bibitem{Mendell}
{\sc Mendell, H.}
\newblock Aristotle and mathematics.
\newblock In {\em The Stanford Encyclopedia of Philosophy}, E.~N. Zalta, Ed.,
  spring 2017~ed. Metaphysics Research Lab, Stanford University, 2017.
\newblock
  \url{https://plato.stanford.edu/archives/spr2017/entries/aristotle-mathematics/}.

\bibitem{Mermin}
{\sc Mermin, N.~D.}
\newblock Is the moon there when nobody looks? reality and the quantum theory.
\newblock {\em Physics Today 38\/} (April 1985), 38--47.
\newblock \url{https://physicstoday.scitation.org/doi/10.1063/1.880968}.

\bibitem{Meyer}
{\sc Meyer, R.~K.}
\newblock Entailment.
\newblock {\em Journal of Philosophy 68\/} (1971), 808–818.
\newblock
  \url{https://www.pdcnet.org/jphil/content/jphil_1971_0068_0021_0808_0818}.

\bibitem{Milne}
{\sc Milne, J.~S.}
\newblock {\em Arithmetic Duality Theorems}.
\newblock Academic Press, 1986.

\bibitem{Milne4}
{\sc Milne, J.~S.}
\newblock The tate conjecture over finite fields (aim talk), 2007.
\newblock \url{http://www.jmilne.org/math/articles/2007e.pdf}.

\bibitem{Milne3}
{\sc Milne, J.~S.}
\newblock The riemann hypothesis over finite fields: From weil to the present
  day.
\newblock In {\em The Legacy of Bernhard Riemann after One Hundred and Fifty
  Years}, S.-T.~Y. Lizhen~Ji, Frans~Oort, Ed. International Press, 2015,
  pp.~487--565.
\newblock \url{https://www.jmilne.org/math/xnotes/pRH.html}.

\bibitem{Milne2}
{\sc Milne, P.}
\newblock Frege's folly: bearerless names and basic law v.
\newblock In {\em The Cambridge Companion to Frege}, M.~Potter and T.~Ricketts,
  Eds. Cambridge University Press, 2010, pp.~473--475.
\newblock
  \url{https://books.google.com/books?id=LasbqbKdrjQC&printsec=frontcover#v=onepage&q&f=false},
  \url{http://hdl.handle.net/1893/6820}.

\bibitem{Minto}
{\sc Minto, W.}
\newblock {\em Logic, Inductive and Deductive}.
\newblock Charles Scribner's Sons, 1915.
\newblock \url{http://www.gutenberg.org/files/31796/31796-h/31796-h.htm},
  \url{https://books.google.com/books?id=vQ9VAAAAMAAJ}.

\bibitem{Mints}
{\sc Mints, G.}
\newblock Notes on constructive negation.
\newblock {\em Synthese 148}, 3 (2006), 701–717.
\newblock \url{https://link.springer.com/article/10.1007/s11229-004-6294-3}.

\bibitem{Mortansen}
{\sc Mortensen, C.}
\newblock Inconsistent mathematics.
\newblock In {\em The Stanford Encyclopedia of Philosophy}, E.~N. Zalta, Ed.,
  winter 2017~ed. Metaphysics Research Lab, Stanford University, 2017.
\newblock
  \url{https://plato.stanford.edu/archives/fall2017/entries/mathematics-inconsistent/}.

\bibitem{Moschovakis}
{\sc Moschovakis, J.}
\newblock Intuitionistic logic.
\newblock In {\em The Stanford Encyclopedia of Philosophy}, E.~N. Zalta, Ed.,
  spring 2015~ed. Metaphysics Research Lab, Stanford University, 2015.
\newblock
  \url{https://plato.stanford.edu/archives/spr2015/entries/logic-intuitionistic/}.

\bibitem{Mumford}
{\sc Mumford, D.}
\newblock A note of shimura's paper "discontinuous groups and abelian
  varieties".
\newblock {\em Mathematische Annalen 181\/} (1969), 345--351.
\newblock \url{https://link.springer.com/article/10.1007/BF01350672}.

\bibitem{Nedo}
{\sc Nedo, M., and Ranchetti, M.}
\newblock {\em Ludwig Wittgenstein: sein Leben in Bildern und Texten}.
\newblock Suhrkamp, 1983.

\bibitem{Norton}
{\sc Norton, J.~D.}
\newblock Special theory of relativity: Relativity of simultaneity.
\newblock
  \url{https://web.archive.org/web/20190331152943/https://www.pitt.edu/~jdnorton/teaching/HPS_0410/chapters/Special_relativity_rel_sim/index.html}.
\newblock Lecture notes for HPS 0410: Einstein for Everyone. [Online; Last
  Updated 14-January-2015.].

\bibitem{Nunez}
{\sc Nunez, C.}
\newblock Introduction to bosonic string theory.
\newblock In {\em INIS IAEA Jorge Andre Swieca Summer School on Particles and
  Fields}. International Atomic Energy Agency (IAEA), Jan 2009.
\newblock \url{https://inis.iaea.org/search/search.aspx?orig_q=RN:40084717}.

\bibitem{Odintsov}
{\sc Odintsov, S.~P., and Wansing, H.}
\newblock The logic of generalized truth values and the logic of bilattices.
\newblock {\em Studia Logica 103}, 1 (02 2015), 91–112.
\newblock \url{https://doi.org/10.1007/s11225-014-9546-3}.

\bibitem{Odlyzko}
{\sc Odlyzko, A.~M., and Schönhage, A.}
\newblock Fast algorithms for multiple evaluations of the {Riemann} zeta
  function.
\newblock {\em Trans. Amer. Math. Soc. 309\/} (1988), 797--809.
\newblock
  \url{https://www.ams.org/journals/tran/1988-309-02/S0002-9947-1988-0961614-2/S0002-9947-1988-0961614-2.pdf}.

\bibitem{Oort}
{\sc Oort, F., and Schappacher, N.}
\newblock Early history of the riemann hypothesis in positive characteristic.
\newblock In {\em The Legacy of Bernhard Riemann after One Hundred and Fifty
  Years}, S.-T.~Y. Lizhen~Ji, Frans~Oort, Ed. International Press, 2015,
  p.~595–631.
\newblock
  \url{http://www.staff.science.uu.nl/~oort0109/EigArt-HistpRH-2016.pdf}.

\bibitem{Overholt}
{\sc Overholt, M.}
\newblock {\em A Course in Analytic Number Theory}.
\newblock American Mathematical Society, 2014.
\newblock
  \url{https://books.google.com/books?id=kBsHBgAAQBAJ&source=gbs_navlinks_s}.

\bibitem{Pais}
{\sc Pais, A.}
\newblock Einstein and the quantum theory.
\newblock {\em Reviews of Modern Physics 51\/} (October 1979), 863--914.
\newblock \url{https://link.aps.org/doi/10.1103/RevModPhys.51.863}.

\bibitem{Panti1998}
{\sc Panti, G.}
\newblock Multi-valued logics.
\newblock In {\em Quantified Representation of Uncertainty and Imprecision},
  P.~Smets, Ed., vol.~1. Springer Netherlands, 1998, pp.~25--74.
\newblock \url{https://doi.org/10.1007/978-94-017-1735-9_2}.

\bibitem{Parsons}
{\sc Parsons, T.}
\newblock The traditional square of opposition.
\newblock In {\em The Stanford Encyclopedia of Philosophy}, E.~N. Zalta, Ed.,
  summer 2017~ed. Metaphysics Research Lab, Stanford University, 2017.
\newblock \url{https://plato.stanford.edu/archives/sum2017/entries/square/}.

\bibitem{Pelletier}
{\sc Pelletier, F., and Linsky, B.}
\newblock Russell vs. frege on definite descriptions as singular terms.
\newblock In {\em Russell Vs. Meinong: The Legacy of "On Denoting"}, N.~Griffin
  and D.~Jacquette, Eds. Routledge, 11 2008, ch.~3, pp.~40--64.
\newblock
  \url{http://citeseerx.ist.psu.edu/viewdoc/summary?doi=10.1.1.158.701}.

\bibitem{Pelletier2}
{\sc Pelletier, F.~J.}
\newblock The law of non-contradiction: New philosophical essays, edited by
  graham priest, j.c. beall, and bradley armour-garb.
\newblock {\em The Bulletin of Symbolic Logic 12}, 1 (03 2006), 131--135.
\newblock \url{https://www.jstor.org/stable/3515889}.

\bibitem{Pelletier3}
{\sc Pelletier, F.~J., and Stainton, R.~J.}
\newblock On "the denial of bivalence is absurd".
\newblock {\em Australasian Journal of Philosophy 81}, 3 (09 2003), 369--382.
\newblock
  \url{https://sites.ualberta.ca/~francisp/Phil428.526/PellStaintonBivalencePublished03.pdf}.

\bibitem{Penrose}
{\sc Penrose, R.}
\newblock {\em The Road to Reality}.
\newblock Vintage Books, 2007.

\bibitem{Perzanowski}
{\sc Perzanowski, J.}
\newblock Fifty years of parainconsistent logics.
\newblock {\em Logic and Logical Philosophy 7\/} (1999), 21--24.
\newblock
  \url{https://web.archive.org/web/20060404050604/http://www.logika.uni.torun.pl/llp/07/50l.pdf}.

\bibitem{Perzanowski2}
{\sc Perzanowski, J.}
\newblock Parainconsistency, or inconsistency tamed, investigated and
  exploited.
\newblock {\em Logic and Logical Philosophy 9\/} (2001), 5--24.
\newblock
  \url{http://apcz.umk.pl/czasopisma/index.php/LLP/article/download/LLP.2001.001/2694}.

\bibitem{Borwein}
{\sc Peter~Borwein, Stephen~Choi, B. R. a. A.~W.}
\newblock {\em The Riemann Hypothesis: A Resource for the Afficionado and
  Virtuoso Alike}.
\newblock Springer-Verlag, 2007.
\newblock \url{http://wayback.cecm.sfu.ca/~pborwein/TEMP_PROTECTED/book.pdf}.

\bibitem{PhysicsOverflow}
{\sc {Physics Overflow contributors}}.
\newblock Where and how exactly does string theory and q.e.d. use zeta function
  regularization? --- physics overflow.
\newblock
  \url{https://www.physicsoverflow.org/10092/where-exactly-string-theory-zeta-function-regularization}.
\newblock [Online; Last edited 26-March-2014].

\bibitem{Piatetskii-Shapiro}
{\sc Piatetskii-Shapiro, I.}
\newblock Interrelations between the tate and hodge conjectures for abelian
  varieties.
\newblock {\em Matematiceskij Sbornik (Recueil Math\'ematique de la Soci\'et\'e
  Math\'ematique de Moscou) 85(127)}, 4 (1971), 610–620.
\newblock (English translation in Math. USSR Sbornik 14 (1971), 615–625).

\bibitem{Plato}
{\sc Plato}.
\newblock Euthyphro.
\newblock In {\em Essential Dialogues of Plato}, P.~de~Blas, Ed. Barnes \&
  Noble Books, 2005.

\bibitem{Plisko}
{\sc Plisko, V.}
\newblock Law of the {Excluded} {Middle}.
\newblock In {\em Encyclopedia of Mathematics}, {Last} modified on 17 {March}
  2014~ed. Springer-Verlag, 2014.
\newblock
  \url{http://www.encyclopediaofmath.org/index.php?title=Law_of_the_excluded_middle&oldid=31386}.

\bibitem{Plisko2}
{\sc Plisko, V.~E.}
\newblock The kolmogorov calculus as a part of minimal calculus.
\newblock {\em Russian Mathematical Surveys 43}, 6 (1988), 95–110.
\newblock
  \url{https://iopscience.iop.org/article/10.1070/RM1988v043n06ABEH001993}.

\bibitem{Pohlmann}
{\sc Pohlmann, H.}
\newblock Algebraic cycles on abelian varieties of complex multiplication type.
\newblock {\em Annals of Mathematics 88}, 2 (1968), 161–180.
\newblock \url{https://www.jstor.org/stable/1970570}.

\bibitem{Poincare2}
{\sc Poincar\'e, H.}
\newblock Papers on topology: Analysis situs and its five supplements.
\newblock \url{https://www.maths.ed.ac.uk/~v1ranick/papers/poincare2009.pdf}.

\bibitem{Poincare}
{\sc Poincar\'e, H.}
\newblock l'analysis situs.
\newblock {\em Journal de l’\'Ecole Polytechnique 1\/} (1895), 1–121.
\newblock \url{
  https://gallica.bnf.fr/ark:/12148/bpt6k4337198.image.r=langFR.f7.pagination}.

\bibitem{Priest5}
{\sc Priest, G.}
\newblock The logic of paradox.
\newblock {\em Journal of Philosophical Logic 8}, 1 (1979), 219–241.
\newblock \url{https://www.jstor.org/stable/30227165}.

\bibitem{Priest7}
{\sc Priest, G.}
\newblock Logic of paradox revisited.
\newblock {\em Journal of Philosophical Logic 13}, 2 (1984), 153--179.
\newblock \url{https://www.jstor.org/stable/30227027}.

\bibitem{Priest10}
{\sc Priest, G.}
\newblock Paraconsistent logic.
\newblock In {\em Handbook of Philosophical Logic}, D.~M. Gabbay and
  F.~Guenthner, Eds., 2nd~ed., vol.~6. Springer-Verlag, Dordrecht, 2002,
  pp.~287--393.
\newblock \url{https://doi.org/10.1007/978-94-017-0460-1_4}.

\bibitem{Priest1}
{\sc Priest, G.}
\newblock Paraconsistency and dialetheism.
\newblock In {\em Handbook of the History of Logic: The Many Valued and
  Nonmonotonic Turn in Logic}, D.~Gabbay and J.~Woods, Eds., vol.~8. Elsevier,
  2007, pp.~137--139.
\newblock
  \url{https://books.google.com/books?id=3TNj1ZkP3qEC&source=gbs_navlinks_s}.

\bibitem{Priest4}
{\sc Priest, G.}
\newblock {\em An Introduction to Non-Classical Logic: From If to Is}, 2nd~ed.
\newblock Cambridge University Press, 2008.

\bibitem{Priest6}
{\sc Priest, G.}
\newblock The logic of the {Catuskoti}.
\newblock {\em Comparative Philosophy 1}, 2 (2010), 24--54.
\newblock
  \url{https://scholarworks.sjsu.edu/comparativephilosophy/vol1/iss2/6/}.

\bibitem{Priest9}
{\sc Priest, G.}
\newblock What's so bad about contradictions?
\newblock In {\em The Law of Non-Contradiction: New Philosophical Essays},
  J.~C.~B. Graham~Priest and B.~Armour-Garb, Eds. Clarendon Press, 2011,
  pp.~23--38.

\bibitem{Priest8}
{\sc Priest, G.}
\newblock {\em Logic: A Very Short Introduction}, 2nd~ed.
\newblock Oxford University Press, 2017.

\bibitem{Priest2}
{\sc Priest, G., Berto, F., and Weber, Z.}
\newblock Dialetheism.
\newblock In {\em The Stanford Encyclopedia of Philosophy}, E.~N. Zalta, Ed.,
  fall 2018~ed. Metaphysics Research Lab, Stanford University, 2018.
\newblock
  \url{https://plato.stanford.edu/archives/fall2018/entries/dialetheism/}.

\bibitem{Priest3}
{\sc Priest, G., Tanaka, K., and Weber, Z.}
\newblock Paraconsistent logic.
\newblock In {\em The Stanford Encyclopedia of Philosophy}, E.~N. Zalta, Ed.,
  summer 2018~ed. Metaphysics Research Lab, Stanford University, 2018.
\newblock
  \url{https://plato.stanford.edu/archives/sum2018/entries/logic-paraconsistent/}.

\bibitem{Raatikainen}
{\sc Raatikainen, P.}
\newblock Gödel's incompleteness theorems.
\newblock In {\em The Stanford Encyclopedia of Philosophy}, E.~N. Zalta, Ed.,
  fall 2018~ed. Metaphysics Research Lab, Stanford University, 2018.
\newblock
  \url{https://plato.stanford.edu/archives/fall2018/entries/goedel-incompleteness/}.

\bibitem{Ramsey}
{\sc Ramsey, F.~P.}
\newblock The foundations of mathematics.
\newblock {\em Proc. London Math. Soc. 25}, 428 (1926), 338--384.
\newblock
  \url{https://londmathsoc.onlinelibrary.wiley.com/doi/abs/10.1112/plms/s2-25.1.338}.

\bibitem{Rescher}
{\sc Rescher, N.}
\newblock {\em Many-valued logic}.
\newblock McGraw-Hill, 1969.
\newblock \url{https://books.google.com/books?id=ZyTXAAAAMAAJ}.

\bibitem{Ribet}
{\sc Ribet, K.~A.}
\newblock Division fields of abelian varieties with complex multiplication.
\newblock {\em M\'emoires 2\/} (1980), 75--94.
\newblock \url{http://www.numdam.org/item/MSMF_1980_2_2__75_0/}.

\bibitem{Routley}
{\sc Richard~Routley, Val~Plumwood, R. K.~M., and Brady, R.~T.}
\newblock {\em Relevant Logics and Their Rivals}.
\newblock Ridgeview, 1982.
\newblock \url{https://philpapers.org/rec/ROUVT}. [Online; accessed
  16-October-2018].

\bibitem{riemann1859number}
{\sc Riemann, B.}
\newblock On the number of prime numbers less than a given quantity (ueber die
  anzahl der primzahlen unter einer gegebenen grösse).
\newblock {\em Monatsberichte der Berliner Akademie\/} (November 1859).
\newblock Translated in 1998.
  \url{http://www.claymath.org/sites/default/files/ezeta.pdf}. [Online;
  accessed 28-July-2017].

\bibitem{Rowland}
{\sc Rowland, T., and Weisstein, E.~W.}
\newblock Pole.
\newblock In {\em MathWorld--A Wolfram Web Resource}. Wolfram Research, Inc.,
  2018.
\newblock
  \url{https://web.archive.org/web/20181106205407/http://mathworld.wolfram.com/Pole.html}.

\bibitem{Russell3}
{\sc Russell, B.}
\newblock {\em The Principles of Mathematics}, vol.~I.
\newblock Cambridge University Press, 1903.
\newblock \url{https://people.umass.edu/klement/pom/pom-portrait.pdf}.

\bibitem{Russell2}
{\sc Russell, B.}
\newblock On denoting.
\newblock {\em Mind 14}, 56 (1905), 479--493.
\newblock \url{https://www.uvm.edu/~lderosse/courses/lang/Russell(1905).pdf}.

\bibitem{Russell}
{\sc Russell, B.}
\newblock {\em The Problems of Philosophy}.
\newblock Henry Holt and Co., 1912.
\newblock \url{http://www.gutenberg.org/files/5827/5827-h/5827-h.htm}.

\bibitem{Russell6}
{\sc Russell, B.}
\newblock {\em Marriage and Morals}.
\newblock Garden City Publishing Co., Inc., 1929.
\newblock \url{https://archive.org/details/in.ernet.dli.2015.227102}.

\bibitem{Russell5}
{\sc Russell, B.}
\newblock {\em The History of Western Philosophy}.
\newblock Simon \& Schuster, Inc., 1945.

\bibitem{Russell4}
{\sc Russell, B.}
\newblock {\em The Philosophy of Logical Atomism}.
\newblock Taylor \& Francis Group, 2010.
\newblock
  \url{https://sites.ualberta.ca/~francisp/NewPhil448/RussellPhilLogicalAtomismPears.pdf}.

\bibitem{Sadegh-Zadeh}
{\sc Sadegh-Zadeh, K.}
\newblock {\em Handbook of Analytic Philosophy of Medicine}.
\newblock Springer, Dordrecht, 2015.
\newblock
  \url{https://link.springer.com/chapter/10.1007%2F978-94-017-9579-1_32}.

\bibitem{First-Order-Logic}
{\sc Sakharov, A.}
\newblock First-order logic.
\newblock In {\em MathWorld--A Wolfram Web Resource}. Wolfram Research, Inc.,
  2014.
\newblock
  \url{https://web.archive.org/web/20180624154257/http://mathworld.wolfram.com/First-OrderLogic.html}.
  [Online; accessed 24-June-2018].

\bibitem{Sarnak}
{\sc Sarnak, P.}
\newblock Problems of the millennium: The riemann hypothesis (2004).
\newblock {\em Chi Annual Report\/} (2004).
\newblock
  \url{http://www.claymath.org/library/annual_report/ar2004/04report_sarnak.pdf}.

\bibitem{Schroedinger}
{\sc Schrödinger, E.}
\newblock Die gegenwärtige situation in der quantenmechanik.
\newblock {\em Naturwissenschaften 23\/} (1935), 807–812, 823–828,
  844–849.
\newblock \url{https://plus.maths.org/content/not-carrot}.

\bibitem{Scruton}
{\sc Scruton, R.}
\newblock {\em Modern Philosophy}.
\newblock Penguin Books, 1994.

\bibitem{Serre}
{\sc Serre, J.-P.}
\newblock Facteurs locaux des fonctions z\^eta des vari\'et\'es alg\'ebriques
  (d\'efinitions et conjectures).
\newblock {\em S\'eminaire Delange-Pisot-Poitou. Th\'eorie des nombres 11}, 2
  (1969-1970), 1--15.
\newblock \url{http://www.numdam.org/item/SDPP_1969-1970__11_2_A4_0}.

\bibitem{Serre2}
{\sc Serre, J.-P.}
\newblock Groupes alg\'ebriques associ\'e aux modules de hodge-tate.
\newblock {\em Ast\'erisque 65\/} (1979), 155–188.
\newblock \url{http://www.numdam.org/book-part/AST_1979__65__155_0}.

\bibitem{Sharon}
{\sc Sharon, A.}
\newblock Applying the law of identity (loi), the law of non-contradiction
  (lnc), and the law of the excluded middle (lem) to the dirichlet and riemann
  versions of the zeta function, 28 Jun 2018.
\newblock \url{https://arxiv.org/pdf/1802.08062v4.pdf }.

\bibitem{Shioda}
{\sc Shioda, T.}
\newblock What is known about the hodge conjecture?
\newblock In {\em Algebraic Varieties and Analytic Varieties}, S.~Iitaka, Ed.
  Mathematical Society of Japan, 1983, pp.~55--68.
\newblock \url{https://projecteuclid.org/euclid.aspm/1524598012}.

\bibitem{Siegel}
{\sc Siegel, C.}
\newblock Über riemanns nachlaß zur analytischen zahlentheorie.
\newblock {\em Quellen Studien zur Geschichte der Math. Astron. und Phys. Abt.
  B: Studien 2\/} (1932), 45--80.

\bibitem{Silverman}
{\sc Silverman, J.~H.}
\newblock {\em The Arithmetic of Elliptic Curves}, 1st~ed., vol.~106 of {\em
  Graduate Texts in Mathematics}.
\newblock Springer-Verlag, 1992.

\bibitem{Smiley}
{\sc Smiley, T.}
\newblock Sense without detonation.
\newblock {\em Analysis 20\/} (1960), 125--35.

\bibitem{Smith2}
{\sc Smith, D.}
\newblock {\em Philosophy in 50 Milestone Moments}.
\newblock Metro Books, 2017.

\bibitem{Smith3}
{\sc Smith, N.~J.}
\newblock Many-valued logics.
\newblock In {\em Routledge Companion to the Philosophy of Language},
  G.~Russell and D.~G. Fara, Eds., 1st~ed. Routledge, 2010.
\newblock
  \url{https://web.archive.org/web/20180408200831/https://www-personal.usyd.edu.au/~njjsmith/papers/smith-many-valued-logics.pdf}.

\bibitem{Smith}
{\sc Smith, R.}
\newblock Aristotle's {Logic}.
\newblock In {\em The Stanford Encyclopedia of Philosophy}, E.~N. Zalta, Ed.,
  spring 2018~ed. Metaphysics Research Lab, Stanford University, 2018.
\newblock
  \url{https://plato.stanford.edu/archives/spr2018/entries/aristotle-logic/}.

\bibitem{Hamilton}
{\sc Socrates}.
\newblock Republic and {Parmenides}.
\newblock In {\em The Collected Dialogues of Plato}, E.~Hamilton and H.~Cairns,
  Eds. Princeton University Press, 1961, p.~436b.

\bibitem{Speranza}
{\sc Speranza, J., and Horn, L.~R.}
\newblock History of negation.
\newblock In {\em Logic: A History of its Central Concepts}, D.~Gabbay,
  F.~Pelletier, and J.~Woods, Eds. North Holland, 2012.
\newblock
  \url{https://books.google.com/books?id=9mwtRDXVM2wC&lpg=PA148&ots=kj9DmMRVlC&dq=%22vacuous%20subjects%22&pg=PA149#v=onepage&q=%22vacuous%20subjects%22&f=false}.

\bibitem{Stewart}
{\sc Stewart, I.}
\newblock {\em Significant Figures}.
\newblock Hachette Book Group, 2017.

\bibitem{Stillwell}
{\sc Stillwell, J.}
\newblock {\em Elements of Mathematics: From Euclid to Gödel}.
\newblock Princeton University Press, 2016.

\bibitem{Stover}
{\sc Stover, C., and Weisstein, E.~W.}
\newblock Function.
\newblock In {\em MathWorld--A Wolfram Web Resource}. Wolfram Research, Inc.,
  2018.
\newblock
  \url{https://web.archive.org/web/20180903044145/http://mathworld.wolfram.com/Function.html}.
  [Online; accessed 8-October-2018].

\bibitem{Strawson}
{\sc Strawson, P.~F.}
\newblock On referring.
\newblock {\em Mind 59}, 235 (July 1950), 320--344.
\newblock
  \url{http://semantics.uchicago.edu/kennedy/classes/f09/semprag1/strawson50.pdf}.

\bibitem{Suszko}
{\sc Suszko, R.}
\newblock The fregean axiom and polish mathematical logic in the 1920s.
\newblock {\em Studia Logica 36}, 4 (December 1977), 377–380.
\newblock Summary of the talk given to the 22nd Conference on the History of
  Logic, Cracow (Poland), July 5–9, 1976.
  \url{https://link.springer.com/article/10.1007%2FBF02120672}.

\bibitem{Sutherland}
{\sc Sutherland, A.}
\newblock 18.783 elliptic curves lecture \#25.
\newblock
  \url{https://web.archive.org/web/20190709164241/https://math.mit.edu/classes/18.783/2017/LectureNotes25.pdf}.
\newblock [Online; Dated 15-May-2017. Posted Online 17-May-2017].

\bibitem{Tankeev}
{\sc Tankeev, S.~G.}
\newblock Cycles on simple abelian varieties of prime dimension.
\newblock {\em zv. Akad. Nauk SSSR Ser. Mat. 46\/} (1982), 155–170.
\newblock
  \url{http://www.mathnet.ru/php/archive.phtml?wshow=paper&jrnid=im&paperid=1611&option_lang=eng}.

\bibitem{Tarski3}
{\sc Tarski, A.}
\newblock The concept of truth in formalized languages.
\newblock In {\em Logic, Semantics, Metamathematics}, J.~Corcoran, Ed. Hackett,
  1983, pp.~152--278.
\newblock English translation of Tarski's 1936 article.
  \url{http://www.thatmarcusfamily.org/philosophy/Course_Websites/Readings/Tarski%20-%20The%20Concept%20of%20Truth%20in%20Formalized%20Languages.pdf}.

\bibitem{Tarski4}
{\sc Tarski, A.}
\newblock Some observations on the concepts of $\omega$-consistency and
  $\omega$-completeness.
\newblock In {\em Logic, Semantics, Metamathematics}, J.~Corcoran, Ed. Hackett,
  1983, p.~279.
\newblock English translation of Tarski's 1933 article.
  \url{http://www.thatmarcusfamily.org/philosophy/Course_Websites/Readings/Tarski%20-%20The%20Concept%20of%20Truth%20in%20Formalized%20Languages.pdf}.

\bibitem{Tarski}
{\sc Tarski, A.}
\newblock {\em Introduction to Logic}.
\newblock Dover Publications, 1995.

\bibitem{Tarski2}
{\sc Tarski, A., Mostowski, A., and Robinson, R.~M.}
\newblock {\em Undecidable Theories}.
\newblock Dover Publications, 1953.

\bibitem{Tate}
{\sc Tate, J.}
\newblock On the conjectures of birch and swinnerton-dyer and a geometric
  analog.
\newblock {\em Seminaire Bourbaki 66}, 306 (1965/66).
\newblock \url{http://www.numdam.org/article/SB_1964-1966__9__415_0.pdf}.

\bibitem{Tate2}
{\sc Tate, J.~T.}
\newblock Algebraic cycles and poles of zeta functions.
\newblock In {\em Arithmetical Algebraic Geometry}, O.~F.~C. Schilling, Ed.
  Harper \& Row, 1965, p.~93–110.
\newblock \url{(Proc. Conf. Purdue Univ., Dec. 5-7, 1963).}

\bibitem{Wiles3}
{\sc Taylor, R., and Wiles, A.}
\newblock Ring-theoretic properties of certain hecke algebras.
\newblock {\em Annals of Mathematics 141}, 3 (1995), 553–572.
\newblock \url{https://www.jstor.org/stable/2118560}.

\bibitem{Titchmarsh}
{\sc Titchmarsh, E.~C., and Heath-Brown, D.~R.}
\newblock {\em The Theory of the Riemann Zeta-function}, 2nd~ed.
\newblock Oxford University Press, 1999.

\bibitem{Tolstova}
{\sc Tolstova, Y.~N.}
\newblock A weakening of intuitionistic logic.
\newblock {\em Journal of Soviet Mathematics 1}, 1 (1939), 132–138.
\newblock \url{https://link.springer.com/article/10.1007%2FBF01117480 }.

\bibitem{Tong}
{\sc Tong, D.}
\newblock 2. the quantum string.
\newblock In {\em Lectures on String Theory}. University of Cambridge, 2009.
\newblock
  \url{https://web.archive.org/web/20180614034832/http://www.damtp.cam.ac.uk/user/tong/string/two.pdf}.

\bibitem{Toppan}
{\sc Toppan, F.}
\newblock String theory and zeta-function, Oct. 2000.
\newblock
  \url{https://web.archive.org/web/20060924190201/http://www.cbpf.br/~toppan/mo00200.pdf}.

\bibitem{Totaro}
{\sc Totaro, B.}
\newblock {Why} believe the {Hodge} {Conjecture}?, March 2012.
\newblock
  \url{https://burttotaro.wordpress.com/2012/03/18/why-believe-the-hodge-conjecture/}.

\bibitem{Totaro2}
{\sc Totaro, B.}
\newblock {Recent} {Progress} on the {Tate} {Conjecture}.
\newblock {\em Bulletin of the American Mathematical Society 54}, 4 (October
  2017), 575–590.
\newblock
  \url{http://www.ams.org/journals/bull/2017-54-04/S0273-0979-2017-01588-1/S0273-0979-2017-01588-1.pdf}.

\bibitem{Troelstra}
{\sc Troelstra, A., and Schwichtenberg, H.}
\newblock {\em Basic Proof Theory}.
\newblock Cambridge University Press, 2000.

\bibitem{Turco}
{\sc Turco, R., Colonnese, M., and Nardelli, M.}
\newblock Links between string theory and the riemann’s zeta function, Nov.
  2009.
\newblock
  \url{http://empslocal.ex.ac.uk/people/staff/mrwatkin/zeta/nardelli2010a.pdf}.

\bibitem{Turing2}
{\sc Turing, A.~M.}
\newblock A method for the calculation of the zeta-function.
\newblock {\em Proc. London Math. Soc. 48}, 2 (1943), 180--197.
\newblock
  \url{https://londmathsoc.onlinelibrary.wiley.com/doi/10.1112/plms/s2-48.1.180}.

\bibitem{Turing}
{\sc Turing, A.~M.}
\newblock Systems of logic based on ordinals.
\newblock In {\em The Undecidable}, M.~Davis, Ed. Dover Publications, Inc.,
  1993.

\bibitem{Ulmer2}
{\sc Ulmer, D.}
\newblock Park city lectures on elliptic curves over function fields.
\newblock \url{https://arxiv.org/abs/1101.1939}. Published 10 Jan 2011.

\bibitem{Ulmer}
{\sc Ulmer, D.}
\newblock Curves and jacobians over function fields.
\newblock In {\em Arithmetic geometry over global function fields}, Adv.
  Courses Math. CRM Barcelona. Birkhäuser/Springer, 2014, p.~283–337.
\newblock \url{https://doi.org/10.1007/978-3-0348-0853-8_5}.

\bibitem{Urquhart}
{\sc Urquhart, A.}
\newblock Many-valued logic.
\newblock In {\em Handbook of Philosophical Logic}, F.~G. D.~Gabbay, Ed.,
  vol.~166 of {\em Synthese Library book series (SYLI)}. Springer-Verlag, 1986,
  p.~71–116.
\newblock \url{https://link.springer.com/chapter/10.1007/978-94-009-5203-4_2}.

\bibitem{Urquhart2}
{\sc Urquhart, A.}
\newblock Basic many-valued logic.
\newblock In {\em Handbook of Philosophical Logic}, D.~M. Gabbay and
  F.~Guenthner, Eds., 2nd~ed., vol.~2. Springer Netherlands, 2001,
  pp.~249--295.
\newblock \url{https://doi.org/10.1007/978-94-017-0452-6_4}.

\bibitem{Vafeiadou}
{\sc Vafeiadou, G., and Moschovakis, J.~R.}
\newblock Intuitionistic mathematics and logic.
\newblock
  \url{https://web.archive.org/web/20190602013826/http://www.math.ucla.edu/~joan/gvfjrmadobereader.pdf}.

\bibitem{VanGulick}
{\sc Van~Gulick, R.}
\newblock Consciousness.
\newblock In {\em The Stanford Encyclopedia of Philosophy}, E.~N. Zalta, Ed.,
  spring 2018~ed. Metaphysics Research Lab, Stanford University, 2018.

\bibitem{Veneziano2}
{\sc Veneziano, G.}
\newblock Construction of a crossing-symmetric, regge-behaved amplitude for
  linearly rising trajectories.
\newblock {\em Nuovo Cimento A 57\/} (1968), 190--7.
\newblock \url{https://link.springer.com/article/10.1007%2FBF02824451 }.

\bibitem{Visser}
{\sc Visser, A.}
\newblock Four valued semantics and the liar.
\newblock {\em Journal of Philosophical Logic 13}, 2 (1984), 181--212.
\newblock \url{https://www.jstor.org/stable/30227028}.

\bibitem{Wajsberg}
{\sc Wajsberg, M.}
\newblock Untersuchungen über den aussagenkalkül von a. heyting.
\newblock {\em Wiadomości matamatyczne 46\/} (1938), 45–101.
\newblock \url{https://doi.org/10.2307/2267799}.

\bibitem{Weil}
{\sc Weil, A.}
\newblock {\em Oeuvres Scientifiques / Collected Papers}, vol.~II (1951-1964).
\newblock Springer-Verlag, 1979.

\bibitem{Weisstein}
{\sc Weisstein, E.~W.}
\newblock Riemann series theorem.
\newblock In {\em MathWorld--A Wolfram Web Resource}. Wolfram Research, Inc.,
  2005.
\newblock
  \url{https://web.archive.org/web/20051127221618/http://mathworld.wolfram.com/RiemannSeriesTheorem.html}.

\bibitem{Weisstein2}
{\sc Weisstein, E.~W.}
\newblock Taniyama-shimura conjecture.
\newblock In {\em MathWorld--A Wolfram Web Resource}. Wolfram Research, Inc.,
  2010.
\newblock
  \url{https://web.archive.org/web/20100209063657/http://mathworld.wolfram.com/Taniyama-ShimuraConjecture.html}.

\bibitem{Westermann}
{\sc Westermann, C.}
\newblock {\em Argumentationen und Begründungen in der Ethik und Rechtslehre}.
\newblock Duncker \& Humblot, 1977.

\bibitem{Weyl}
{\sc Weyl, H.}
\newblock {\em The Concept of a Riemann Surface}.
\newblock Dover Publications, 2009.

\bibitem{Wheeler}
{\sc Wheeler, J.~A.}
\newblock Physics in copenhagen in 1934 and 1935.
\newblock In {\em Niels Bohr: A Centenary Volume}, A.~P. French and P.~J.
  Kennedy, Eds. Harvard University Press, 1985, p.~221–226.

\bibitem{Whitehead1}
{\sc Whitehead, A.~N.}
\newblock {\em Science and the Modern World}.
\newblock Simon and Schuster, 1925.
\newblock \url{https://books.google.com/books?id=L6kZPLbCrScC}.

\bibitem{Whitehead2}
{\sc Whitehead, A.~N., and Russell, B.}
\newblock {\em Principia mathematica}, 2nd~ed., vol.~1.
\newblock Cambridge University Press, 1925.
\newblock \url{https://archive.org/details/PrincipiaMathematicaVolumeI}.

\bibitem{Whittaker}
{\sc Whittaker, E.~T., and Watson, G.~N.}
\newblock {\em A Course of Modern Analysis}, 4th~ed.
\newblock Cambridge University Press, 1920.
\newblock
  \url{https://archive.org/stream/courseofmodernan00whit#page/n0/mode/2up}.
  [Online; accessed 15-October-2017].

\bibitem{Aristotle3}
{\sc {Wikipedia contributors}}.
\newblock Aristotle --- {Wikipedia}{,} the free encyclopedia.
\newblock
  \url{https://en.wikipedia.org/w/index.php?title=Aristotle&oldid=900073268}.
\newblock [Online; Last edited 3-June-2019].

\bibitem{Classical}
{\sc {Wikipedia contributors}}.
\newblock Classical logic --- {Wikipedia}{,} the free encyclopedia.
\newblock
  \url{https://en.wikipedia.org/w/index.php?title=Classical_logic&oldid=822435045}.
\newblock [Online; Last edited 26-January-2018].

\bibitem{DecisionProblem}
{\sc {Wikipedia contributors}}.
\newblock Decision problem --- {Wikipedia}{,} the free encyclopedia.
\newblock
  \url{https://en.wikipedia.org/w/index.php?title=Decision_problem&oldid=873678788}.
\newblock [Online; Last edited 14-December-2018].

\bibitem{Dirichlet_L-function}
{\sc {Wikipedia contributors}}.
\newblock Dirichlet l-function --- {Wikipedia}{,} the free encyclopedia.
\newblock
  \url{https://en.wikipedia.org/w/index.php?title=Dirichlet_L-function&oldid=862708763}.
\newblock [Online; Last edited 6-October-2018].

\bibitem{Dirichlet_series}
{\sc {Wikipedia contributors}}.
\newblock Dirichlet series --- {Wikipedia}{,} the free encyclopedia.
\newblock
  \url{https://en.wikipedia.org/w/index.php?title=Dirichlet_series&oldid=872492273}.
\newblock [Online; Last edited 7-December-2018].

\bibitem{Galilean_inv}
{\sc {Wikipedia contributors}}.
\newblock Galilean invariance --- {Wikipedia}{,} the free encyclopedia.
\newblock
  \url{https://en.wikipedia.org/w/index.php?title=Galilean_invariance&oldid=866155140}.
\newblock [Online; Last edited 28-October-2018].

\bibitem{GRH}
{\sc {Wikipedia contributors}}.
\newblock Generalized riemann hypothesis --- {Wikipedia}{,} the free
  encyclopedia.
\newblock \url{
  https://en.wikipedia.org/w/index.php?title=Generalized_Riemann_hypothesis&oldid=883901968}.
\newblock [Online; Last edited 18-February-2019].

\bibitem{IncompletenessTheorems}
{\sc {Wikipedia contributors}}.
\newblock Gödel's incompleteness theorems --- {Wikipedia}{,} the free
  encyclopedia.
\newblock
  \url{https://en.wikipedia.org/w/index.php?title=G%C3%B6del%27s_incompleteness_theorems&oldid=871124813
  }.
\newblock [Online; Last edited 29-November-2018].

\bibitem{HasseWeil}
{\sc {Wikipedia contributors}}.
\newblock Hasse–weil zeta function --- {Wikipedia}{,} the free encyclopedia.
\newblock
  \url{https://en.wikipedia.org/w/index.php?title=Hasse%E2%80%93Weil_zeta_function&oldid=862709306}.
\newblock [Online; Last edited 6-October-2018].

\bibitem{Tree}
{\sc {Wikipedia contributors}}.
\newblock If a tree falls in a forest --- {Wikipedia}{,} the free encyclopedia.
\newblock
  \url{https://en.wikipedia.org/w/index.php?title=If_a_tree_falls_in_a_forest&oldid=896534541}.
\newblock [Online; Last edited 11-May-2019].

\bibitem{Popper}
{\sc {Wikipedia contributors}}.
\newblock Karl popper --- {Wikipedia}{,} the free encyclopedia.
\newblock
  \url{https://en.wikipedia.org/w/index.php?title=Karl_Popper&oldid=899287706}.
\newblock [Online; Last edited 29-May-2019].

\bibitem{L-function}
{\sc {Wikipedia contributors}}.
\newblock L-function --- {Wikipedia}{,} the free encyclopedia.
\newblock
  \url{https://en.wikipedia.org/w/index.php?title=L-function&oldid=904595602}.
\newblock [Online; Last edited 3-July-2019].

\bibitem{LiarParadox}
{\sc {Wikipedia contributors}}.
\newblock Liar paradox --- {Wikipedia}{,} the free encyclopedia.
\newblock
  \url{https://en.wikipedia.org/w/index.php?title=Liar_paradox&oldid=871230217}.
\newblock [Online; Last edited 29-November-2018].

\bibitem{Lindelof}
{\sc {Wikipedia contributors}}.
\newblock Lindelöf hypothesis --- {Wikipedia}{,} the free encyclopedia.
\newblock
  \url{https://en.wikipedia.org/w/index.php?title=Lindel%C3%B6f_hypothesis&oldid=887266061
  }.
\newblock [Online; Last edited 11-March-2019].

\bibitem{LogicalTruth}
{\sc {Wikipedia contributors}}.
\newblock Logical truth --- {Wikipedia}{,} the free encyclopedia.
\newblock
  \url{https://en.wikipedia.org/w/index.php?title=Logical_truth&oldid=886686295}.
\newblock [Online; Last edited 7-March-2019].

\bibitem{Material_implication}
{\sc {Wikipedia contributors}}.
\newblock Material implication (rule of inference) --- {Wikipedia}{,} the free
  encyclopedia.
\newblock
  \url{https://en.wikipedia.org/w/index.php?title=Material_implication_(rule_of_inference)&oldid=884836188}.
\newblock [Online; Last edited 24-February-2019].

\bibitem{Minimal}
{\sc {Wikipedia contributors}}.
\newblock Minimal logic --- {Wikipedia}{,} the free encyclopedia.
\newblock
  \url{https://en.wikipedia.org/w/index.php?title=Minimal_logic&oldid=822881725}.
\newblock [Online; Last edited 29-January-2018].

\bibitem{Muenchhausen-Trilemma}
{\sc {Wikipedia contributors}}.
\newblock Münchhausen trilemma --- {Wikipedia}{,} the free encyclopedia.
\newblock
  \url{https://en.wikipedia.org/w/index.php?title=M%C3%BCnchhausen_trilemma&oldid=873254462}.
\newblock [Online; Last edited 12-December-2018].

\bibitem{Bohr}
{\sc {Wikipedia contributors}}.
\newblock Niels bohr --- {Wikipedia}{,} the free encyclopedia.
\newblock
  \url{https://en.wikipedia.org/w/index.php?title=Niels_Bohr&oldid=901863245}.
\newblock [Online; Last edited 19-June-2019].

\bibitem{Nyaya}
{\sc {Wikipedia contributors}}.
\newblock Nyaya --- {Wikipedia}{,} the free encyclopedia.
\newblock
  \url{https://en.wikipedia.org/w/index.php?title=Nyaya&oldid=897539143}.
\newblock [Online; Last edited 17-May-2019].

\bibitem{Explosion}
{\sc {Wikipedia contributors}}.
\newblock Principle of explosion --- {Wikipedia}{,} the free encyclopedia.
\newblock
  \url{https://en.wikipedia.org/w/index.php?title=Principle_of_explosion&oldid=898404042}.
\newblock [Online; Last edited 23-May-2019].

\bibitem{Euler}
{\sc {Wikipedia contributors}}.
\newblock Proof of the euler product formula for the riemann zeta function ---
  {Wikipedia}{,} the free encyclopedia.
\newblock
  \url{https://en.wikipedia.org/w/index.php?title=Proof_of_the_Euler_product_formula_for_the_Riemann_zeta_function&oldid=852102893}.
\newblock [Online; Last edited 26-July-2018].

\bibitem{Relevance}
{\sc {Wikipedia contributors}}.
\newblock Relevance logic --- {Wikipedia}{,} the free encyclopedia.
\newblock
  \url{https://en.wikipedia.org/w/index.php?title=Relevance_logic&oldid=836099696}.
\newblock [Online; Last edited 12-April-2018].

\bibitem{Riemann-Siegel}
{\sc {Wikipedia contributors}}.
\newblock Riemann–siegel formula --- {Wikipedia}{,} the free encyclopedia.
\newblock
  \url{https://en.wikipedia.org/w/index.php?title=Riemann%E2%80%93Siegel_formula&oldid=855414500
  }.
\newblock [Online; Last edited 18-August-2018].

\bibitem{UndefinabilityTheorem}
{\sc {Wikipedia contributors}}.
\newblock Tarski's undefinability theorem --- {Wikipedia}{,} the free
  encyclopedia.
\newblock
  \url{https://en.wikipedia.org/w/index.php?title=Tarski%27s_undefinability_theorem&oldid=862710629
  }.
\newblock [Online; Last edited 6-October-2018].

\bibitem{3VL}
{\sc {Wikipedia contributors}}.
\newblock Three-valued logic --- {Wikipedia}{,} the free encyclopedia.
\newblock
  \url{https://en.wikipedia.org/w/index.php?title=Three-valued_logic&oldid=863963861}.
\newblock [Online; Last edited 14-October-2018].

\bibitem{Veneziano1}
{\sc {Wikipedia contributors}}.
\newblock Veneziano amplitude --- {Wikipedia}{,} the free encyclopedia.
\newblock
  \url{https://en.wikipedia.org/w/index.php?title=Veneziano_amplitude&oldid=829103935}.
\newblock [Online; Last edited 6-March-2018].

\bibitem{Wittgenstein}
{\sc {Wikipedia contributors}}.
\newblock Wittgenstein --- {Wikipedia}{,} the free encyclopedia.
\newblock
  \url{https://en.wikipedia.org/w/index.php?title=Ludwig_Wittgenstein&oldid=867922333}.
\newblock [Online; Last edited 8-November-2018].

\bibitem{ZFC}
{\sc {Wikipedia contributors}}.
\newblock Zermelo–fraenkel set theory --- {Wikipedia}{,} the free
  encyclopedia.
\newblock
  \url{https://en.wikipedia.org/w/index.php?title=Zermelo%E2%80%93Fraenkel_set_theory&oldid=896993185
  }.
\newblock [Online; Last edited 14-May-2019].

\bibitem{Poles}
{\sc {Wikipedia contributors}}.
\newblock Zeros and poles --- {Wikipedia}{,} the free encyclopedia.
\newblock \url{https://en.wikipedia.org/w/index.php?oldid=878048355}.
\newblock [Online; Last edited 12-January-2019].

\bibitem{Analysis_Situs}
{\sc {Wikipedia contributors}}.
\newblock Analysis situs (paper) --- {Wikipedia}{,} the free encyclopedia.
\newblock
  \url{https://en.wikipedia.org/w/index.php?title=Analysis_Situs_(paper)&oldid=904974336},
  2019.
\newblock [Online; last edited 5-July-2019].

\bibitem{Relativity}
{\sc {Wikipedia contributors}}.
\newblock Relativity of simultaneity --- {Wikipedia}{,} the free encyclopedia.
\newblock
  \url{https://en.wikipedia.org/w/index.php?title=Relativity_of_simultaneity&oldid=904754466},
  2019.
\newblock [Online; last edited 4-July-2019].

\bibitem{Wiles}
{\sc Wiles, A.}
\newblock The {Birch} and {Swinnerton-Dyer} {Conjecture}.
\newblock \url{http://www.claymath.org/sites/default/files/birchswin.pdf}.
  [Online; accessed 25-January-2018].

\bibitem{Wiles2}
{\sc Wiles, A.}
\newblock Modular elliptic curves and fermat’s last theorem.
\newblock {\em Annals of Mathematics 141}, 3 (1995), 443–551.
\newblock \url{https://www.jstor.org/stable/2118559}.

\bibitem{witten1991}
{\sc Witten, E.}
\newblock On quantum gauge theories in two dimensions.
\newblock {\em Communications in Mathematical Physics 141}, 1 (1991), 153--209.
\newblock \url{https://projecteuclid.org:443/euclid.cmp/1104248198}.

\bibitem{Wolenski}
{\sc Woleński, J.}
\newblock Lvov-warsaw school.
\newblock In {\em The Stanford Encyclopedia of Philosophy}, E.~N. Zalta, Ed.,
  fall 2017~ed. Metaphysics Research Lab, Stanford University, 2017.
\newblock
  \url{https://plato.stanford.edu/archives/fall2017/entries/lvov-warsaw/}.

\bibitem{Woods2003}
{\sc Woods, J.}
\newblock {\em Paradox and Paraconsistency}.
\newblock Cambridge University Press, 2003.

\bibitem{Woods2005}
{\sc Woods, J.}
\newblock Dialectical considerations on the logic of contradiction: Part i.
\newblock {\em Logic Journal of the IGPL 13\/} (2005), 231–260.

\bibitem{Lukasiewicz}
{\sc Łukasiewicz, J.}
\newblock O zasadzie sprzeczności u arystotelesa.
\newblock {\em Bulletin International de l’Acad\'emie des Sciences de
  Cracovie 1-2\/} (1910), 15--38.

\bibitem{Lukasiewicz4}
{\sc Łukasiewicz, J.}
\newblock O logice trójwartościowej.
\newblock {\em Ruch Filozoficny 5\/} (1920), 170--171.

\bibitem{Lukasiewicz3}
{\sc Łukasiewicz, J.}
\newblock On the principle of contradiction in aristotle.
\newblock {\em The Review of Metaphysics 24}, 3 (March 1971), 485--509.
\newblock \url{https://www.jstor.org/stable/20125812}.

\bibitem{Lukasiewicz2}
{\sc Łukasiewicz, J.}
\newblock Zasada sprzeczności u arystotelesa.
\newblock {\em Filozofia Nauki 5}, 1 (1997), 147--164.
\newblock (Reprint of original).
  \url{http://bazhum.muzhp.pl/media//files/Filozofia_Nauki/Filozofia_Nauki-r1997-t5-n1/Filozofia_Nauki-r1997-t5-n1-s147-164/Filozofia_Nauki-r1997-t5-n1-s147-164.pdf}.

\end{thebibliography}

\end{document}